\documentclass[review]{elsarticle}

\usepackage{geometry}
\usepackage{float}
\usepackage{amsmath,amssymb,amsfonts,amsthm}    

\usepackage{graphicx}                           
\usepackage[font=small,skip=0pt]{caption}
\usepackage{subcaption}
\usepackage{mathtools}
\usepackage{multirow,makecell}
\usepackage{array}
\usepackage{accents}
\usepackage{bbm}
\usepackage[hidelinks]{hyperref}
\hypersetup{
	colorlinks = true,
	linkcolor = blue,
	citecolor = blue,
	urlcolor = blue
}
\usepackage{leftidx}
\usepackage{chngcntr}
\usepackage[page]{appendix}
\usepackage{xcolor}
\usepackage{relsize}
\usepackage{bigstrut}
\usepackage{booktabs}
\usepackage{natbib}
\usepackage[makeroom]{cancel}
\theoremstyle{plain}
\newtheorem{thm}{Theorem}[section]

\theoremstyle{definition}
\newtheorem{zremark}[thm]{Remark}
\newenvironment{remark}{\par\small\zremark}{\endzremark}

\theoremstyle{definition}
\newtheorem{zobservation}[thm]{Observation}
\newenvironment{observation}{\par\small\zobservation}{\endzobservation}

\makeatletter
\DeclareRobustCommand{\rvdots}{%
  \vbox{
    \baselineskip4\p@\lineskiplimit\z@
    \kern-\p@
    \hbox{.}\hbox{.}\hbox{.}
  }}
\makeatother
\usepackage{cleveref}

\journal{Computer Methods in Applied Mechanics and Engineering}

\usepackage{etoolbox}
\patchcmd{\pprintMaketitle}
  {\fi\hrule}
  {\fi\ifvoid\extrainfobox\else\unvbox\extrainfobox\par\vskip5pt\fi\hrule}
  {}{}

\newenvironment{extrainfo}
  {\global\setbox\extrainfobox=\vbox\bgroup\parindent=0pt }
  {\egroup}
\newsavebox\extrainfobox

\usepackage{numcompress}\biboptions{authoryear}

\date{}
\begin{document}
\emergencystretch 3em

\begin{frontmatter}

\title{A selectively reduced degree basis for efficient mixed nonlinear isogeometric beam formulations with extensible directors}

\author[LBB]{Myung-Jin Choi\corref{mycorrespondingauthor}}
\ead{choi@lbb.rwth-aachen.de}

\author[AICES,Gdansk,IIT]{Roger A. Sauer}
\ead{sauer@aices.rwth-aachen.de}

\author[LBB]{Sven Klinkel}
\ead{klinkel@lbb.rwth-aachen.de}

\cortext[mycorrespondingauthor]{Corresponding author}

\address[LBB]{Chair of Structural Analysis and Dynamics, RWTH Aachen University, Mies-van-der-Rohe Str.\,1,\,52074 Aachen, Germany}
\address[AICES]{Aachen Institute for Advanced Study in Computational Engineering Science (AICES), RWTH Aachen University, Templergraben 55, 52062 Aachen, Germany}
\address[Gdansk]{Faculty of Civil and Environmental Engineering, Gda\'{n}sk University of Technology, ul. Narutowicza 11/12, 80-233 Gda\'{n}sk, Poland}
\address[IIT]{Department of Mechanical Engineering, Indian Institute of Technology Guwahati, Assam 781039, India}

\begin{extrainfo}
	\begin{center}
		{\small To be published\footnote{This PDF is the personal version of an article whose final publication is available at \href{www.sciencedirect.com}{www.sciencedirect.com}} in \textit{Computer Methods in Applied Mechanics and Engineering}}\\ 
	{\small Submitted on 19 June 2023, Revised on 15 August 2023, Accepted on 15 August 2023}
\end{center}
\end{extrainfo}

\begin{abstract}
The effect of higher order continuity in the solution field by using NURBS basis function in isogeometric analysis (IGA) is investigated for an efficient mixed finite element formulation for elastostatic beams. It is based on the Hu-Washizu variational principle considering geometrical and material nonlinearities. Here we present a reduced degree of basis functions for the additional fields of the stress resultants and strains of the beam, which are allowed to be discontinuous across elements. This approach turns out to significantly improve the computational efficiency and the accuracy of the results. We consider a beam formulation with extensible directors, where cross-sectional strains are enriched to avoid Poisson locking by an enhanced assumed strain method. In numerical examples, we show the superior per degree-of-freedom accuracy of IGA over conventional finite element analysis, due to the higher order continuity in the displacement field. We further verify the efficient rotational coupling between beams, as well as the path-independence of the results.\\
\end{abstract}


\end{frontmatter}


%
\section{Introduction}
Beam models have been widely utilized for an efficient and accurate mechanical simulation of slender rods, rod-like bodies, and their assemblies, across many application areas, see e.g., simulations of Brownian dynamics of microstructures \citep{cyron2012numerical}, atomistic structures \citep{schmidt2015continuum}, and entangled materials \citep{durville2005numerical}. Accordingly, many beam models have been developed, whose applicability is, however, often limited by a certain range of geometrical characteristics of the beam, e.g., the slenderness ratio, and the cross-sectional shape. For an extensive classification of linear and nonlinear beam theories, we refer readers to the recent review in \citet{meier2019geometrically}. In this paper, we aim at expanding the applicability of our nonlinear beam formulation from \textcolor{black}{the} previous work \citep{choi2021isogeometric} to the thin beam limit, with much improved robustness and efficiency. In this paper, we particularly investigate the effect of higher order continuity in \textcolor{black}{the displacement fields of the center axis and directors}, in the framework of isogeometric analysis (IGA), using non-uniform rational B-splines (NURBS) basis functions.

Isogeometric analysis was presented by \citet{hughes2005isogeometric} to incorporate the exact geometrical description, inherent in computer-aided design (CAD) model, into analysis, by employing the same spline basis functions utilized in CAD system. This approach enables not only the exact description of the initial geometry, but also the higher order continuity in the solution field. In this study, we focus on the superiority of the $k$-version of mesh refinement (smooth degree elevation) in IGA in terms of per degree-of-freedom (DOF) accuracy in elastostatic nonlinear beam problems. We compare the results with conventional FEA with Lagrange shape functions. For further extension to dynamic problems like structural vibration and wave propagation, one may refer to \citet{hughes2008duality}. The smoothness property in IGA enables us to use much less DOFs per element, compared with conventional FEA. To simply illustrate this attractive property, we consider a mesh in a one-dimensional domain. For $C^0$-continuous finite elements, the number of nodes (or basis functions) is expressed by
\begin{equation}
	n_\mathrm{nd}={p}\cdot{{n_\mathrm{el}}}+1,
\end{equation}
where $p$ and $n_\mathrm{el}$ denote the degree of basis functions, and the number of elements, respectively. In contrast, for $C^{p-1}$-continuous IGA, the number of control points (or basis functions) is 
\begin{equation}
	n_\mathrm{cp}={{n_\mathrm{el}}}+p.
\end{equation}
That is, as we increase the number of elements ($n_\mathrm{el}$), the increase in the number of nodes in FEA is much larger than that in IGA, and this gap increases as we choose higher $p$. This property simply extends to multi-dimensional cases, see, e.g., \citet{cottrell2006isogeometric}. However, it turns out that IGA also suffers from numerical locking \citep{echter2010numerical} in constrained problems like bending-dominated slender beams. It was also shown that the higher order continuity in the displacement and rotation fields may accentuate the locking \citep{adam2014improved}. 

Locking can be attributed to a \textit{field-inconsistency} in the approximated strains in the conventional displacement-based finite element formulation \citep{prathap2013finite}. In order to alleviate locking, many approaches have been developed, which include reduced integration methods, and mixed variational formulations. In the framework of IGA, they can be divided into \textit{local} (element-wise) and \textit{global} (patch-wise) approaches. An \textit{element-wise} selective and reduced integration (SRI) method falls into the first category. A SRI rule for quadratic NURBS basis functions was presented in \citet{bouclier2012locking}, and it was generalized by \citet{adam2014improved} in terms of the degree of basis functions and the inter-element continuity. A key observation of the latter work is that additional constraints in the stiffness matrix may arise due to higher inter-element continuity in the membrane and transverse shear strain fields, which can also lead to a significant spurious increase of the bending stiffness. Thus, for IGA, a local approach always needs to be combined with an additional treatment, e.g., further reduction of the number of quadrature points in the SRI, see \citet{adam2014improved}. A local approach can be also found in the context of $\bar B$ projection methods for IGA. For example, \citet{hu2016order} presented a reduction of the degree of basis in the local (element-wise) projection space, which was extended to shell problems in \citet{hu2020isogeometric}. On the other hand, a \textit{patch-wise} quadrature rule, which falls into the global approach, was initiated by an observation that the conventional element-wise quadrature rule for FEA may not be optimal for IGA \citep{hughes2010efficient}. Patch-wise reduced integration approaches to alleviate locking were also developed, e.g., SRI in \citet{adam2015selective}, and a Greville quadrature rule in \citet{zou2021galerkin}. A global approach in mixed formulations suffers from a significant increase of the computational cost due to the inversion of the global Gram matrix, and the resulting \textit{dense} global stiffness matrix, see, e.g., the $\bar B$ projection method in \citet{bouclier2012locking}. To overcome this, a local least square (LLSQ) method \citep{govindjee2012convergence} was employed in \citet{bouclier2013efficient}. \textcolor{black}{The LLSQ method converts the given least square problem into an approximated one, which solves a set of independent \textit{element-wise} equations. This is much more efficient than solving the original least square problem, since it inverts the element Gram matrices, not the global one. However, the approximation error may increase, as degree of basis $p$ and the inter-element continuity increase.} A similar concept was employed in the mixed formulation for alleviating membrane locking in plane Kirchhoff rods by \citet{greco2017efficient}. This approach was extended to shell problems in \citet{kikis2022two}, based on the Hellinger-Reissner variational principle, where the test functions are allowed to have \textcolor{black}{inter-element} discontinuities, so that the finite element equations for the test functions of the additional solution fields can be treated \textcolor{black}{element-wisely}. This enables an efficient element-wise static condensation. However, this selection of different function spaces for the test functions and \textcolor{black}{solution} leads to unsymmetric stiffness matrix. A mixed formulation was also developed in the context of isogeometric collocation method for nonlinear beams, e.g., \cite{marino2017locking} and \cite{weeger2017isogeometric}, where the control variables of the additional fields are not condensed out, which may significantly increase the number of DOFs, but the stiffness matrix is still sparse. \textcolor{black}{Isogeometric collocation methods are extended to deal with elasto-visco-plasticity and visco-elasticity problems in \citet{weeger2022mixed} and \citet{ferri2023efficient}, respectively.} A more intuitive way to achieve field-consistency, in geometrically linear problems, is to utilize one degree lower bases for the rotation field, e.g., see \citet{da2012isogeometric} for plate problems and \citet{kikis2019adjusted} for plate and shell problems.

A beam can be regarded as a spatial curve with attached deformable director vectors (or \textit{directors}). This curve is then also called \textit{directed} or \textit{Cosserat} curve. In this paper, the mixed formulation is based on the first order beam kinematics in \cite{choi2021isogeometric}, where it was shown that the displacement-based beam formulation suffers from transverse shear, membrane, and curvature-thickness locking. For an illustration of the curvature-thickness locking, one may refer to \cite{betsch1995assumed}. It was also observed that the tangent stiffness matrix becomes ill-conditioned in the thin beam limit, which leads to instability in the Newton-Raphson solution process. In our present study, a Hu-Washizu variational principle is employed, where we introduce additional solution fields for the stress resultants and strains of the beam. This provides the following advantages:
\begin{itemize}
	\item Alleviation of transverse shear, membrane, and curvature-thickness locking,
	\item \textcolor{black}{Improved convergence of the Newton-Raphson iteration for larger load increments \citep{klinkel2006robust, betsch2016energy}}, and in the thin beam limit,
	\item In contrast to the Hellinger-Reissner principle, the stress resultants and strain of the beam are independent from each other, so that we can use nonlinear constitutive laws \citep{santos2010hybrid}.
\end{itemize}
Our mixed formulation is basically a \textit{local approach}, and we alleviate the additional artificial constraints due to the higher order continuity, by reducing the degree of basis functions in the additional solution fields, i.e., the beam's stress resultants and strains. This approach has the following novelties over previous local and global approaches:
\begin{itemize}
	\item It enables an element-wise static condensation. 
	\item It reduces the number of internal variables, and the size of relevant stiffness (sub-)matrices, which also makes the matrix operations in the condensation process computationally more efficient. 
	\item It uses the same function space for the test and solution functions of the additional fields, so that the resulting stiffness matrix is always symmetric \textcolor{black}{for conservative loads}.\footnote{\textcolor{black}{Here we assume no rotational coupling conditions. The relevant discussion is given in Section \ref{rcont_betw_b}.}}
\end{itemize}

The remainder of this paper is organized as follows. In Section \ref{bkin}, we briefly review the beam kinematics with extensible directors, and recall the expressions of the strains, and the stress resultants of the beam. In Section \ref{3fmixed}, we present a mixed finite element formulation, based on the Hu-Washizu variational principle. In Section \ref{ig_fe_disc_sec}, we present the isogeometric finite element discretization. In Section \ref{rcont_betw_b}, we present a formulation for imposing the rotational continuity between beams. In Section \ref{num_example}, several numerical examples are presented. Section \ref{conclu} concludes the paper.

\section{Beam kinematics}
\label{bkin}
\begin{figure}[H]	
	\centering
	\begin{subfigure}[b] {0.55\textwidth} \centering
		\includegraphics[width=\linewidth]{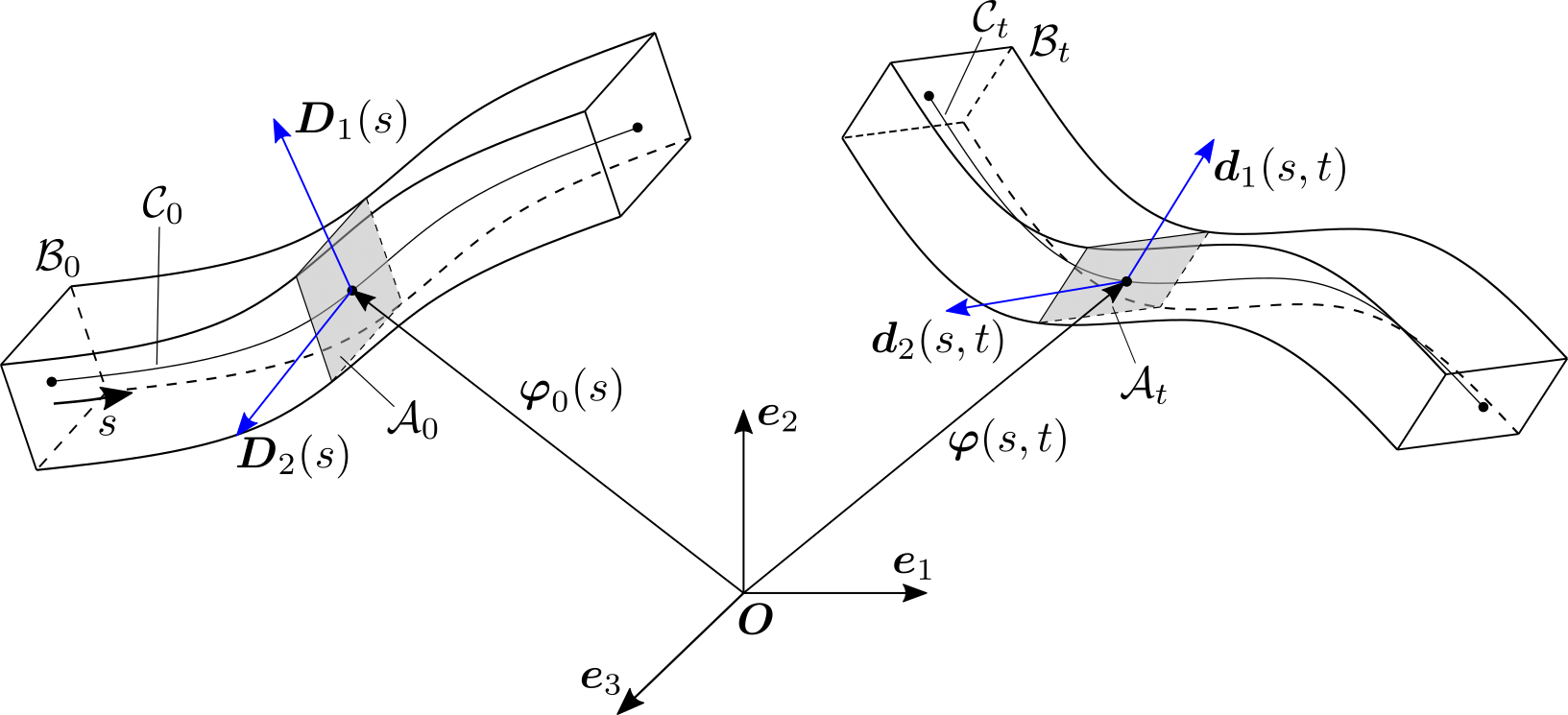}	
	\end{subfigure}	
	\caption{An illustration of the beam kinematics in the initial and current configurations, whose domains are denoted by $\mathcal{B}_0$ and $\mathcal{B}_t$, respectively. $\mathcal{C}_0$ and $\mathcal{C}_t$ indicate the initial and current center axes, respectively. $\mathcal{A}_0$ and $\mathcal{A}_t$ indicate the initial and current cross-sections, respectively. $\left\{{\boldsymbol{e}}_1,{\boldsymbol{e}}_2,{\boldsymbol{e}}_3\right\}$ denotes the standard Cartesian basis in $\Bbb{R}^3$. This figure is redrawn with modifications from \citet{choi2021isogeometric}.}
	\label{bkin_conf_init_cur}
\end{figure}

\begin{figure}[H]	
	\centering
	\begin{subfigure}[b] {0.35\textwidth} \centering
		\includegraphics[width=\linewidth]{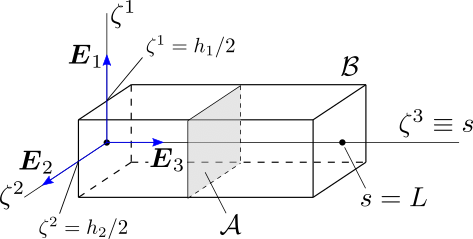}	
	\end{subfigure}		
	\caption{An example of the reference domain $\mathcal{B}$ having a rectangular cross-section with dimension ${h_1}\times{h_2}$. $\mathcal{A}$ denotes the reference domain of the cross-section. This figure is redrawn with modifications from \citet{choi2021isogeometric}.}
	\label{bkin_ref_conf}
\end{figure}

Two transverse directions of the beam are defined by the principal directions of the second moment of inertia tensor in the initial (undeformed) cross-section. The origin of the transverse coordinates $\zeta^\alpha\,(\alpha\in\left\{1,2\right\})$ is defined by the geometrical center of the initial cross-section, which coincides with the mass center under the assumption of constant mass density in the initial configuration. The line connecting these center points of the cross-sections is called \textcolor{black}{a \textit{center axis}}, whose position is denoted by ${\boldsymbol{\varphi }}({s},t)$. Here, ${s}$ denotes the arc-length coordinate along the initial center axis, and $t$ denotes time. In this paper, the argument $t$ is often omitted for brevity. The initial beam configuration is expressed by
\begin{equation}
	{{\boldsymbol{x}}}_0({\zeta ^1},{\zeta ^2},s) = {\boldsymbol{\varphi }}_0(s) + {\zeta ^\gamma }{{\boldsymbol{D}}_\gamma }(s),
\end{equation}
where ${\boldsymbol{\varphi }}_0(s)$ denotes the position of the beam's initial center axis, and the initial cross-sectional plane is spanned by two \textit{initial directors} ${\boldsymbol{D}}_\gamma(s) \in \Bbb{R}^3$ $(\gamma\in\left\{1,2\right\})$. Here and hereafter, unless stated otherwise, repeated Greek indices like $\alpha$, $\beta$, and $\gamma$ imply summation over $1$ to $2$, and repeated Latin indices like $i$ and $j$ imply summation over $1$ to $3$. We define a covariant basis ${\boldsymbol{G}}_i\coloneqq{\partial{\boldsymbol{x}}_0}/{\partial{\zeta^i}}$ $(i\in\left\{1,2,3\right\})$ with $\zeta^3\equiv{s}$, and a contravariant basis $\left\{{\boldsymbol{G}}^1,{\boldsymbol{G}}^2,{\boldsymbol{G}}^3\right\}$ is defined by the orthogonality condition $\boldsymbol{G}_i\cdot\boldsymbol{G}^j=\delta_i^j$, where $\delta_i^j$ denotes the Kronecker-delta. According to the first order beam kinematics \citep{rhim1998vectorial,durville2012contact,choi2021isogeometric}, the position vector is a linear function of the coordinates ${\zeta ^\gamma }$ along the transverse directions, i.e.,
\begin{equation}
	\label{bkin_x_d_zta_ph}
	{{\boldsymbol{x}}}({\zeta ^1},{\zeta ^2},s,t) = {\boldsymbol{\varphi }}(s,t) + {\zeta ^\gamma }{{\boldsymbol{d}}_\gamma }(s,t),
\end{equation}
where ${\boldsymbol{\varphi }}(s,t)$ denotes the current position of the center axis, and the current cross-sectional plane is spanned by two \textit{current directors} ${\boldsymbol{d}}_\gamma(s,t) \in \Bbb{R}^3$ $(\gamma\in\left\{1,2\right\})$. Fig.\,\ref{bkin_conf_init_cur} schematically illustrates this beam kinematics. Fig.\,\ref{bkin_ref_conf} shows a reference domain in the case of a rectangular cross-section, where $\left\{{\boldsymbol{E}}_1,{\boldsymbol{E}}_2,{\boldsymbol{E}}_3\right\}$ denotes the standard Cartesian basis in $\Bbb{R}^3$, with ${\boldsymbol{E}}_3$ along the axial direction, and ${\boldsymbol{E}}_1$ and ${\boldsymbol{E}}_2$ along two transverse directions.
\noindent The covariant components of the Green-Lagrange strain tensor $\boldsymbol{E}=E_{ij}\,\boldsymbol{G}^i\otimes\boldsymbol{G}^j$ can be expressed in terms of the beam strains, as \citep{choi2021isogeometric}
\begin{equation}
	\label{gl_strn_mata_beam_strn}
	\underline {\boldsymbol{E}}= \left[ {\begin{array}{*{20}{c}}
			0&0&0&0&0&0&0&0&0&0&0&0&1&0&0\\
			0&0&0&0&0&0&0&0&0&0&0&0&0&1&0\\
			1&{{\zeta ^1}}&{{\zeta ^2}}&{{\zeta ^1}{\zeta ^1}}&{{\zeta ^2}{\zeta ^2}}&{{\zeta ^1}{\zeta ^2}}&0&0&0&0&0&0&0&0&0\\
			0&0&0&0&0&0&0&0&0&0&0&0&0&0&1\\
			0&0&0&0&0&0&1&0&{{\zeta ^1}}&{{\zeta ^2}}&0&0&0&0&0\\
			0&0&0&0&0&0&0&1&0&0&{{\zeta ^1}}&{{\zeta ^2}}&0&0&0
	\end{array}} \right]\left\{ {\begin{array}{*{20}{c}}
			\varepsilon \\
			{\boldsymbol{\rho }}\\
			{\boldsymbol{\kappa }}\\
			{\boldsymbol{\delta }}\\
			{\boldsymbol{\gamma }}\\
			{\boldsymbol{\chi }}
	\end{array}} \right\} \eqqcolon {\boldsymbol{A}\,\boldsymbol{\varepsilon}},
\end{equation}
with $\underline {\boldsymbol{E}} \coloneqq {\left[ {{E_{11}},{E_{22}},{E_{33}},2{E_{12}},2{E_{13}},2{E_{23}}} \right]^{\rm{T}}}$, where we have defined the beam strain arrays
\begin{equation*}
	{\boldsymbol{\rho }} \coloneqq \left\{ {\begin{array}{*{20}{c}}
			{{\rho _1}}\\
			{{\rho _2}}
	\end{array}} \right\},\,\,{\boldsymbol{\kappa }} \coloneqq \left\{ {\begin{array}{*{20}{c}}
			{{\kappa _{11}}}\\
			{{\kappa _{22}}}\\
			{2{\kappa _{12}}}
	\end{array}} \right\},\,\,{\boldsymbol{\delta }} \coloneqq \left\{ {\begin{array}{*{20}{c}}
			{{\delta _1}}\\
			{{\delta _2}}
	\end{array}} \right\},\,\,{\boldsymbol{\gamma }} \coloneqq \left\{ {\begin{array}{*{20}{c}}
			{{\gamma _{11}}}\\
			{{\gamma _{12}}}\\
			{{\gamma _{21}}}\\
			{{\gamma _{22}}}
	\end{array}} \right\},\,\,\mathrm{and}\,\,{\boldsymbol{\chi }}\coloneqq \left\{ {\begin{array}{*{20}{c}}
			{{\chi _{11}}}\\
			{{\chi _{22}}}\\
			{2{\chi _{12}}}
	\end{array}} \right\},
\end{equation*}
and
\begin{subequations}
	\label{beam_th_strn_comp_basis_d123}
	\begin{alignat}{2}
			\varepsilon  &\coloneqq \frac{1}{2}({\left\| {{{\boldsymbol{\varphi }}_{\!,s}}} \right\|^2} - 1)\,\,&&({\text{axial stretching strain}}),\\
			{{\rho }_\alpha } &\coloneqq {{\boldsymbol{\varphi}}_{\!,s}}\!\cdot\!{{\boldsymbol{d}}_{\alpha ,s}} - {{\boldsymbol{\varphi}}_{0,s}}\!\cdot\!{{\boldsymbol{D}}_{\alpha ,s}}\,\,&&({\text{bending strain}}),\\
			{{\delta }_\alpha } &\coloneqq {{\boldsymbol{\varphi}}_{\!,s}}\!\cdot\!{{\boldsymbol{d}}_\alpha }- {{\boldsymbol{\varphi }}_{0,s}}\!\cdot\!{{\boldsymbol{D}}_\alpha }\,\,&&({\text{transverse shear strain}}),\label{def_compat_tshear}\\
			{{\gamma }_{\alpha \beta }} &\coloneqq {{\boldsymbol{d}}_\alpha }\!\cdot\!{{\boldsymbol{d}}_{\beta ,s}} - {{\boldsymbol{D}}_\alpha }\!\cdot\!{{\boldsymbol{D}}_{\beta ,s}}\,\,&&({\text{couple shear strain}}),\label{def_b_strn_coup_sh}\\
			{{\chi }_{\alpha \beta }} &\coloneqq \frac{1}{2}({{\boldsymbol{d}}_\alpha }\!\cdot\!{{\boldsymbol{d}}_\beta } - \boldsymbol{D}_{\alpha}\!\cdot\!\boldsymbol{D}_{\beta})\,\,&&({\text{cross-section stretching and in-plane shear strains}}),\label{strn_comp_chi}\\
			{{{\kappa}} _{\alpha \beta }} &\coloneqq \frac{1}{2}\left( {{{\boldsymbol{d}}_{\alpha ,s}}\!\cdot\!{{\boldsymbol{d}}_{\beta ,s}} - {{\boldsymbol{D}}_{\alpha ,s}}\!\cdot\!{{\boldsymbol{D}}_{\beta ,s}}} \right)\,\,&&({\text{high-order bending strain}}).\label{strn_high_bend}
		\end{alignat}
	\end{subequations}
	\begin{remark} \label{rem_inp_cs_strn}
		\textit{Constant in-plane cross-sectional strains in the beam kinematics.} Two extensible directors can represent constant in-plane cross-sectional strains according to Eq.\,(\ref{strn_comp_chi}), and the work-conjugate stress resultants may not vanish as well. That is, the conventional zero-stress condition is not imposed here. The inability to represent linear in-plane (cross-sectional) strains may artificially increase the bending stiffness, which is so called \textit{Poisson locking}. In Section \ref{3fmixed_var_form}, we discuss the further enrichment of the cross-sectional strain, based on the enhanced assumed strain (EAS) method.
	\end{remark}
	\section{Mixed finite element formulation}
	\label{3fmixed}
	\subsection{Variational formulation}
	\label{3fmixed_var_form}
	We assume that the \textit{strain energy density} (defined as the strain energy per unit undeformed volume) is expressed in terms of the Green-Lagrange strain tensor $\boldsymbol{E}=\boldsymbol{E}(\boldsymbol{u})$ with ${\boldsymbol{u}}\coloneqq{\boldsymbol{x}}-{\boldsymbol{x}}_0$, as
	\begin{equation}
		\Psi = \Psi(\boldsymbol{E}).
	\end{equation}
	Then the total strain energy of the beam is obtained from
	\begin{align}\label{tot_strn_energy_beam}
		U \coloneqq \int_{{\mathcal{B}_0}} {\Psi \,{\mathrm{d}}{\mathcal{B}_0}} = {\int_0^L {\int_{\mathcal{A}} {{{\Psi}}\,{j_0}\,{\mathrm{d}}\mathcal{A}}\,{\mathrm{d}}s}},
	\end{align}
	where $j_0=({\boldsymbol{G}}_1\times{\boldsymbol{G}}_2)\cdot{\boldsymbol{G}}_3$ denotes the Jacobian of the mapping ${\boldsymbol{x}_0}(\zeta^1,\zeta^2,s):\, \mathcal{B} \to {\mathcal{B}_0}$, such that $\mathrm{d}\mathcal{B}_0={j_0}\,\mathrm{d}\mathcal{B}$ \citep{choi2021isogeometric}. In the Hu-Washizu variational principle, the total strain energy is expressed by
	\begin{equation}
		\label{hw_var_strn_e}
		{U_{{\mathrm{HW}}}} \coloneqq \int_0^L {\int_\mathcal{A} {{\Psi _{{\mathrm{HW}}}}\,{j_0}\,{\mathrm{d}}\mathcal{A}\,} {\mathrm{d}}s},
	\end{equation}
	with the strain energy density,
	\begin{equation}
		\label{hw_var_strn_ed}
		{\Psi _{{\rm{HW}}}} = {\Psi _{{\rm{HW}}}}(\boldsymbol{E}(\boldsymbol{u}),\boldsymbol{E}_\mathrm{p},\boldsymbol{S}_\mathrm{p}) = {\Psi}({{\boldsymbol{E}}_{\mathrm{p}}}) + {{\boldsymbol{S}}_{\mathrm{p}}}:\left\{ {{\boldsymbol{E}}({\boldsymbol{u}}) - {{\boldsymbol{E}}_{\mathrm{p}}}} \right\},
	\end{equation}
	where we have three independent solution fields, the displacement vector $\boldsymbol{u}$, the \textit{physical} Green-Lagrange strain tensor $\boldsymbol{E}_\mathrm{p}$, and the \textit{physical} second Piola-Kirchhoff stress tensor $\boldsymbol{S}_\mathrm{p}$. $\boldsymbol{E}{(\boldsymbol{u})}$ denotes the \textit{geometrical} (or compatible) Green-Lagrange strain tensor. The physical Green-Lagrange strain can be decomposed into the \textit{physical kinematic} part, and the \textit{enhanced} part, as \citep{simo1990class} 
	\begin{equation}
		\label{dec_phy_str_k_eh}
		{{\boldsymbol{E}}_{\rm{p}}} = \underbrace{{\boldsymbol{E}}_{\rm{p}}^{\rm{c}}}_{\textrm{kinematic}} + \underbrace {{\boldsymbol{\tilde E}}}_{{\rm{enhanced}}},
	\end{equation}
	where the additional strain part ${{\boldsymbol{\tilde E}}}$ is intended to enrich higher order cross-sectional strains. It is noted that, in contrast to the beam formulation of \citet{wackerfuss2009mixed}, which considers no cross-sectional strains (rigid cross-section) in the kinematic assumption, the kinematic part in Eq.\,(\ref{dec_phy_str_k_eh}) includes a constant in-plane strain field, see Remark \ref{rem_inp_cs_strn}. \textcolor{black}{Here, we employ the following orthogonality condition \citep{simo1990class,bischoff1997shear}
	\begin{equation}
		\label{orthog_cond_contin_SpE}
		\int_{{\mathcal{B}_0}} {{{\boldsymbol{S}}_{\rm{p}}}:{\boldsymbol{\tilde E}}\,{\rm{d}}{\mathcal{B}_0}}  = 0.
	\end{equation}
	}
	From Eq.\,(\ref{gl_strn_mata_beam_strn}), the physical kinematic part can be further decomposed into
	\begin{equation}
		\label{phy_kin_dec}
		\underline{{\boldsymbol{E}}^\mathrm{c}_{\mathrm{p}}}(\zeta^1,\zeta^2,s) = {\boldsymbol{A}}({\zeta ^1},{\zeta ^2})\,{{\boldsymbol{\varepsilon }}_{\rm{p}}}(s),
	\end{equation}
	where $\boldsymbol{\varepsilon}_\mathrm{p}(s)$ denotes the array of physical (kinematic) beam strains, i.e., ${{\boldsymbol{\varepsilon }}_{\mathrm{p}}} \coloneqq {\left[ {{\varepsilon _{\mathrm{p}}},{\boldsymbol{\rho }}_{\mathrm{p}}^{\mathrm{T}},{\boldsymbol{\kappa }}_{\mathrm{p}}^{\mathrm{T}},{\boldsymbol{\delta }}_{\mathrm{p}}^{\mathrm{T}},{\boldsymbol{\gamma }}_{\mathrm{p}}^{\mathrm{T}},{\boldsymbol{\chi }}_{\mathrm{p}}^{\mathrm{T}}} \right]^{\mathrm{T}}}$. Note that we use $\underline{(\bullet)}$ to denote the Voigt notation of a \textcolor{black}{symmetric} second order tensor. The enhanced part is decomposed into
	\begin{equation}
		\label{phy_enh_dec}
		{\underline{\boldsymbol{\tilde E}}}(\zeta^1,\zeta^2,s) = {\boldsymbol{\Gamma }}({\zeta ^1},{\zeta ^2})\,{\boldsymbol{\alpha }}(s).
	\end{equation}
	In order to alleviate Poisson locking, we enrich the linear and bilinear strains in the cross-section, using the polynomial basis functions
	\begin{equation}
		{\boldsymbol{\Gamma }}({\zeta ^1},{\zeta ^2}) \coloneqq \left[ {\begin{array}{*{20}{c}}
				{{\zeta ^1}}&{{\zeta ^2}}&{{\zeta ^1}{\zeta ^2}}&0&0&0&0&0&0\\
				0&0&0&{{\zeta ^1}}&{{\zeta ^2}}&{{\zeta ^1}{\zeta ^2}}&0&0&0\\
				0&0&0&0&0&0&0&0&0\\
				0&0&0&0&0&0&{{\zeta ^1}}&{{\zeta ^2}}&{{\zeta ^1}{\zeta ^2}}\\
				0&0&0&0&0&0&0&0&0\\
				0&0&0&0&0&0&0&0&0
		\end{array}} \right],
	\end{equation}
	with nine coefficient functions $\alpha_i(s)$, $i\in\left\{1,2,...,9\right\}$ \citep{choi2021isogeometric}. Here, we enrich only the in-plane cross-sectional strains, i.e., ${\tilde E}_{33}={\tilde E}_{13}={\tilde E}_{23}=0$. For further enrichment of higher order strains including the out-of-plane ones, one may refer to \citet{wackerfuss2009mixed} and \citet{wackerfuss2011nonlinear}. 
	\begin{remark} It was shown in \citet{wriggers1996note} that the EAS method may suffer from numerical instability under large compression, for example, for the bi-linear quadrilateral element in plane strain problems with a constant stress field. However, the hour-glass mode may not appear in the present beam formulation, since the deformation mode does not exist in the solution space of in-plane cross-sectional deformations, represented by two extensible directors.
	\end{remark} 
	\noindent Using Eqs.\,(\ref{phy_kin_dec}) and (\ref{phy_enh_dec}), Eq.\,(\ref{dec_phy_str_k_eh}) can be rewritten, as
	\begin{equation}
		\underline{{\boldsymbol{E}}_{\rm{p}}} = \left[{\begin{array}{*{20}{c}}
				{{\boldsymbol{A}}({\zeta ^1},{\zeta ^2})}&{{\boldsymbol{\Gamma }}({\zeta ^1},{\zeta ^2})}
		\end{array}}\right]\,\left\{ {\begin{array}{*{20}{c}}
				{{{\boldsymbol{\varepsilon }}_{\rm{p}}}(s)}\\
				{\,{\boldsymbol{\alpha }}(s)}
		\end{array}} \right\}.
	\end{equation}
	Then, from Eq.\,(\ref{hw_var_strn_e}), we obtain the total strain energy for beams, as
	\begin{equation}
		\label{mod_uhw_b}
		{{\tilde U}_{{\rm{HW}}}} = \int_0^L {\left\{\psi ({{\boldsymbol{\varepsilon }}_{\rm{p}}},{\boldsymbol{\alpha }}) + {\boldsymbol{\varepsilon }}({\boldsymbol{y}}) \cdot {{\boldsymbol{r}}_{\rm{p}}} - {{\boldsymbol{\varepsilon }}_{\rm{p}}} \cdot {{\boldsymbol{r}}_{\rm{p}}} - \cancel{{{\boldsymbol{\alpha}}}\cdot{{{\boldsymbol{\tilde r}}}_{\rm{p}}}}\right\}\,{\rm{d}}s},
	\end{equation}
	where we have defined the strain energy density per unit undeformed length (i.e., line energy density)
	\begin{equation}
		\psi ({{\boldsymbol{\varepsilon }}_{\rm{p}}},{\boldsymbol{\alpha }}) \coloneqq \int_\mathcal{A} {\Psi ({{\boldsymbol{E}}_\mathrm{p}}({{\boldsymbol{\varepsilon }}_{\rm{p}}},{\boldsymbol{\alpha }}))\,{j_0}\,{\rm{d}}\mathcal{A}},
	\end{equation}
	and the \textcolor{black}{array of physical stress resultants}
	\begin{equation}
		\label{comp_strs_res_rp_ar}
		{{\boldsymbol{r}}_{\rm{p}}} \coloneqq \int_\mathcal{A} {{{\boldsymbol{A}}^{\rm{T}}}{{{\boldsymbol{S}}}_{\rm{p}}}\,{j_0}\,{\rm{d}}\mathcal{A}}=\left[{\tilde n}_\mathrm{p},{\tilde m}^{1}_\mathrm{p},{\tilde m}^{2}_\mathrm{p},{\tilde h}^{\!11}_\mathrm{p},{\tilde h}^{\!22}_\mathrm{p},{\tilde h}^{\!12}_\mathrm{p},{\tilde q}^{1}_\mathrm{p},{\tilde q}^{2}_\mathrm{p},{\tilde m}^{\!11}_\mathrm{p},{\tilde m}^{\!12}_\mathrm{p},{\tilde m}^{\!21}_\mathrm{p},{\tilde m}^{\!22}_\mathrm{p},{\tilde \ell}^{11}_\mathrm{p},{\tilde \ell}^{22}_\mathrm{p},{\tilde \ell}^{12}_\mathrm{p}\right]^\mathrm{T},
	\end{equation}
	and
	\begin{equation}
		{{\boldsymbol{\tilde r}}_{\rm{p}}} \coloneqq \int_\mathcal{A} {{{\boldsymbol{\Gamma }}^{\rm{T}}}{{{\boldsymbol{S}}}_{\rm{p}}}\,{j_0}\,{\rm{d}}\mathcal{A}}.
	\end{equation}
	\textcolor{black}{The stress resultant field ${{\boldsymbol{\tilde r}}}_\mathrm{p}$ is removed from Eq.\,(\ref{mod_uhw_b}) by the orthogonality condition in Eq.\,(\ref{orthog_cond_contin_SpE}), as 
	\begin{equation}
		\int_0^L {{\boldsymbol{\alpha }} \cdot \left(\int_\mathcal{A} {{{\boldsymbol{\Gamma }}^{\rm{T}}}{{\boldsymbol{S}}_{\rm{p}}}}\,{j_0}\,{\rm{d}}\mathcal{A}\right)\,{\rm{d}}s}  = \int_0^L {{\boldsymbol{\alpha }} \cdot {{{\boldsymbol{\tilde r}}}_{\rm{p}}}\,{\rm{d}}s}  = 0.
	\end{equation}
	This condition should be satisfied, regardless of the beam's initial geometry.} Therefore, we finally have the following five independent solution fields along the center axis:
	\begin{itemize}
		\item center axis position ${\boldsymbol{\varphi}}(s)$,
		\item directors ${{\boldsymbol{d}}_\alpha}(s)$ $(\alpha\in\left\{1,2\right\})$,
		\item physical stress resultant ${\boldsymbol{r}}_\mathrm{p}(s)$,
		\item physical (kinematic) strain ${\boldsymbol{\varepsilon}}_\mathrm{p}(s)$,
		\item (physical) enhanced strain ${\boldsymbol{\alpha}}(s)$. 
	\end{itemize}
	Hereafter, for brevity, we use the notation ${\boldsymbol{y}}\coloneqq{{\left[ {{\boldsymbol{\varphi }}^\mathrm{T},{{\boldsymbol{d}}_1^\mathrm{T}},{{\boldsymbol{d}}_2^\mathrm{T}}} \right]^{\rm{T}}}}$. The first variation of the beam strains can be expressed by \citep[Section A.4.2]{choi2021isogeometric}
	\begin{equation}
		\label{var_strn_oper_b}
		\delta {\boldsymbol{\varepsilon }}({\boldsymbol{y}}) = {{\Bbb{B}}_{{\rm{total}}}}\,\delta {\boldsymbol{y}}.
	\end{equation}
	Taking the first variation of Eq.\,(\ref{mod_uhw_b}) and substituting Eq.\,(\ref{var_strn_oper_b}) yields
	\begin{equation}
		\label{vform_hw}
		G_{{\mathop{\mathrm {int}}} }^{{\mathrm{HW}}} \equiv \delta {{\tilde U}_{{\mathrm{HW}}}} = \int_0^L {\left\{ {\begin{array}{*{20}{c}}
					{\delta {\boldsymbol{y}}}\\
					{\delta {{\boldsymbol{r}}_{\rm{p}}}}\\
					{\delta {{\boldsymbol{\varepsilon }}_{\rm{p}}}}\\
					{\delta {\boldsymbol{\alpha }}}
			\end{array}} \right\} \cdot \left\{ {\begin{array}{*{20}{c}}
					{{\Bbb{B}}_{{\rm{total}}}^{\rm{T}}{{\boldsymbol{r}}_{\rm{p}}}}\\
					{{\boldsymbol{\varepsilon }}({\boldsymbol{y}}) - {{\boldsymbol{\varepsilon }}_{\rm{p}}}}\\
					{{\partial _{{{\boldsymbol{\varepsilon }}_{\rm{p}}}}}\psi ({{\boldsymbol{\varepsilon }}_{\rm{p}}},{\boldsymbol{\alpha }}) - {{\boldsymbol{r}}_{\rm{p}}}}\\
					{{\partial _{\boldsymbol{\alpha }}}\psi ({{\boldsymbol{\varepsilon }}_{\rm{p}}},{\boldsymbol{\alpha }})}
			\end{array}} \right\}{\rm{d}}s},
	\end{equation}
	where we have defined \textcolor{black}{the stress resultants, calculated from the strain energy density function, as}
	\begin{equation}
		{\partial _{{{\boldsymbol{\varepsilon }}_{\rm{p}}}}}\psi({{\boldsymbol{\varepsilon }}_{\rm{p}}},{\boldsymbol{\alpha }})  \coloneqq \int_\mathcal{A} {{{\boldsymbol{A}}^{\rm{T}}}\underline{\partial _{{{\boldsymbol{E}}_{\rm{p}}}}\Psi} \,{j_0}\,{\rm{d}}\mathcal{A}},
	\end{equation}
	\begin{equation}
		{\partial _{\boldsymbol{\alpha }}}\psi({{\boldsymbol{\varepsilon }}_{\rm{p}}},{\boldsymbol{\alpha }})  \coloneqq \int_\mathcal{A} {{{\boldsymbol{\Gamma }}^{\rm{T}}}\underline{{\partial _{{{\boldsymbol{E}}_{\rm{p}}}}}\Psi} \,{j_0}\,{\rm{d}}\mathcal{A}},
	\end{equation}
	with ${\partial _{{{\boldsymbol{E}}_{\rm{p}}}}}\Psi  \coloneqq \partial \Psi /\partial {{\boldsymbol{E}}_{\rm{p}}}$.
	The external virtual work is given by \citep{choi2021isogeometric}
	\begin{equation}
		\label{var_ex_vwork}
		{G_{{\rm{ext}}}}\left({\delta {\boldsymbol{y}}} \right) = \int_0^L {\delta {{\boldsymbol{y}}^{\rm{T}}}{\boldsymbol{\bar R}}\,{\mathrm{d}}s} + {\left[ {\delta {{\boldsymbol{y}}^{\rm{T}}}{{{\boldsymbol{\bar R}}}_0}} \right]_{s\in{\Gamma _{\mathrm{N}}}}},
	\end{equation}
	with ${\boldsymbol{\bar R}} \coloneqq {\left[ {{{{\boldsymbol{\bar n}}}^{\rm{T}}},{{{\boldsymbol{\bar {\tilde m}}}}{{_{}}^{1}}}^\mathrm{T},{{\boldsymbol{\bar {\tilde m}}}}{{_{}}^{2}}^\mathrm{T}} \right]^{\rm{T}}}$ and ${{\boldsymbol{\bar R}}_0} \coloneqq {\left[ {{{{\boldsymbol{\bar n}}}_0}^{\rm{T}},{\boldsymbol{\bar {\tilde m}}}{{_0^1}^{\rm{T}}},{\boldsymbol{\bar {\tilde m}}}{{_0^2}^{\rm{T}}}} \right]^{\rm{T}}}$, where ${{\boldsymbol{\bar n}}}$ and ${{\boldsymbol{\bar {\tilde m}}}}{{_{}}^{\alpha}}$ ($\alpha\in\left\{1,2\right\}$) denote the external stress resultants from the external load on the lateral surface of the beam. ${{\boldsymbol{\bar n}}}_0$ and ${{\boldsymbol{\bar {\tilde m}}}}{{_{}}_0^{\alpha}}$ denote the prescribed stress resultants at the boundary $s\in\Gamma_\mathrm{N}$. Here, $\Gamma_\mathrm{N}\,\textcolor{black}{\subset}\,{{\emptyset}\cup\left\{0,L\right\}}$ denotes the Neumann boundary, where $\emptyset$ denotes the empty set. Then, we obtain the following variational equation: Find ${\boldsymbol{y}} \in \mathcal{V}$, ${{\boldsymbol{r}}_{\rm{p}}},{{\boldsymbol{\varepsilon }}_{\rm{p}}} \in {\mathcal{V}}_\mathrm{p}$, and ${\boldsymbol{\alpha}}\in\mathcal{V}_\mathrm{a}$ such that
	\begin{equation}
		\label{vareq_hw_mix_space}
		G_{{\rm{int}}}^{{\rm{HW}}}\left({\boldsymbol{y}},{{\boldsymbol{r}}_{\rm{p}}},{{\boldsymbol{\varepsilon }}_{\rm{p}}},\textcolor{black}{\boldsymbol{\alpha}},\delta {\boldsymbol{y}},\delta {{\boldsymbol{r}}_{\rm{p}}},\delta {{\boldsymbol{\varepsilon }}_{\rm{p}}},\textcolor{black}{\delta\boldsymbol{\alpha}}\right) = {G_{{\rm{ext}}}}\left( {\delta {\boldsymbol{y}}} \right),\,\,\forall\,\delta {\boldsymbol{y}} \in {\bar{\mathcal{V}}},\,\,\delta {{\boldsymbol{r}}_{\mathrm{p}}},\delta {{\boldsymbol{\varepsilon}}_{\rm{p}}} \in {\mathcal{V}}_\mathrm{p},\,\,\mathrm{and}\,\,\delta {{\boldsymbol{\alpha}}} \in {\mathcal{V}}_\mathrm{a},
	\end{equation}
	with the solution space given by
	\begin{equation}
		\mathcal{V} \coloneqq \left\{ {\left. {\boldsymbol{y} \in {{\left[ {{H^1}(0,L)} \right]}^d}} \right|{\boldsymbol{\varphi }} = {{{\boldsymbol{\bar \varphi }}}_0},\,\,{{\boldsymbol{d}}_1} = {{{\boldsymbol{d}}}_{10}},\,\,{\rm{and}}\,\,{{\boldsymbol{d}}_2} = {{{\boldsymbol{d}}}_{20}}\,\,{\rm{on}}\,\,s \in {\Gamma _{\rm{D}}}\,} \right\},
	\end{equation}
	where $\Gamma_\mathrm{D}$ denotes the Dirichlet boundary, such that ${\Gamma _{\rm{D}}} \cup {\Gamma _{\rm{N}}} = \left\{ {0,L} \right\}$ and ${\Gamma _{\rm{D}}} \cap {\Gamma _{\rm{N}}} = \emptyset$. $d=9$ is the number of independent components in $\boldsymbol{y}$.	${\bar{\boldsymbol{\varphi}}}_0$ and ${\boldsymbol{d}}_{\alpha0}$ ($\alpha\in\left\{1,2\right\}$) denote the prescribed center axis position and director vectors, respectively, and the variational space is given by
	\begin{equation}
		{\bar {\mathcal{V}}} \coloneqq \left\{ {\left. {\delta {\boldsymbol{y}} \in {{\left[ {{H^1}(0,L)} \right]}^d}} \right|\delta {\boldsymbol{\varphi }} = \delta {{\boldsymbol{d}}_1} = \delta {{\boldsymbol{d}}_2} = {\boldsymbol{0}}\,\,{\rm{on}}\,\,s \in {\Gamma _{\rm{D}}}\,} \right\}.
	\end{equation}
	It is chosen such that the physical stress resultants and physical strains belong to ${\mathcal{V}}_\mathrm{p}\coloneqq\left[{L^2(0,L)}\right]^{d_\mathrm{p}}$, and the enhanced strains belong to ${\mathcal{V}}_\mathrm{a}\coloneqq\left[{L^2(0,L)}\right]^{d_{\mathrm{a}}}$, where the number of independent functions are $d_\mathrm{p}=15$, and $d_{\mathrm{a}} = 9$. This means that no inter-element continuity is required for those additional fields, so that we can locally (element-wisely) condense out those corresponding nodal degrees-of-freedom after the finite element approximation. Those condensed coefficients are kept as \textit{internal variables}, which requires additional computer storage. It is noted that we simply follow the conventional treatment to consider the physical stress and strain fields in the open domain $(0,L)\ni{s}$, without any extra boundary conditions, see, e.g., \citet{wackerfuss2009mixed}. A further investigation on this aspect remains future work. 
	\textcolor{black}{
	\begin{remark}
		\label{rem_strongf_eq}
		From the stationarity condition ${G^\mathrm{HW}_\mathrm{int}}-{G^\mathrm{HW}_\mathrm{ext}}=0$, we obtain the following Euler-Lagrange equations\footnote{We refer this procedure to \citet[Section 9.3]{santos2010hybrid}.}, which states the strong form equations of the given problem: Find ${\boldsymbol{y}} \in {{\Bbb{R}}^d}$, satisfying
		\begin{equation}
			\label{dbdc_str_feq}
			{\boldsymbol{\varphi }} = {{\boldsymbol{\bar \varphi }}_0},\,\,{{\boldsymbol{d}}_1} = {{\boldsymbol{d}}_{10}},\,\,\text{and}\,\,{{\boldsymbol{d}}_2} = {{\boldsymbol{d}}_{20}}\,\,\mathrm{at}\,\,s\in\Gamma_\mathrm{D},		
		\end{equation}
		and ${{\boldsymbol{r}}_{\rm{p}}} \in {{\Bbb{R}}^{{d_{\rm{p}}}}}$, ${{\boldsymbol{\varepsilon }}_{\rm{p}}} \in {{\Bbb{R}}^{{d_{\rm{p}}}}}$, and ${\boldsymbol{\alpha }} \in {{\Bbb{R}}^{{d_{\rm{a}}}}}$, such that
		\begin{subequations}
			\label{elag_eq}
			\begin{alignat}{3}
				{{\boldsymbol{n}}_{{\rm{p}},s}} + {\boldsymbol{\bar n}} &= {\boldsymbol{0}}&&{}\quad\text{(linear momentum balance),}\label{elag_lmnt_bal}\\
				{\boldsymbol{\tilde m}}_{{\rm{p}},s}^\alpha  - {\boldsymbol{l}}_{\rm{p}}^\alpha  + {{\boldsymbol{\bar {\tilde m}}}^\alpha } &= {\boldsymbol{0}},\,\,\alpha \in \left\{ {1,2} \right\}&&{}\quad\text{(director momentum balance),}\label{elag_dir_mnt_bal}\\
				{{\boldsymbol{\varepsilon }}({\boldsymbol{y}}) - {{\boldsymbol{\varepsilon }}_{\rm{p}}}}&=\boldsymbol{0}&&{}\quad\text{(compatibility),}\label{elag_compat}\\
				{{\partial _{{{\boldsymbol{\varepsilon }}_{\rm{p}}}}}\psi ({{\boldsymbol{\varepsilon }}_{\rm{p}}},{\boldsymbol{\alpha }}) - {{\boldsymbol{r}}_{\rm{p}}}}&=\boldsymbol{0}&&{}\quad\text{(constitutive equation),}\label{constituti_eq}\\
				{{\partial _{\boldsymbol{\alpha }}}\psi ({{\boldsymbol{\varepsilon }}_{\rm{p}}},{\boldsymbol{\alpha }})}&=\boldsymbol{0}&&{}\quad\text{(zero higher order stress resultants),}\label{zw_ehn}
			\end{alignat}
		\end{subequations}
		at $s\in\left(0,L\right)$, and the Neumann boundary conditions
		\begin{subequations}
			\label{elag_eq_neum_bdc}
			\begin{alignat}{3}
				{{\boldsymbol{n}}_{\rm{p}}} - {{{\boldsymbol{\bar n}}}_0} &= {\boldsymbol{0}},\label{nmbdc_f}\\
				{\boldsymbol{\tilde m}}_{\rm{p}}^\alpha  - {\boldsymbol{\bar {\tilde m}}}_0^\alpha  &= {\boldsymbol{0}},\,\,\alpha  \in \left\{ {1,2} \right\},\label{nmbdc_m}
			\end{alignat}
		\end{subequations}
		at $s\in\Gamma_\mathrm{N}$, where those physical stress resultant components in the array $\boldsymbol{r}_\mathrm{p}$ of Eq.\,(\ref{comp_strs_res_rp_ar}) are related to
		\begin{subequations}
			\label{elag_eq_def_nml}
			\begin{alignat}{3}
				{{\boldsymbol{n}}_{\rm{p}}} &\equiv {\tilde n_{\rm{p}}}{{\boldsymbol{\varphi }}_{,s}} + \tilde m_{\rm{p}}^\alpha {{\boldsymbol{d}}_{\alpha,s}} + \tilde q_{\rm{p}}^\alpha {{\boldsymbol{d}}_\alpha}&&{}\quad\text{(stress resultant),}\label{elag_strs_res_phy}\\
				{\boldsymbol{\tilde m}}_{\rm{p}}^\alpha  &\equiv \tilde m_{\rm{p}}^\alpha {{\boldsymbol{\varphi }}_{,s}} + \tilde h_{\rm{p}}^{\alpha \beta }{{\boldsymbol{d}}_{\beta ,s}} + \tilde m_{\rm{p}}^{\beta \alpha}{{\boldsymbol{d}}_\beta },\,\,\alpha \in \left\{ {1,2} \right\}&&{}\quad\text{(director stress couple),}\label{elag_dir_strs_coup_phy}\\
				{\boldsymbol{l}}_{\rm{p}}^\alpha  &\equiv \tilde \ell _{\rm{p}}^{\alpha \beta }{{\boldsymbol{d}}_\beta } + \tilde q_{\rm{p}}^\alpha {{\boldsymbol{\varphi }}_{,s}} + \tilde m_{\rm{p}}^{\alpha \beta }{{\boldsymbol{d}}_{\beta ,s}},\,\,\alpha \in \left\{ {1,2} \right\}&&{}\quad\text{(through-the-thickness stress resultant).}\label{elag_strs_res_phy}
			\end{alignat}
		\end{subequations}		
	\end{remark}
	}
	\subsection{Linearization}
	The internal virtual work $G^\mathrm{HW}_\mathrm{int}$ of Eq.\,(\ref{vform_hw}) contains geometrical or material nonlinearities. Therefore, we employ a Newton-Raphson method to solve the variational equation of Eq.\,(\ref{vareq_hw_mix_space}). An external load is incrementally applied, and using the equilibrium configuration at the previous $n$th load step as an initial guess, the solution at the next, $(n+1)$th load step is found. The iterative scheme to find the solution is stated as follows: For a given solution ${}^{n + 1}{\boldsymbol{\eta}}^{(i - 1)}\coloneqq\left\{ {{}^{n + 1}{{\boldsymbol{y}}^{(i - 1)}},{}^{n + 1}{\boldsymbol{r}}_{\rm{p}}^{(i - 1)},{}^{n + 1}{\boldsymbol{\varepsilon }}_{\rm{p}}^{(i - 1)},{}^{n + 1}{\boldsymbol{\alpha}}^{(i - 1)}} \right\} \in \mathcal{V} \times {\mathcal{V}_\mathrm{p}} \times {\mathcal{V}_\mathrm{p}}\times {\mathcal{V}_{\mathrm{a}}}$ at the $(i-1)$th iteration in the $(n+1)$th load step, find the solution increment $\Delta{\boldsymbol{\eta}}\coloneqq\left\{ {\Delta {\boldsymbol{y}},\Delta {{\boldsymbol{r}}_{\rm{p}}},\Delta {{\boldsymbol{\varepsilon }}_{\rm{p}}},\Delta {{\boldsymbol{\alpha }}}} \right\} \in \bar {\mathcal{V}} \times {\mathcal{V}_\mathrm{p}} \times {\mathcal{V}_\mathrm{p}}\times{{\mathcal{V}_{\mathrm{a}}}}$ such that
	\begin{align}
		\label{lin_var_eq_g_hw}
		&\Delta G_{{\rm{int}}}^{{\rm{HW}}}\left({}^{n + 1}{{\boldsymbol{\eta}}^{(i - 1)}},\delta {\boldsymbol{\eta}},\Delta {\boldsymbol{\eta}}\right)= {G_{{\rm{ext}}}}\left( {\delta {\boldsymbol{y}}} \right) - G_{{\rm{int}}}^{{\rm{HW}}}\left( {{}^{n + 1}{{\boldsymbol{\eta}}^{(i - 1)}}},\delta {\boldsymbol{\eta}} \right),\nonumber\\
		&\textcolor{black}{\forall\,\delta {\boldsymbol{\eta}}\coloneqq\left\{{\delta {\boldsymbol{y}},\delta {{\boldsymbol{r}}_{\rm{p}}},\delta {{\boldsymbol{\varepsilon }}_{\rm{p}}},\delta {{\boldsymbol{\alpha }}}}\right\}\in \bar {\mathcal{V}} \times {\mathcal{V}_\mathrm{p}} \times {\mathcal{V}_\mathrm{p}}\times{{\mathcal{V}_{\mathrm{a}}}}}.
	\end{align}
	The solution is then updated by
	\begin{align}
		\left. \arraycolsep=1.4pt\renewcommand{\arraystretch}{1.15}\begin{array}{rclccrl}
			{}^{n + 1}{{\boldsymbol{y}}^{(i)}} &=& {}^{n + 1}{{\boldsymbol{y}}^{(i - 1)}} &+& \Delta {\boldsymbol{y}},\,\,{}^{n + 1}{{\boldsymbol{y}}^{(0)}} &=& {}^n{\boldsymbol{y}},\\
			{}^{n + 1}{\boldsymbol{r}}_{\rm{p}}^{(i)} &=& {}^{n + 1}{\boldsymbol{r}}_{\rm{p}}^{(i - 1)} &+& \Delta {{\boldsymbol{r}}_{\rm{p}}},\,\,{}^{n + 1}{\boldsymbol{r}}_{\rm{p}}^{(0)} &=& {}^n{{\boldsymbol{r}}_{\rm{p}}},\\
			{}^{n + 1}{\boldsymbol{\varepsilon }}_{\rm{p}}^{(i)} &=& {}^{n + 1}{\boldsymbol{\varepsilon }}_{\rm{p}}^{(i - 1)} &+& \Delta {{\boldsymbol{\varepsilon }}_{\rm{p}}},\,\,{}^{n + 1}{\boldsymbol{\varepsilon }}_{\rm{p}}^{(0)} &=& {}^n{{\boldsymbol{\varepsilon }}_{\rm{p}}},\\
			{}^{n + 1}{\boldsymbol{\alpha}}^{(i)} &=& {}^{n + 1}{\boldsymbol{\alpha}}^{(i - 1)} &+& \Delta {{\boldsymbol{\alpha}}},\,\,{}^{n + 1}{\boldsymbol{\alpha}}^{(0)} &=& {}^n{{\boldsymbol{\alpha}}},
		\end{array} \right\}
	\end{align}
	where the initial guess is given by the converged solution in the previous ($n$th) load step. We obtain the increment of the internal virtual work $\Delta G_{{\rm{int}}}^{{\rm{HW}}}$ by taking the directional derivative of Eq.\,(\ref{vform_hw}), as
	\begin{align}
		\Delta {G}_{{\mathop{\mathrm{int}}} }^{{\mathrm{HW}}} = \int_0^L {{{\left\{ {\begin{array}{*{20}{c}}
							{\delta {\boldsymbol{y}}}\\
							{\delta {{\boldsymbol{r}}_{\rm{p}}}}\\
							{\delta {{\boldsymbol{\varepsilon }}_{\rm{p}}}}\\
							{\delta {\boldsymbol{\alpha }}}
					\end{array}} \right\}}^{\rm{T}}}\left[ {\begin{array}{*{20}{c}}
					{{{\boldsymbol{Y}}^{\rm{T}}}{{\boldsymbol{k}}_\mathrm{G}}{\boldsymbol{Y}}}&{{\Bbb{B}}_{{\rm{total}}}^{\rm{T}}}&{{{\bf{0}}_{d \times {d_\mathrm{p}}}}}&{{{\bf{0}}_{d \times {d_{\rm{a}}}}}}\\
					{}&{{{\bf{0}}_{{d_\mathrm{p}} \times {d_\mathrm{p}}}}}&{ - {{\bf{1}}_{{d_\mathrm{p}} \times {d_\mathrm{p}}}}}&{{{\bf{0}}_{{d_\mathrm{p}} \times {d_{\rm{a}}}}}}\\
					{}&{}&{{{\Bbb{C}}^{\varepsilon\varepsilon}_{{\rm{p}}}}}&{{\Bbb{C}}}^{\textrm{a}{\varepsilon}\,\mathrm{T}}_{\rm{p}}\\
					{{\rm{sym}}{\rm{.}}}&{}&{}&{{{\Bbb{C}}_{\rm{p}}^{{\rm{aa}}}}}
			\end{array}} \right]\left\{ {\begin{array}{*{20}{c}}
					{\Delta {\boldsymbol{y}}}\\
					{\Delta {{\boldsymbol{r}}_{\rm{p}}}}\\
					{\Delta {{\boldsymbol{\varepsilon}}}_{\rm{p}}}\\
					{\Delta {\boldsymbol{\alpha }}}
			\end{array}} \right\}{\rm{d}}s},
	\end{align}
	with the operator \citep[Section A.4.4]{choi2021isogeometric}
	%
	\begin{equation}\label{def_matY_oper}
		{\boldsymbol{Y}} \coloneqq {\left[ {\begin{array}{*{20}{c}}
					{{{(\bullet)}_{,s}}{{\boldsymbol{1}}_{3\times3}}}&{{{\bf{0}}_{3 \times 3}}}&{{{\bf{0}}_{3 \times 3}}}\\
					{{{\bf{0}}_{3 \times 3}}}&{{{(\bullet)}_{,s}}{{\boldsymbol{1}}_{3\times3}}}&{{{\bf{0}}_{3 \times 3}}}\\
					{{{\bf{0}}_{3 \times 3}}}&{{{\bf{0}}_{3 \times 3}}}&{{{(\bullet)}_{,s}}{{\boldsymbol{1}}_{3\times3}}}\\
					{{{\bf{0}}_{3 \times 3}}}&{{\boldsymbol{1}}_{3\times3}}&{{{\bf{0}}_{3 \times 3}}}\\
					{{{\bf{0}}_{3 \times 3}}}&{{{\bf{0}}_{3 \times 3}}}&{{\boldsymbol{1}}_{3\times3}}
			\end{array}} \right]_{15 \times 9}},
	\end{equation}
	\textcolor{black}{where ${\bf{0}}_{m\times{n}}$ and ${\bf{1}}_{n\times{n}}$ denote the $m\times{n}$ null matrix, and $n\times{n}$ identity matrix, respectively.} Here we have defined the constitutive matrices
	\begin{equation}
		{\Bbb{C}_{\rm{p}}^{\varepsilon\varepsilon}} \coloneqq \int_\mathcal{A} {{{\boldsymbol{A}}^{\rm{T}}}{{\underline{\underline{{\boldsymbol{\mathcal{C}}}_{\rm{p}}}}}}{\boldsymbol{A}}\,{j_0}}\,{\rm{d}}\mathcal{A},
	\end{equation}
	\begin{equation}
		{{\Bbb{C}}_\mathrm{p}^{{\rm{aa}}}} \coloneqq \int_\mathcal{A} {{\boldsymbol{\Gamma}}^{\rm{T}}}\underline{\underline{\boldsymbol{\mathcal{C}}_{\rm{p}}}}\boldsymbol{\Gamma}\,{j_0}\,{\rm{d}}\mathcal{A},
	\end{equation}
	\begin{equation}
		{{\Bbb{C}}_\mathrm{p}^{{\rm{a}}\varepsilon}} \coloneqq \int_\mathcal{A} {{\boldsymbol{\Gamma}}^{\rm{T}}}\underline{\underline{\boldsymbol{\mathcal{C}}_{\rm{p}}}}\boldsymbol{A}\,{j_0}\,{\rm{d}}\mathcal{A},
	\end{equation}
	where
	\begin{equation}
		{{\boldsymbol{\mathcal{C}}}_{\rm{p}}} \coloneqq \dfrac{{{\partial ^2}\Psi }}{{\partial {{\boldsymbol{E}}_{\mathrm{p}}}\,\partial {{\boldsymbol{E}}_{\mathrm{p}}}}},
	\end{equation}
and ${{\underline{\underline{(\bullet)}}}}$ denotes the Voigt notation of a fourth-order tensor \textcolor{black}{with both major and minor symmetries}.
\section{An isogeometric finite element discretization}
\label{ig_fe_disc_sec}
\subsection{NURBS curve: An exact representation of the initial geometry}
The initial geometry of the beam's center axis can be represented by a spline curve. In the framework of IGA, we use the same spline basis functions utilized in the CAD model. Here we summarize the construction of NURBS curves. More detailed explanation on the properties of NURBS, and algorithms for the knot insertion and degree elevation can be found in \citet{piegl1996nurbs}. Those two operations are not commutative. Degree elevation followed by knot insertion maintains the maximum continuity of the original curve, which is so called $k$-refinement. The other way around, increasing the degree of each curve segment after the knot insertion, corresponds to the classical $p$-refinement. For a specific example of these mesh refinement processes, see \citet{hughes2005isogeometric}. For a given patch of a NURBS curve, we have the knot vector ${\tilde \varXi}={\left\{{\xi_1},{\xi_2},...,{\xi_{{n_{\mathrm{cp}}}+p+1}}\right\}}$, where ${\xi_i}\in{\Bbb R}$ is the $i$th knot, $p$ is the degree of basis functions, and ${n_{\mathrm{cp}}}$ is the total number of basis functions (or control points). B-spline basis functions are recursively defined \citep{piegl1996nurbs}. For $p=0$, they are defined by
\begin{equation} \label{Bspline_basis_0}
	B_I^0(\xi ) = 
	\begin{cases} 
		1&{{\rm{if~~ }}{\xi _I} \le \xi  < {\xi _{I + 1}}},\\
		0&{{\text{otherwise,  }}}
	\end{cases}
\end{equation}
and for $p=1,2,3,...,$ they are defined by
\begin{equation} \label{Bspline_basis_p}
	B_I^p(\xi ) = \frac{{\xi  - {\xi _I}}}{{{\xi _{I + p}} - {\xi _I}}}B_I^{p - 1}(\xi ) + \frac{{{\xi _{I + p + 1}} - \xi }}{{{\xi _{I + p + 1}} - {\xi _{I + 1}}}}B_{I + 1}^{p - 1}(\xi ),
\end{equation}
where $\xi\in\varXi\subset{\Bbb R}$ denotes the parametric coordinate, and $\varXi\coloneqq\left[\xi_{1},\xi_{{n_{\mathrm{cp}}}+p+1}\right]$ represents the parametric domain. We employ NURBS to exactly represent the initial geometries from a conic section like circle and ellipse. From the B-spline basis functions, the NURBS basis functions are defined by 
\begin{equation}\label{nurbs_basis_1d_def}
	{N^p_I}(\xi ) = \frac{{B_I^p(\xi )\,{w_I}}}{{\sum\limits_{J = 1}^{n_{\mathrm{cp}}} {B_J^p(\xi )\,{w_J}} }},
\end{equation}
where ${w_I}$ denotes the given weight of the $I$th control point. If the weights are equal, the NURBS becomes a B-spline. The geometry of the beam's initial center axis can be represented by a NURBS curve, as
\begin{equation}\label{beam_curve_pos_nurbs}
	{\boldsymbol{X}}(\xi ) = \sum\limits_{I = 1}^{n_{\mathrm{cp}}} {{N^p_I}(\xi )\,{{\boldsymbol{X}}_{\!I}}},
\end{equation}
where $\boldsymbol{X}_{\!I}$ are the control point positions. The arc-length parameter along the initial center axis can be expressed by the mapping $s(\xi):{\varXi}\to{\left[0,L\right]}$, defined by
\begin{equation}\label{beam_curve_pos_nurbs_alen_map}
	s(\xi )\coloneqq \int_{{\xi _1}}^{\eta  = \xi } {\left\| {{{\boldsymbol{X}}_{\!,\eta }}(\eta )} \right\|{\mathrm{d}}\eta }.
\end{equation}
Then the Jacobian of the mapping is derived as
\begin{align}\label{beam_curve_pos_nurbs_alen_map_jcb}
	\tilde j\coloneqq \frac{{{\mathrm{d}}s}}{{{\mathrm{d}}\xi }}= \left\| {{{\boldsymbol{X}}_{\!,\xi }}(\xi )} \right\|.
\end{align}
In the discretization of the variational form, we often use the notation ${N^p_{I,s}}$ for brevity, which is defined by
\begin{align}\label{beam_curve_pos_nurbs_alen_map_jcb}
	{N^p_{I,s}} \coloneqq {N^p_{I,\xi }}\frac{{{\mathrm{d}}\xi }}{{{\mathrm{d}}s}} = \frac{1}{{\tilde j}}{N^p_{I,\xi }},
\end{align}
where ${N^p_{I,\xi }}$ denotes the differentiation of the basis function ${N^p_{I}(\xi)}$ with respect to $\xi$.
\subsection{Discretization of the variational form}
We may choose different degrees of basis functions for the displacements of the center axis, and directors, which are denoted by $p$ and $p_\mathrm{d}$, respectively. In this paper, we propose to use $p_\mathrm{d}=p-1$ for the field consistency in the finite element approximation of the transverse shear strain of Eq.\,(\ref{def_compat_tshear}). In the entire domain of a curve patch, we have the approximated current center axis position
\begin{equation}
	\label{cdisp_entire_patch0}
	{\boldsymbol{\varphi}}^{\,h}(s(\xi)) = \sum\limits_{I = 1}^{n_\mathrm{cp}} {N_I^{{p}}(\xi )\,{{{\boldsymbol{\varphi}}}_{{I}}}},\,\,\xi\in{\varXi},
\end{equation}
and total director displacement ${{\boldsymbol{\bar d}}_\alpha } \coloneqq {\boldsymbol{d}}_\alpha - {{\boldsymbol{D}}_\alpha }$ ($\alpha\in\left\{1,2\right\}$)
\begin{equation}
	\label{tot_dirv_disp_entire_patch0}
	{\boldsymbol{\,\bar d}}_\alpha ^{\,h}(s(\xi)) = \sum\limits_{J = 1}^{n_\mathrm{cp}^\mathrm{d}} {N_J^{{p_\mathrm{d}}}(\xi )\,{{{\bf{\bar d}}}_{\alpha{J}}}},\,\,\alpha\in\left\{1,2\right\},\,\,\xi\in{\varXi},
\end{equation}
where ${n_\mathrm{cp}^\mathrm{d}}$ denotes the total number of control coefficients for the director displacement field in the patch. In IGA, we define an element by a half-open interval between two distinct knots, i.e., a non-zero knot span. Let ${\varXi}_e\coloneqq{\left[ {\xi _1^e,\xi _2^e} \right)}\ni{\xi}$ denote the $e$th element, such that $\varXi={\varXi_1}\cup{\varXi_2}\cup\cdots\cup{\varXi}_{n_\mathrm{el}}$, where $n_\mathrm{el}$ denotes the total number of elements. For the last element of an open curve (i.e., $e=n_\mathrm{el}$), we consider a closed interval ${\varXi}_{n_\mathrm{el}}\coloneqq{\left[ {\xi _1^{n_\mathrm{el}},\xi _2^{n_\mathrm{el}}} \right]}$ to include the end point. \textcolor{black}{We use the same mesh (non-zero knot spans) in the parametric domain, for the displacement fields of the axis, and directors. Therefore, we have different control net for those two fields, with $n^\mathrm{d}_\mathrm{cp} = n_\mathrm{cp} - 1$.} We have $n_e=p+1$, and $n^\mathrm{d}_e=p_\mathrm{d}+1$ local support basis functions in every element, for the axis and director displacement fields, respectively. Then, we can rewrite Eqs.\,(\ref{cdisp_entire_patch0}) and (\ref{tot_dirv_disp_entire_patch0}) in each element, as
\begin{equation}
	\label{approx_caxis_pos}
	{{\boldsymbol{\varphi }}^{h}}(s(\xi )) = \sum\limits_{I = 1}^{{n_e}} {N_I^p(\xi )\,{{\boldsymbol{\varphi }}^e_I}},\,\,\xi\in{\varXi_e},
\end{equation}
and
\begin{equation}
	\label{tot_dirv_disp}
	{\boldsymbol{\,\bar d}}_\alpha ^{\,h}(s(\xi)) = \sum\limits_{J = 1}^{n_e^\mathrm{d}} {N_J^{{p_\mathrm{d}}}(\xi )\,{{{\bf{\bar d}}}^e_{\alpha{J}}}},\,\,\xi\in{\varXi_e}.
\end{equation}
Combining Eqs.\,(\ref{approx_caxis_pos}) and (\ref{tot_dirv_disp}), we obtain
\begin{equation}
	{{\boldsymbol{y}}^h} = \left\{ {\begin{array}{*{20}{c}}
			{{{\boldsymbol{\varphi }}^{h}}}\\
			{\begin{array}{*{20}{c}}
					{{\boldsymbol{\bar d}}_1^{\,h}}\\
					{{\boldsymbol{\bar d}}_2^{\,h}}
			\end{array}}
	\end{array}} \right\} = {\Bbb{N}_e}(\xi )\,{{{\bf{y}}^e}},\,\,\xi\in\varXi_e,
\end{equation}
with 
\begin{equation}
	{\Bbb{N}_e}\coloneqq\left[ {\begin{array}{*{20}{c}}
			{N_1^p{{\bf{1}}_{3 \times 3}}}& \cdots &{N_{{n_e}}^p{{\bf{1}}_{3 \times 3}}}&\vline& {}&{{{\bf{0}}_{3 \times 6n_e^\mathrm{d}}}}&{}\\
			{}&{{{\bf{0}}_{6 \times 3{n_e}}}}&{}&\vline& {N_1^{{p_d}}{{\bf{1}}_{6 \times 6}}}& \cdots &{N_{n_e^\mathrm{d}}^{{p_\mathrm{d}}}{{\bf{1}}_{6 \times 6}}}
	\end{array}} \right],
\end{equation}
and ${{{\bf{y}}^e}}\coloneqq\left[{{{\boldsymbol{\varphi}}_1^{e\,\mathrm{T}}}},\cdots,{{{\boldsymbol{\varphi}}_{{n_e}}^{e\,\mathrm{T}}}},{{\bf{\bar d}}_1^{e\,\mathrm{T}}},\cdots,{{\bf{\bar d}}_{n_e^\mathrm{d}}^{e\,\mathrm{T}}}
\right]^\mathrm{T}$, where \textcolor{black}{${{\boldsymbol{\varphi }}^e_I} \in {{\Bbb{R}}^3}$, and ${{\bf{\bar d}}^e_{J}} \coloneqq \left[{\bf{\bar d}}_{1J}^{e\,\mathrm{T}},{\bf{\bar d}}_{2J}^{e\,\mathrm{T}}\right]^\mathrm{T}\in {{\Bbb{R}}^6}$} denote the control coefficient vectors for the current axis position, and the total displacement of directors, respectively. \textcolor{black}{Note that, due to $n_e\ne{n_e^\mathrm{d}}$, we need to separate the arrangement of control coefficients for the center axis position and director displacement parts.} 
\textcolor{black}{
	\begin{remark} \textit{An exact construction of the initial director field}. It is noted that, in Eq.\,(\ref{tot_dirv_disp}), we approximate the \textit{total displacement} of the directors, i.e., the difference vector between the current and initial directors. This is to avoid the approximation error in representing the initial geometry, since the NURBS-based approximation may not preserve the orthonormality of the initial directors. In Section \ref{num_ex_45deg_arc_cant}, we consider a numerical example of an initially curved beam, where we employ the smallest rotation method (see \citet{choi2019isogeometric} and references therein) to construct the initial orthonormal director field. 	
	\end{remark}
	\begin{remark} 
		\label{rem_red_pd_curv}
		\textit{Minimum required degree $p_\mathrm{d}$ to exactly represent rigid body rotations}. The degree elevation of a NURBS curve does not alter the curve either geometrically or parametrically \citep{piegl1996nurbs}. If the center axis is initially modeled by a NURBS curve of the minimum required degree $p_\mathrm{g}$ to represent the \textit{exact} geometry, and $p,p_\mathrm{d}\ge{p_\mathrm{g}}$, the parameterization of the displacement fields of the center axis, and directors, with basis functions of degrees $p$ and $p_\mathrm{d}$, respectively, is \textit{consistent} with that of the initial geometry. In order to represent the rigid body rotations exactly for initially curved beams, it is required to use $p_\mathrm{d}\ge{p_\mathrm{g}}$, see Section \ref{cant45_case1_subsec_lin} for a relevant example.
	\end{remark}
}
\noindent In the approximation of physical stress resultants, and strains, we allow for inter-element discontinuity for an element-wise static condensation process, and use Lagrange interpolation functions. We first define a mapping from a parametric domain $\left[-1,1\right]\ni{\bar \xi}$, where the Lagrange functions are defined, to the closed interval ${\bar \varXi}_e\coloneqq{\left[ {\xi _1^e,\xi _2^e} \right]}\ni{\xi}$, as
\begin{equation}	
	\bar \xi  = 1 - 2\left( {\frac{{\xi _2^e - \xi }}{{\xi _2^e - \xi _1^e}}} \right).
\end{equation}
The physical stress resultants are approximated, using the Lagrange polynomial functions of degree $p_\mathrm{p}$, as
\begin{equation}
	\label{approx_phy_rp}
	{\boldsymbol{r}}_{\rm{p}}^h(s(\xi )) = \left[ {\begin{array}{*{20}{c}}
			{{L}^{p_\mathrm{p}}_1(\bar \xi ){{\boldsymbol{1}}_{15 \times 15}}}& \cdots &{{L}^{p_\mathrm{p}}_{{n}^\mathrm{p}_e}(\bar \xi ){{\boldsymbol{1}}_{15 \times 15}}}
	\end{array}} \right]\left\{ {\begin{array}{*{20}{c}}
			{{{\bf{r}}^e_{1}}}\\
			\vdots \\
			{{{\bf{r}}^e_{{{n}}^\mathrm{p}_e}}}
	\end{array}} \right\} \eqqcolon {{{\Bbb{L}}}}_e(\bar \xi)\,{\bf{r}}^e,\,\,{\bar \xi}\in\left[-1,1\right],
\end{equation}
where ${{\bf{r}}^e_{I}} \in {\Bbb{R}^{15}}$ denotes the coefficient array for the physical stress resultants. $L^{p_\mathrm{p}}_I$ denotes the $I$th Lagrange polynomial function of degree $p_\mathrm{p}$ in each element, and $n^\mathrm{p}_e={p_\mathrm{p}}+1$ denotes the number of basis functions in $e$th element. Similarly, the physical strains are approximated by 
\begin{equation}
	\label{disc_pstrn}
	{\boldsymbol{\varepsilon }}_{\rm{p}}^h(s(\xi)) = {{{\Bbb{L}}}_e}(\bar \xi)\,{\bf{e}}^e,\,\,\mathrm{with}\,\,{\bf{e}}^e\coloneqq{\left\{ {\begin{array}{*{20}{c}}
				{{{\bf{e}}^e_{1}}}\\
				\vdots \\
				{{{\bf{e}}^e_{{{n}}^\mathrm{p}_e}}}
		\end{array}} \right\}},\,\,{\bar \xi}\in\left[-1,1\right],
\end{equation}
where ${{\bf{e}}^e_{I}} \in {{\Bbb{R}}^{15}}$ denotes the coefficient vector for the physical (kinematic) strains, $I\in\left\{1,2,...,{{n}^\mathrm{p}_e}\right\}$. The physical enhanced strain parameters are approximated by 
\begin{equation}\label{interp_var_alpha_beta}
	{\boldsymbol{\alpha }^h}(s(\xi)) = \left[ {\begin{array}{*{20}{c}}
			{{{L}^{{p}_\mathrm{a}}_1}(\bar \xi )\,{{\bf{1}}_{9 \times 9}}}&\cdots&{{{L}^{{p}_\mathrm{a}}_{{{n}}^\mathrm{a}_e}}(\bar \xi )\,{{\bf{1}}_{9 \times 9}}}
	\end{array}} \right]\left\{ {\begin{array}{*{20}{c}}
			{{\boldsymbol{\alpha }}^e_1}\\
			\vdots\\
			{{\boldsymbol{\alpha }}^e_{{{n}}^\mathrm{a}_e}}
	\end{array}} \right\} \eqqcolon \,{{\tilde{\Bbb{L}}}_{e}(\bar \xi)}\,{{\boldsymbol{\alpha }}^e},
\end{equation}
where ${{\boldsymbol{\alpha }}^e}$ is \textcolor{black}{the array of nodal coefficients} of the enhanced strain parameters, and $n^\mathrm{a}_e={p_\mathrm{a}}+1$ denotes the number of basis functions in $e$th element. In this paper, we use $p_\mathrm{a}=1$, as in \citet{choi2021isogeometric}. Substituting Eqs.\,(\ref{approx_phy_rp}) and (\ref{disc_pstrn}) into Eq.\,(\ref{vform_hw}), the internal virtual work is approximated by
\begin{equation}
	\label{spat_disc_int_vw_hw}
	G_{{\mathop{\rm int}} }^{{\rm{HW}}} \approx \sum\limits_{e = 1}^{{n_{{\rm{el}}}}} {\int_{{\varXi _e}} {\left\{ {\begin{array}{*{20}{c}}
					{\delta {{\bf{y}}^e}}\\
					{\delta {\bf{r}}^e}\\				
					{\delta {\bf{e}}^e}\\
					{\delta {{\boldsymbol{\alpha}}^e}}
			\end{array}} \right\} \cdot {\bf{F}}_{{\mathop{\rm int}} }^{e}\,\tilde j\,{\rm{d}}\xi}}= {\delta {{\bf{y}}^{*\,\mathrm{T}}}}\left(\mathop {\mathlarger{\mathlarger{\bf{A}}}}\limits_{e = 1}^{{n_{{\rm{el}}}}} {\bf{F}}_{{\mathop{\rm int}} }^{e}\right),
\end{equation}
with the elemental internal load vector, ${\bf{F}}_{{\mathop{\rm int}} }^{e} \coloneqq {\left[ {{\bf{f}}_\mathrm{y}^{e\,{\rm{T}}},{\bf{f}}_\mathrm{r}^{e\,{\rm{T}}},{\bf{f}}_\varepsilon ^{e\,{\rm{T}}},{\bf{f}}_\mathrm{a}^{e\,{\rm{T}}}} \right]^{\rm{T}}}$, defined by 
\begin{subequations}
	\label{sub_mat_int_f_vw}
	\begin{alignat}{3}
		\left[{\bf{f}}_\mathrm{y}^e\right]_{{m_e}\times{1}} &\coloneqq \int_{{\varXi _e}}{\Bbb{B}}_{{\rm{total}}}^{e\,{\rm{T}}}\,{{\boldsymbol{r}}_{\rm{p}}^h}\,\tilde j\,{\rm{d}}\xi,\label{submat_intf_disp}\\
		\left[{\bf{f}}_\mathrm{r}^{e}\right]_{{{m^\mathrm{p}_e}}\times{1}} &\coloneqq \int_{{\varXi _e}}{\Bbb{L}}^{\mathrm{T}}_e\left\{ {{\boldsymbol{\varepsilon }}({{\boldsymbol{y}}^h}) - {\boldsymbol{\varepsilon }}_{\rm{p}}^h} \right\}\tilde j\,{\rm{d}}\xi,\\
		\left[{\bf{f}}_\varepsilon ^e\right]_{{{m^\mathrm{p}_e}}\times{1}} &\coloneqq \int_{{\varXi _e}}{\Bbb{L}}^{\mathrm{T}}_e\left\{ {{{{\partial _{{{\boldsymbol{\varepsilon }}_{\rm{p}}}}}\psi ({{\boldsymbol{\varepsilon}}^h_{\rm{p}}},{\boldsymbol{\alpha }}^h)}} - {\boldsymbol{r}}_{\rm{p}}^h} \right\}\tilde j\,{\rm{d}}\xi,\\
		\left[{\bf{f}}^e_\mathrm{a}\right]_{{{m^\mathrm{a}_e}}\times{1}} &\coloneqq \int_{{\varXi _e}}{\tilde{\Bbb{L}}}^{\mathrm{T}}_e\,{{\partial _{\boldsymbol{\alpha }}}\psi ({{\boldsymbol{\varepsilon }}^h_{\rm{p}}},{\boldsymbol{\alpha}}^h)}\,\tilde j\,{\rm{d}}\xi.
	\end{alignat}
\end{subequations}
In Eq.\,(\ref{spat_disc_int_vw_hw}), $\bf{A}$ denotes the finite element assembly operator, and ${\bf{y}}^*$ denotes the global array of the coefficients in the center axis position, director displacement, and additional stress and strain fields. In Eq.\,(\ref{sub_mat_int_f_vw}), $m_e=3{n_e}+6{n^\mathrm{d}_e}$ denotes the DOF number of the approximated center axis position, and director displacement fields in each element, and $m_e^\mathrm{p}=d_\mathrm{p}\cdot{n^\mathrm{p}_e}$ denotes the DOF number of the approximated physical stress resultants or (kinematic) strain field in each element, and $m_e^\mathrm{a}=d_\mathrm{a}\cdot{n^\mathrm{a}_e}$ denotes the DOF number of the approximated enhanced strain field in each element. The detailed expression of matrix ${\Bbb{B}}_{{\rm{total}}}^{e}$ in Eq.\,(\ref{submat_intf_disp}) can be found in Section \ref{app_oper_intvw}. The increment of internal virtual work is also approximated by
\begin{equation}	
	\label{spat_disc_tstiff_hw_int}
	\Delta G_{{\mathop{\rm int}}}^{{\rm{HW}}} \approx \sum\limits_{e = 1}^{{n_{{\rm{el}}}}} {{\left\{ {\begin{array}{*{20}{c}}
					{\delta {{\bf{y}}^e}}\\
					{\delta {\bf{r}}^e}\\
					{\delta {\bf{e}}^e}\\
					{\delta {\boldsymbol{\alpha}}^e}		
			\end{array}} \right\} \cdot {\bf{K}}_{{\mathop{\rm int}} }^{e}\left\{ {\begin{array}{*{20}{c}}
					{\Delta {{\boldsymbol{y}}^e}}\\
					{\Delta {\bf{r}}^e}\\
					{\Delta {\bf{e}}^e}\\
					{\Delta {\boldsymbol{\alpha}}^e}
			\end{array}} \right\}}} =\delta {{\bf{y}}^{* {\rm{T}}}}\left({\mathop {\mathlarger{{\mathlarger{\bf{A}}}}}\limits_{e = 1}^{{n_{{\rm{el}}}}} {\bf{K}}_{{\mathop{\rm int}} }^e}\right)\Delta {{\bf{y}}^*},
\end{equation}
and the element tangent stiffness matrix
\begin{equation}
	{\bf{K}}_{{\mathop{\rm int}} }^e \coloneqq {\left[ {\renewcommand{\arraystretch}{1.5}\begin{array}{*{20}{c}}
				{{\bf{k}}_\mathrm{yy}^e}&{{\bf{k}}_\mathrm{ry}^{e\,\mathrm{T}}}&{{{\bf{0}}_{{m_e} \times {{m}^\mathrm{p}_e}}}}&{{{\bf{0}}_{{m_e} \times {{m}^\mathrm{a}_e}}}}\\
				{}&{{{\bf{0}}_{{{m}^\mathrm{p}_e} \times {{m}^\mathrm{p}_e}}}}&{{\bf{k}}_\mathrm{r\varepsilon}^e} & {{\bf{0}}_{{{m}^\mathrm{p}_e}\times{{m}^\mathrm{a}_e}}}\\
				&{}&{{\bf{k}}_{\varepsilon \varepsilon }^e} & {\bf{k}}_\mathrm{a\varepsilon}^{e\,\mathrm{T}}\\
				{{\rm{sym}}{\rm{.}}} & & & {\bf{k}}_\mathrm{aa}^{e}
		\end{array}} \right]},
\end{equation}
with
\begin{subequations}
	\begin{alignat}{3}
		{\bf{k}}_\mathrm{yy}^e&\coloneqq{\int_{{\varXi _e}}{{\Bbb{Y}}_e^{\rm{T}}}{{\boldsymbol{k}}_{\rm{G}}}{{\Bbb{Y}}_e}\,\tilde j\,{\rm{d}}\xi},\label{gtan_uti_from_prev}\\
		{\bf{k}}_\mathrm{ry}^e&\coloneqq{\int_{{\varXi _e}}{{\Bbb{L}}_{e}^{\mathrm{T}}}{\Bbb{B}}_{{\rm{total}}}^{e}\,\tilde j\,{\rm{d}}\xi},\\
		{\bf{k}}_\mathrm{r\varepsilon}^e&\coloneqq - {\int_{{\varXi _e}}{{\Bbb{L}}_e^{\rm{T}}}{{\Bbb{L}}_e}\,\tilde j\,{\rm{d}}\xi},\label{la_submat_k_re}\\
		{\bf{k}}_{\varepsilon \varepsilon }^e &\coloneqq{\int_{{\varXi _e}}{{\Bbb{L}}_e^{\rm{T}}}{{\Bbb{C}}^{{\varepsilon \varepsilon }}_{\rm{p}}}{{\Bbb{L}}_e}\,\tilde j\,{\rm{d}}\xi},\\
		{\bf{k}}_\mathrm{a\varepsilon}^e&\coloneqq {\int_{{\varXi _e}}{{\tilde{\Bbb{L}}}_e^{\rm{T}}}{{\Bbb{C}}^\mathrm{a\varepsilon}_{\rm{p}}}{{\Bbb{L}}_e}\,\tilde j\,{\rm{d}}\xi},\\
		{\bf{k}}_\mathrm{aa}^e &\coloneqq{\int_{{\varXi _e}}{{\tilde{\Bbb{L}}}_e^{\rm{T}}}{{\Bbb{C}}^\mathrm{aa}_{\rm{p}}}{{\tilde{\Bbb{L}}}_e}\,\tilde j\,{\rm{d}}\xi}.
	\end{alignat}
\end{subequations}
Detailed expressions of ${{\boldsymbol{k}}_{\rm{G}}}$ and ${{\Bbb{Y}}_e}$ in Eq.\,(\ref{gtan_uti_from_prev}) can be found in Appendix \ref{app_oper_gtan}. The external virtual work of Eq.\,(\ref{var_ex_vwork}) is also approximated by
\begin{equation}
	\label{fedisc_nbs_ext_v}
	G_{{\mathop{\rm ext}} }(\delta{\boldsymbol{y}}) \approx {\delta {{\bf{y}}^{\mathrm{T}}}}{\bf{F}}_{{\mathop{\rm ext}} },\,\,\mathrm{with}\,\,{\bf{F}}_{{\mathop{\rm ext}} }\coloneqq{\mathop{\mathlarger{\mathlarger{\bf{A}}}}\limits_{e = 1}^{{n_{{\rm{el}}}}}}{{\bf{F}}_{{\mathop{\rm ext}} }^{e}} + {\mathlarger{\mathlarger{{\bf{A}}}}}{\left[ {{{{\boldsymbol{\bar R}}}_0}} \right]_{s\in{\Gamma _{\text{N}}}}},
\end{equation}
where $\bf{y}$ denotes the global array of the control coefficients for the center axis position, and director displacement fields. We have also defined
\begin{equation}
	{{\bf{F}}_{{\mathop{\rm ext}} }^{e}}\coloneqq\int_{{\varXi _e}}{\Bbb{N}}_e^{\mathrm{T}}\,{{\boldsymbol{\bar R}}}\,\,\tilde j\,{\rm{d}}\xi.
\end{equation}
Substituting Eqs.\,(\ref{spat_disc_int_vw_hw}), (\ref{spat_disc_tstiff_hw_int}), and (\ref{fedisc_nbs_ext_v}) into Eq.\,(\ref{lin_var_eq_g_hw}) gives
%
\begin{equation}
	\label{disc_var_eq_lin}
	\delta {{\bf{y}}^{ *{\rm{T}}}}\,\!\left({\mathop{\mathlarger{\mathlarger{\bf{A}}}}\limits_{e = 1}^{{n_{{\rm{el}}}}}}\,{\bf{K}}_{{{\rm int}} }^e\right)\Delta {{\bf{y}}^*} = \delta{\bf{y}}^\mathrm{T}\,{\bf{F}}_\mathrm{ext}-\delta {{\bf{y}}^{ *{\rm{T}}}}\,\!\left(\mathop{\mathlarger{\mathlarger{\bf{A}}}}\limits_{e = 1}^{{n_{{\rm{el}}}}} {\bf{F}}_{{{\rm int}} }^{e}\right).
\end{equation}
\subsubsection{Element-wise static condensation}
Since the physical strain and stress resultants may have discontinuities between adjacent elements, Eq.\,(\ref{disc_var_eq_lin}) can be rewritten as
\begin{subequations}
	\begin{equation}
		\sum\limits_{e = 1}^{{n_{{\rm{el}}}}} \delta {{\bf{y}}^{e\,{\rm{T}}}}\left( {{\bf{k}}_\mathrm{yy}^e\Delta {{\bf{y}}^e} + {\bf{k}}_\mathrm{ry}^{e\,\mathrm{T}}\Delta {\bf{r}}^e} \right) = \delta{\bf{y}}^\mathrm{T}\,{\bf{F}}_\mathrm{ext} - \sum\limits_{e = 1}^{{n_{{\rm{el}}}}} {\delta {{\bf{y}}^{e\,{\rm{T}}}}{\bf{f}}_\mathrm{y}^e},\label{split_disc_y}		
	\end{equation}
	\begin{alignat}{3}
		{\bf{k}}_\mathrm{ry}^{e}\Delta {{\bf{y}}^e} &+ {\bf{k}}_\mathrm{r\varepsilon }^e\Delta {{\bf{e}}^e} &=  - {\bf{f}}_\mathrm{r}^e,\,\,e\in\left\{1,\cdots,{n_\mathrm{el}}\right\},\label{split_disc_r}\\
		{\bf{k}}_\mathrm{r\varepsilon }^{e\,{\rm{T}}}\Delta {\bf{r}}^e + {\bf{k}}_\mathrm{\varepsilon \varepsilon }^e\Delta {{\bf{e}}^e}  &+ {\bf{k}}_\mathrm{a \varepsilon}^{e\,\mathrm{T}}\Delta {{\boldsymbol{\alpha}}^e}&=  - {\bf{f}}_\varepsilon ^e,\,\,e\in\left\{1,\cdots,{n_\mathrm{el}}\right\}\label{split_disc_e},\\
		{\bf{k}}_\mathrm{a\varepsilon }^{e}\Delta {{\bf{e}}^e} &+ {\bf{k}}_\mathrm{aa}^{e}\Delta {{\boldsymbol{\alpha}}^e}&=  - {\bf{f}}_\mathrm{a}^e,\,\,e\in\left\{1,\cdots,{n_\mathrm{el}}\right\} \label{split_disc_a}.
	\end{alignat}
\end{subequations}
Since the matrices ${\bf{k}}_\mathrm{r\varepsilon }^{e}$ and ${\bf{k}}_\mathrm{aa}^{e}$ are invertible, we obtain from Eqs.\,(\ref{split_disc_r})-(\ref{split_disc_a})
\begin{subequations}
	\label{del_elem_coeff_e_a_r}
	\begin{align}
		\Delta {{\bf{e}}^e} &=  - {\bf{k}}_\mathrm{r\varepsilon }^{e\,\, - 1}\left( {{\bf{f}}_\mathrm{r}^e + {\bf{k}}_\mathrm{ry}^{e}\Delta {{\bf{y}}^e}} \right),\label{del_e_elem_vec}\\
		\Delta {\bf{r}}^e &=   -{\bf{k}}_\mathrm{r\varepsilon }^{e\,\, - {\rm{T}}}\left({\bf{f}}_\varepsilon ^e+{\bf{k}}_{\varepsilon \varepsilon }^e\Delta{\bf{e}}^e + {\bf{k}}_\mathrm{a \varepsilon}^{e\,\mathrm{T}}\Delta {{\boldsymbol{\alpha}}^e}\right),\label{del_r_p_e_static_c}\\
		\Delta {{\boldsymbol{\alpha}}^e} &=  - {\bf{k}}_\mathrm{aa}^{e\,\, - 1}\left( {{\bf{f}}_\mathrm{a}^e + {\bf{k}}_\mathrm{a\varepsilon }^e\Delta {{\bf{e}}^e}} \right).\label{del_a_elem_vec}
	\end{align}
\end{subequations}
Then, Eq.\,(\ref{split_disc_y}) can be rewritten as
\begin{equation}
	\label{eq_bef_app_dbdc_kd_f}
	\delta {{\bf{y}}^{{\rm{T}}}}{\bf{\bar K}}\,\Delta {{\bf{y}}} = \delta {{\bf{y}}^{{\rm{T}}}}{\bf{\bar F}},
\end{equation}
where we define
\begin{equation}
	\label{glob_tstiff_cond}
	{\bf{\bar K}} \coloneqq \mathop {\mathlarger{\mathlarger{\bf{A}}}}\limits_{e = 1}^{{n_{{\rm{el}}}}} {\bf{\bar K}}^e_\mathrm{int},
\end{equation}
and
\begin{equation}
	{\bf{\bar F}} \coloneqq \mathop {\mathlarger{\mathlarger{{\bf{A}}}}}\limits_{e = 1}^{{n_{{\rm{el}}}}} \left({{\bf{F}}^e_\mathrm{ext}-{\bf{\bar F}}^e_\mathrm{int}}\right) + {\mathlarger{\mathlarger{{\bf{A}}}}}{\left[ {{{{\boldsymbol{\bar R}}}_0}} \right]_{s\in{\Gamma _{\text{N}}}}},
\end{equation}
with
\begin{equation}
	\label{glob_tstiff_cond_elem}
	{\bf{\bar K}}^e_\mathrm{int} \coloneqq {{{\bf{k}}_\mathrm{yy}^e + {\bf{k}}_\mathrm{ry}^{e\,{\rm{T}}}{\bf{k}}_{\mathrm{r}\varepsilon }^{e\,\, - {\rm{T}}}\underbrace{\left({\bf{k}}_{\varepsilon \varepsilon }^e - {\bf{k}}_{\mathrm{a}\varepsilon }^{e\,{\rm{T}}}{\bf{k}}_\mathrm{aa}^{e\, - 1}{\bf{k}}_{\mathrm{a}\varepsilon }^e \right)}_{\eqqcolon{{\bf{\bar k}}_{\varepsilon \varepsilon }^e}}{\bf{k}}_{r\varepsilon }^{e\, - 1}{\bf{k}}_\mathrm{ry}^e}},
\end{equation}
and
\begin{equation}
	\label{glob_resid_vec_fe_int}
	{\bf{\bar F}}^e_\mathrm{int} \coloneqq {{\bf{f}}_\mathrm{y}^e - {\bf{k}}_\mathrm{ry}^{e\,{\rm{T}}}{\bf{k}}_{\mathrm{r}\varepsilon }^{e\, - {\rm{T}}}\left\{{\bf{f}}_\varepsilon ^e - {{\bf{k}}_{\varepsilon \varepsilon }^e{\bf{k}}_{\mathrm{r}\varepsilon }^{e\, - 1}}{\bf{f}}_\mathrm{r}^e			
		-{\bf{k}}_{\mathrm{a}\varepsilon }^{e\,{\rm{T}}}{\bf{k}}_\mathrm{aa}^{e\, - 1}\left({\bf{f}}_\mathrm{a}^e-{\bf{k}}_{\mathrm{a}\varepsilon }^e{\bf{k}}_{\mathrm{r}\varepsilon }^{e\, - 1} {\bf{f}}_\mathrm{r}^e\right)\right\}}.
\end{equation}
The only matrices which need to be inverted are ${\bf{k}}^e_{\mathrm{r}\varepsilon}$ and ${\bf{k}}^e_\mathrm{aa}$. By applying the displacement boundary conditions, we finally have the following reduced system of linear equations at $i$th iteration in the $(n+1)$th load step
\begin{equation}
	\label{kyf_sys_lin_eqs}
	{}^{n+1}{{\bf{\bar K}}^{(i-1)}_{\rm{r}}}\,\Delta {\bf{y}}_{\rm{r}} = {}^{n+1}{\bf{\bar F}}^{(i-1)}_{\rm{r}},
\end{equation}
where $(\bullet)_\mathrm{r}$ denotes the reduced vector or matrix.
\textcolor{black}{
\begin{remark} Since we use different control nets for the spatial discretization of the center axis displacement, and the director displacement fields (i.e., $n_\mathrm{cp}\ne{n_\mathrm{cp}^\mathrm{d}}$), the arrangement of control coefficients for those two fields should be separated in $\delta {\bf{y}}$, and $\Delta {\bf{y}}$, as well as in $\Delta {\bf{y}}^e$ and $\delta {\bf{y}}^e$. This eventually makes the tangent stiffness matrix $\bf{\bar K}$ banded locally, not globally, see the examples of sparsity patterns in Figs.\,\ref{sparsity_iga_loc_red} and \ref{sparsity_iga_loc_red_nel20}.
\end{remark}
\begin{remark} \label{remk_gloc_app_term}\textit{Global and local approaches to approximate the physical stress resultant and (kinematic) strain fields}. We consider the following two approaches, combined with the proposed mixed formulation within the framework of IGA:
\begin{itemize}
	\item A \textbf{global} approach (``\textbf{glo}'') using NURBS basis functions of degree $p_\mathrm{p} = {p-1}$ for the physical stress resultant and strain fields having $C^{{p}_\mathrm{p}-1}$ inter-element continuity, which is combined with a patch-wise static condensation,
	\item A \textbf{local} approach allowing inter-element discontinuity of the physical stress resultant and strain fields, which is combined with an element-wise static condensation. For IGA, depending on the selection of degree $p_\mathrm{p}$, the local approach is subdivided into
	\begin{itemize}
		\item ``\textbf{loc}'': $p_\mathrm{p}=p-1$,
		\item ``\textbf{loc-ur}'': Uniformly reduced degree $p_\mathrm{p}=1$,
		\item ``\textbf{loc-sr}'': Selectively reduced degree $p_\mathrm{p}$, given by Table \ref{tab_selec_deg_r}.
	\end{itemize}
\end{itemize}
The further reduction of $p_\mathrm{p}$ in ``loc-ur'' and ``loc-sr'' aims at the alleviation of locking due to the higher order inter-element continuity in the displacement field of IGA. In FEA, we use only the local approach, ``loc''. In all the approaches, the physical enhanced strain field is allowed to be discontinuous across the elements, which allows an element-wise static condensation, and we always use $p_\mathrm{a}=1$.
\end{remark} 
\begin{remark} \label{time_complexity_stat_cond} \textit{Counting operations in the static condensation process}. Here we discuss the time complexity of the element-wise static condensation process in the calculation of the element tangent stiffness matrix ${\bf{\bar K}}_{{\mathop{\rm int}} }^e$ in Eq.\,(\ref{glob_tstiff_cond_elem}). The inversion of a $n\times{n}$ matrix using the Gauss elimination process requires $O(n^3)$ operations \citep{moin2010fundamentals}. For a slightly more efficient algorithm, one may refer to the relevant comments in \citet{kikis2022two} and references therein. The static condensation in Eq.\,(\ref{glob_tstiff_cond_elem}) is performed in two steps. In the \textit{first step}, we condense out the DOFs of the \textit{enhanced} strains. The calculation of ${{\bf{\bar k}}_{\varepsilon \varepsilon }^e}$ consists of the matrix inversion ${{\bf{k}}_\mathrm{aa}^{e\,-1}}$, which requires $O({m_e^{\mathrm{a}\,3}})$ operations, where $m^\mathrm{a}_e\coloneqq{n^\mathrm{a}_e\cdot{d_\mathrm{a}}}={(p_\mathrm{a}+1)\cdot{d_\mathrm{a}}}$ denotes the DOF number of the approximated enhanced strain field in each element, i.e., the dimension of the square matrix ${\bf{k}}^e_\mathrm{aa}$. The subsequent two matrix--matrix multiplications, and the addition of the resulting matrix to ${{\bf{k}}_{\varepsilon \varepsilon }^e}$ require $O(m^{\mathrm{a}\,2}_{e}\cdot{m^\mathrm{p}_e} + m^{\mathrm{a}}_{e}\cdot{m^{\mathrm{p}\,2}_e})$, and $O(m_e^{\mathrm{p}\,2})$ operations, respectively, where $m^\mathrm{p}_e\coloneqq{n^\mathrm{p}_e}\cdot{d_\mathrm{p}}={(p_\mathrm{p}+1)\cdot{d_\mathrm{p}}}$ denotes the DOF number of the approximated physical stress resultants or (kinematic) strain field in each element. In this paper, we use $p_\mathrm{a}=1$, so that the time complexity for the first step is $O({p_\mathrm{p}}^2)$ with $n_\mathrm{el}$ fixed, and is $O(n_\mathrm{el})$ with $p_\mathrm{a}$ and $p_\mathrm{p}$ fixed, due to the element-wise operations. In the \textit{second step}, we condense out the DOFs of the physical stress resultants and strains, where we need the matrix inversion ${{\bf{k}}_{\mathrm{r}\varepsilon}^{e\,-1}}$, which requires $O({m^{\mathrm{p}\,3}_e})$ operations. The subsequent matrix--matrix multiplications, and the addition of the resulting matrix to ${{\bf{k}}}^e_\mathrm{yy}$ require $O(m_e\cdot{m^{\mathrm{p}\,2}_e}+{m_e}^2\cdot{m^\mathrm{p}_e})$, and $O({m_e}^{2})$ operations, respectively, where $m_e\coloneqq{3\cdot{n_e}+6\cdot{n^\mathrm{d}_e}={3\cdot(p+1)+6\cdot{(p_\mathrm{d}+1)}}}$ denotes the DOF number of the approximated center axis, and director displacement fields in each element. Note that, in the global approach, the dimension of the inverted matrix is proportional to the total number of elements ($n_\mathrm{el}$). Table\,\ref{tab_comp_time_complex} summarizes the time complexity of the static condensation in each approach. It is clear that the element-wise condensation is much more efficient than the patch-wise one in the global approach (``glo'') for the same $p$, which is more pronounced, as $n_\mathrm{el}$ increases. Further, among the local approaches, the operations in ``loc-ur'' and ``loc-sr'' with reduced $p_\mathrm{p}$ are much less expensive than that in ``loc'', since $p_\mathrm{p}$ is much smaller than $p$, which is more pronounced, as $p$ or $n_\mathrm{el}$ increases. 
\end{remark}
}
\begin{table}[H]
	\centering
	\caption{\textcolor{black}{Time complexity of the static condensation in the global and local approaches. Note that the degrees $p_\mathrm{p}$ in the approach ``loc-ur'' and ``loc-sr'' of IGA are much smaller than $p$.}}
	\label{tab_comp_time_complex}	
	\begin{tabular}{llllll}
		\toprule
		& \multicolumn{2}{c}{First step (condensation of $\boldsymbol{\alpha}$)}  &       & \multicolumn{2}{c}{Second step (condensation of $\boldsymbol{r}_\mathrm{p}$ and $\boldsymbol{\varepsilon}_\mathrm{p}$)} \\
		\cmidrule{1-3}\cmidrule{5-6}      & {\makecell[l]{$n_\mathrm{el}$ changing,\\with $p$ fixed}} & {\makecell[l]{$p_\mathrm{a}=1$, $p$ changing,\\with $n_\mathrm{el}$ fixed}} &       & {\makecell[l]{$n_\mathrm{el}$ changing,\\with $p$ fixed}} & {\makecell[l]{$p$ changing,\\with $n_\mathrm{el}$ fixed}} \\
		\midrule
		glo   & $O({n_\mathrm{el}})$ & $O(p^2)$ &       & $O({n_\mathrm{el}}^3)$ & $O(p^3)$ \\
		loc   & $O(n_\mathrm{el})$ & $O(p^2)$ &       & $O(n_\mathrm{el})$ & $O(p^3)$ \\
		loc-ur & $O(n_\mathrm{el})$ & $O({p_\mathrm{p}}^2)$ &       & $O(n_\mathrm{el})$ & $O({{p_\mathrm{p}}}\cdot{p^2})$ \\
		loc-sr & $O(n_\mathrm{el})$ & $O({p_\mathrm{p}}^2)$ &       & $O(n_\mathrm{el})$ & $O({{p_\mathrm{p}}}\cdot{p^2})$ \\
		\bottomrule
	\end{tabular}%
\end{table}
\subsubsection{Update of the configuration and internal variables}
Using the solution increment $\Delta {\bf{y}}_\mathrm{r}$ obtained by solving Eq.\,(\ref{kyf_sys_lin_eqs}), we update the control coefficients $\boldsymbol{\varphi}_I$ and ${\bf{\bar d}}_J$ by the following procedure: At the $i$th iteration in the $(n+1)$th load step, we update
\begin{subequations}
	\begin{align}
		{}^{n + 1}{\boldsymbol{\varphi}}_I^{(i)} &= {}^{n + 1}{\boldsymbol{\varphi}}_I^{(i - 1)} + \Delta {\boldsymbol{\varphi}}_I,\,\,I\in\left\{1,2,\cdots,n_\mathrm{cp}\right\},\\ 
		{}^{n + 1}{\bf{\bar d}}_J^{(i)} &= {}^{n + 1}{\bf{\bar d}}_J^{(i - 1)} + \Delta {\bf{d}}_J,\,\,J\in\left\{1,2,\cdots,n^\mathrm{d}_\mathrm{cp}\right\},\label{additive_update_dir_coef}
	\end{align}
\end{subequations}
with the initial guesses ${}^{n + 1}{\boldsymbol{\varphi}}_I^{(0)}\equiv{}^{n}{\boldsymbol{\varphi}}_I$, and ${}^{n + 1}{\bf{\bar d}}_J^{(0)}\equiv{}^{n}{\bf{\bar d}}_J$, from the equilibrium configuration in the previous ($n$th) load step. Further, from $\Delta {\bf{y}}_\mathrm{r}$, $\Delta {\bf{y}}^e$ can be simply extracted for each element, and is substituted into Eq.\,(\ref{del_elem_coeff_e_a_r}) to obtain the increment of internal variables $\Delta {\bf{e}}^e$, $\Delta {\bf{r}}^e$, and $\Delta {\boldsymbol{\alpha}}^e$. Then, for every element $e\in\left\{1,\cdots,{n_\mathrm{el}}\right\}$, we update
\begin{subequations}
	\begin{alignat}{3}
		{}^{n + 1}{{\bf{e}}^{e\,(i)}} &= {}^{n + 1}{{\bf{e}}^{e\,(i - 1)}} &&+ \Delta {{\bf{e}}^e},\\
		{}^{n + 1}{{\bf{r}}^{e\,(i)}} &= {}^{n + 1}{{\bf{r}}^{e\,(i - 1)}} &&+ \Delta {{\bf{r}}^e},\\
		{}^{n + 1}{{\boldsymbol{\alpha}}^{e\,(i)}} &= {}^{n + 1}{{\boldsymbol{\alpha }}^{e\,(i - 1)}} &&+ \Delta {{\boldsymbol{\alpha }}^e},
	\end{alignat}
\end{subequations}
with the initial guesses ${}^{n + 1}{{\bf{e}}^{e\,(0)}}\equiv{}^{n}{{\bf{e}}^{e}}$, ${}^{n + 1}{{\bf{r}}^{e\,(0)}}\equiv{}^{n}{{\bf{r}}^{e}}$, and ${}^{n + 1}{{\boldsymbol{\alpha}}^{e\,(0)}}\equiv{}^{n}{{\boldsymbol{\alpha}}^{e}}$, from the converged solution in the previous ($n$th) load step. 

\section{Imposition of rotational continuity between beams}
\label{rcont_betw_b}
In order to facilitate the connection of multiple beams with an arbitrary initial intersection angle and rotational continuity, we introduce rotational degrees-of-freedom for directors at the ends of beam. Here we are only concerned with the rigid joint condition, that is, all the connected beams' cross-sections may have only three rotational degrees-of-freedom. However, in contrast to the conventional rigid joint, our formulation may still allow stretching along the two directors in each cross-section. The presented formulation is limited to the finite element basis functions having Kronecker-delta property at the end points of the beam. Therefore, for IGA, we only use clamped knot vectors. This enables us to simply \textcolor{black}{apply the necessary transformation operations after the finite element approximation, considering} only the director coefficients corresponding to the end points, as in the conventional finite element formulation, e.g., see \citet{romero2002objective}.
\subsection{Multiplicative decomposition of the directors}
The director vectors can be decomposed into a unit director and a scalar stretch, as \citep{simo1990stress}
\begin{equation}
	\label{mdec_dir_alph}
	{{\boldsymbol{d}}_\alpha }(s) = {\lambda _\alpha }(s)\,{{\boldsymbol{t}}_\alpha }(s),
\end{equation} 
with the stretch \textcolor{black}{ratio} ${\lambda _\alpha }(s) \coloneqq \left\| {{{\boldsymbol{d}}_\alpha }}(s) \right\|>0$, where no summation is implied on $\alpha\in\left\{1,2\right\}$. Further, the unit vectors can be expressed by rotational transformation as
\begin{equation}
	\label{rot_tr_tvec_a}
	\boldsymbol{t}_\alpha(s) = {\boldsymbol{\Lambda}_\alpha}(s)\,{\boldsymbol{E}_\alpha}\,\,\mathrm{with}\,\,{\boldsymbol{\Lambda}_\alpha}(s)\in{\mathrm{S}_\alpha^2},\,\,\alpha\in\left\{1,2\right\},
\end{equation}
where we define \citep{simo1989stress}
\begin{equation}
	\mathrm{S}_{\alpha}^2 \coloneqq \left\{ {\left. {{\boldsymbol{\Lambda }} \in {\rm{SO}}(3)} \right|{\boldsymbol{\Lambda}\boldsymbol{\psi}} = {\boldsymbol{\psi }},\,\,\text{and}\,\,{\boldsymbol{\psi }} \in {{\Bbb{R}}^3}\,\,\text{satisfies}\,\,{\boldsymbol{\psi }} \cdot {\boldsymbol{E}}_\alpha = 0} \right\}.
\end{equation}
$\mathrm{SO(3)}$ defines the set of proper orthogonal tensors in three-dimensional space, and $\mathrm{S}_{\alpha}^2$ is a subset of $\mathrm{SO(3)}$ such that the axis of rotation is orthogonal to ${{\boldsymbol{E}}_\alpha}\in\Bbb{R}^3$. \textcolor{black}{That is, each unit director ${\boldsymbol{t}}_\alpha$ can be parameterized by two rotational DOFs. Here, for simplicity, we restrict our discussion by removing the DOF of changing the angle between two unit directors ${\boldsymbol{t}}_1$ and ${\boldsymbol{t}}_2$. Therefore, the rotational motion of two unit directors is described by three rotational DOFs, which constitutes the orthogonal tensor $\boldsymbol{\Lambda}$, such that
\begin{align}
	\label{rot_tr_tvec_lambda_3}
	\boldsymbol{t}_\alpha(s) = {\boldsymbol{\Lambda}}(s)\,{\boldsymbol{E}_\alpha}\,\,\mathrm{with}\,\,{\boldsymbol{\Lambda}}(s)\in{\mathrm{SO(3)}},\,\,\alpha\in\left\{1,2\right\}.
\end{align}
At an admissible perturbed configuration, Eq.\,(\ref{rot_tr_tvec_lambda_3}) can be rewritten as
\begin{equation}
	\label{rot_tr_tvec_a_ptb}	
	{{\boldsymbol{t}}_{\alpha \varepsilon }}(s) = {{\boldsymbol{\Lambda }}_{\varepsilon }}(s)\,{{\boldsymbol{E}}_\alpha },\,\,\varepsilon\in\Bbb{R},
\end{equation}
where the subscript $(\bullet)_\varepsilon$ indicates the dependence of $(\bullet)$ to the perturbation amount $\varepsilon$. Hereafter, we often omit the argument $s$ for brevity. Taking the directional derivative of Eq.\,(\ref{rot_tr_tvec_a_ptb}) yields
\begin{equation}
	\label{del_unit_vec_t_alph}
	\delta {{\boldsymbol{t}}_\alpha } \coloneqq \frac{\mathrm{d}}{{\mathrm{d}\varepsilon}}
	{\left. {{{\boldsymbol{t}}_{\alpha \varepsilon }}} \right|_{\varepsilon  = 0}} = \widehat {\delta {{\boldsymbol{\theta }}}}\,{{\boldsymbol{t}}_\alpha },\,\,\alpha\in\left\{1,2\right\},
\end{equation}
where $\widehat {\delta {{\boldsymbol{\theta }}}} \coloneqq \delta {{\boldsymbol{\Lambda }}}{\boldsymbol{\Lambda }}^{\rm{T}}$ is a skew-symmetric tensor, with the notation $\widehat {(\bullet)}$ representing the skew-symmetric tensor associated for a given vector $(\bullet)\in{{\Bbb R}^3}$, such that $\widehat {(\bullet)}{\boldsymbol{h}} = (\bullet) \times {\boldsymbol{h}},\,\forall {\boldsymbol{h}} \in {{\Bbb{R}}^3}$. Thus, in the cross-section at a selected end point, we eventually have three rotational DOFs, i.e., the components of $\delta{\boldsymbol{\theta}}\in\Bbb{R}^3$. Note that the stretches along the directors, $\lambda_\alpha>0$ ($\alpha\in\left\{1,2\right\}$) are kinematically not constrained.} Taking the first variation of the director in Eq.\,(\ref{mdec_dir_alph}), and substituting Eq.\,(\ref{del_unit_vec_t_alph}) gives
\begin{equation}
	\label{transf_gen_dely_red}
	\left\{ {\begin{array}{*{20}{c}}
		{\delta {{\boldsymbol{d}}_1}}\\
		{\delta {{\boldsymbol{d}}_2}}
		\end{array}} \right\} = {\boldsymbol{\Xi }}\left\{ {\begin{array}{*{20}{c}}
			{\delta {\boldsymbol{\theta }}}\\
			{\delta {\mu _1}}\\
			{\delta {\mu _2}}
	\end{array}} \right\},\,\,\mathrm{with}\,\,{\boldsymbol{\Xi }} \coloneqq \left[ {\begin{array}{*{20}{c}}
	{ - \widehat {{{\boldsymbol{d}}_1}}}&{{{\boldsymbol{d}}_1}}&{{{\boldsymbol{0}}_{3\times{1}}}}\\
	{ - \widehat {{{\boldsymbol{d}}_2}}}&{{{\boldsymbol{0}}_{3\times{1}}}}&{{{\boldsymbol{d}}_2}}
\end{array}} \right],
\end{equation}
where \textcolor{black}{${\mu _\alpha } \coloneqq \ln {\lambda _\alpha }$} $(\alpha\in\left\{1,2\right\})$ defines the logarithmic stretch ratio.
\subsection{Isogeometric finite element formulation}
\begin{figure}[H]	
	\centering
	\begin{subfigure}[b] {0.5\textwidth} \centering
		\includegraphics[width=\linewidth]{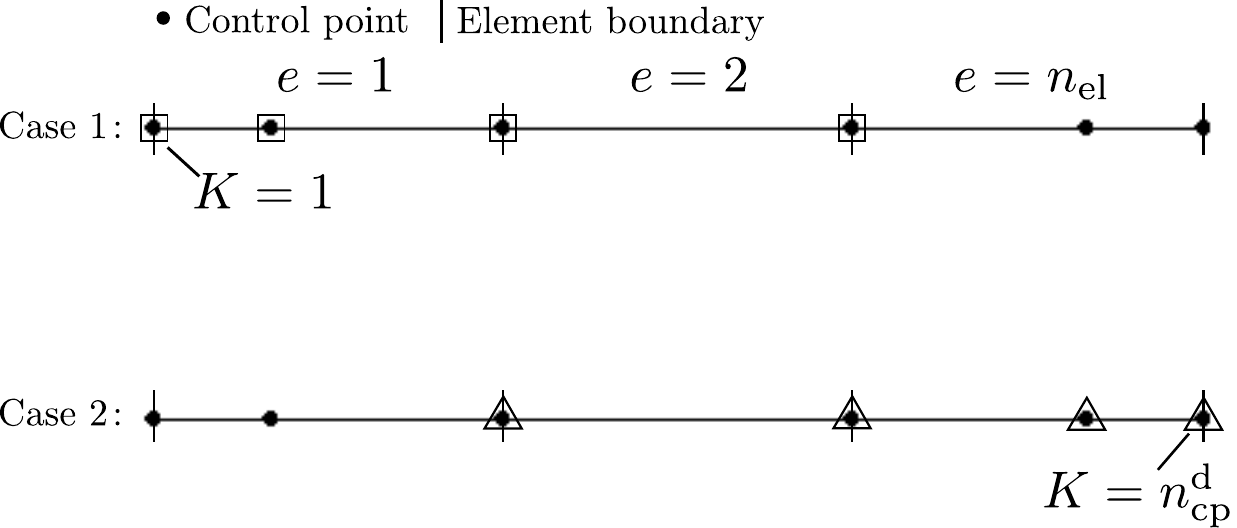}	
	\end{subfigure}		
	\caption{Illustration of two cases of joint positions. The curve patch (solid line) of $p_\mathrm{d}=3$ consists of three elements, and $\square$ and $\triangle$ indicate the control points having local support at the first ($e=1$) and the last element ($e=n_\mathrm{el}$) of the patch, respectively. A joint condition can be imposed at the left end (the first control point, i.e., $K=1$) (top figure) or the right end (the last control point, i.e., $K=n^\mathrm{d}_\mathrm{cp}$) (bottom figure).}
	\label{end_cp_loc_sup}
\end{figure}

Fig.\,\ref{end_cp_loc_sup} illustrates two cases of junction positions in a curve patch, where we can apply the joint condition to 
\begin{itemize}
	\item Case 1 (joint at the left-end): the first control point ($K=1$) having local support at the first element $e=1$, or
	\item Case 2 (joint at the right-end): the last control point ($K={n^\mathrm{d}_\mathrm{cp}}$) having local support at the last element $n_\mathrm{el}$,
\end{itemize}
where $K$ denotes the index of control points in a given curve patch. Here, for simplicity, we consider that an end element has a joint condition only at one end. 
\noindent We also denote a set of indices of elements having a joint condition by $\mathcal{J}$. That is, if $e$th element has a joint condition, $e\in\mathcal{J}$, otherwise, $e\not\in\mathcal{J}$. At the end control point of the patch, using Eq.\,(\ref{transf_gen_dely_red}), we reduce the six director DOFs to five, three rotation parameters and two logarithmic stretch ratios, such that
\begin{equation}
	\delta{{\bf{d}}_{K}}\coloneqq\left\{ {\begin{array}{*{20}{c}}
			{\delta {{\bf{d}}_{1K}}}\\
			{\delta {{\bf{d}}_{2K}}}
	\end{array}} \right\} = {{\bf{\Xi }}_K}\left\{ {\begin{array}{*{20}{c}}
	{\delta {{\bf{\Theta }}_K}}\\
	{\delta {\mu _{1K}}}\\
	{\delta {\mu _{2K}}}
\end{array}} \right\},
\end{equation}
with
\begin{equation}
	\label{del_y_red_lda_tht}
	{{\bf{\Xi }}_K} \coloneqq \left[ {\begin{array}{*{20}{c}}
		{ - \widehat {{{\bf{d}}_{1K}}}}&{{{\bf{d}}_{1K}}}&{{{\bf{0}}_{3 \times 1}}}\\
		{ - \widehat {{{\bf{d}}_{2K}}}}&{{{\bf{0}}_{3 \times 1}}}&{{{\bf{d}}_{2K}}}
	\end{array}} \right]_{6\times5},
\end{equation}
where $K=1$ for $e=1$ (i.e., Case 1), and $K=n^\mathrm{d}_\mathrm{cp}$ for $e=n_\mathrm{el}$ (i.e., Case 2), see also Table\,\ref{tab_corr_indic_jct} for the pair of indices in each case. Here, we also define \textcolor{black}{${{\bf{d}}_{\alpha K}} \coloneqq {\boldsymbol{d}}_\alpha ^h|_{s=0}$} for $K=1$, and \textcolor{black}{${{\bf{d}}_{\alpha K}} \coloneqq {\boldsymbol{d}}_\alpha ^h|_{s=L}$} for $K=n^\mathrm{d}_\mathrm{cp}$.
\begin{remark}
	\label{fict_dof_trans_rem}
	In our computer implementation, for convenience, we introduce a fictitious DOF\,($\delta{\Theta}_K^*$) to have the same number of control coefficients at the selected boundary control points as we have at internal control points. The matrix ${{\boldsymbol{\Xi }}_K}$ is also modified to have an additional sixth column of zeros accordingly, as
	\begin{equation}
		\label{trans_star_del_y}
		\delta {{\bf{d}}_K} = {\bf{\Xi }}_K^*\delta {{\bf{d}}^*_K},\,\,\mathrm{with}\,\,\,{\boldsymbol{\Xi }}_K^* \coloneqq \left[ {\begin{array}{*{20}{c}}
				{{{{\boldsymbol{\Xi }}}_K}}&{{{\boldsymbol{0}}_{6 \times 1}}}
		\end{array}} \right]_{6\times6},\,\,\mathrm{and}\,\,\delta {{\bf{d}}^*_K}\coloneqq\left\{ {\begin{array}{*{20}{c}}
			{\delta {{\bf{\Theta }}_K}}\\
		{\delta {\mu _{1K}}}\\
		{\delta {\mu _{2K}}}\\
		{\delta {\Theta_K ^ * }}
	\end{array}} \right\},\,\,K\in\left\{1,n^\mathrm{d}_\mathrm{cp}\right\}.
	\end{equation}
	Note that the fictitious DOF ($\delta{\Theta}_K^*$) is associated with the change of the angle between two directors, which is not allowed in the rigid joint condition. The artificial DOFs are removed, together with the constrained displacement DOFs, so that they are not considered in the final reduced system of linear equations.
\end{remark}
\noindent Thus, if the first or last element of the given NURBS patch contains a junction, we apply the transformation
\begin{equation}
	\label{del_trans_y_y_td}
	\delta {{\bf{d}}^e} = {{\boldsymbol{\Xi }}^{e}}\delta {{\bf{d}}^{*e}},\,\,e\in\mathcal{J},
\end{equation}
with
\begin{equation}
	{{\boldsymbol{\Xi }}^{e}} \coloneqq \begin{cases}\renewcommand\arraystretch{1.5}\left[ {\begin{array}{*{20}{c}}
				{{{{\boldsymbol{\Xi }}}^{*}_1}}&{{{{\boldsymbol{0}}_{6 \times 6({n^\mathrm{d}_e} - 1)}}}}\\
				{{{{\boldsymbol{0}}_{6({n^\mathrm{d}_e} - 1) \times 6}}}}&{{{\boldsymbol{1}}_{6({n^\mathrm{d}_e} - 1) \times 6({n^\mathrm{d}_e} - 1)}}}
		\end{array}} \right],\,\,\mathrm{if}\,\,e=1,\,\,&\text{(Case 1)}\\\\
		\left[\renewcommand\arraystretch{1.5}{\begin{array}{*{20}{c}}
				{{{\boldsymbol{1}}_{6({n^\mathrm{d}_e} - 1) \times 6({n^\mathrm{d}_e} - 1)}}}&{{{{\boldsymbol{0}}_{6({n^\mathrm{d}_e} - 1) \times 6}}}}\\
				{{{{\boldsymbol{0}}_{6 \times 6({n^\mathrm{d}_e} - 1)}}}}&{{{{\boldsymbol{\Xi }}}^{*}_{{n^\mathrm{d}_\mathrm{cp}}}}}
		\end{array}} \right],\,\,\mathrm{if}\,\,e=n_\mathrm{el},\,\,&\text{(Case 2)}
	\end{cases}
\end{equation}
and we define
\begin{equation}
	\delta {{\bf{d}}^{*e}}\coloneqq\begin{cases}\renewcommand\arraystretch{1.5}\left[{
			{\delta {\bf{d}}_1^{*\,\mathrm{T}}},
			{\delta {{\bf{d}}_2^{e\,\mathrm{T}}}},
			\cdots,
			{\delta {{\bf{d}}_{{n^\mathrm{d}_e}-1}^{e\,\mathrm{T}}}},
			{\delta {{\bf{d}}_{{n^\mathrm{d}_e}}^{e\,\mathrm{T}}}}} \right]^\mathrm{T},\,\,\mathrm{if}\,\,e=1,\,\,&\text{(Case 1)}\\\\
		\left[{{\delta {{\bf{d}}_1^{e\,\mathrm{T}}}},
			{\delta {{\bf{d}}_2^{e\,\mathrm{T}}}},
			\cdots,
			{\delta {\bf{d}}_{{n^\mathrm{d}_{{e}}}-1}^{e\,\mathrm{T}}},
			{\delta {\bf{d}}_{{n^\mathrm{d}_{\mathrm{cp}}}}^{*\,\mathrm{T}}}} \right]^\mathrm{T},\,\,\mathrm{if}\,\,e=n_\mathrm{el}.\,\,&\text{(Case 2)}
	\end{cases}
\end{equation}
For convenience, we further define
\begin{equation}
	{{\bf{\tilde \Xi }}^e} \coloneqq \left[ {\begin{array}{*{20}{c}}
			{{{\bf{1}}_{3{n_e} \times 3{n_e}}}}&{{{\bf{0}}_{3{n_e} \times 6n_e^\mathrm{d}}}}\\
			{{{\bf{0}}_{6n_e^\mathrm{d} \times 3{n_e}}}}&{{{\bf{\Xi }}^{e}}}
	\end{array}} \right],\,\,\mathrm{and}\,\,\delta {{\bf{\tilde y}}^e} \coloneqq \left\{ {\begin{array}{*{20}{c}}
	{\delta {{\boldsymbol{\varphi }}^e}}\\
	{\delta {{\bf{d}}^{ * e}}}
\end{array}} \right\},
\end{equation}
such that
\begin{equation}
	\label{del_dir_delta_y_e}
	\delta {{\bf{y}}^e} = {{\bf{\tilde \Xi }}^e}\delta {{\bf{\tilde y}}^e}.
\end{equation}
The same transformation applies to the increment,
\begin{equation}
	\label{inc_dir_delta_y_e}
	\Delta {{\bf{y}}^e} = {{\bf{\tilde \Xi }}^e}\Delta {{\bf{\tilde y}}^e}.
\end{equation}
For an element that contains a junction at the end ($e\in\mathcal{J}$), the internal load vector of the element is transformed, as
\begin{equation}
	\label{trans_elem_int_f_vec}
	\delta {{\bf{y}}^e}^{\rm{T}}{\bf{\bar F}}_{{\mathop{\rm int}} }^e = \delta {{\bf{\tilde y}}^{e}}^{\rm{T}}{{\bf{\tilde F}}}_{{\mathop{\rm int}} }^{e}\,\,\mathrm{with}\,\,{{\bf{\tilde F}}}_{{\mathop{\rm int}} }^{e}\coloneqq{{\bf{\tilde \Xi }}^{e\mathrm{T}}}{\bf{\bar F}}_{{\mathop{\rm int}} }^e,
\end{equation}
where ${{\bf{\bar F}}}_{{\mathop{\rm int}} }^{e}$ is given in Eq.\,(\ref{glob_resid_vec_fe_int}). In the same way, the elemental external load vector is also transformed by
\begin{equation}
	\label{trans_elem_ext_f_vec}
	\delta {{\bf{y}}^e}^{\rm{T}}{\bf{F}}_{{\mathop{\rm ext}} }^e = \delta {{\bf{\tilde y}}^{e}}^{\rm{T}}{{\bf{\tilde F}}}_{{\mathop{\rm ext}} }^{e}\,\,\mathrm{with}\,\,{{\bf{\tilde F}}}_{{\mathop{\rm ext}} }^{e}\coloneqq{{\bf{\tilde \Xi }}^{e\mathrm{T}}}{\bf{F}}_{{\mathop{\rm ext}} }^e.
\end{equation}
\subsubsection{Linearization}
The transformation of Eq.\,(\ref{del_trans_y_y_td}) also contributes to the geometric tangent stiffness, due to its dependence on the current directors at the junction. Applying a perturbation to the transformation matrix, the term $\delta{\bf{y}}^{e\,\mathrm{T}}{\bf{f}}^e_\mathrm{y}$ in Eq.\,(\ref{spat_disc_int_vw_hw}) can be rewritten, as
\begin{equation}
	\label{ptb_int_f_trans}
	\delta {\bf{y}}{_\varepsilon ^{e\,{\rm{T}}}}{{\bf{f}}_{{\mathrm{y}}}^e} = \delta {{\bf{\tilde y}}^{e\,\rm{T}}}{\bf{\tilde \Xi }}_\varepsilon ^{e{\rm{T}}}{{\bf{f}}_{{\mathrm{y}}}^e},\,\,\varepsilon\in\Bbb{R},\,\,e\in\mathcal{J}.
\end{equation}
Then, the directional derivative is obtained by
\begin{equation}
	\label{dir_deriv_gstiff_trans}
	\dfrac{\mathrm{d}}{{\mathrm{d}\varepsilon }}{\left. {\left( {\delta {{{\bf{\tilde y}}}^{e}}^{\rm{T}}{\bf{\tilde \Xi }}_\varepsilon^{e\rm{T}}{{\bf{f}}_{{{y}}}^e}} \right)} \right|_{\varepsilon  = 0}} = \left\{ {\begin{array}{*{20}{c}}
			{\delta {{\bf{\Theta }}_K}}\\
			{\delta {\mu _{1K}}}\\
			{\delta {\mu _{2K}}}
	\end{array}} \right\}^{\!\rm{T}}{{\Bbb{K}}_K}\left\{ {\begin{array}{*{20}{c}}
	{\Delta {{\bf{\Theta }}_K}}\\
	{\Delta {\mu _{1K}}}\\
	{\Delta {\mu _{2K}}}
\end{array}} \right\}= \delta{\bf{d}}_K^{*\,\rm{T}}{{\Bbb{K}}^*_K}\Delta{\bf{d}}_K^{*},
\end{equation}
where
\begin{equation}
	\label{tstiff_trans_k_td_i}
	{\Bbb{K}}_K = {\left[\renewcommand\arraystretch{1.2} {\begin{array}{*{20}{c}}
				{{{{\bf{d}}_{\alpha K}} \otimes {\bf{\tilde m}}^\alpha_{K} - \left( {{{\bf{d}}_{\alpha K}} \cdot {\bf{\tilde m}}^{\alpha}_{K}} \right){{\bf{1}}_{3 \times 3}}}}&{{\bf{m}}_K^1}&{{\bf{m}}_K^2}\\
				{{\bf{m}}{{_K^{1\,{\rm{T}}}}}}&{{{\bf{d}}_{1K}} \cdot {\bf{\tilde m}}_K^1}&0\\
				{{\bf{m}}{{_K^{2\,\mathrm{T}}}}}&0&{{{\bf{d}}_{2K}} \cdot {\bf{\tilde m}}_K^2}
		\end{array}} \right]_{5 \times 5}},
\end{equation}
with a repeated index $\alpha\in\left\{1,2\right\}$, and no summation implied on $K\in\left\{1,n^\mathrm{d}_\mathrm{cp}\right\}$. ${\bf{\tilde m}}_K^\alpha$ denotes the director stress couple at the junction, which is simply obtained by the components, conjugate to $\delta{\bf{d}}_{\alpha K}$, in the elemental array ${\bf{f}}_y^e$. Note that $\delta{\bf{d}}_{\alpha K}\equiv \delta{\bf{d}}^e_{\alpha I}$, based on the correspondence between the indices $e$, $I$, and $K$ in each of the Cases 1 and 2, see Table\,\ref{tab_corr_indic_jct}.
We also define the moment at the junction, ${{\bf{m}}^\alpha_K}\coloneqq{{\bf{d}}_{\alpha K}}\times{\bf{\tilde m}}^\alpha_K$ (no sum on $\alpha\in\left\{1,2\right\}$). In the second equality of Eq.\,(\ref{dir_deriv_gstiff_trans}), we also introduce the fictitious \textcolor{black}{sixth} row and column of zeros, 
\begin{equation}
	{\Bbb{K}}_K^* = {\left[\renewcommand\arraystretch{1.2} {\begin{array}{*{20}{c}}
				{{{{\Bbb{K}}}_K}}&{{{\bf{0}}_{5 \times 1}}}\\
				{{{\bf{0}}_{1 \times 5}}}&0
		\end{array}} \right]_{6 \times 6}}.
\end{equation}
Then, the geometrical element tangent stiffness matrix due to the joint condition for the $e$th element ($e\in\mathcal{J}$) is obtained by
\begin{equation}
	\label{egtstiff_aug_jct}
	{{\Bbb{\tilde K}}^{e}} = \left[ {\begin{array}{*{20}{c}}
			{{{\bf{0}}_{3{n_e} \times 3{n_e}}}}&{{{\bf{0}}_{3{n_e} \times 6n_e^\mathrm{d}}}}\\
			{{{\bf{0}}_{6n_e^\mathrm{d} \times 3{n_e}}}}&{{{\Bbb{K}}^{e}}}
	\end{array}} \right],
\end{equation}
with
\begin{align}
	{\Bbb{K}}^{e} \coloneqq \begin{cases}\left[\renewcommand\arraystretch{1.5} {\begin{array}{*{20}{c}}
				{{\Bbb{K}}_1^*}&{{{\bf{0}}_{6 \times 6({n^\mathrm{d}_e} - 1)}}}\\
				{{{\bf{0}}_{6({n^\mathrm{d}_e} - 1) \times 6}}}&{{{\bf{0}}_{6({n^\mathrm{d}_e} - 1) \times 6({n^\mathrm{d}_e} - 1)}}}
		\end{array}} \right],\,\,\mathrm{if}\,\,e=1,\,\,&\text{(Case 1)}\\\\
		\left[\renewcommand\arraystretch{1.5} {\begin{array}{*{20}{c}}
				{{{\bf{0}}_{6({n^\mathrm{d}_e} - 1) \times 6({n^\mathrm{d}_e} - 1)}}}&{{{\bf{0}}_{6({n^\mathrm{d}_e} - 1) \times 6}}}\\
				{{{\bf{0}}_{6 \times 6({n^\mathrm{d}_e} - 1)}}}&{{\Bbb{K}}_{{n^\mathrm{d}_\mathrm{cp}}}^*}
		\end{array}} \right],\,\,\mathrm{if}\,\,e=n_\mathrm{el}.\,\,&\text{(Case 2)}
	\end{cases}
\end{align}  
Applying the transformation of Eqs.\,(\ref{del_dir_delta_y_e}) and (\ref{inc_dir_delta_y_e}), and adding the additional geometric tangent stiffness of Eq.\,(\ref{egtstiff_aug_jct}), we obtain the tangent stiffness matrix of the element containing the joint, as
\begin{equation}
	\label{elem_K_jct_end_elem}
	{\Bbb{\bar K}}_{{\mathop{\rm int}} }^{e} \coloneqq {{\boldsymbol{\tilde \Xi}}^{e}}^{\rm{T}}{\bf{\bar K}}_{{\mathop{\rm int}} }^e\,{\boldsymbol{\tilde \Xi} ^{e}} + {{\tilde{\Bbb{K}}}^{e}},\,\,e\in\mathcal{J},
\end{equation}
where ${\bf{\bar K}}_{{\mathop{\rm int}} }^e$ is given in Eq.\,(\ref{glob_tstiff_cond_elem}).
Then, considering the joint condition, Eq.\,(\ref{eq_bef_app_dbdc_kd_f}) can be rewritten as
\begin{equation}
	\label{eq_bef_app_dbdc_kd_f_jct}
	\delta {{\bf{\tilde y}}^{{\rm{T}}}}{\Bbb{\bar K}}\,\Delta {{\bf{\tilde y}}} = \delta {{\bf{\tilde y}}^{{\rm{T}}}}{\Bbb{\bar F}},
\end{equation}
where
\begin{equation}
	{{\Bbb{\bar K}}} \coloneqq \mathop {\mathlarger{\mathlarger{{\bf{A}}}}}\limits_{e \not \in {{\mathcal{J}}}} {{\bf{\bar K}}^e_\mathrm{int}} + \mathop {\mathlarger{\mathlarger{{\bf{A}}}}}\limits_{e \in {{\mathcal{J}}}} {{\Bbb{\bar K}}^{e}_\mathrm{int}},
\end{equation}
and
\begin{equation}
	{\Bbb{\bar F}}\coloneqq \mathop {\mathlarger{\mathlarger{{\bf{A}}}}}\limits_{e \not \in {{\mathcal{J}}}} \left( {{\bf{F}}^e_{\rm{ext}} - {\bf{\bar F}}^e_{{\mathop{\rm int}} }} \right) + \mathop {\mathlarger{\mathlarger{{\bf{A}}}}}\limits_{e \in {{\mathcal{J}}}} \left( {{\bf{\tilde F}}^{e}_{\rm{ext}} - {\bf{\tilde F}}^{e}_{{\mathop{\rm int}} }} \right) + {\mathlarger{\mathlarger{{\bf{A}}}}}{\left[ {{{{\boldsymbol{\bar R}}}_0}} \right]_{s\in{\Gamma _{\text{N}}}}}.
\end{equation}
After removing the rows and columns corresponding to the constrained DOFs from the \textcolor{black}{Dirichlet (displacement) boundary conditions}, and the fictitious DOFs (see Remark \ref{fict_dof_trans_rem}) in the joint conditions, we obtain the system of linear equations 
\begin{equation}
	\label{red_kd_r_eq_trans}
	{{{{\Bbb{\bar K}}}_{\mathrm{r}}}}\,\Delta {{{\bf{\tilde y}}}_{\mathrm{r}}} = {{{{\Bbb{\bar F}}}_{\mathrm{r}}}},
\end{equation}
where the subscript $(\bullet)_\mathrm{r}$ in Eq.\,(\ref{red_kd_r_eq_trans}) denotes the \textit{reduced} vector or matrix.   
\begin{remark} \textit{Symmetry of the tangent stiffness matrix at equilibrium}. It appears that the additional geometric tangent stiffness matrix due to the joint condition is generally non-symmetric due to the term ${{\bf{d}}_{\alpha K}} \otimes {\bf{\tilde m}}^\alpha_K$. It has been discussed in \citet[Sections 4.1 and 4.2]{simo1986three} that the presence of a rotation group, which is a nonlinear manifold, in configuration space leads to the non-symmetry; however, they also proved that the skew-symmetric part of the geometric tangent stiffness operator vanishes at the equilibrium configuration. In a similar manner, we further examine the skew-symmetric part of the tangent stiffness matrix, i.e., \textcolor{black}{${\Bbb{K}}_K^{*\rm{A}} \coloneqq \frac{1}{2}\left( {{{{\Bbb{K}}}^*_K} - {{{\Bbb{K}}}_K^{*\rm{T}}}} \right)$}, which gives
	\begin{align}
		\label{skew_part_k_a}
		{{\delta {{\bf{d}}^*_K} \cdot {{{\Bbb{K}}}_K^{*\rm{A}}}\Delta {{\bf{d}}^*_K}}} = \delta {{\bf{\Theta }}_K} \cdot {\dfrac{1}{2}}\left( {{{\bf{d}}_{\alpha K}} \otimes {\bf{\tilde m}}^\alpha_K - {\bf{\tilde m}}^\alpha_K \otimes {{\bf{d}}_{\alpha K}}} \right)\Delta {{\bf{\Theta }}_K}= {\dfrac{1}{2}}{{\bf{m}}^\alpha_K \cdot\delta {{\bf{\Theta }}_K} \times \Delta {{\bf{\Theta }}_K}}.
	\end{align}
	For multiple beams (NURBS patches) connected to a joint with rotational continuity condition ${\left( {\delta {{\bf{\Theta }}_K}} \right)_1} = \cdots = {\left( {\delta {{\bf{\Theta }}_K}} \right)_N} \equiv \delta {\bf{\Theta }}$ and ${\left( {\Delta {{\bf{\Theta }}_K}} \right)_1} = \cdots = {\left( {\Delta {{\bf{\Theta }}_K}} \right)_N} \equiv \Delta {\bf{\Theta }}$, where the subscript $(\bullet)_i$ denotes that the quantity belongs to $i$-th constituent beam (patch), and $N$ denotes the total number of connected beams. Then, we can rewrite Eq.\,(\ref{skew_part_k_a}) as
	\begin{equation}
		\sum\limits_{i = 1}^N {{{\left( {\delta {{\bf{d}}^*_K} \cdot {\Bbb{K}}_K^{*\rm{A}}\Delta {{\bf{d}}^*_K}} \right)}_i}}  = \delta {\bf{\Theta }} \times\Delta {\bf{\Theta }}\cdot \dfrac{1}{2}\sum\limits_{i = 1}^N {{{\left( {{\bf{m}}^\alpha_K} \right)}_i}}.
	\end{equation}
	By the moment equilibrium at the joint,
	\begin{equation}
		\sum\limits_{i = 1}^N {{{\left( {{\bf{m}}_K^\alpha} \right)}_i}}  = {\boldsymbol{0}}.
	\end{equation}
	Therefore, the skew-symmetric part of the tangent stiffness matrix vanishes at equilibrium. That is, in the \textit{equilibrium configuration}, the symmetry of the global tangent stiffness matrix solely depends on whether the external loading is conservative. In \textit{non-equilibrium configurations}, the additional geometrical tangent stiffness due to the rotational continuity conditions leads to an unsymmetric tangent stiffness matrix.
\end{remark}
\subsubsection{Configuration update procedure}
For an end control point under the joint condition, the update process in Eq.\,(\ref{additive_update_dir_coef}) should be replaced by a \textit{multiplicative} one for an \textit{exact} update. Using Eq.\,(\ref{mdec_dir_alph}), the end directors are decomposed into
\begin{equation}
	{}^{n + 1}{\bf{d}}_{\alpha K}^{(i)} = {}^{n + 1}\lambda _{\alpha K}^{(i)}\,{}^{n + 1}{\bf{t}}_{\alpha K}^{(i)}.	
\end{equation}
The finite rotation of the unit director, due to the (finite) incremental rotation ${\Delta\boldsymbol{\Theta }}$, can be expressed by \citep{argyris1982excursion}
\begin{align}
	{}^{n + 1}{\bf{t}}_{{\alpha}K}^{(i)} = \exp \left[ {\widehat {\Delta {{\boldsymbol{\Theta }}}}} \right]{}^{n + 1}{\bf{t}}_{{\alpha}K} ^{(i - 1)},\,\,{}^{n + 1}{{\bf{t}}_{{\alpha}K}^{(0)}} \equiv {}^n{{\bf{t}}_{{\alpha}K}},\label{upd_cfg_tv}
\end{align}
where
\begin{equation}
	\label{erod_formula_expmap}
	\exp \left[ {\widehat {\Delta \boldsymbol{\Theta }}} \right] = {\boldsymbol{1}} + \dfrac{{\sin \Theta }}{\Theta }\,\widehat {\Delta\boldsymbol{\Theta }} + \dfrac{{1 - \cos{\Theta }}}{{{\Theta ^2}}}\,{\widehat {\Delta \boldsymbol{\Theta }}^2},\,\,\mathrm{if}\,\,\Theta\coloneqq{\left\| \Delta{\boldsymbol{\Theta }} \right\|}\ne0,
\end{equation}
and it becomes the identity tensor, if $\Theta = 0$. The stretch ratio is also updated by the (finite) increment of the logarithmic stretch \textcolor{black}{ratio} $\Delta \mu _{\alpha K}$, as
\begin{align}
	{}^{n + 1}\lambda _{\alpha K}^{(i)} = \exp \left[ {{}^{n + 1}\mu _{\alpha K}^{(i)}} \right], \label{upd_cfg_mu}	
\end{align}
with
\begin{equation}
	{}^{n + 1}\mu _{\alpha K}^{(i)} = {}^{n + 1}\mu _{\alpha K}^{(i - 1)} + \Delta \mu _{\alpha K},\,\,{}^{n + 1}\mu _{\alpha K}^{(0)} \equiv {}^{n}\mu _{\alpha K}.
\end{equation}

\begin{remark} It should be noted that the presented configuration update procedure does not require any secondary storage for the cross-sectional orientation and stretch at the previous iteration step, i.e., ${}^{n + 1}{\bf{t}}_{{\alpha}K} ^{(i - 1)}$, and ${}^{n + 1}{\mu}_{{\alpha}K} ^{(i - 1)}$, respectively. This is due to the fact that they can be simply calculated using the primary control coefficients of the director displacement ${\bf{\bar d}}_K\coloneqq\left[{\bf{\bar d}}_{1K}^\mathrm{T},{\bf{\bar d}}_{2K}^\mathrm{T}\right]^\mathrm{T}$ ($K\in\left\{1,n^\mathrm{d}_\mathrm{cp}\right\}$) in Eq.\,(\ref{tot_dirv_disp_entire_patch0}), that is, the displacement of directors at the joint. 
\end{remark}
\section{Numerical examples}
\label{num_example}
\textcolor{black}{In numerical examples, we compare several approaches of approximating the physical stress resultant and strain fields, see Remark \ref{remk_gloc_app_term}.} In all the considered approaches, we use \textcolor{black}{the element-wise} full Gauss integration along the axis (``FI''). For the integration over the cross-section, we also use Gauss integration, with $3\times{3}$ quadrature points. 
\subsection{Cantilever beam under bending moment}
\begin{figure}[H]
	\centering
	\begin{subfigure}[b] {0.55\textwidth} \centering
		\includegraphics[width=\linewidth]{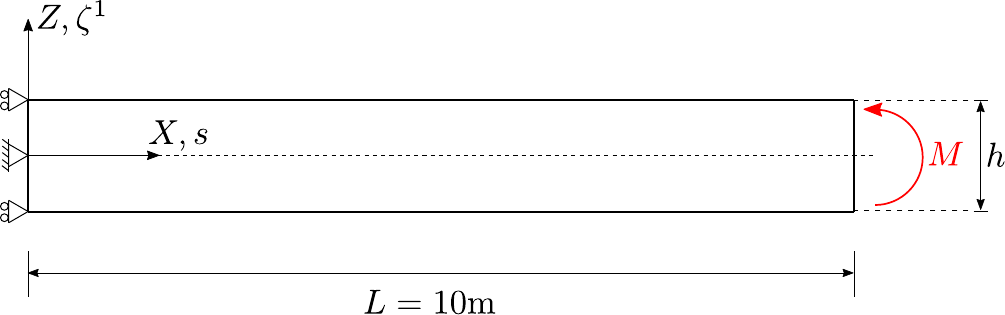}	
	\end{subfigure}			
	\caption{Cantilever beam under bending moment. Geometry and boundary conditions. Note that the thickness stretch ($\Delta {d_{13}}$) is free at $s=0$.}
	\label{ex_end_bend_mnt}
\end{figure}

We consider an initially straight beam of length $L=10\,\mathrm{m}$ with a rectangular cross-section of width $w=1\,\mathrm{m}$ and thickness $h$. A St.\,Venant--Kirchhoff type material is considered, with Young's modulus $E=1.2\times{10^7}\,\mathrm{Pa}$, and Poisson's ratio $\nu=0$. The beam's left-end is kinematically constrained, and the other end is subject to a moment load $M=2\pi{EI}/L$, where $I=wh^3/12$ denotes the second area moment of inertia of the cross-section, see Fig.\,\ref{ex_end_bend_mnt}. The moment load is applied by a distributed follower load on the cross-section at the right-end, whose expressions of the external virtual work and its increment can be found in \citet{choi2021isogeometric}. Here, we apply the total moment load in 10 load steps, with uniform increments. At the left-end, the cross-section's translation and rotation are constrained, but the transverse normal (through-the-thickness) stretching is allowed, so that 
\begin{equation}
	{\boldsymbol{\varphi }} = {\boldsymbol{\varphi }}_0,\,\,\Delta {d_{11}} = \Delta {d_{12}} = 0,\,\,\mathrm{and}\,\,\Delta {d_{13}}\,\,\text{is free},
\end{equation} 
where $\Delta {d_{\alpha i}} \coloneqq \Delta {{\boldsymbol{d}}_\alpha} \cdot {{\boldsymbol{e}}_i}$. For the given problem, we have an analytical solution for the thin beam limit (pure bending condition), such that the moment $M$ deforms the beam's axis into a circle with radius $R=EI/M$. 
%
%
\subsubsection{Alleviation of locking for $p=2$ in IGA by reducing $p_\mathrm{d}$}
\label{allev_p2_red_pd_iga_lock}
\begin{figure}[H]
	\centering
	\begin{subfigure}[b] {0.425\textwidth} \centering
		\includegraphics[width=\linewidth]{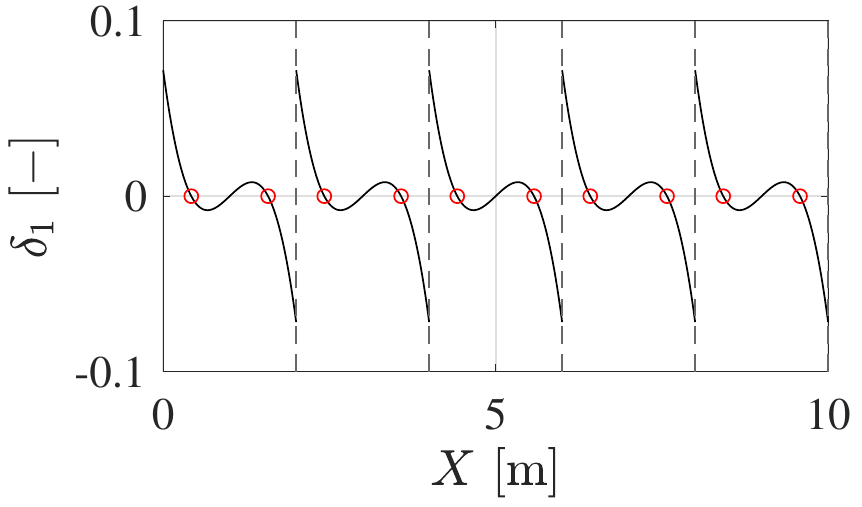}	
		\caption{FEA, $p=p_d=2$, URI}		
		\label{endmnt_tshear_fea_p2_uri}
	\end{subfigure}			
	\begin{subfigure}[b] {0.425\textwidth} \centering
		\includegraphics[width=\linewidth]{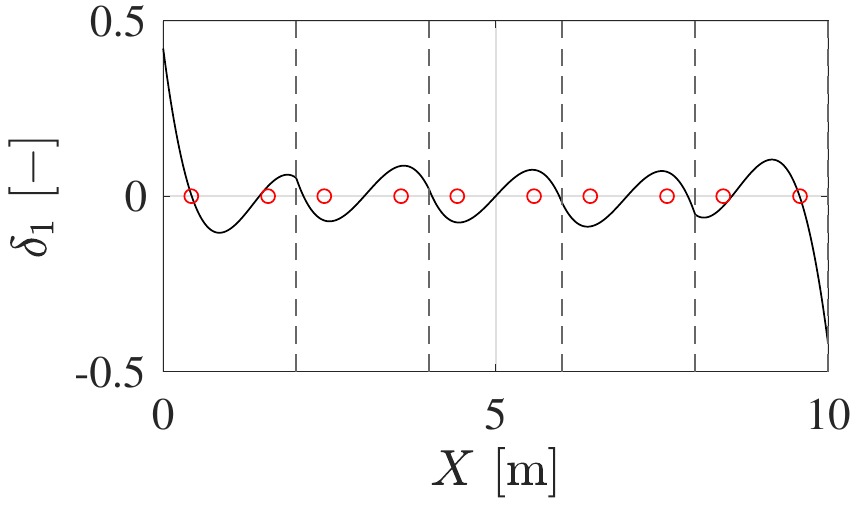}	
		\caption{IGA, $p=p_d=2$, URI}		
		\label{endmnt_tshear_iga_p2_uri}		
	\end{subfigure}		
	\caption{Cantilever beam under bending moment: Comparison of the transverse shear strain distribution along the beam's axis according to FEA and IGA. The slenderness ratio is $L/h=10^{2}$. In both cases, we use quadratic basis functions, and five elements. In every element, we use two Gauss quadrature points ($n_\mathrm{G}=2$), i.e., uniformly reduced integration (URI). The hollow circles represent the positions of the Gauss quadrature points, and the dashed vertical lines represent the element boundaries. All results are obtained from the displacement-based formulation.}
	\label{ex_end_bend_mnt_shear_strn}
\end{figure}

\begin{figure}[H]
	\centering
	\begin{subfigure}[b] {0.4875\textwidth} \centering
		\includegraphics[width=\linewidth]{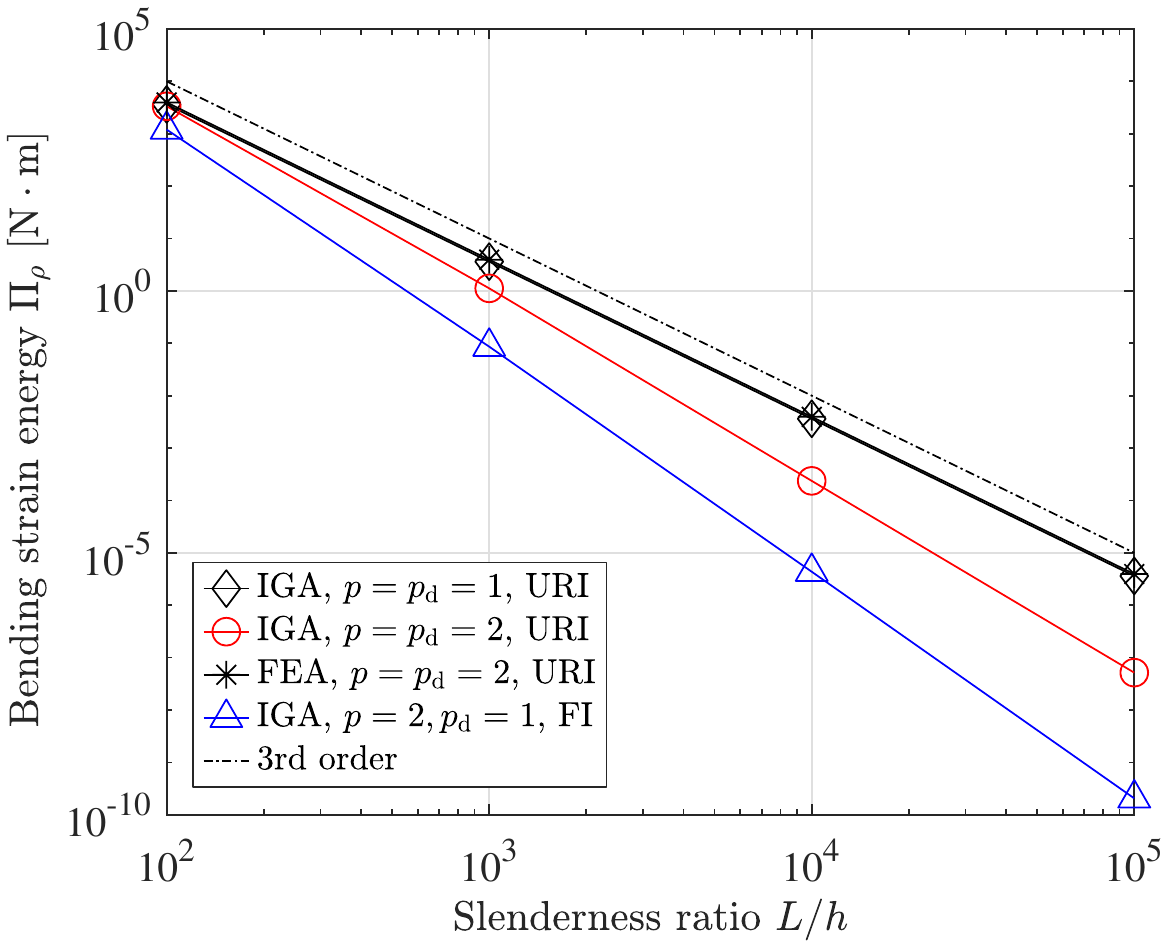}	
		\caption{Bending strain}		
		\label{pbend_strn_e_bend}		
	\end{subfigure}			
	\begin{subfigure}[b] {0.4875\textwidth} \centering
		\includegraphics[width=\linewidth]{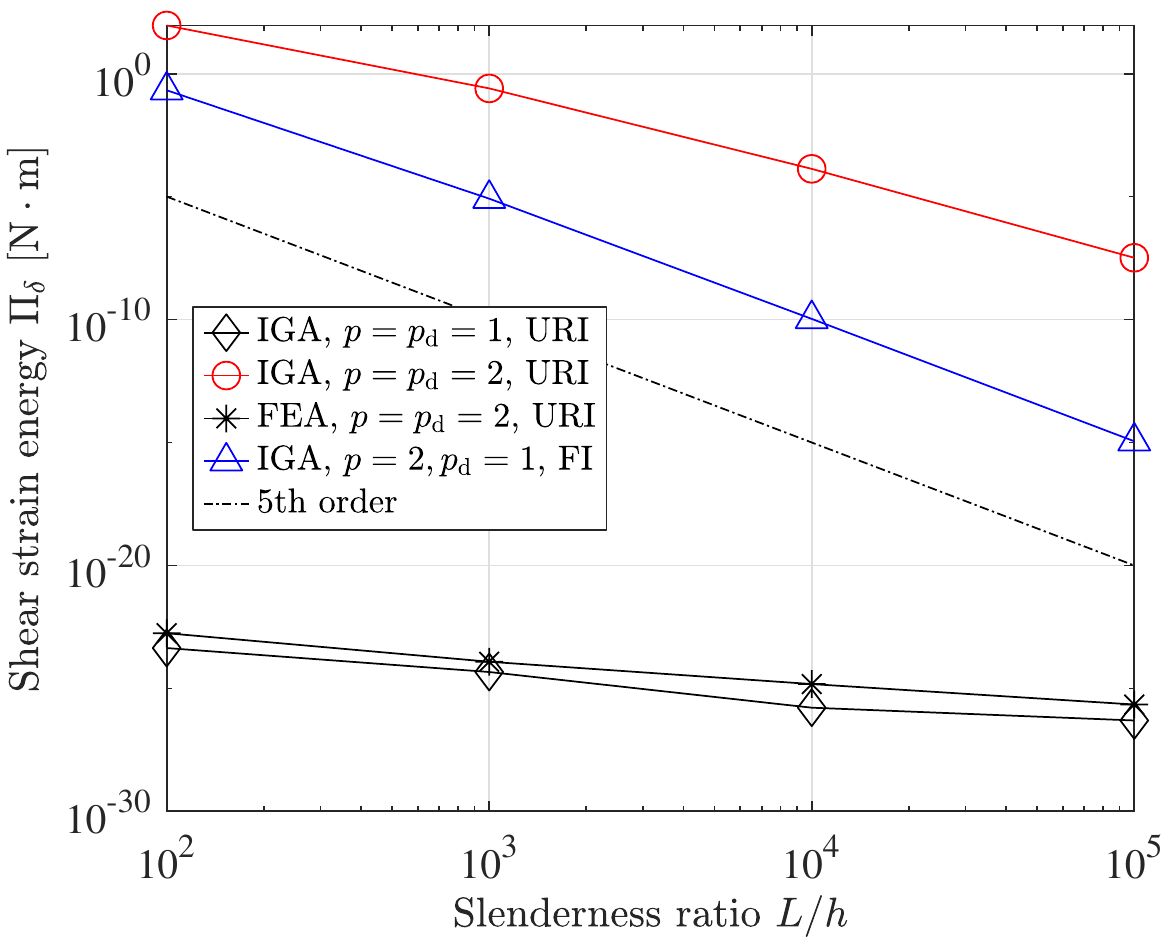}	
		\caption{Transverse shear strain}				
		\label{pbend_strn_tshear}				
	\end{subfigure}	
	\begin{subfigure}[b] {0.4875\textwidth} \centering
		\includegraphics[width=\linewidth]{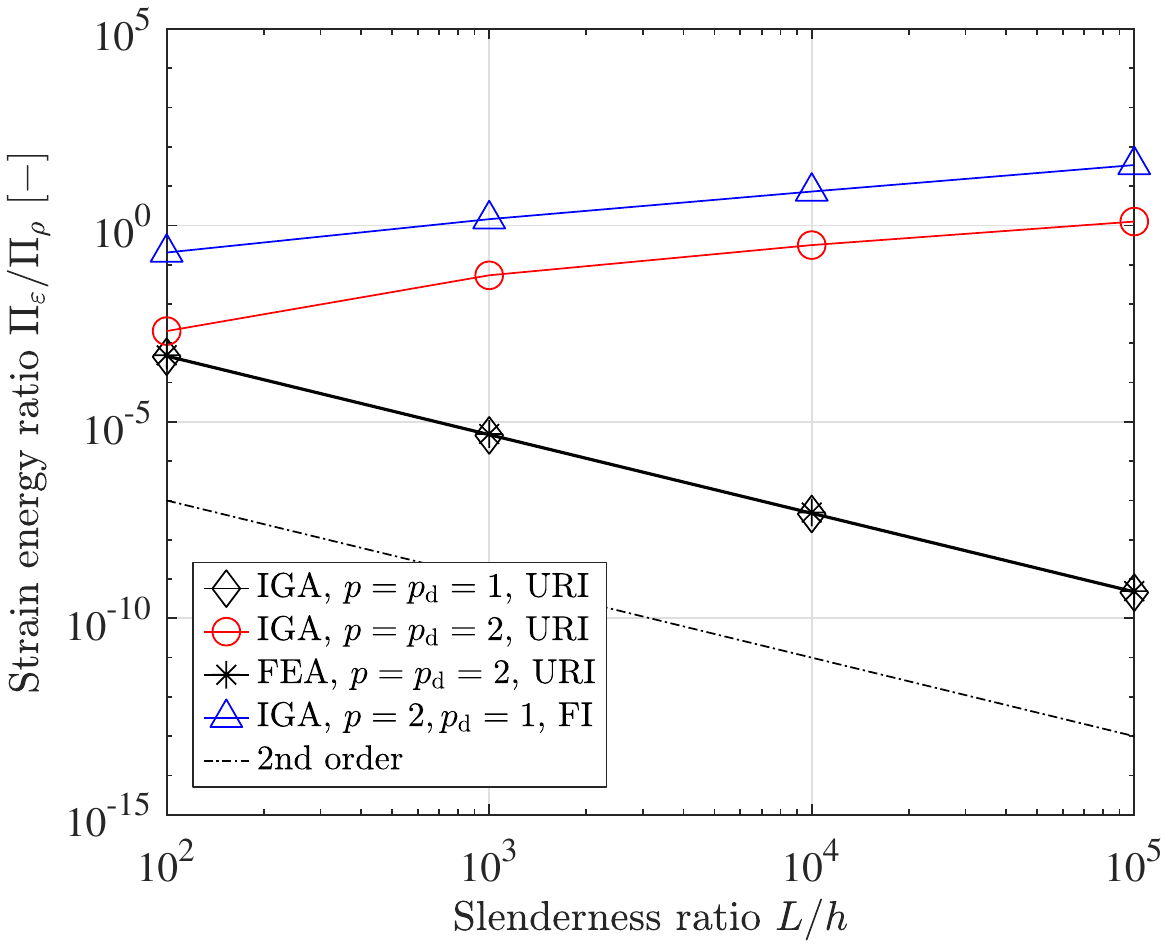}	
		\caption{Membrane strain}		
		\label{pbend_strn_mem}		
	\end{subfigure}				
	\begin{subfigure}[b] {0.4875\textwidth} \centering
		\includegraphics[width=\linewidth]{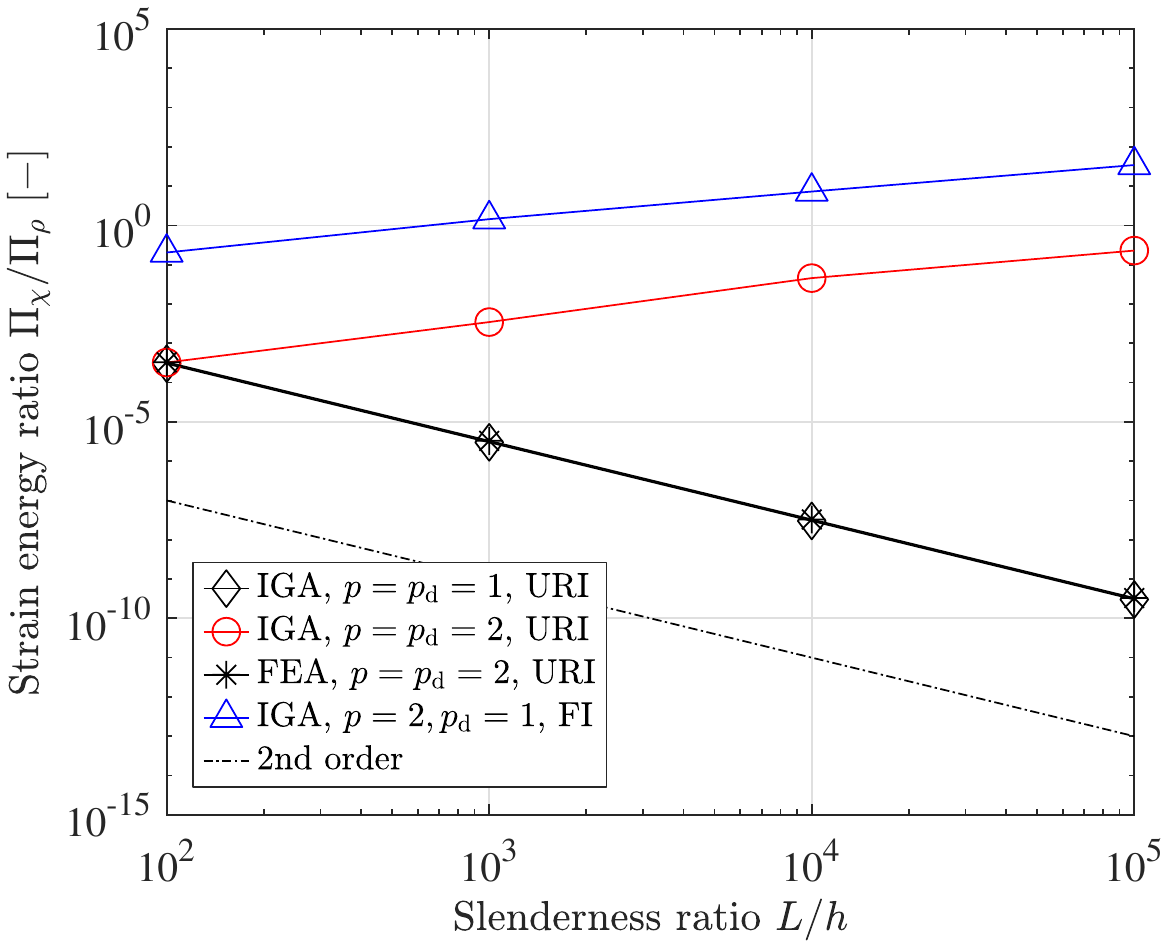}	
		\caption{Transverse normal strain}		
		\label{pbend_strn_inp}		
	\end{subfigure}	
	\caption{Cantilever beam under bending moment: Comparison of the strain energies resulting from different selection of $p_\mathrm{d}$, including the cases of using (i) the same degree of basis $p={p_\mathrm{d}}=1,2$ with the uniformly reduced integration (URI), and (ii) different degree of basis ${p}=2$ and ${p_\mathrm{d}}=1$ with three Gauss quadrature points per element (${n_\mathrm{G}}=3$), i.e., full integration (FI). Note that IGA using $p=1$ gives the same result with that from FEA using $p=1$. All results are from the displacement-based formulation, and $n_\mathrm{el}=10$.}
	\label{ex_end_bend_mnt_compare_strne_ratio_uri}
\end{figure}

Figs.\,\ref{endmnt_tshear_fea_p2_uri} and \ref{endmnt_tshear_iga_p2_uri} show the distribution of the transverse shear strain along the beam's axis resulting from FEA and IGA, respectively. In both cases, we use \textcolor{black}{the element-wise} uniformly reduced integration (URI). The FEA result shows that the positions of the Gauss quadrature points in every element coincide with the zeros of the transverse shear strain, which is discontinuous across the elements. Thus, transverse shear locking can be alleviated by URI in FEA. However, in the IGA results, it is seen that the continuous strain distribution vanishes in other locations than the Gauss quadrature points, \textcolor{black}{except at the first and last quadrature points}. This means that an element-wise \textcolor{black}{uniformly} reduced integration is not an effective tool for IGA to alleviate transverse shear locking\textcolor{black}{, see \citet{adam2014improved} for a more detailed discussion}. Transverse shear locking can be explained by the field-inconsistency paradigm \citep{prathap2013finite}. To circumvent this inconsistency, we may consider to use one degree lower basis functions for the director field. In Fig.\,\ref{ex_end_bend_mnt_compare_strne_ratio_uri}, we investigate the effectiveness of using $p_\mathrm{d} = p - 1$ to alleviate the transverse shear locking. Note that it is combined with \textcolor{black}{the element-wise} full Gauss integration along the axis (``FI'').

\noindent We first recall the expressions of the bending, membrane (axial), and transverse normal (through-the-thickness) strain energies in \citet{choi2021isogeometric},
\begin{equation}
	\label{theo_bend_strn_e}
	{\Pi _\rho} \coloneqq \int_0^L {{{\tilde m}^1}{\rho_1}\,{\rm{d}}s}\sim{h^3},
\end{equation}
\begin{equation}
	\label{theo_memb_strn_e}
	{\Pi _\varepsilon } \coloneqq \int_0^L {\tilde n\,\varepsilon \,{\rm{d}}s}\sim{h^5},
\end{equation}
and
\begin{align}
	\label{the_strn_e_inp_cs}
	{\Pi _\chi } \coloneqq \int_0^L {{{\tilde l}^{11}}{\chi _{11}}\,{\rm{d}}s}\sim{h^5},
\end{align}
respectively. Therefore, the ratios $\Pi_{\varepsilon}/{\Pi_\rho}$ and $\Pi_{\chi}/{\Pi_\rho}$ should decrease quadratically with decreasing the initial cross-sectional thickness $h$. \textcolor{black}{Further, we investigate the transverse shear strain energy, defined by
\begin{equation}
	\label{the_strn_e_tshear}
	{\Pi _\delta } \coloneqq \int_0^L {{{\tilde q}^{1}}{\delta _{1}}\,{\rm{d}}s}.	
\end{equation} 
}
\noindent In the calculation of the strain energies in Eqs.\,(\ref{theo_bend_strn_e})-(\ref{the_strn_e_tshear}), we use 
\begin{itemize}
	\item compatible strains in Eq.\,(\ref{beam_th_strn_comp_basis_d123}), with the approximated displacements, for the \textbf{displacement-based} formulation,
	\item approximated physical stress resultants in Eq.\,(\ref{approx_phy_rp}), and physical strains in Eq.\,(\ref{disc_pstrn}) for the \textbf{mixed} formulation.
\end{itemize}

\begin{figure}[htb]
	\centering
	\begin{subfigure}[b] {0.4875\textwidth} \centering
		\includegraphics[width=\linewidth]{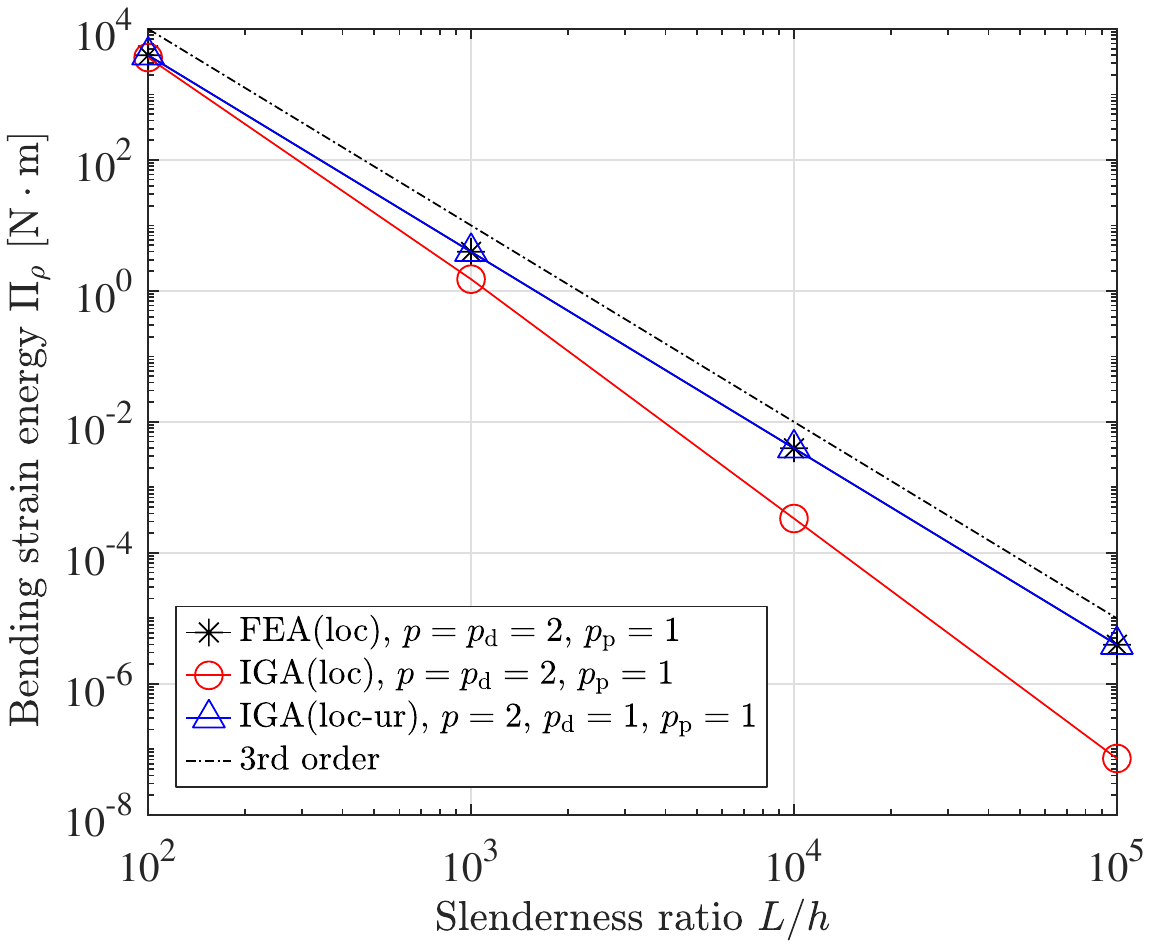}	
		\caption{Bending strain}		
		\label{pbend_strn_e_p2_mix_bend}		
	\end{subfigure}			
	\begin{subfigure}[b] {0.4875\textwidth} \centering
		\includegraphics[width=\linewidth]{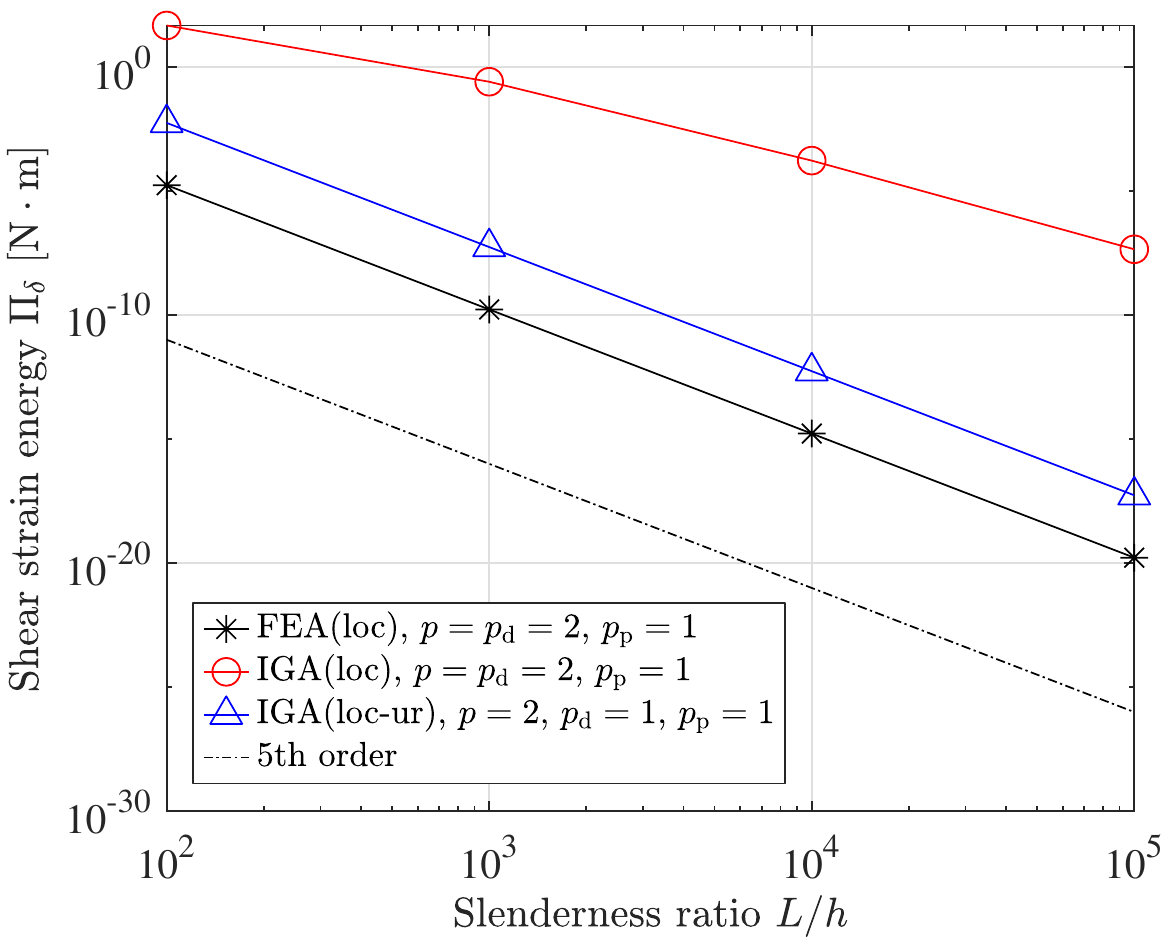}	
		\caption{Transverse shear strain}				
		\label{pbend_strn_e_p2_mix_tshear}				
	\end{subfigure}	
	\begin{subfigure}[b] {0.4875\textwidth} \centering
		\includegraphics[width=\linewidth]{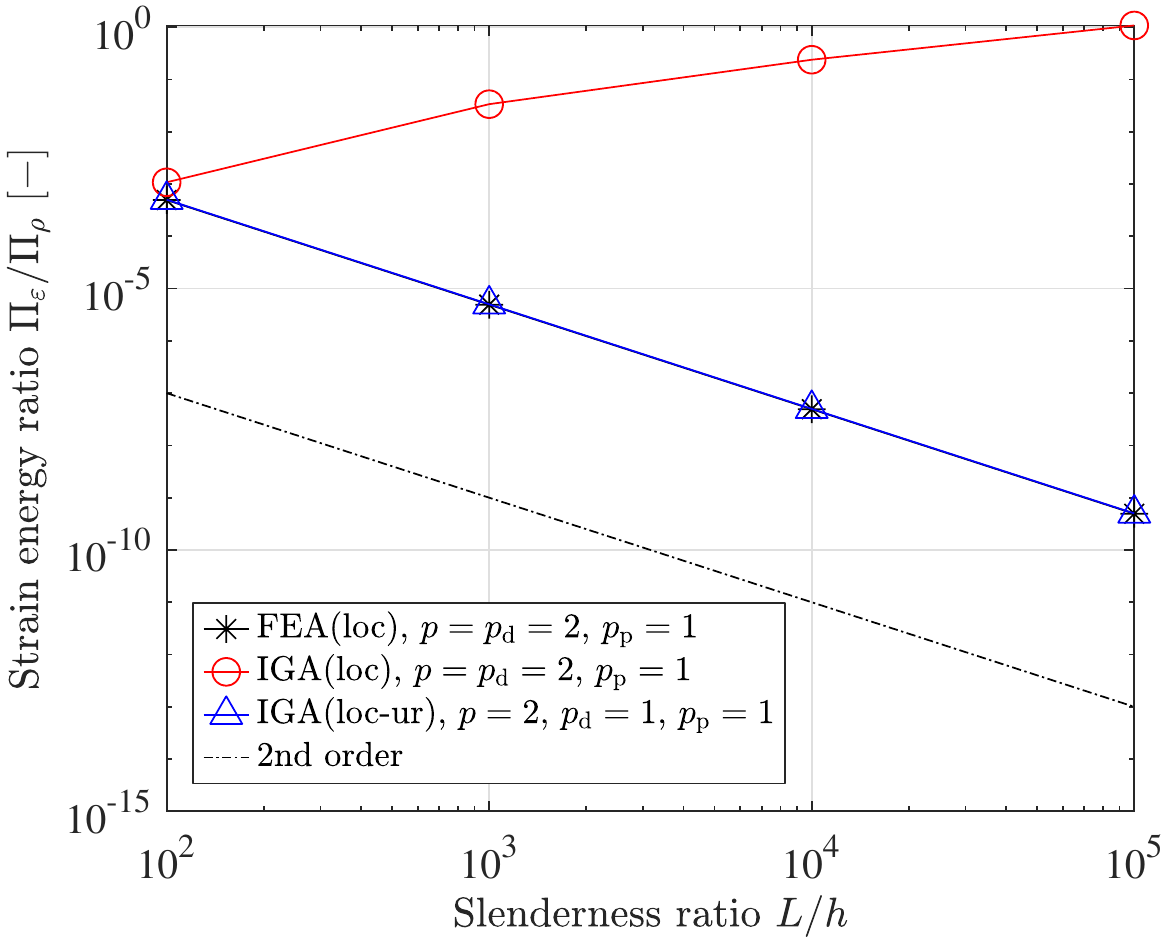}	
		\caption{Membrane strain}		
		\label{pbend_strn_e_p2_mix_memb}				
	\end{subfigure}				
	\begin{subfigure}[b] {0.4875\textwidth} \centering
		\includegraphics[width=\linewidth]{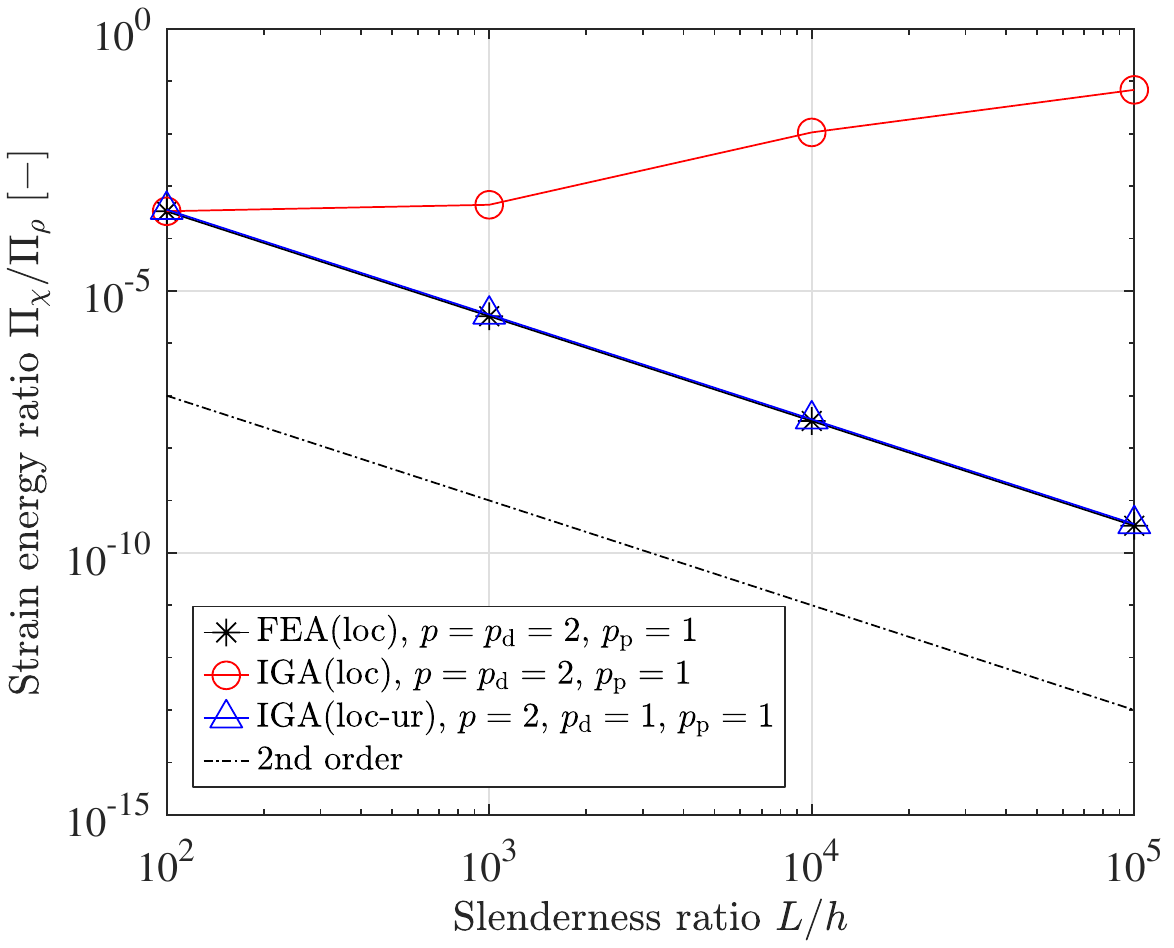}	
		\caption{Transverse normal strain}		
		\label{pbend_strn_e_p2_mix_inp_cs}		
	\end{subfigure}	
	\caption{Cantilever beam under bending moment: Comparison of the strain energies resulting from different selection of $p_\mathrm{d}$, including the cases of using (i) the same degrees of bases $p={p_\mathrm{d}}=2$, and (ii) different degrees of bases ${p}=2$ and ${p_\mathrm{d}}=1$, and $p_\mathrm{p} = 1$. All results are from the mixed formulation with $p_\mathrm{p} = 1$, combined with the local approach, full integration along the axis ($n_\mathrm{G}=3$), and $n_\mathrm{el}=10$.}
	\label{ex_end_bend_mnt_compare_strne_ratio_mix}
\end{figure}

\noindent Fig.\,\ref{pbend_strn_e_bend} shows that the FEA result with $p=2$ and URI (black curve with cross markers) exhibits a cubic rate of decrease of the bending strain energy, as expected by Eq.\,(\ref{theo_bend_strn_e}). However, the result from IGA with $p=2$ and URI (red curve) shows, as expected, deviation from the analytical rate of decrease due to locking. In IGA, only in the case of $p=1$, shows the correct decrease rate, since it has no inter-element continuity in the transverse shear and membrane strain fields, as in FEA. It is also seen in Fig.\,\ref{pbend_strn_e_bend} that IGA with ${p}=2$ and ${p_\mathrm{d}}=1$ (blue curve) suffers from severe locking. However, in Fig.\,\ref{pbend_strn_tshear}, \textcolor{black}{the reduction of ${p_\mathrm{d}}$ to 1 (blue curve) yields much smaller spurious transverse shear strain energy than that from IGA ($p=p_\mathrm{d}=2$) with URI (red curve)}. In Figs.\,\ref{pbend_strn_mem} and \ref{pbend_strn_inp}, it is seen that the IGA results show spurious increase in both of the membrane, and transverse normal strains. In contrast, the results having no higher inter-element continuity in the strain fields (black curves) show the quadratic rate of decrease in both cases. This adverse effect of higher order continuity in IGA can be also observed in results from the mixed formulation. In Fig.\,\ref{ex_end_bend_mnt_compare_strne_ratio_mix}, we investigate the decrease rate of the strain energy, obtained by using the mixed formulation, for different cases of $p_\mathrm{d}$ and the inter-element continuity of the strains. 

\noindent In Fig.\,\ref{ex_end_bend_mnt_compare_strne_ratio_mix}, it is seen that the IGA using the mixed formulation, combined with a local (element-wise) static condensation (red curve), suffers from severe locking, see spurious strain energy in Figs.\,\ref{pbend_strn_e_p2_mix_tshear}-\ref{pbend_strn_e_p2_mix_inp_cs}, and the resulting spuriously higher decrease rate of the bending strain energy in Fig.\,\ref{pbend_strn_e_p2_mix_bend}. This is due to the fact that it does not consider the locking arising from the continuity condition of the strain fields. Surprisingly, if we reduce $p_\mathrm{d}$ to 1 \textcolor{black}{(blue curve)}, the transverse shear strain vanishes at fifth order rate, as in the FEA result, and all the other spurious strain energies also vanish. 

\begin{figure}[H]
	\centering
	\begin{subfigure}[b] {0.49\textwidth} \centering
		\includegraphics[width=\linewidth]{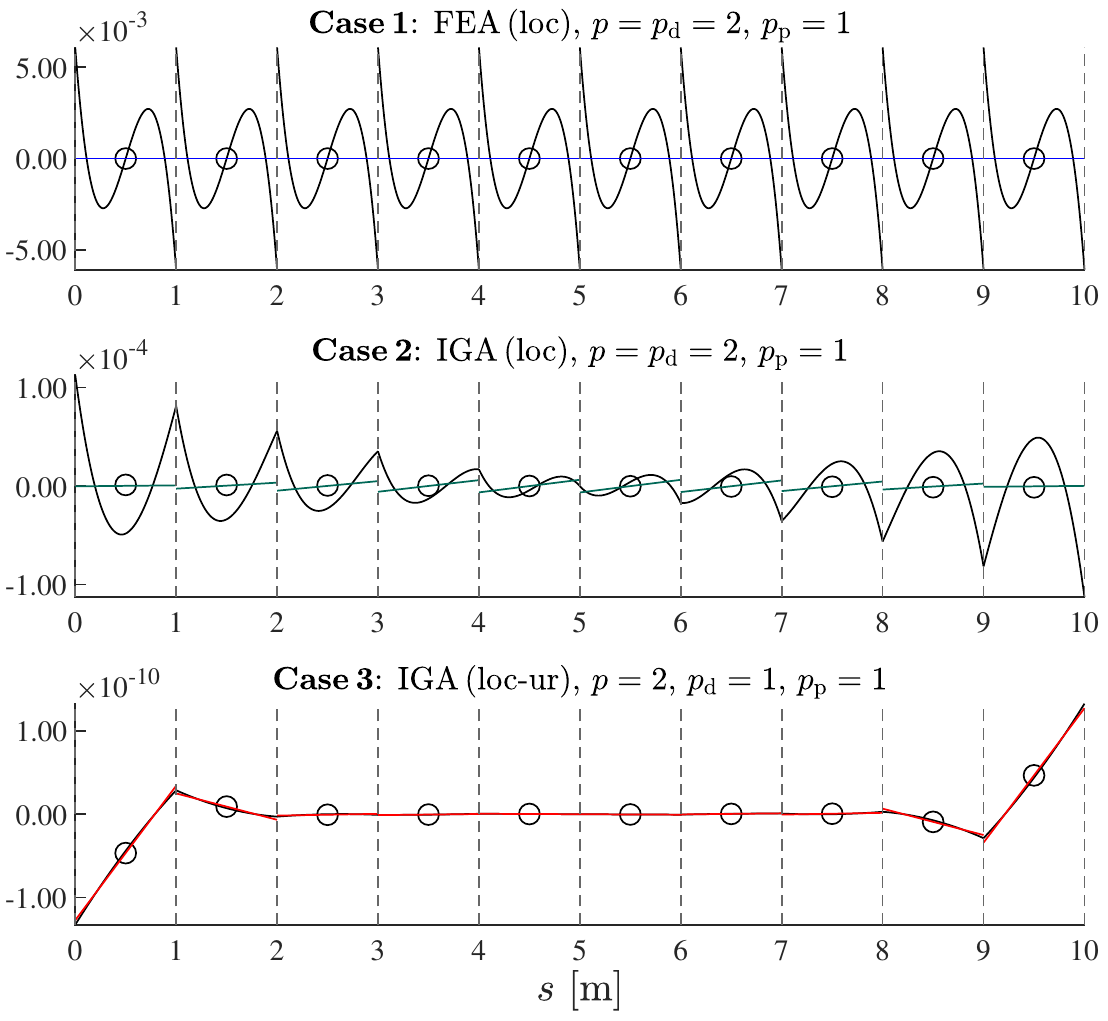}	
		\caption{Comparison of the geometrical and physical strains}	
		\label{pbend_strn_dist_tshear_geom_phy}			
	\end{subfigure}			
	\begin{subfigure}[b] {0.49\textwidth} \centering
		\includegraphics[width=\linewidth]{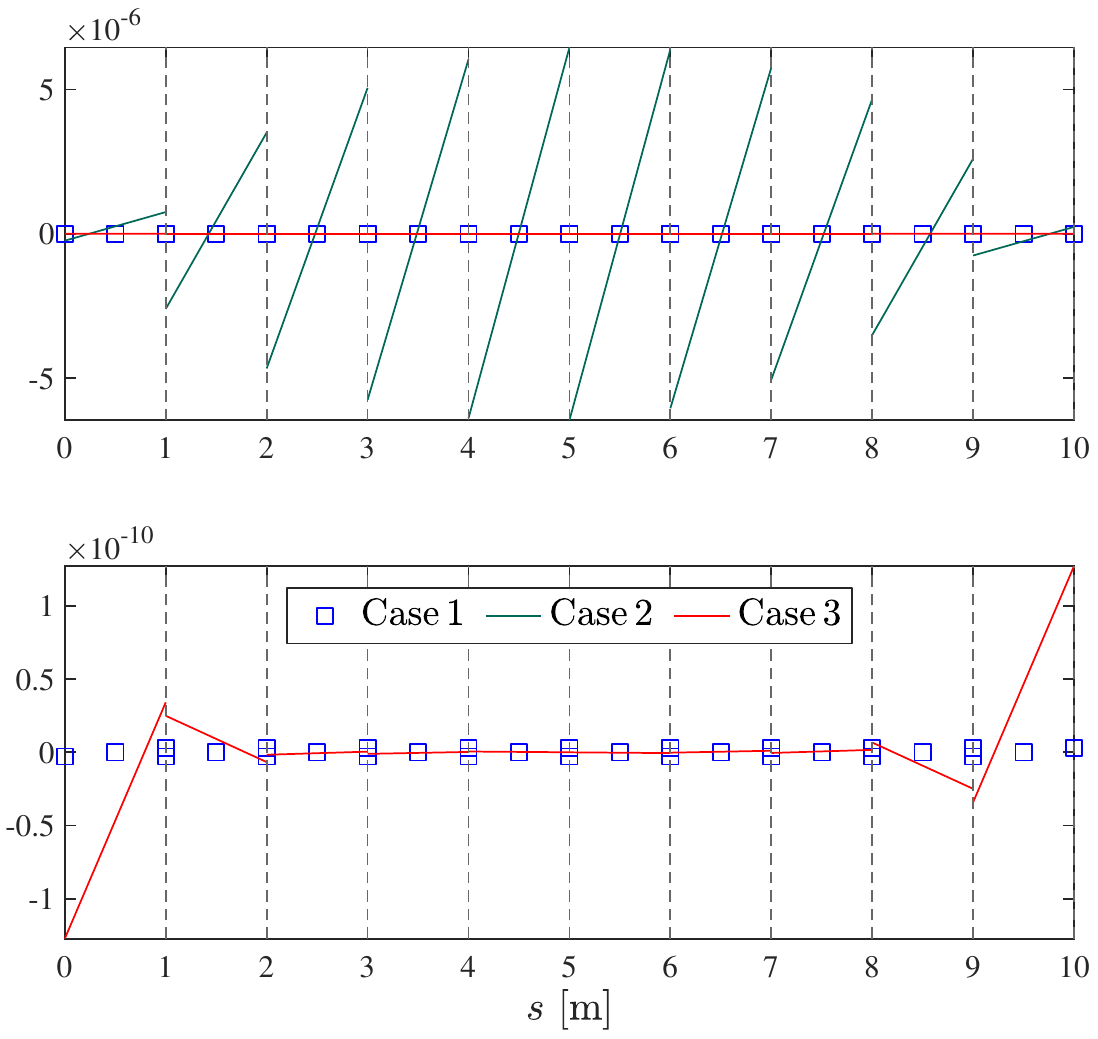}	
		\caption{Physical strains (bottom: Cases 1 and 3 only)}		
		\label{pbend_strn_dist_tshear_phy}
	\end{subfigure}	
	\caption{Cantilever beam under bending moment: Distribution of the \textbf{transverse shear} strain $\delta_1$ along the axis. (a) Comparison of the geometrical (black curves) and physical strains (colored curves) along the axis, obtained by three different formulations (Cases 1-3). In the result of Case 3, the black and red curves are overlapping. (b) Those physical strains in (a) are re-plotted in the common ordinate. In all results, we use ten elements ($n_\mathrm{el}=10$), and the selected slenderness ratio is $L/h=10^{5}$. The dashed lines represent the element boundaries, and the hollow circles in (a) represent element-wise average of the geometric strain.}
	\label{ex_end_bend_mnt_strn_dist_tshear}
\end{figure}

In order to further investigate the alleviation of locking, we compare the distribution of the axial ($\varepsilon$), transverse shear ($\delta_1$), transverse normal ($\chi_{11}$), and couple shear ($\gamma_{11}$) strains along the axis, resulting from three different finite element approximations in the mixed formulation using $p=2$. Fig.\,\ref{ex_end_bend_mnt_strn_dist_tshear} compares the distribution of the transverse shear strain. In every case shown in Fig.\,\ref{ex_end_bend_mnt_strn_dist_tshear}, the physical strain agrees very well with the element-wise average of the geometrical strain (black hollow circles), which means the compatibility condition is weakly satisfied well. For the polynomial basis functions of $p={p_\mathrm{d}}=2$, the transverse shear strain $\delta_1$ of Eq.\,(\ref{def_compat_tshear}) is a cubic polynomial. In Fig.\,\ref{pbend_strn_dist_tshear_geom_phy}, the FEA result (Case 1) clearly shows the cubic (black) curve within each element; however, in the IGA result (Case 2), the geometric strain field in each element (black curve) is very close to a quadratic function. This is attributed to the additional constraints in IGA for the inter-element $C^0$-continuity of $\delta_1$, and it eventually leads to oscillatory distribution of the physical strains to satisfy the compatibility condition. In the FEA result (Case 1), the geometric strain does not vanish, but the corresponding physical strain, using the degree of basis $p_\mathrm{p} = 1$, vanishes properly, see also the physical strain distribution (\textcolor{black}{blue} square markers) in Fig.\,\ref{pbend_strn_dist_tshear_phy}. In contrast, the IGA results with the same degrees of bases $p=p_\mathrm{d}=2$, and ${p_\mathrm{p}}=1$ (Case 2) shows a significant amount of the transverse shear strain, which leads to a severe (artificial) increase of the bending stiffness. By decreasing the degree $p_\mathrm{d}$ in Case 3, the transverse shear locking is effectively alleviated. It is noticeable that the geometrical strain also vanishes, in contrast to the FEA result using $p_\mathrm{d}=p=2$ (Case 1). 

\begin{figure}[H]
	\centering
	\begin{subfigure}[b] {0.49\textwidth} \centering
		\includegraphics[width=\linewidth]{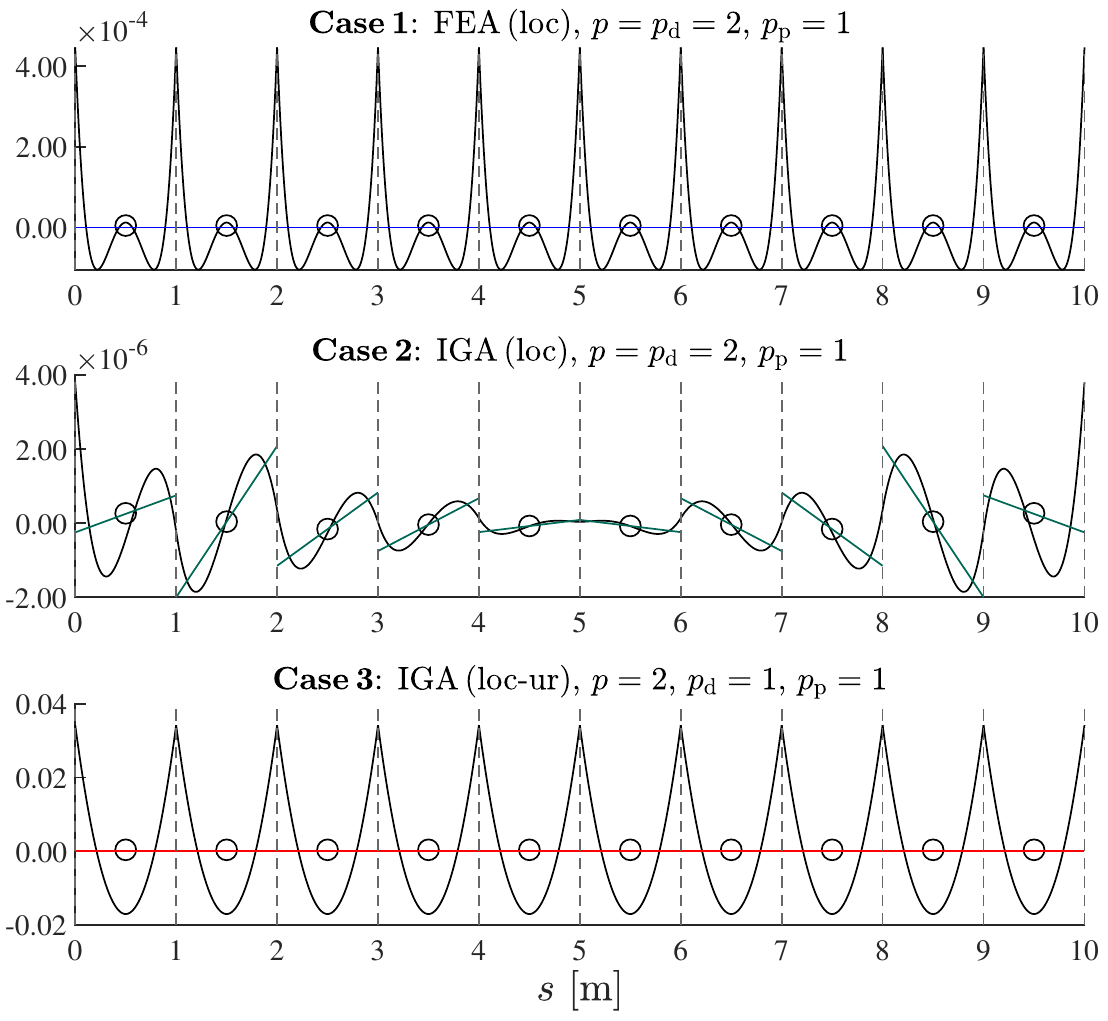}	
		\caption{Transverse normal strain}		
		\label{pbend_strn_dist_chi_geom_phy}					
	\end{subfigure}			
	\begin{subfigure}[b] {0.49\textwidth} \centering
		\includegraphics[width=\linewidth]{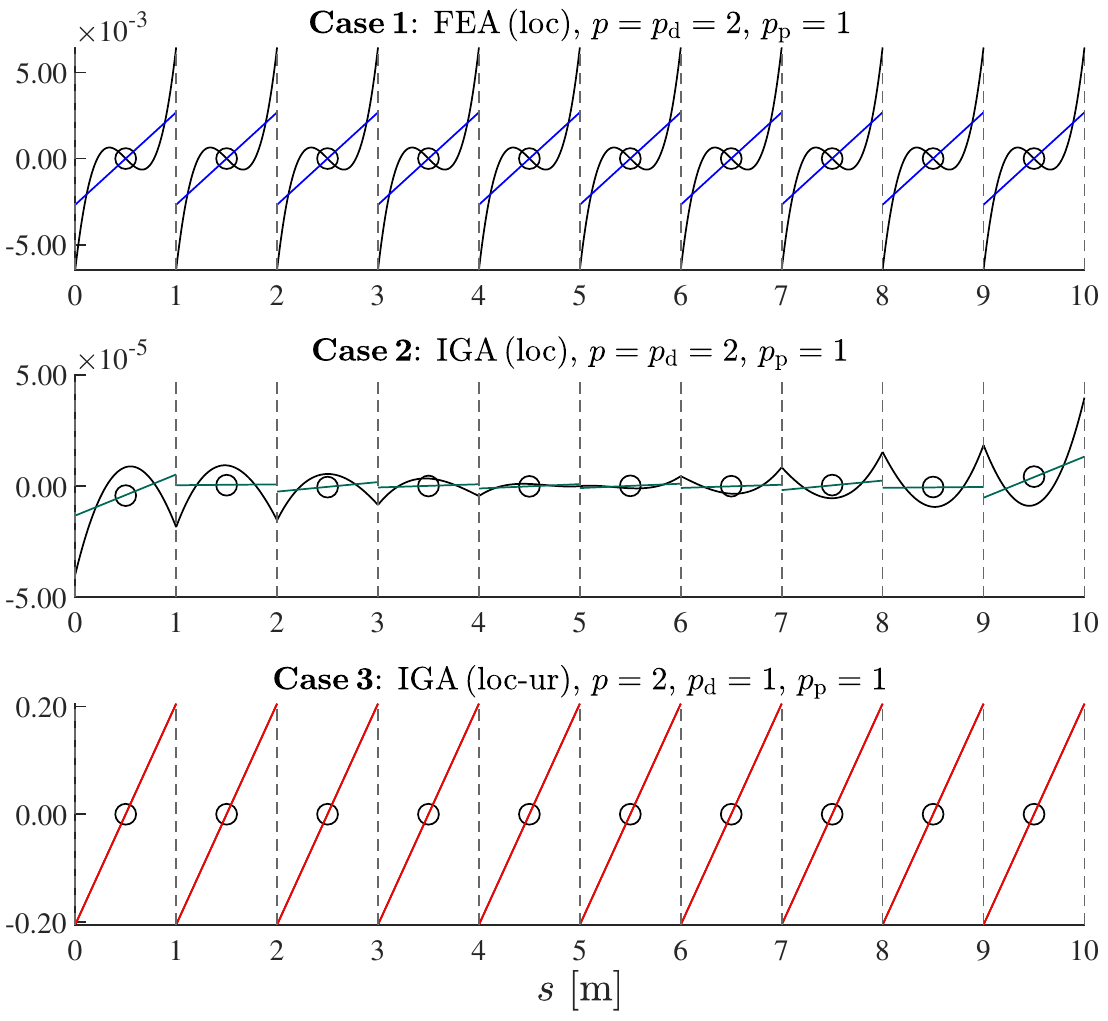}	
		\caption{Couple shear strain}		
		\label{pbend_strn_dist_gamma_geom_phy}
	\end{subfigure}	
	\caption{Cantilever beam under bending moment: Distribution of the \textbf{transverse normal} \textcolor{black}{($\chi_{11}$)}, and \textbf{couple shear} \textcolor{black}{($\gamma_{11}$)} strains along the axis. The geometrical (black curves) and physical (colored curves) strains in the three cases of the mixed formulation are compared. The hollow circles represent the element-wise average of the geometrical strain. \textcolor{black}{In the result of Case 3 in (b), the black and red curves are overlapping.} In all results, $n_\mathrm{el}=10$, and $L/h=10^{5}$. The dashed lines represent the element boundaries.}	
	\label{ex_end_bend_mnt_strn_dist_chi_gamma}
\end{figure}

\noindent Fig.\,\ref{ex_end_bend_mnt_strn_dist_chi_gamma} compares the distribution of the transverse normal strain
(through-the-thickness stretch), and the couple shear strain along the axis. For degree ${p_\mathrm{d}}=p=2$ polynomial bases of director displacements, the transverse normal strain $\chi_{11}$ in Eq.\,(\ref{strn_comp_chi}) is a quartic polynomial. Similar to the previous case of transverse shear strain distribution, it is observed in Fig.\,\ref{pbend_strn_dist_chi_geom_phy} that the geometric strain field in the result of Case 2 is closer to a cubic function within each element, but the FEA result shows a quartic distribution of the geometric strain in each element. In Case 3, we use $p_\mathrm{d}=p-1=2$, so that $\chi_{11}$ is a quadratic polynomial, see the parabolic (black) curve in the result of Case 3 in Fig.\,\ref{pbend_strn_dist_chi_geom_phy}. The geometrical couple shear strain $\gamma_{11}$ is one degree lower than $\chi_{11}$, due to $\gamma_{11}=\chi_{11,s}$, see the geometric strains (black curves) in Fig.\,\ref{pbend_strn_dist_gamma_geom_phy}. Here we observe an inconsistency between the strain fields ${\gamma _{11}}$ and $\chi _{11,s}$, i.e., ${\gamma _{11}} \ne \chi _{11,s}^{\rm{p}}$, in all the cases. However, this couple shear strain is a higher order strain than the transverse shear strain, in terms of the transverse coordinate $\zeta^{1}$ in the Green-Lagrange strain component $E_{31}$, see Eq.\,(\ref{gl_strn_mata_beam_strn}). Thus, it may not significantly affect the overall response in thin beams. Further alleviation of this inconsistency remains future work.

\begin{figure}[H]
	\centering
	\begin{subfigure}[b] {0.49\textwidth} \centering
		\includegraphics[width=\linewidth]{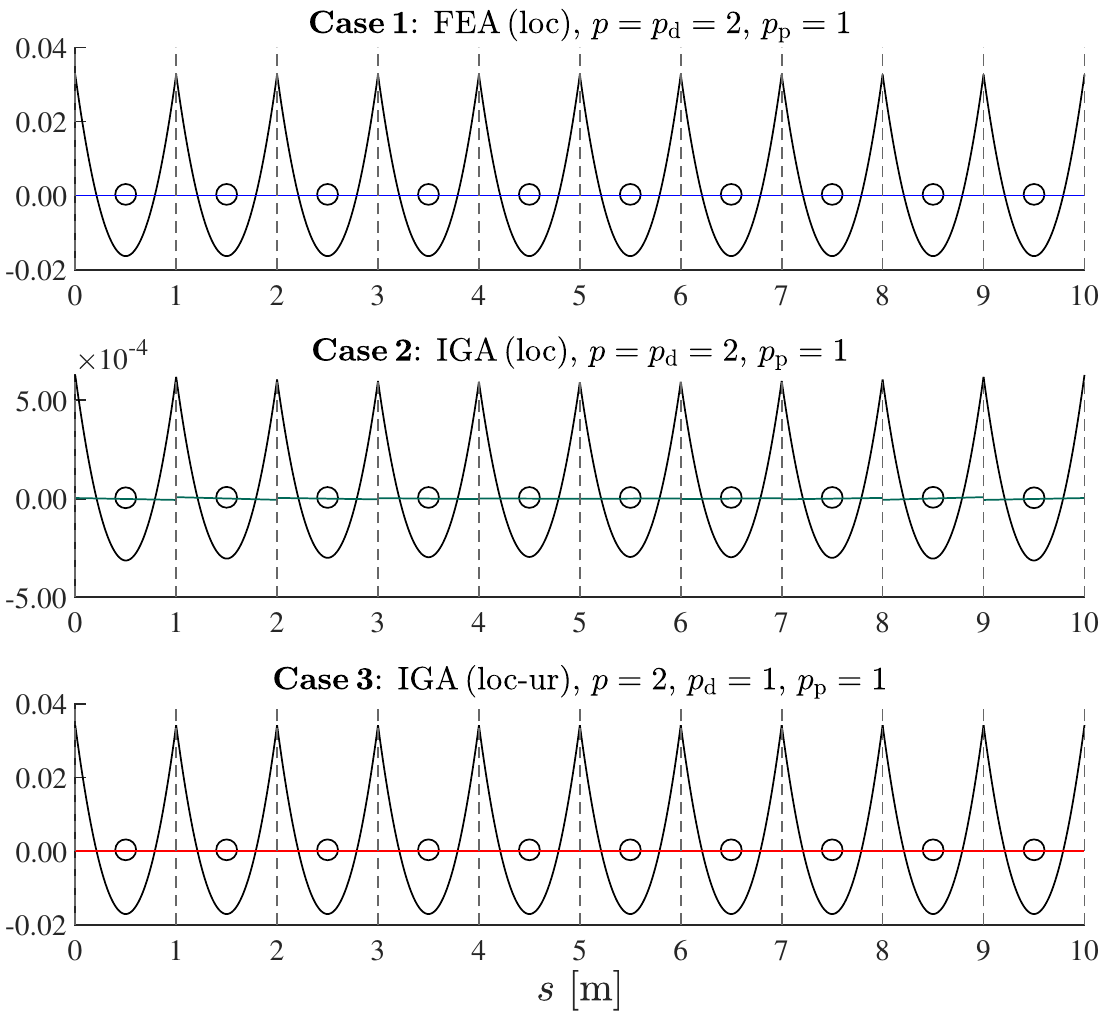}	
		\caption{Comparison of the geometrical and physical strains}		
		\label{pbend_strn_dist_memb_geom_phy}					
	\end{subfigure}			
	\begin{subfigure}[b] {0.49\textwidth} \centering
		\includegraphics[width=\linewidth]{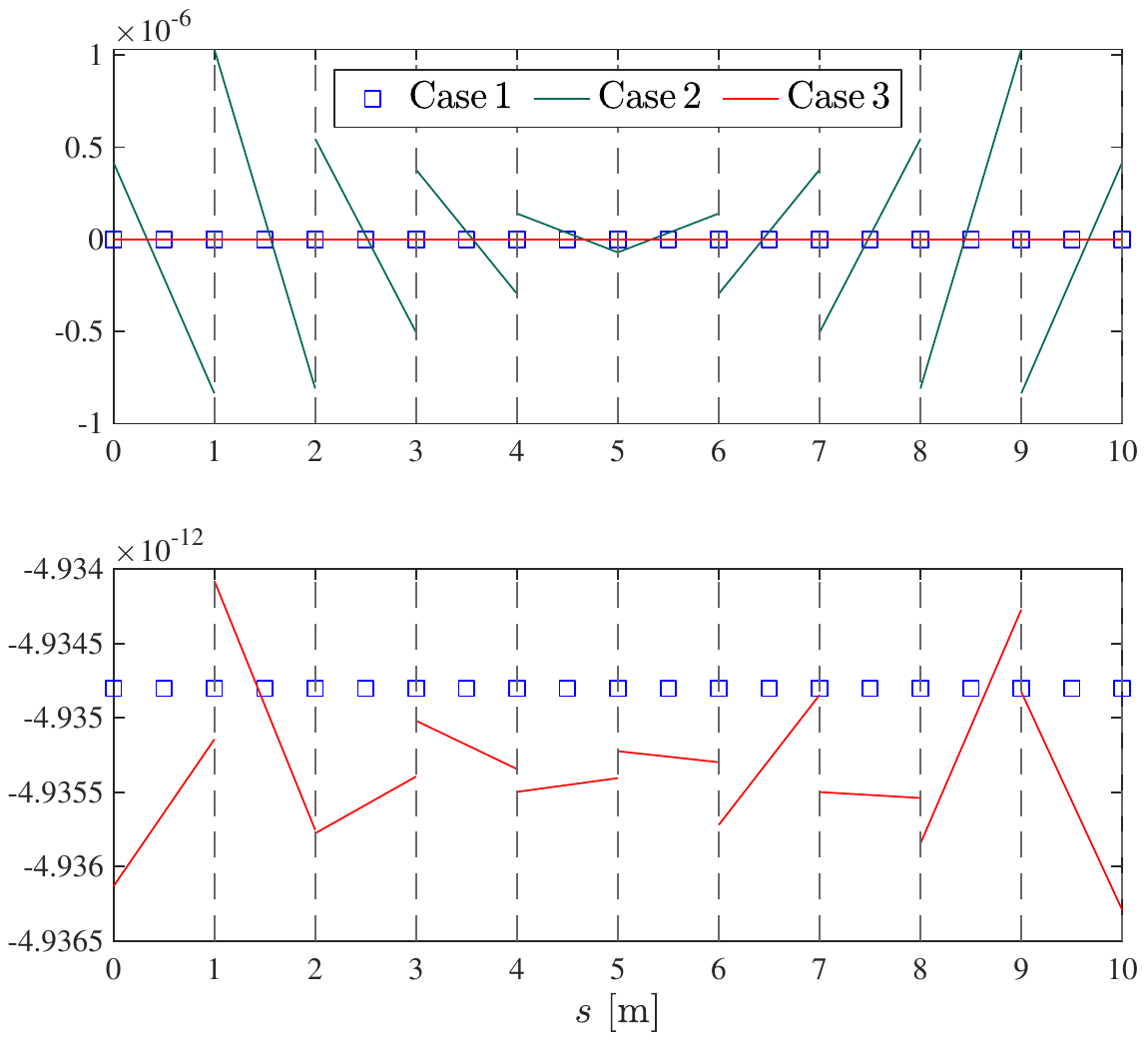}	
		\caption{Comparison of the physical strains only}		
		\label{pbend_strn_dist_memb_phy}
	\end{subfigure}	
	\caption{Cantilever beam under bending moment: Distribution of the \textbf{axial} strain $\varepsilon$ along the axis. (a) Comparison of the distribution of the geometrical (black curves) and physical (colored curves) axial strains along the axis in the three cases. The hollow circles in (a) represent element-wise average of the geometric strain. (b) Those physical strains in (a) are re-plotted in the common ordinate. In all results, $n_\mathrm{el}=10$, and $L/h=10^{5}$. The dashed lines represent the element boundaries.}	
	\label{ex_end_bend_mnt_strn_dist_memb}
\end{figure}

\noindent Fig.\,\ref{ex_end_bend_mnt_strn_dist_memb} compares the distribution of the axial (membrane) strain along the axis. For the degree $p=2$ polynomial bases of the axial displacement, the axial strain is a quadratic polynomial, see the parabolic (black) curves for the geometric strains in Fig.\,\ref{pbend_strn_dist_memb_geom_phy}. In Fig.\,\ref{pbend_strn_dist_memb_geom_phy}, the physical strain in Case 2 shows slight discontinuities across the elements, which can be more clearly seen in Fig.\,\ref{pbend_strn_dist_memb_phy}. Thus, in full Gauss integration along the axis, the spurious membrane strain is evaluated, which also contributes to the artificial increase of the bending stiffness. It is remarkable that, when we reduce the degree $p_\mathrm{d}$ in Case 3, the spurious physical axial strain vanishes. Although the director is not explicitly associated with the axial strain, the alleviation of other coupled spurious strains may affect the axial deformation as well in geometrically nonlinear problems.
\begin{observation} \textcolor{black}{
	\begin{itemize}
		\item Local approaches like the element-wise uniformly reduced integration, and the mixed formulation, with discontinuous physical stress resultants and strains across the elements, is not effective to alleviate numerical locking in IGA, due to the higher order continuity in the displacement field. 
		\item In the displacement-based formulation, the reduced degree $p_\mathrm{d}=p-1$ alleviates transverse shear locking due to the field consistency; however, it still suffers from severe membrane and curvature-thickness locking.
		\item In the mixed formulation (IGA, ``loc-ur'') with $p=2$ and $p_\mathrm{p}=1$, the reduced degree $p_\mathrm{d}=p-1=1$ turns out to alleviate the numerical locking.
	\end{itemize}
}
\end{observation}
\noindent In the following, it will be shown that the system is again significantly over-stiffened, as we increase the order of basis functions, $p$, and $p_\mathrm{d}={p}-1$. To alleviate this, for $p>2$, we additionally adjust the degree $p_\mathrm{p}$ of the basis functions for the additional fields, physical stress resultants and strains. 
\subsubsection{Alleviation of locking for $p=3$ and $p_\mathrm{d}=2$ in IGA by reducing \texorpdfstring{${{p}_\mathrm{p}}$}{TEXT}}
\label{pbend_alleviate_lock_p3}
\begin{figure}[H]
	\centering
	\begin{subfigure}[b] {0.4875\textwidth} \centering
		\includegraphics[width=\linewidth]{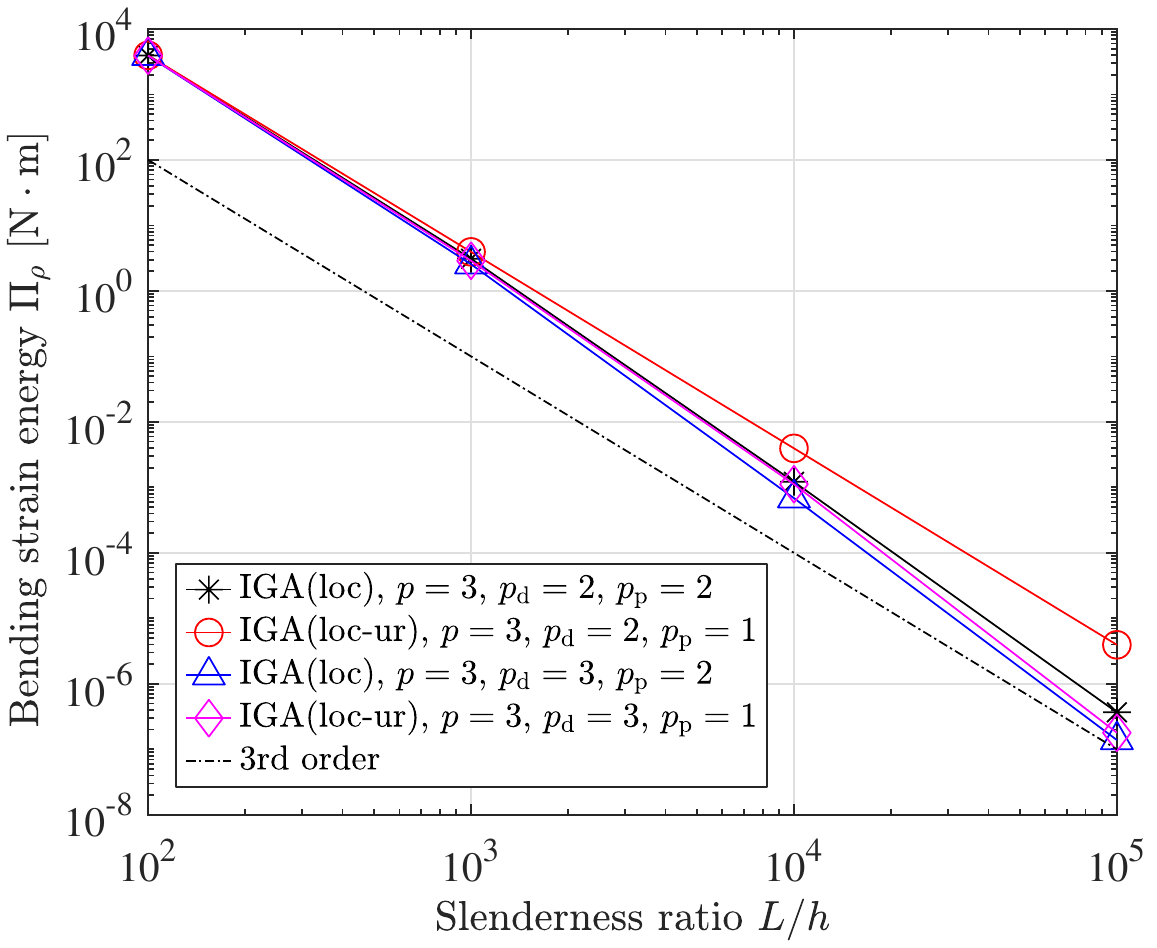}	
		\caption{Bending strain}		
		\label{pbend_mix_strn_e_bend_p3}		
	\end{subfigure}			
	\begin{subfigure}[b] {0.4875\textwidth} \centering
		\includegraphics[width=\linewidth]{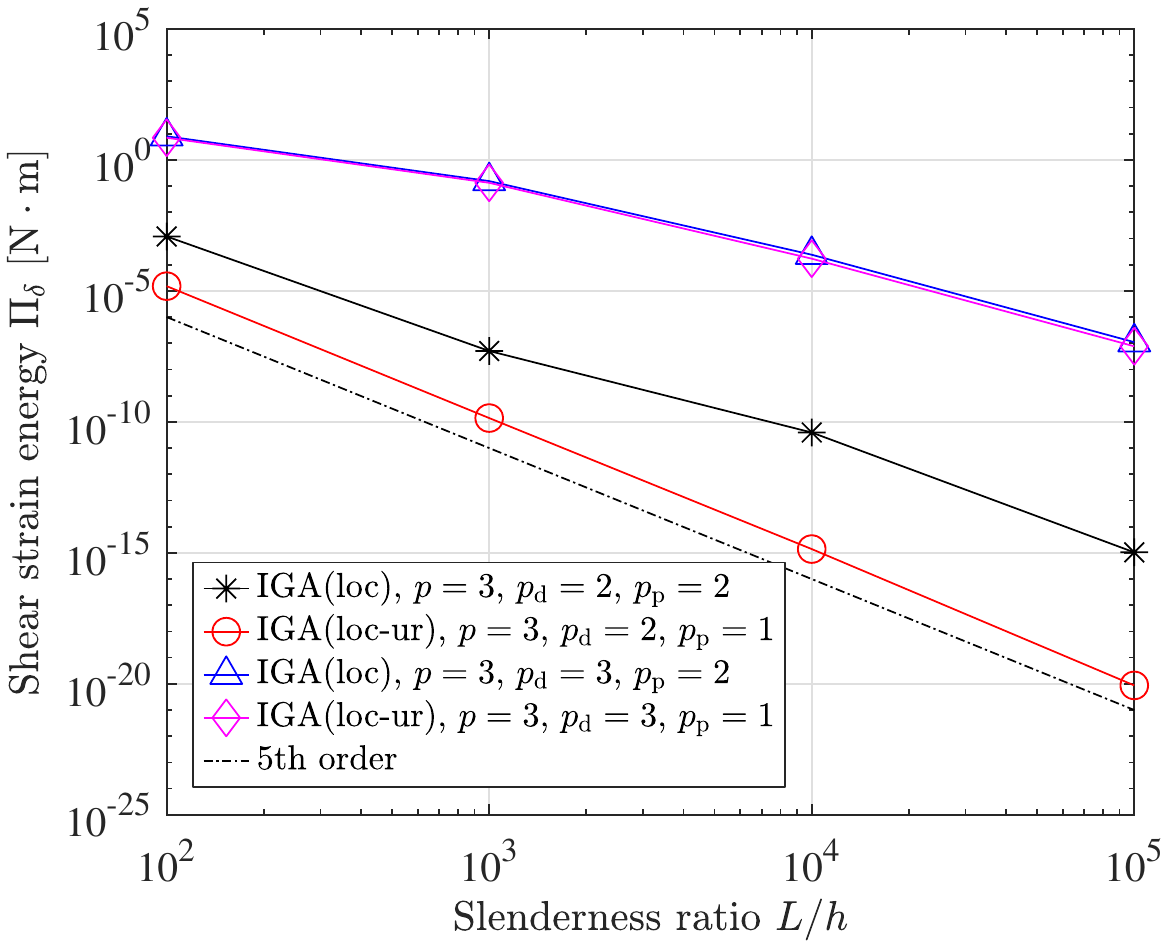}	
		\caption{Transverse shear strain}				
		\label{pbend_mix_strn_tshear_p3}				
	\end{subfigure}	
	\begin{subfigure}[b] {0.4875\textwidth} \centering
		\includegraphics[width=\linewidth]{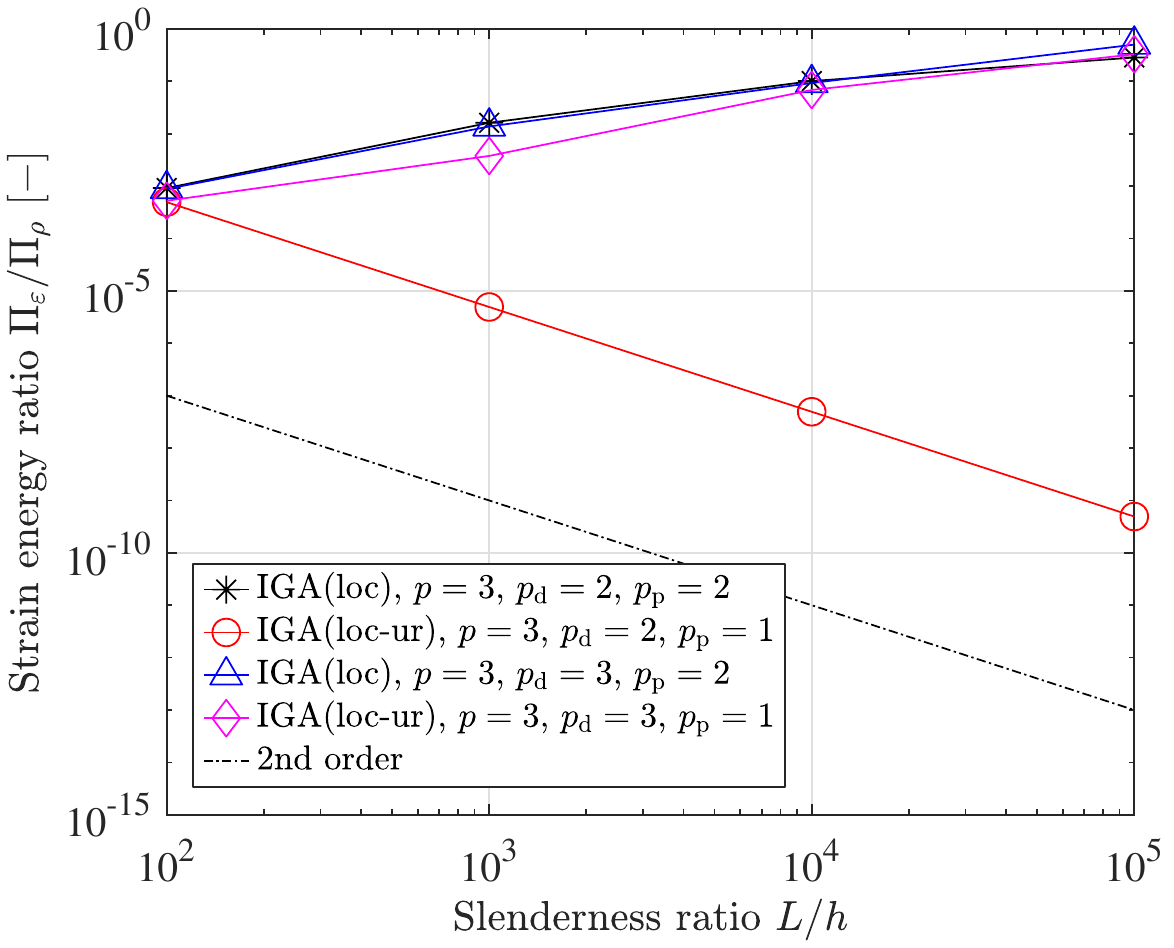}	
		\caption{Membrane strain}		
		\label{pbend_mix_strn_mem_p3}		
	\end{subfigure}				
	\begin{subfigure}[b] {0.4875\textwidth} \centering
		\includegraphics[width=\linewidth]{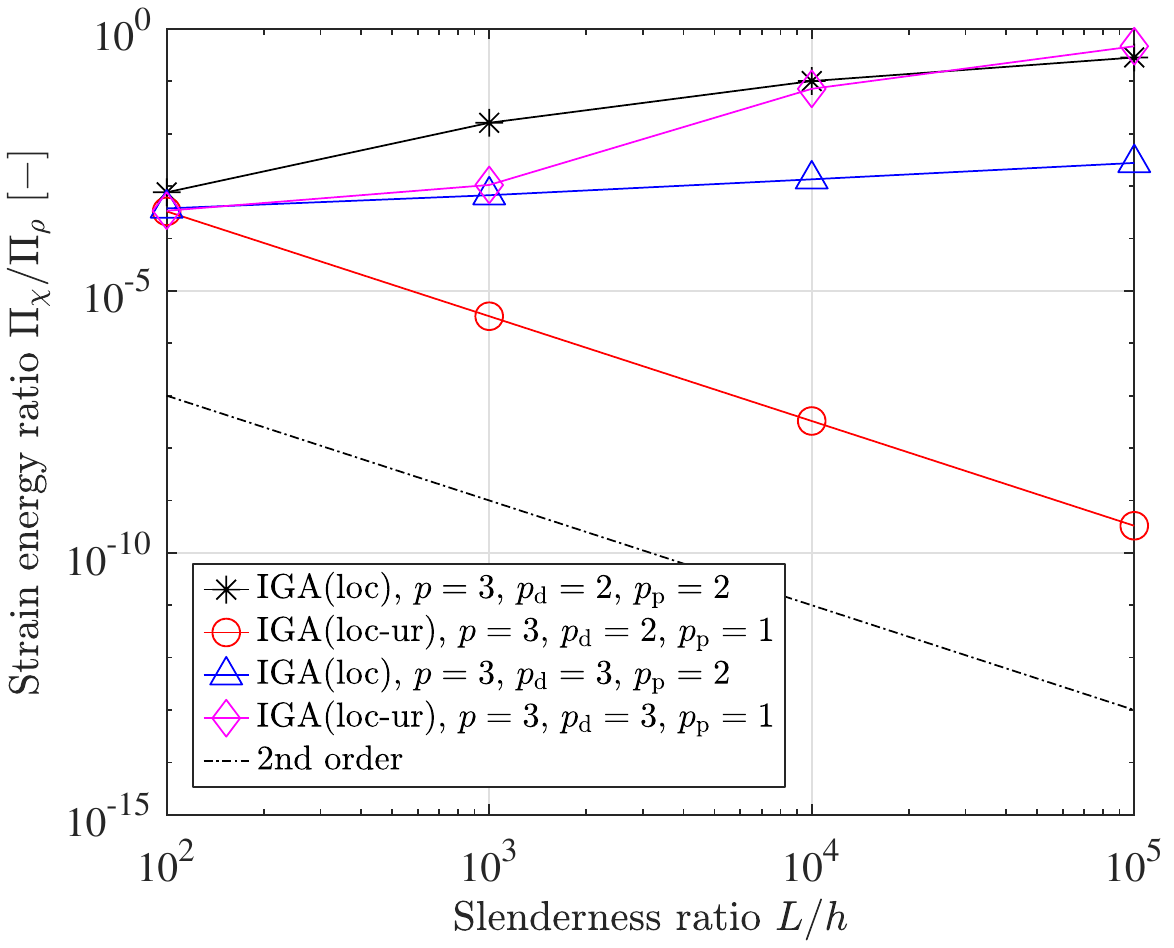}	
		\caption{Transverse normal strain}		
		\label{pbend_mix_strn_inp_p3}		
	\end{subfigure}	
	\caption{Cantilever beam under bending moment: Comparison of the strain energy from the different $p_\mathrm{p}$, and $p_\mathrm{d}$. All results are from the mixed formulation with ${p} = 3$, combined with the local (element-wise) static condensation, full integration ($n_\mathrm{G}=4$), and $n_\mathrm{el}=10$.}
	\label{ex_end_bend_mnt_compare_strne_ratio_mix_p3}
\end{figure}

Fig.\,\ref{ex_end_bend_mnt_compare_strne_ratio_mix_p3} compares the convergence rate of strain energies resulting from using $p_\mathrm{p}={p}-1 = 2$ (``loc'') and $p_\mathrm{p} = 1$ (``loc-ur''). \textcolor{black}{Further, in order to investigate the effect of reducing $p_\mathrm{d}$, we additionally consider two different selections, $p_\mathrm{d}=p=3$, and $p_\mathrm{d}=p-1=2$. In Fig.\,\ref{ex_end_bend_mnt_compare_strne_ratio_mix_p3}, in the results from using $p_\mathrm{p} = p-1 = 2$ (black and blue curves) or $p_\mathrm{d}=p=3$ (magenta curve), deteriorated convergences are clearly observed. In Fig.\,\ref{pbend_mix_strn_tshear_p3}, between the results from using $p_\mathrm{p}=2$, the amount of spurious transverse shear strain energy is significantly lower for $p_\mathrm{d}=p-1=2$ (black curve), than the other from using $p_\mathrm{d}=p=3$ (blue curve).} This also shows the alleviation of transverse shear locking by reducing $p_\mathrm{d} = p - 1$, as we observed in Section \ref{allev_p2_red_pd_iga_lock}. It is remarkable that, when we choose the degrees $p_\mathrm{p} = 1$ and $p_\mathrm{d}=2$, the correct decrease rates are achieved in all cases (red curves). 
%
%
%

\begin{observation}\textcolor{black}{
	\begin{itemize}
		\item It is seen that the conventional selection of $p_\mathrm{p}=p-1$ may lead to significantly over-constrained system in IGA, due to additional artificial constraints arising from the higher order continuity in the displacement field.
		\item It is verified that, by reducing the degree $p_\mathrm{d}$ and $p_\mathrm{p}$, numerical locking can be effectively alleviated.
	\end{itemize}
}
\end{observation}
\noindent In the following two sections, we verify the advantages of IGA \textcolor{black}{combined with the presented approach ``loc-ur''}, in terms of the computational accuracy and efficiency over the conventional FEA, due to the higher order continuity in the NURBS basis functions.
\subsubsection{Verification of the superior convergence behavior of IGA in $h$-refinement}
\label{sec_conv_href_iga}
\begin{figure}[H]
	\begin{subfigure}[b] {0.4975\textwidth} \centering
		\includegraphics[width=\linewidth]{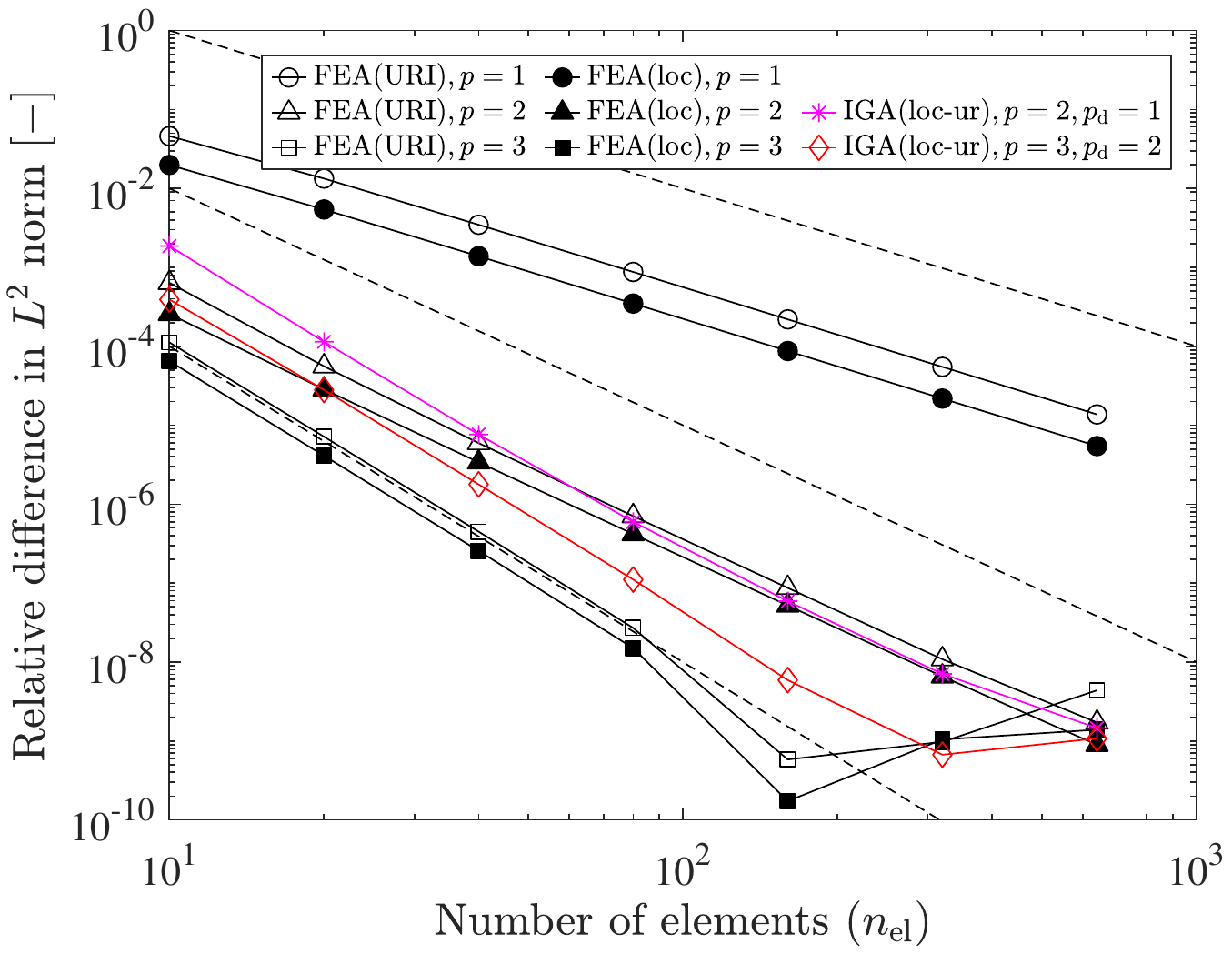}		
		\caption{Comparison of per element accuracy}	
		\label{pbend_per_e_acc}					
	\end{subfigure}
	\begin{subfigure}[b] {0.4975\textwidth} \centering
		\includegraphics[width=\linewidth]{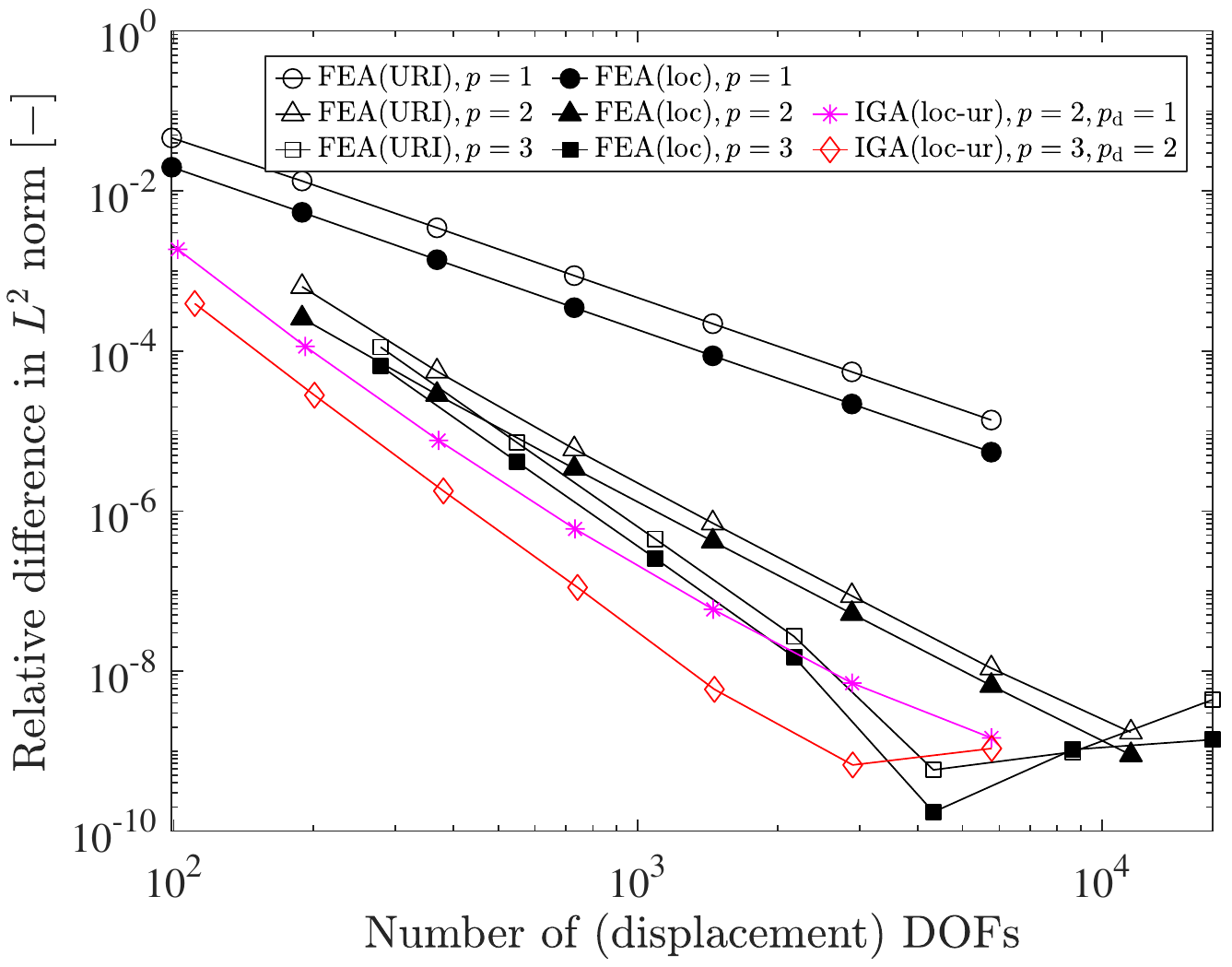}		
		\caption{Comparison of per DOF accuracy}		
		\label{pbend_per_dof_acc}			
	\end{subfigure}
	\caption{Cantilever beam under bending moment: Comparison of the relative difference ${\left\| {{e_u}} \right\|_{{L^2}}}$ in FEA and IGA results. The dashed lines in (a) represent the analytical solution of asymptotic convergence rate, order $p+1=2,3,4$. The selected slenderness ratio is $L/h=10^{5}$. In all IGA results, we use $p_\mathrm{p}=1$ (i.e.,\,``loc-ur''), but for FEA results, we use ${p_\mathrm{p}} = p-1$ (i.e.,\,``loc'').}
	\label{ex_end_bend_mnt_conv_L2_diff_xdisp_per_elem}
\end{figure}

We verify the convergence rate of the presented mixed isogeometric beam formulation, with the reduced degrees of bases, $p_\mathrm{d}$ and $p_\mathrm{p}$. We compare the IGA results with the FEA results using the Lagrange polynomial bases of degree $p$. In FEA, we use $p={p_\mathrm{d}}$, and the following two different formulations to alleviate locking,
\begin{itemize}
	\item the displacement-based formulation with a uniformly reduced integration along the axis (URI) with ${n_\mathrm{G}}=p$,
	\item the mixed formulation with the degree of bases $p_\mathrm{p}={p-1}$ for the additional fields of the physical stress resultants and strains.
\end{itemize}
For all the IGA and FEA results based on the mixed formulation, we use $n_\mathrm{G}=p+1$ Gauss integration points for the exact (full) integration. For the given problem, we have an analytical solution in the thin beam limit (pure bending condition), in which the $X$-displacement of the center axis is expressed by
\begin{align}
	\label{ex_end_mnt_analytic_xdisp}
	{u_{{\rm{ref}}}} = R\sin \frac{X}{R} - L.
\end{align}
For the following verification of the numerical solutions, we utilize the relative $L^2$ norm of the difference in the $X$-displacement ($u$) along the center axis
\begin{equation}\label{def_rel_l2_err}
	{\left\| {{e_u}} \right\|_{{L^2}}} \coloneqq \sqrt {\frac{{\int_0^L {{{\left( {u - {u_{{\rm{ref}}}}} \right)}^2}\,{\mathrm{d}}s } }}{{\int_0^L {{u_{{\rm{ref}}}}^2\,{\mathrm{d}}s } }}}.
\end{equation}
Fig.\,\ref{ex_end_bend_mnt_conv_L2_diff_xdisp_per_elem} compares the convergence of this quantity between the results of FEA and IGA. Figs.\,\ref{pbend_per_e_acc} and \ref{pbend_per_dof_acc} plot the same results with different abscissae: the number of elements, and the number of DOFs, respectively. In Fig.\,\ref{pbend_per_e_acc}, it is seen that the FEA results from using the displacement-based formulation with URI (black curves with hollow markers) shows the same rate of convergence with FEA results from using the mixed formulation (black curves with filled markers) for all degrees of bases, $p=1,2,3$, but the mixed formulation gives more accurate results due to more accurate numerical integration. Further, in cases of very high beam slenderness ratio, the URI requires much more iterations in the Newton-Raphson process, compared with the mixed formulation, which will be shown in Section \ref{pbend_cant_sec_imp_stab}. In the IGA results (red and magenta curves), even though we use $p_\mathrm{d} = {p} - 1$, one degree lower than that of FEA ($p_\mathrm{d}=p$), the convergence rate still agrees very well with the analytical (asymptotic) rate (order $p+1$). Fig.\,\ref{pbend_per_dof_acc} clearly shows that IGA gives superior per DOF accuracy, compared with conventional FEA. It should be noted that the number of DOFs in Fig.\,\ref{pbend_per_dof_acc} only contains the displacement DOFs. In the case of $p=3$, IGA uses much fewer internal DOFs than FEA, due to the lower $p_\mathrm{p}$ for the former, and the difference is proportional to the number of elements. 
\begin{observation} \textcolor{black}{
	\begin{itemize}
		\item In FEA, the element-wise URI effectively alleviates locking, so that it exhibits the same convergence rate with the results of mixed formulation. However, it is less accurate than the mixed formulation, due to the lower accuracy of numerical integration.
		\item In IGA, the solution convergence rate agrees very well with the analytical one, i.e., the decrease rate of $L^2$-difference in order $p+1$, even though the degrees $p_\mathrm{d}$ and $p_\mathrm{p}$ are lowered. Eventually, IGA exhibits superior per DOF accuracy over conventional FEA, and further, it uses much smaller number of internal DOFs for the physical stress resultants and strains.
	\end{itemize}
}
\end{observation}
\subsubsection{Further comparison: $k$-refinement in IGA vs. $p$-refinement in FEA}
\label{ex_end_mnt_kref_vs_p}
\begin{figure}[H]
	\centering
	\begin{subfigure}[b] {0.55\textwidth} \centering
		\includegraphics[width=\linewidth]{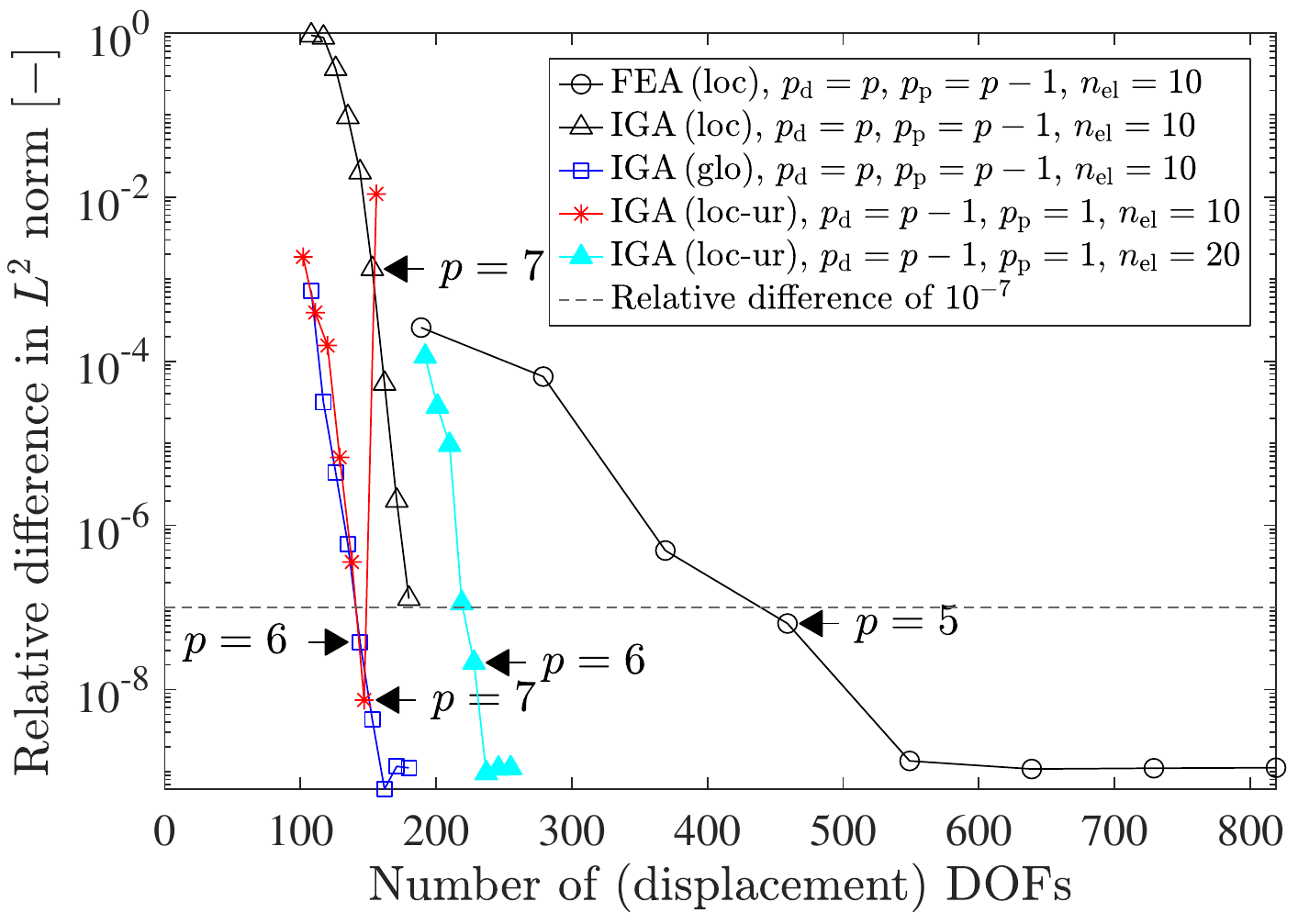}	
	\end{subfigure}			
	\caption{Cantilever beam under bending moment: Comparison of the results from the $k$-refinement in IGA with the $p$-refinement in FEA for the relative difference of the $X-$displacement \textit{per displacement DOF}. In each graph, the markers represent the results of using $p=2,3,\cdots$. In all cases we use the mixed formulation. IGA (``loc-ur'') diverges from $p=9$ and $p=10$ onward for $n_\mathrm{el}=10$ and $n_\mathrm{el}=20$, respectively, so those data are not shown. The slenderness ratio is $L/h=10^5$.}
	\label{ex_end_bend_mnt_conv_L2_diff_xdisp_kref_per_dof}
\end{figure}

\begin{figure}[H]
	\centering
	\begin{subfigure}[b] {0.4\textwidth} \centering
		\includegraphics[width=\linewidth]{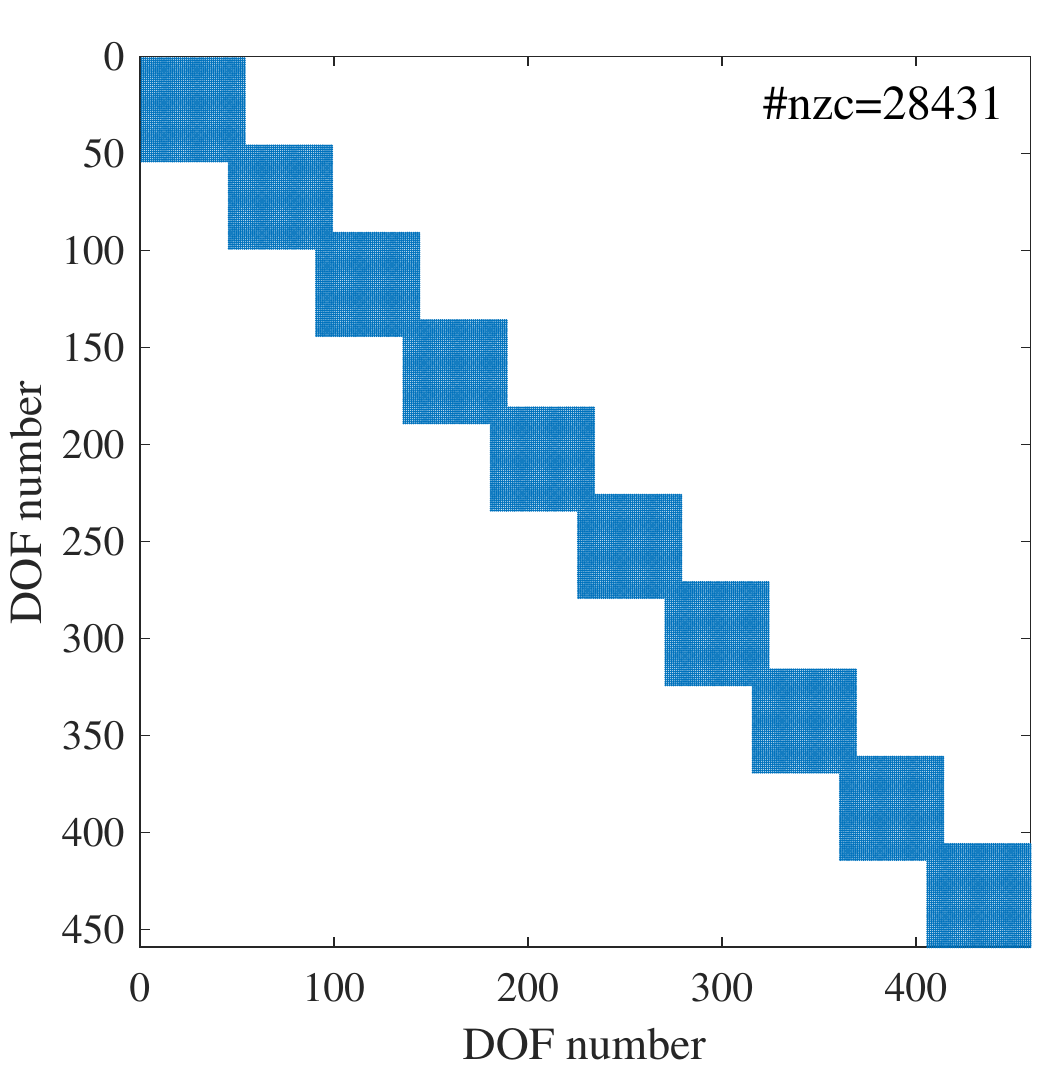}	
		\caption{FEA (loc), $p=p_\mathrm{d}=5$, $p_\mathrm{p}=4$, $n_\mathrm{el}=10$}		
	\end{subfigure}			
	\begin{subfigure}[b] {0.4\textwidth} \centering
		\includegraphics[width=\linewidth]{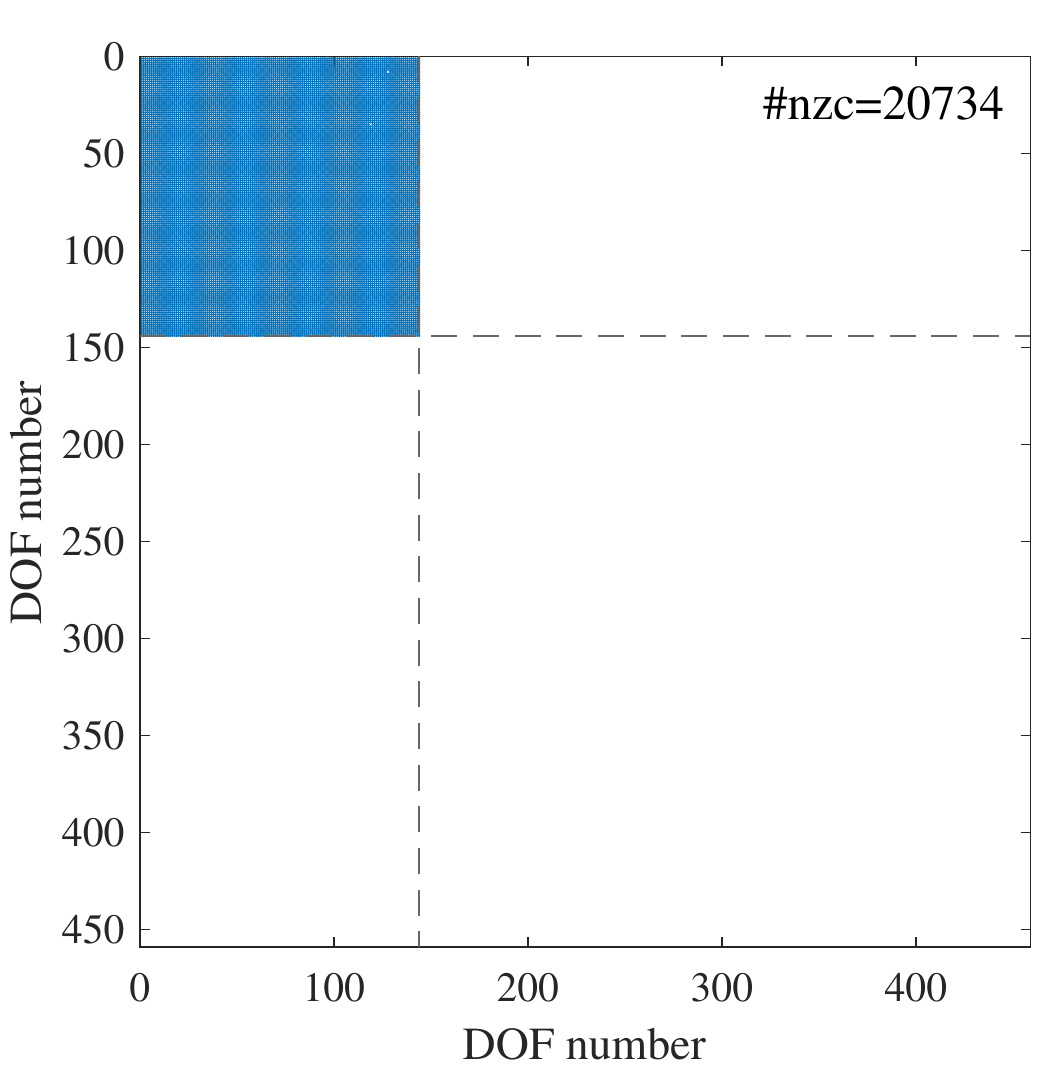}	
		\caption{IGA (glo), $p=p_\mathrm{d}=6$, $p_\mathrm{p} = 5$, $n_\mathrm{el}=10$}		
		\label{sparsity_iga_glo}
	\end{subfigure}	\\
	\begin{subfigure}[b] {0.4\textwidth} \centering
		\includegraphics[width=\linewidth]{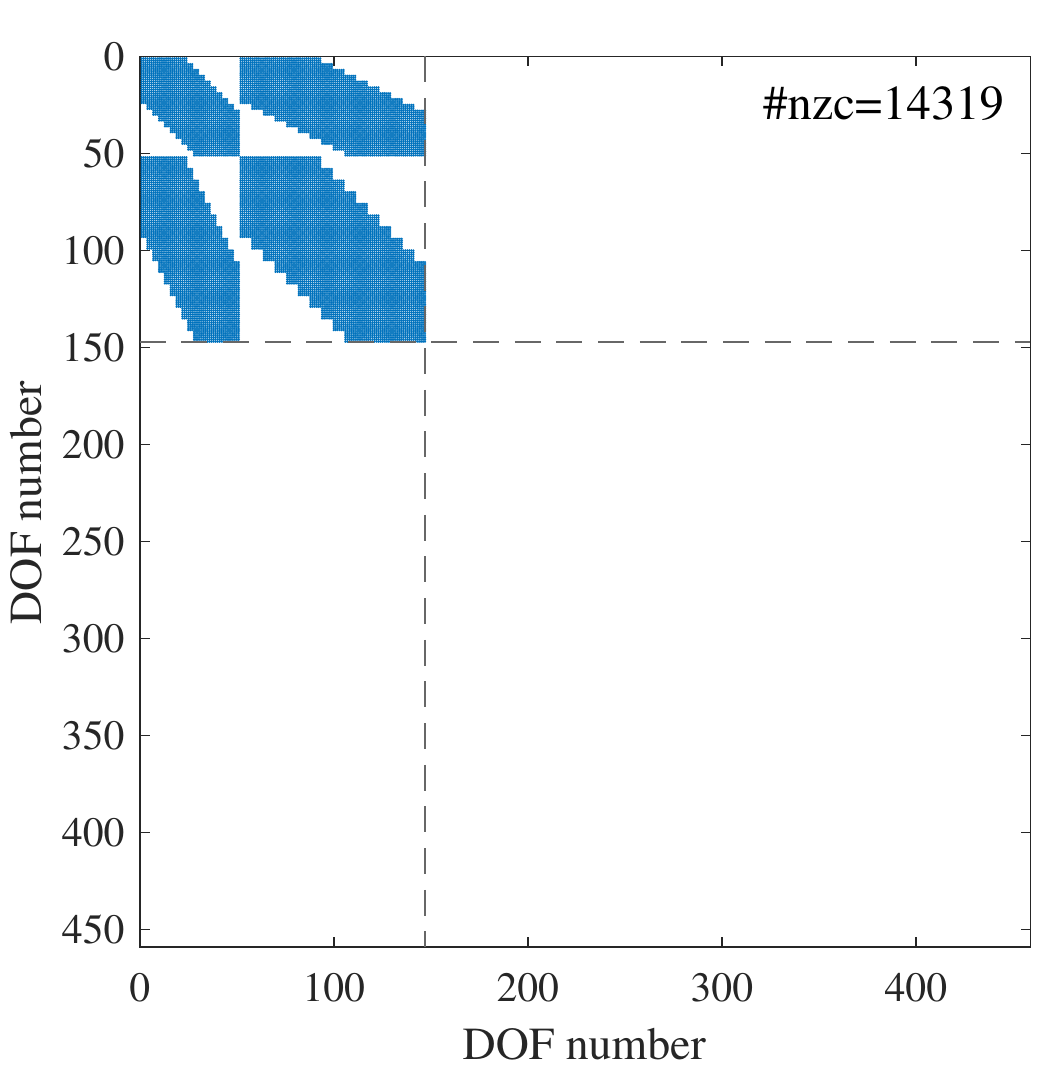}	
		\caption{IGA (loc-ur), $p=7$, $p_\mathrm{d}=6$, $p_\mathrm{p} = 1$, $n_\mathrm{el}=10$}		
		\label{sparsity_iga_loc_red}
	\end{subfigure}	
	\begin{subfigure}[b] {0.4\textwidth} \centering
		\includegraphics[width=\linewidth]{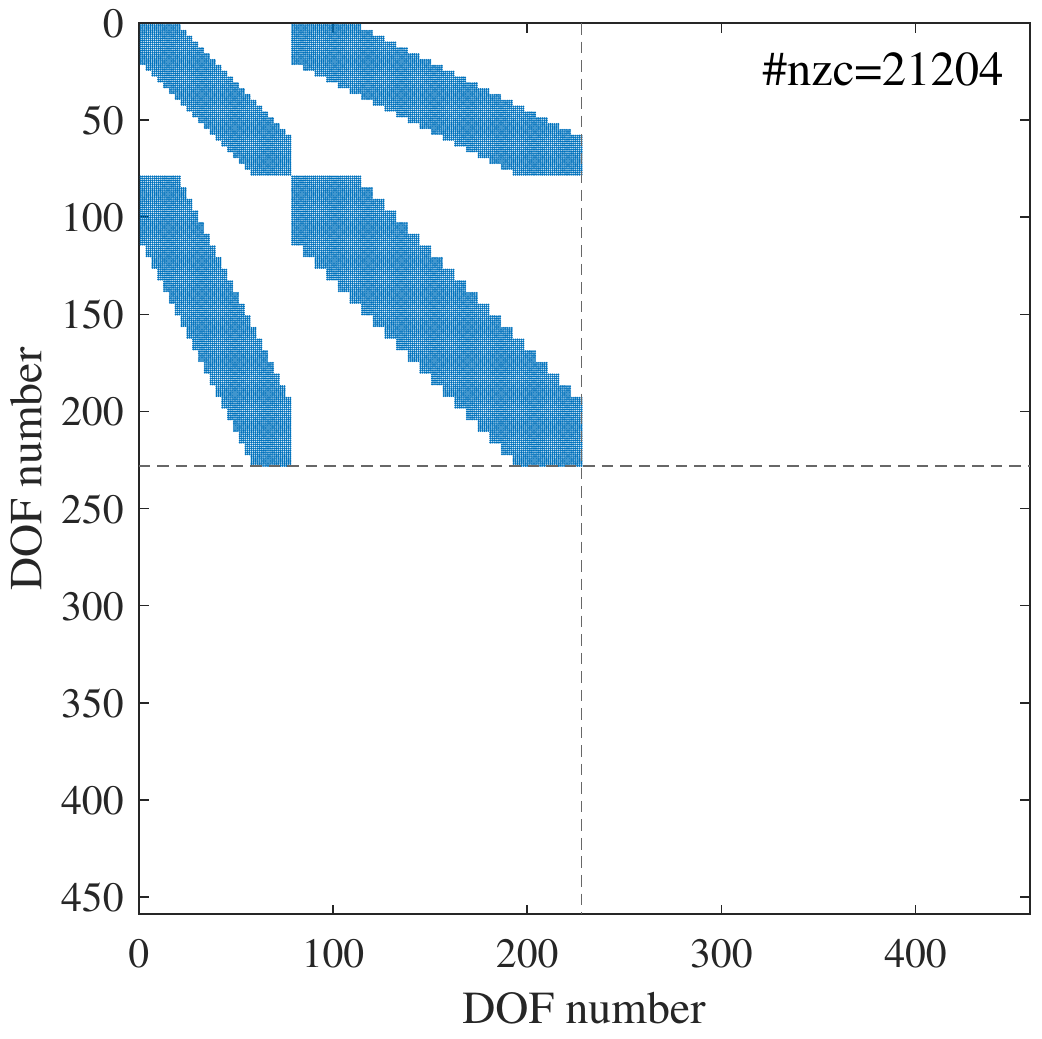}	
		\caption{IGA (loc-ur), $p=6$, $p_\mathrm{d}=5$, $p_\mathrm{p} = 1$, $n_\mathrm{el}=20$}		
		\label{sparsity_iga_loc_red_nel20}
	\end{subfigure}	
	\caption{Cantilever beam under bending moment: comparison of the sparsity pattern in the tangent stiffness matrix at the first iteration of the last (10th) load step. Here {\#}nzc denotes the number of nonzero components. Figures (b-d) are shown in the same scale as (a) for comparison.}
	\label{ex_end_bend_mnt_compare_sparsity_pattern}
\end{figure}

\begin{figure}[H]	
	\centering
	\begin{subfigure}[b] {0.4875\textwidth} \centering
		\includegraphics[width=\linewidth]{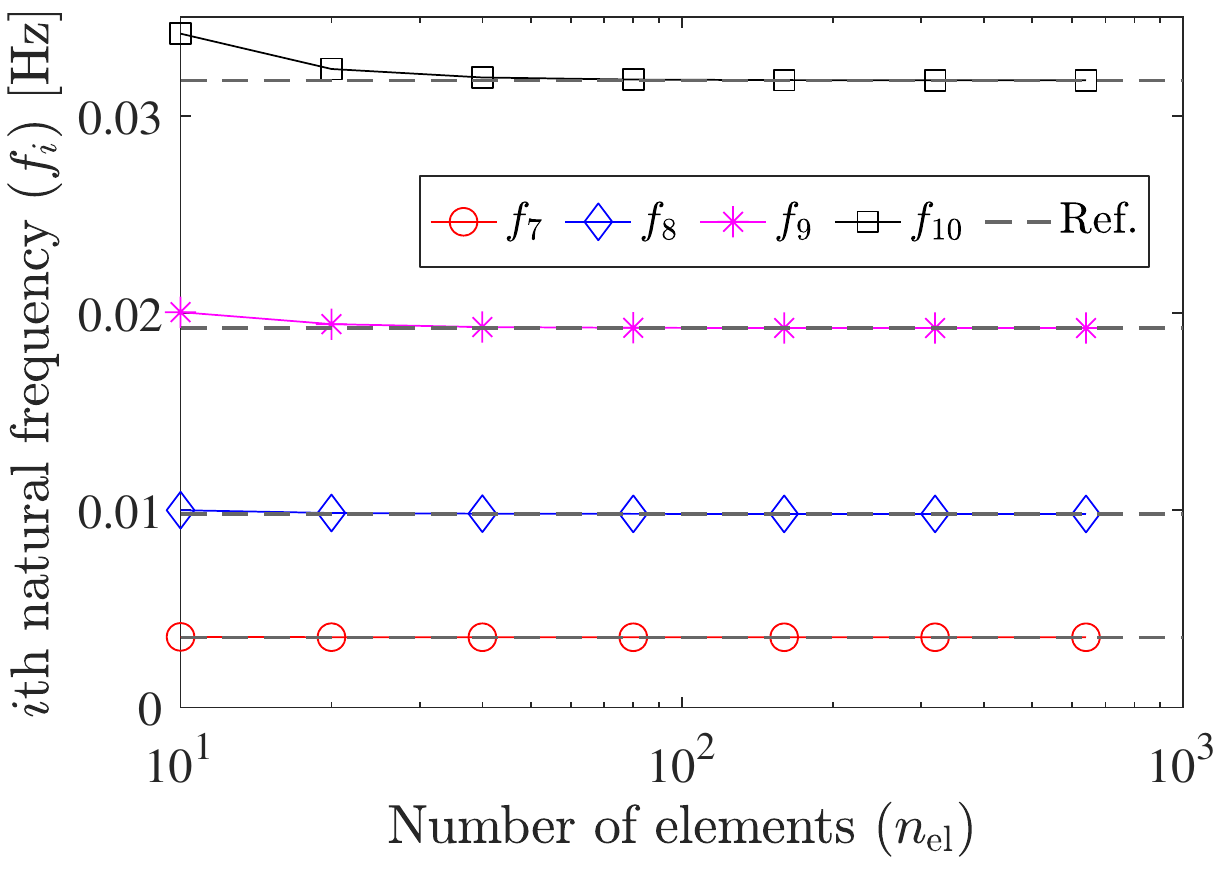}	
		\caption{$p=2$, $p_\mathrm{d}=1$, $p_\mathrm{p}=1$}
	\end{subfigure}		
	\begin{subfigure}[b] {0.4875\textwidth} \centering
		\includegraphics[width=\linewidth]{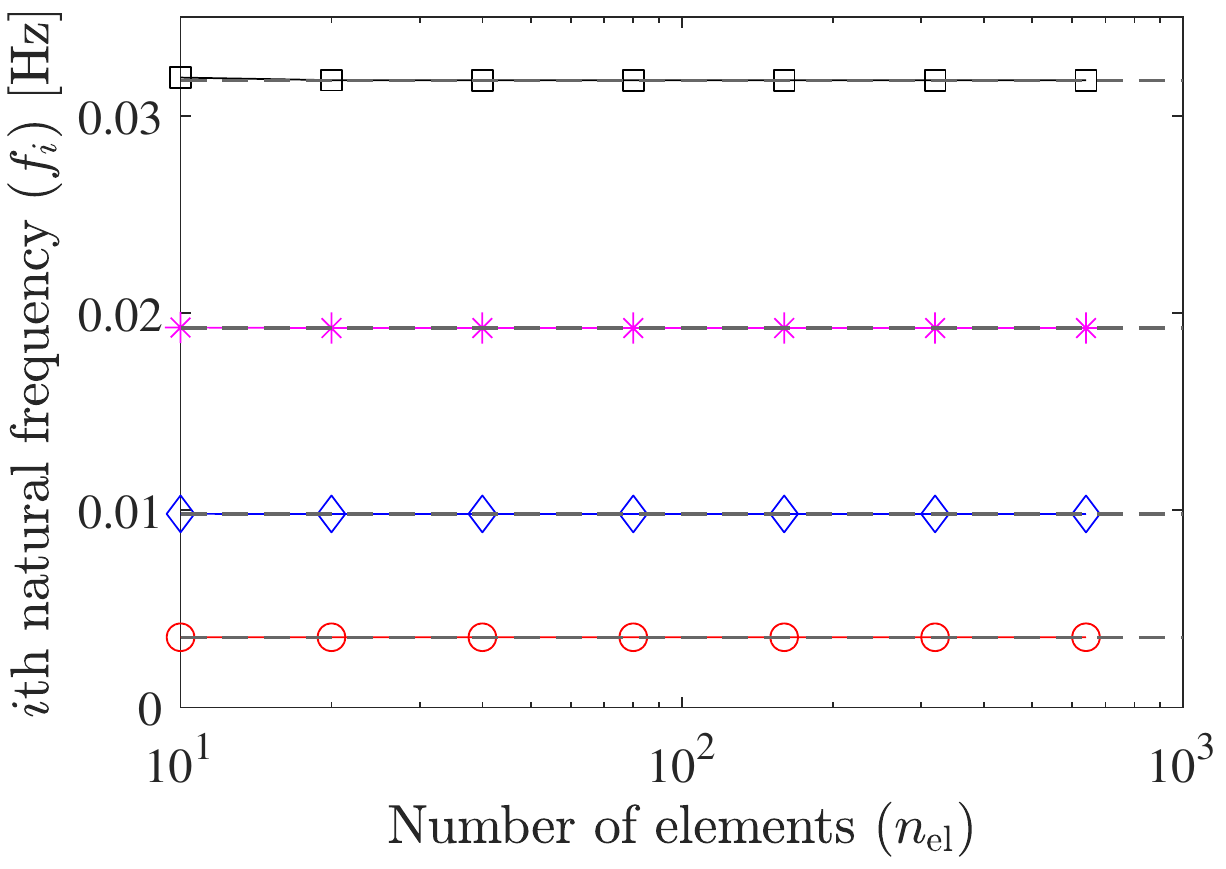}		
		\caption{$p=3$, $p_\mathrm{d}=2$, $p_\mathrm{p}=1$}
	\end{subfigure}		
	\caption{Eigenvalue analysis of a straight beam: Convergence of the first four non-zero eigenfrequencies, ${f_7}\sim{f_{10}}$, with increasing number of elements (i.e., $h$-refinement). The generalized eigenvalue problem is solved at the initial configuration. The dashed lines (``Ref.'') represent the Euler--Bernoulli beam solution for each eigenmode, presented in \citet{han1999dynamics}. \textcolor{black}{The slenderness ratio is $L/h=10^5$}.}
	\label{eval_conv_p2_p3}
\end{figure}

\begin{figure}[H]	
	\centering
	\begin{subfigure}[b] {0.4875\textwidth} \centering
		\includegraphics[width=\linewidth]{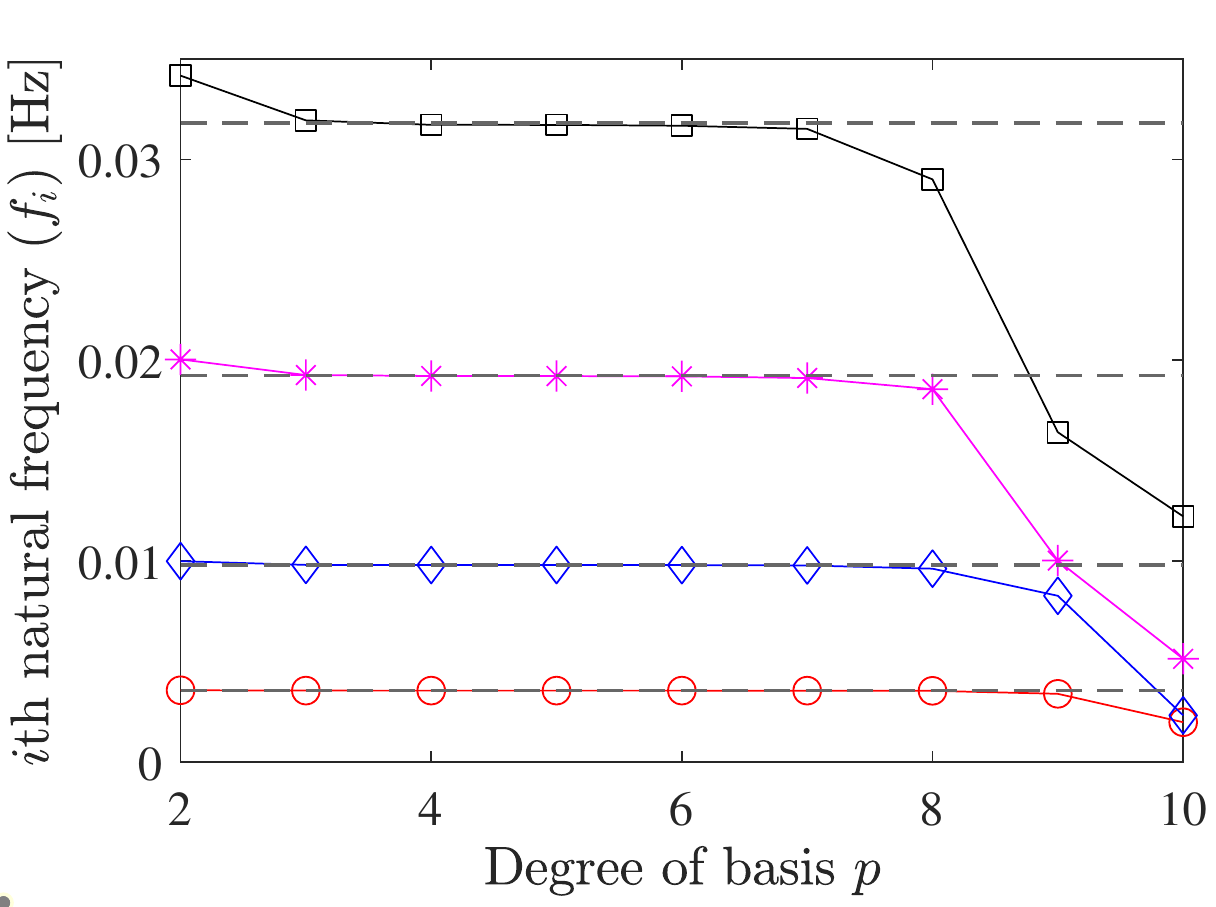}	
		\caption{$n_\mathrm{el}=10$}
		\label{eval_conv_kref_ne10}		
	\end{subfigure}		
	\begin{subfigure}[b] {0.4875\textwidth} \centering
		\includegraphics[width=\linewidth]{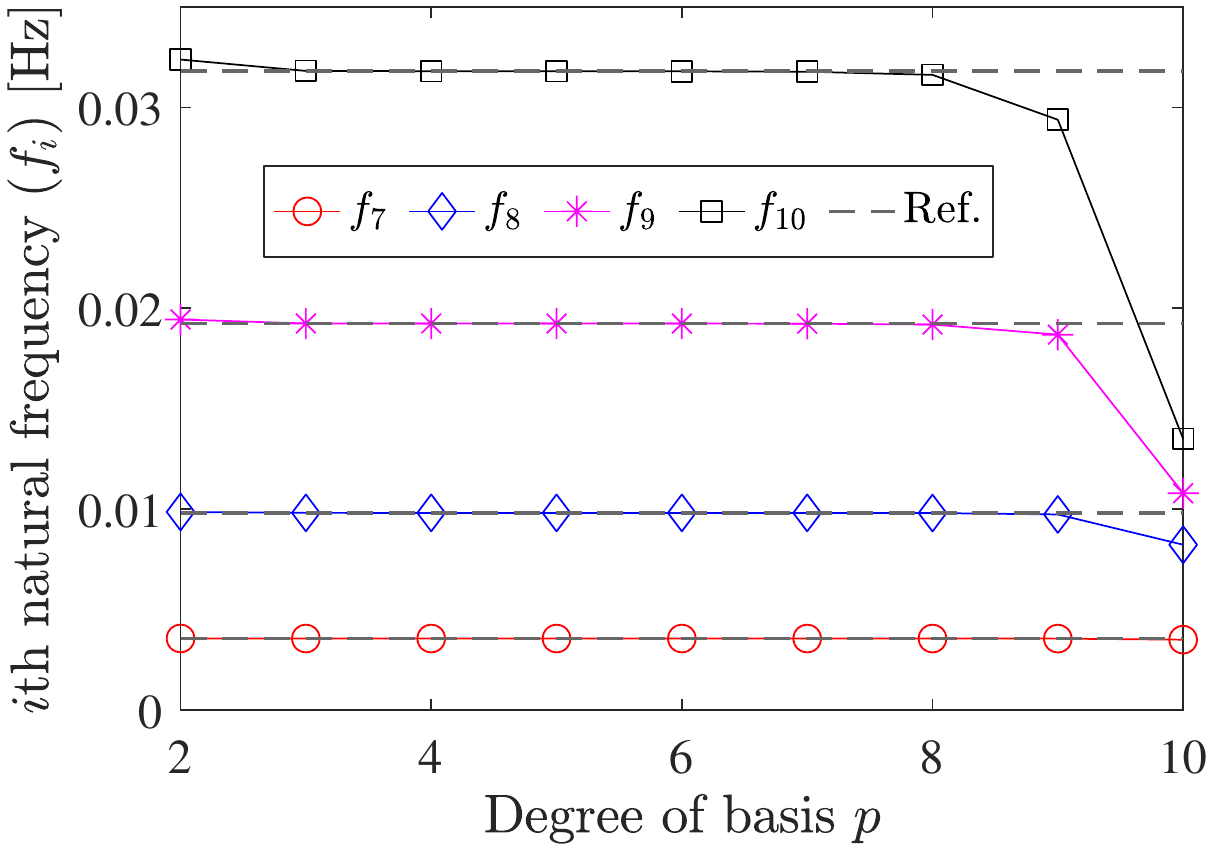}		
		\caption{$n_\mathrm{el}=20$}
		\label{eval_conv_kref_ne20}
	\end{subfigure}			
	\caption{Eigenvalue analysis of a straight beam: Convergence of the first four non-zero eigenfrequencies, ${f_7}\sim{f_{10}}$, with increasing degree $p$ and $p_\mathrm{d}=p-1$ (i.e., $k$-refinement). The generalized eigenvalue problem is solved at the initial configuration. The dashed lines (``Ref.'') represent the Euler--Bernoulli beam solution for each eigenmode, presented in \citet{han1999dynamics}. The slenderness ratio is $L/h=10^5$}
	\label{eval_conv_kref}
\end{figure}

We verify the superior convergence behavior of IGA in $k$-refinement (smooth degree elevation) over classical $p$-refinement (degree elevation) in FEA. The advantage of $k$-refinement is that it enables to maintain the maximum $C^{p-1}$-continuity in the displacement field with bases of degree $p$, in contrast to the $C^0$-continuity from $p$-refinement \citep{hughes2005isogeometric}. We also observe an instability of IGA with uniform reduction of the degree $p_\mathrm{p}$ (``loc-ur'') for very high degree $p\ge8$, but this turns out to be alleviated by increasing the number of elements. In Fig.\,\ref{ex_end_bend_mnt_conv_L2_diff_xdisp_kref_per_dof}, we compare the following three different approaches of the IGA-based mixed formulation with the FEA-based one,
\begin{itemize}
	\item the global approach (``glo'') using $p_\mathrm{p} = p-1$, with $p_\mathrm{d} = p$,
	\item the local approach (``loc'') using $p_\mathrm{p} = p-1$, with $p_\mathrm{d} = p$,
	\item the local approach (``\textcolor{black}{loc-ur}'') using $p_\mathrm{p} = 1$, with $p_\mathrm{d} = p-1$.
\end{itemize}
In Fig.\,\ref{ex_end_bend_mnt_compare_sparsity_pattern}, we further compare the sparsity patterns and the number of nonzero components ({\#}nzc) of the tangent stiffness matrix in the four chosen cases. We choose the cases giving the relative difference around $10^{-7}$ in Fig.\,\ref{ex_end_bend_mnt_conv_L2_diff_xdisp_kref_per_dof} (dashed line): (i) FEA (``loc'') with $p=p_\mathrm{d}=5$ and $p_\mathrm{p}=4$, (ii) IGA (``glo'') with $p=p_\mathrm{d}=6$ and $p_\mathrm{p}=5$, and (iii) IGA (``loc-ur'') with $p=7$, $p_\mathrm{d}=6$ and $p_\mathrm{p}=1$, and (iv) IGA (``loc-ur'') with $p=6$, $p_\mathrm{d}=5$ and $p_\mathrm{p}=1$, where we use $n_\mathrm{el}=10$ for the first three cases, and $n_\mathrm{el}=20$ for the last case. Fig.\,\ref{ex_end_bend_mnt_conv_L2_diff_xdisp_kref_per_dof} shows that the IGA-based mixed formulation (``loc'') without reducing the degrees $p_\mathrm{d}$ and $p_\mathrm{p}$ (black curve with triangle markers) exhibits severe locking, which arises due to the fact that the discontinuous fields of physical stress resultants and strains cannot resolve the parasitic strains from the higher order continuity conditions in the displacement field. Even though the global approach of IGA-based mixed formulation (blue curve with square markers) gives remarkable improvement of the accuracy, it is computationally prohibitive, due to (i) the inversion of the full (global) Gram matrix in the condensation process, and (ii) the resulting dense tangent stiffness matrix, see Fig.\,\ref{sparsity_iga_glo}. In contrast, the presented IGA-based mixed formulation (``loc-ur'') with reduced degrees $p_\mathrm{d}=p-1$ and $p_\mathrm{p} = 1$ (red curve with star markers) gives comparable accuracy as the global approach, at much lower computational cost, due to (i) element-wise (local) condensation, (ii) sparsity of the resulting tangent stiffness matrix, see Figs.\,\ref{sparsity_iga_loc_red} and \ref{sparsity_iga_loc_red_nel20}, and (iii) much fewer number of displacement and internal DOFs due to the reduced degrees $p_\mathrm{d}$ and $p_\mathrm{p}$. However, in the results of IGA (``loc-ur'') with $n_\mathrm{el}=10$, we observe that the performance (i.e., convergence rate with increasing degree $p$) is much better in the range of degree $5\le{p}\le{7}$, compared with other degrees. In cases of very high degrees, $p=9,10$, the solution process even diverges. This result indicates that in the range $2\le{p}\le{4}$, the system\footnote{the system of linear equations in the iterative solution process} is relatively over-constrained, whereas, in the range of $8\le{p}\le{10}$, the system is under-constrained. The instability in the very high range of $p$ can be alleviated by increasing the number of elements ($n_\mathrm{el}$), see, e.g., the results from $n_\mathrm{el}=20$ (cyan curve), which exhibit remarkable per DOF accuracy. 
%
%
%
In order to further investigate the stability of the proposed mixed formulation (``loc-ur'') using $p_\mathrm{p}=1$, with $p_\mathrm{d}=p-1$, we calculate the natural frequencies of the straight beam by solving an eigenvalue problem at the initial configuration. We solve a generalized eigenvalue problem \textcolor{black}{with both ends free}
\begin{equation}
	\label{geigv_prob}
	{\bf{\bar K}}\,\Delta {\bf{y}} = {\omega^2}\,{\bf{M}}\,\Delta {\bf{y}},
\end{equation}
where $\omega$ denotes the angular (natural) frequency, and $\bar {\bf{K}}$ denotes the global tangent stiffness matrix in Eq.\,(\ref{glob_tstiff_cond}), and $\bf{M}$ denotes the global mass matrix (see Appendix \ref{app_cons_mass_mat} for its detailed expression). Here, we consider the initial mass density $\rho_0=1\,\mathrm{kg}/\mathrm{m}^3$. As we consider no displacement boundary conditions, i.e., free-free ends of the beam, the eigenvalue problem of Eq.\,(\ref{geigv_prob}) \textcolor{black}{should give} six zero eigenvalues, associated with the rigid body motions. Here, all the results are calculated by the proposed mixed isogeometric beam formulation (``loc-ur'') using $p_\mathrm{p}=1$, with $p_\mathrm{d}=p-1$. Fig.\,\ref{eval_conv_p2_p3} shows the convergence of the first four non-zero eigenvalues with increasing element number, for two different cases with degrees $p=2$ and $p=3$.  The reference solution (``Ref.'') is obtained by numerically solving\footnote{For this purpose, we use the function \textit{fsolve} in MATLAB.} the frequency equations for Euler-Bernoulli beam model, presented in \citet{han1999dynamics}. In all cases with $p=2$ and $p=3$, the eigenfrequencies $f_i\coloneqq\omega_i/{2\pi}$ $(i=7,8,9,10)$ converge, and the converged values agree very well with the reference solutions. 
\noindent Fig.\,\ref{eval_conv_kref} shows the convergence of the smallest four non-zero eigenfrequencies with increasing the degree of basis $p$ and $p_\mathrm{d}=p-1$, for two different cases of element numbers, $n_\mathrm{el}=10$ and $20$. In both cases, $p_\mathrm{p}=1$. It is seen that the eigenfrequencies drastically decrease in very high degrees $p\ge{8}$, which means that selecting those degrees may suffer from severe numerical instability, as we also observe in the result (red curve) of Fig.\,\ref{ex_end_bend_mnt_conv_L2_diff_xdisp_kref_per_dof}. However, as we increase the number of elements to $n_\mathrm{el}=20$, we obtain much better accuracy of $f_8$ and $f_9$, which eventually leads to very accurate results even for those two degrees $p=8$ and $9$ (cyan curve) in Fig.\,\ref{ex_end_bend_mnt_conv_L2_diff_xdisp_kref_per_dof}. For a further treatment of the locking and instability, one may consider an adjustment of the degree  $p_\mathrm{p}$ \textit{selectively}, in a similar way to the SRI method in \citet{adam2014improved}. \textcolor{black}{This (IGA, ``loc-sr'') will be considered in Section \ref{cant45_case1_subsec_lin}}. Further mathematical investigation on the stability condition in the mixed formulation (see e.g., the generalized inf-sup test in \citet{krischok2019generalized}) remains future work.
\begin{observation} \textcolor{black}{	
	\begin{itemize}
		\item The superior per DOF accuracy of IGA over conventional FEA is more pronounced, as the degree $p$ increases, due to the higher inter-element continuity of the displacement field. 
		\item The uniformly reduced degree $p_\mathrm{p}=1$ in IGA (``loc-ur'') may not be optimal, and exhibits instability in very high degree $p\ge{8}$, which turns out to be alleviated by increasing number of elements. 
	\end{itemize}
}
\end{observation}
\subsubsection{Improved convergence in the thin beam limit}
\label{pbend_cant_sec_imp_stab}
We show that the mixed formulation improves the convergence in the Newton-Raphson iteration. Table \ref{ex_end_bend_mnt_hist_h1e-4} compares the convergence history at the last load step between the displacement-based formulation with uniformly reduced integration (URI) and the mixed formulation. A very high slenderness ratio $L/h=10^5$ is considered. As seen, the mixed formulation needs much fewer iterations for the convergence than the displacement-based one. Further the mixed formulation shows monotonic convergence, but the displacement-based one exhibits severe oscillations during the iteration. \textcolor{black}{Table\,\ref{app_tab_cant_hist_disp-b_fea_uri} shows that, as we increase the number of load steps to 100, the number of iterations in each load step decreases. However, we still observe a serious oscillation in the convergence history, due to the ill-conditioned tangent stiffness matrix in the thin beam limit.} \textcolor{black}{The improved robustness due to the mixed formulation may extend to transient dynamics problems, see e.g., \citet{betsch2016energy}, which remains a future work.}
%
\begin{table}[H]
	\centering
	\caption{Cantilever beam under bending moment: History of the Newton-Raphson iteration at the last (10th) load step, in case of the slenderness ratio $L/h=10^{5}$.}
	\label{ex_end_bend_mnt_hist_h1e-4}	
\begin{tabular}{cccrrr}
\toprule
& \multicolumn{2}{c}{Displacement-based formulation} &       & \multicolumn{2}{c}{Mixed formulation} \\
\cmidrule{2-3}\cmidrule{5-6}      & \multicolumn{2}{c}{FEA, $p={p_\mathrm{d}}=2$, $n_\mathrm{el}=10$, URI} &       & \multicolumn{2}{c}{IGA (loc-ur), $p=2$, $p_\mathrm{d}=1$, $n_\mathrm{el}=10$, FI} \\
\cmidrule{1-3}\cmidrule{5-6}\multirow{2}[2]{*}{Iteration\#} & \multicolumn{1}{c}{\multirow{2}[2]{*}{\makecell{Euclidean norm\\of residual}}} & \multicolumn{1}{c}{\multirow{2}[2]{*}{\makecell{Energy\\norm}}} &       & \multicolumn{1}{c}{\multirow{2}[2]{*}{\makecell{Euclidean norm\\of residual}}} & \multicolumn{1}{c}{\multirow{2}[2]{*}{\makecell{Energy\\norm}}} \\
&       &       &       &       &  \\
\midrule
1     & 6.3E-08 & 3.8E-08 &       & \multicolumn{1}{c}{6.3E-08} & \multicolumn{1}{c}{4.1E-08} \\
2     & 4.1E+02 & 1.5E+02 &       & \multicolumn{1}{c}{4.8E+02} & \multicolumn{1}{c}{1.9E+02} \\
3     & 5.6E+01 & 3.2E+00 &       & \multicolumn{1}{c}{2.3E+02} & \multicolumn{1}{c}{4.4E+01} \\
4     & 2.1E+00 & 3.6E-03 &       & \multicolumn{1}{c}{7.7E+01} & \multicolumn{1}{c}{5.9E+00} \\
5     & 3.9E-03 & 8.7E-09 &       & \multicolumn{1}{c}{1.2E+01} & \multicolumn{1}{c}{1.6E-01} \\
6     & 4.0E-03 & 2.3E-08 &       & \multicolumn{1}{c}{2.8E-01} & \multicolumn{1}{c}{8.0E-05} \\
7     & 2.8E-01 & 3.1E-05 &       & \multicolumn{1}{c}{4.6E-04} & \multicolumn{1}{c}{2.2E-10} \\
8     & 6.5E-03 & 7.1E-08 &       & \multicolumn{1}{c}{1.5E-10} & \multicolumn{1}{c}{2.4E-23} \\
\vdots& \vdots      &   \vdots    &       &       &  \\
299   & 9.9E-05 & 1.8E-11 &       &       &  \\
300   & 1.5E-11 & 8.9E-19 &       &       &  \\
301   & 3.0E-09 & 1.6E-20 &       &       &  \\
302   & 1.5E-11 & 3.2E-26 &       &       &  \\
\bottomrule
\end{tabular}%
\end{table}
\subsection{$45^\circ$-arc cantilever beam}
\label{num_ex_45deg_arc_cant}
\begin{figure}[H]	
	\centering
	\begin{subfigure}[b] {0.6\textwidth} \centering
		\includegraphics[width=\linewidth]{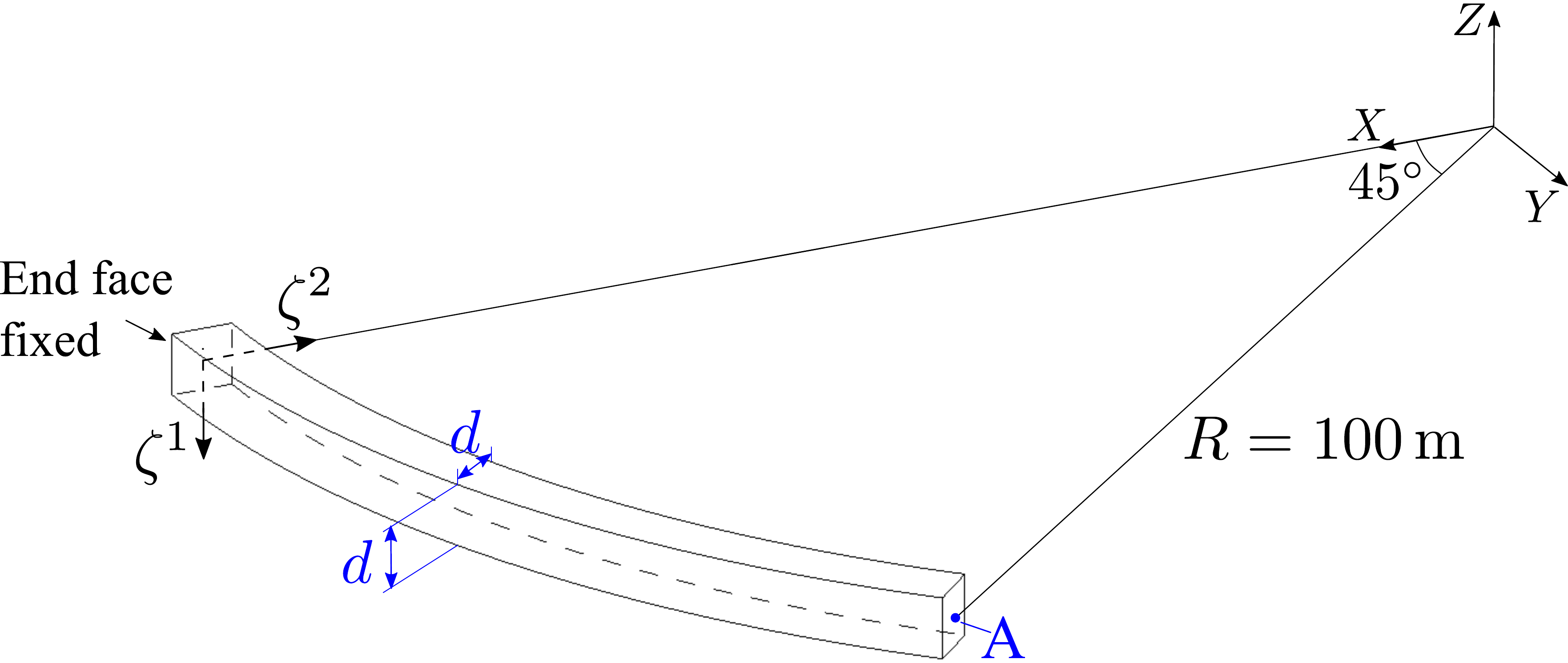}	
	\end{subfigure}		
	\caption{$45^\circ$-arc cantilever beam: Undeformed configuration and boundary conditions. Point A indicates the loaded end (tip) of center axis.}
	\label{45deg_init_config}
\end{figure}

We consider the initial beam center axis lying on the $XY$-plane and describing an $1/8$ of a full circle with radius $R=100\,\mathrm{m}$ and square cross-section of dimension $d\,[\mathrm{m}]$, see Fig.\,\ref{45deg_init_config}. One end face is fixed, and a distributed $Z$-directional force of \textcolor{black}{$\boldsymbol{F}=F\,\boldsymbol{e}_3\,[\mathrm{N}/\mathrm{m}^2]$ is applied on the other end face, where $F$ represents the external force per unit undeformed area, such that the resultant force is 
\begin{equation}
	\label{cant45deg_app_end_f}
	F\cdot{A}=7.2\times{10^3}\dfrac{{E_0}{I}}{{{L_0}^2}},
\end{equation}
with the second area moment of inertia $I=d^{\,4}/12$, and the initial area $A=d^{\,2}$ of the cross-section. $E_0$ and $L_0$ define the unit stress $[{\mathrm{N}/\mathrm{m}^2}]$, and the unit length $[\mathrm{m}]$, respectively. Note that we apply the external (resultant) force $F\cdot{A}=600\,\mathrm{N}$ in the case of the slenderness ratio $R/d=10^2$, as in \citet{bathe1979large}, see also Table \ref{app_num_ex_cant45_ext_f}} \textcolor{black}{for the other slenderness ratios.} Two different cases with isotropic material properties are considered: 
\begin{itemize}
	\item Case 1: St.\,Venant-Kirchhoff type material with zero Poisson's ratio ($\nu=0$),
	\item Case 2: Compressible Neo-Hookean type material with Poisson's ratio $\nu=0.3$.
\end{itemize} 
In both cases the Young's modulus is $E=10\,\mathrm{MPa}$. This example aims at verifying the followings: (i) the improved agreement of the beam solution with the brick element one, due to the enriched cross-sectional strains, (ii) alleviation of locking in IGA-based mixed formulation combined with the local approach, by reducing the degrees $p_\mathrm{d}$ and $p_\mathrm{p}$, (iii) the robustness of the developed mixed finite element formulation in terms of stable Newton-Raphson convergence \textcolor{black}{in the thin beam limit}, (iv) path-independence of the solution from using the presented formulation.
\subsubsection{Case 1: Linear elastic material without Poisson effect}
\label{cant45_case1_subsec_lin}
We compare the displacement at the tip of the center axis (point A in Fig.\,\ref{45deg_init_config}) with the brick element solution\textcolor{black}{, in the case of the slenderness ratio $R/d=10^2$}. For the IGA-based brick element solution, the following three different levels of mesh refinement in the cross-section is considered, with the same 320 NURBS elements (nonzero knot spans) along the center axis:
\begin{itemize}
	\item a single linear element, i.e., $\mathrm{deg.}=(3,1,1)$, and $n_\mathrm{el}=320\times1\times1$,
	\item a single quadratic B-spline element, i.e., $\mathrm{deg.}=(3,2,2)$, and $n_\mathrm{el}=320\times1\times1$,
	\item $8\times8$ cubic B-spline elements, i.e., $\mathrm{deg.}=(3,3,3)$, and $n_\mathrm{el}=320\times8\times8$.
\end{itemize}
Here, the notation $\mathrm{deg.}\!=\!(p_\mathrm{L},p_\mathrm{H},p_\mathrm{W})$, and $n_\mathrm{el}\!=\!{n_\mathrm{L}\times{n_\mathrm{H}}\times{n_\mathrm{W}}}$ denotes the degrees of basis functions, and the number of elements along the axial (L), and two transverse directions (H,W), respectively. Table\,\ref{app_num_ex_cant45_c1_conv} shows the convergence of the brick element solution. In Table\,\ref{ex_cant45_case1_tip_disp_compare}, we compare the tip displacements of the presented beam formulation with reference solutions from literature. It is noticeable that the presented beam formulation (``IGA,\,mixed,\,loc-ur'') gives much better agreement with the brick element solution using a single quadratic element in the cross-section, compared with other reference solutions. However, it still deviates from the brick element solution with multiple cubic elements in the cross-section (i.e.,\,``IGA,\,brick,\,$\mathrm{deg.}\!=\!(3,3,3)$''). This can be further improved by enriching the higher order cross-sectional strains including the out-of-plane ones (cross-sectional warpings), which remains future work. It is seen that the enrichment of the linear and bilinear in-plane cross-sectional strains by the EAS method (i.e.,\,``IGA,\,mixed,\,loc-ur,\,EAS'') improves the agreement slightly. 
\begin{table}[H]
	\centering
	\caption{$45^\circ$-arc cantilever beam (Case 1, \textcolor{black}{$R/d=10^2$}): Comparison of the tip displacements. The reference solutions from literature are obtained from Table 1 of \citet{frischkorn2013solid}}
	\label{ex_cant45_case1_tip_disp_compare}	
	\begin{tabular}{lccc}
		\toprule
		& $u_1$ & $u_2$ & $u_3$ \\
		\midrule
		IGA, brick, $\mathrm{deg.}\!=\!(3,3,3)$, $n_\mathrm{el}\!=\!320\!\times\!8\!\times\!8$ & 13.730631 & -23.825817 & 53.609888 \\
		IGA, brick, $\mathrm{deg.}\!=\!(3,1,1)$, $n_\mathrm{el}\!=\!320\!\times\!1\!\times\!1$ & 13.604502 & -23.568020 & 53.477701 \\
		{IGA, brick, $\mathrm{deg.}\!=\!(3,2,2)$, $n_\mathrm{el}\!=\!320\!\times\!1\!\times\!1$} & {13.604269} & {-23.567644} & {53.477296} \\
		IGA, displacement-based, $p=4$, $n_\mathrm{el}=80$ & 13.604255 & -23.567612 & 53.477268 \\
		IGA, mixed (loc-ur), $p=4$, $n_\mathrm{el}=80$ & 13.604255 & -23.567611 & 53.477267 \\
		{IGA, mixed (loc-ur, EAS), $p=4$, $n_\mathrm{el}=80$} & {13.604256} & {-23.567616} & {53.477273} \\
		\citet{frischkorn2013solid} & 14.11 & -23.38 & 53.50 \\
		\citet{rhim1998vectorial} & 13.70 & -23.64 & 53.46 \\
		\citet{bathe1979large} & 13.40 & -23.50 & 53.40 \\
		\citet{crisfield1990consistent}   & 13.69 & -23.87 & 53.71 \\
		\citet{simo1986three}  & 13.50 & -23.48 & 53.37 \\
		\citet{cardona1988beam}    & 13.74 & -23.67 & 53.50 \\
		\citet{dvorkin1988non}    & 13.60 & -23.50 & 53.30 \\
		\bottomrule
	\end{tabular}%
\end{table}

\begin{figure}[H]	
	\centering
	\begin{subfigure}[b] {0.4455\textwidth} \centering
		\includegraphics[width=\linewidth]{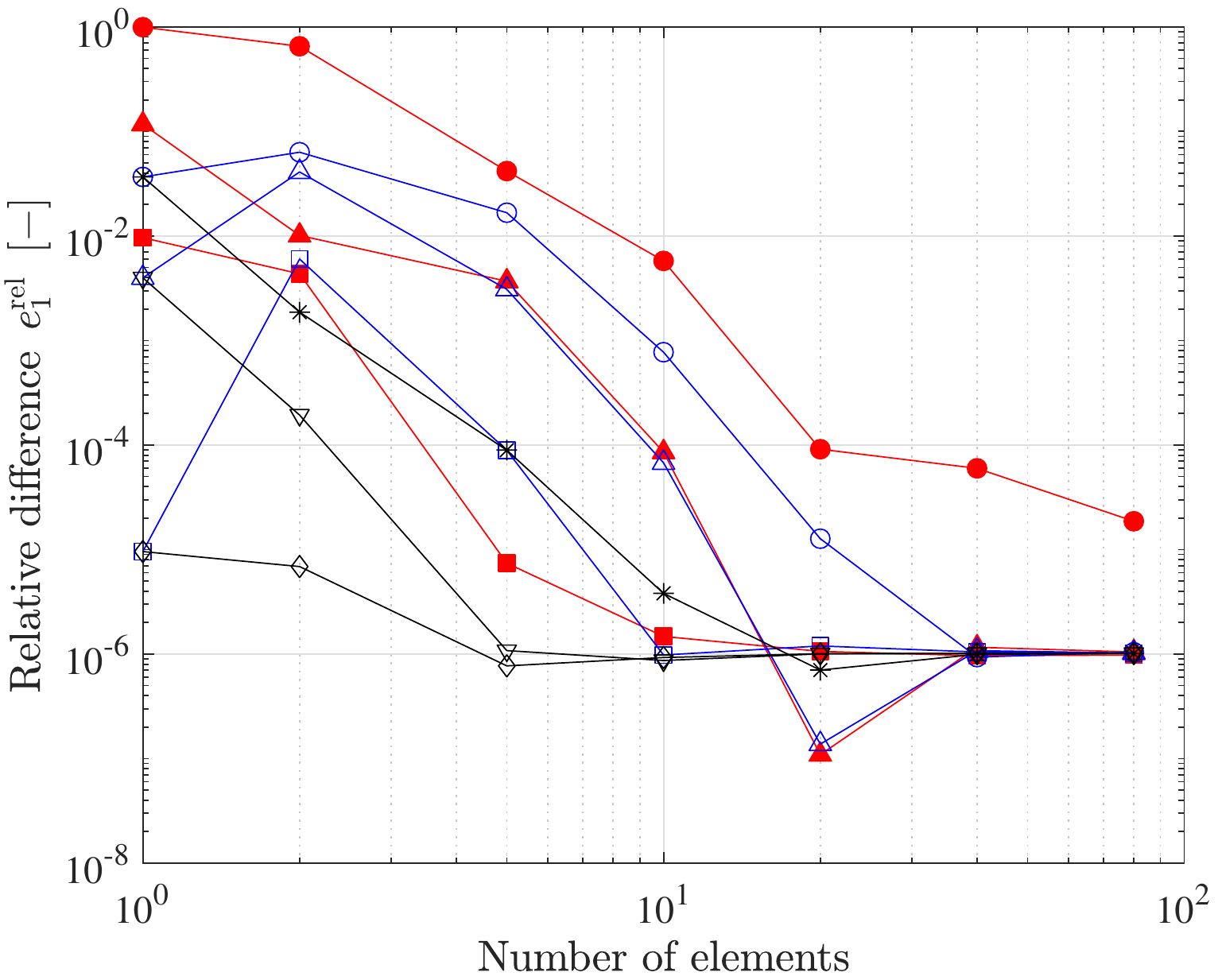}	
		\caption{$X$-displacement (per element accuracy)}
	\end{subfigure}
	\begin{subfigure}[b] {0.4455\textwidth} \centering
	\includegraphics[width=\linewidth]{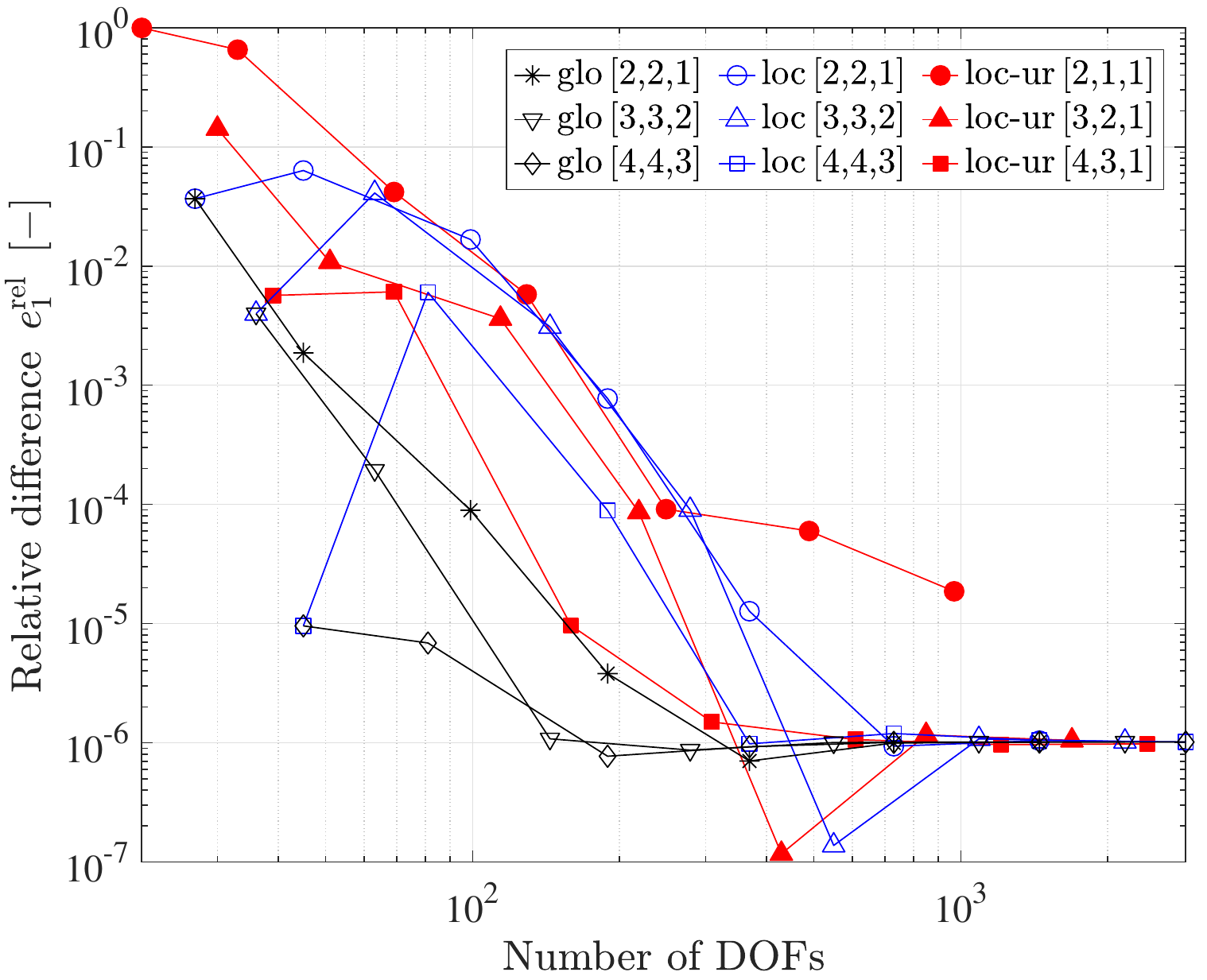}	
	\caption{$X$-displacement (per DOF accuracy)}
	\end{subfigure}
	\begin{subfigure}[b] {0.4455\textwidth} \centering
		\includegraphics[width=\linewidth]{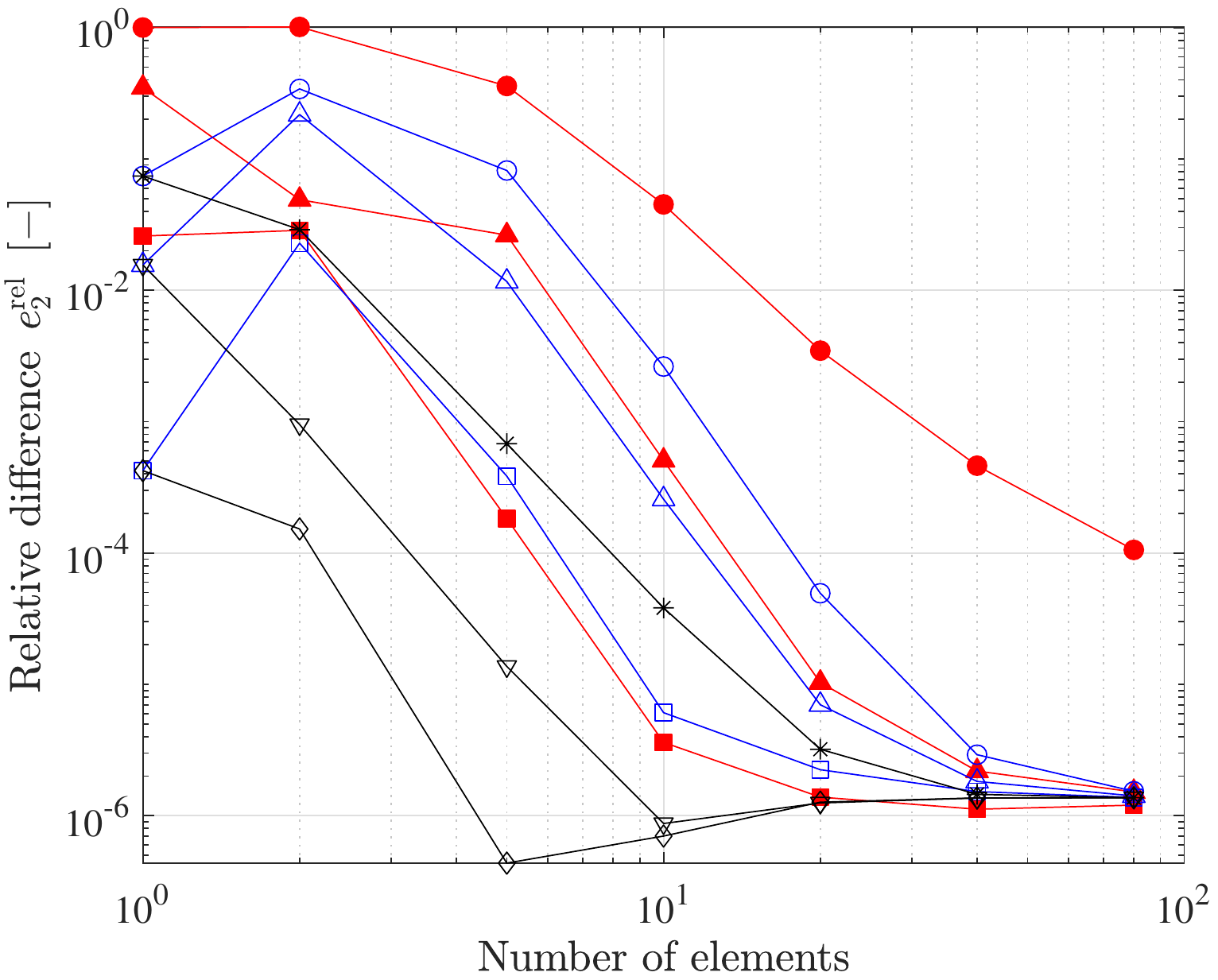}	
		\caption{$Y$-displacement (per element accuracy)}
	\end{subfigure}
	\begin{subfigure}[b] {0.4455\textwidth} \centering
		\includegraphics[width=\linewidth]{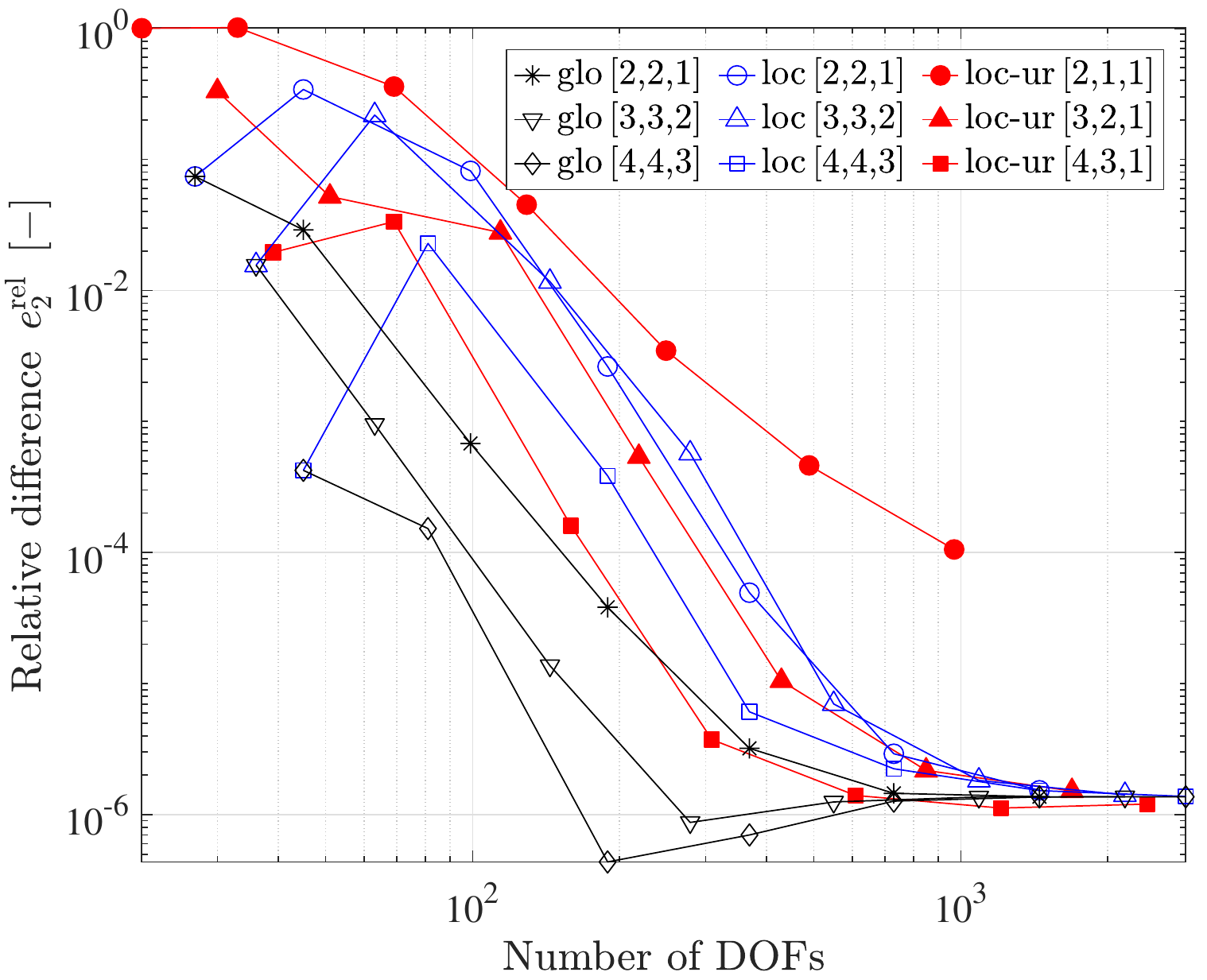}	
		\caption{$Y$-displacement (per DOF accuracy)}
	\end{subfigure}		
	\begin{subfigure}[b] {0.4455\textwidth} \centering
		\includegraphics[width=\linewidth]{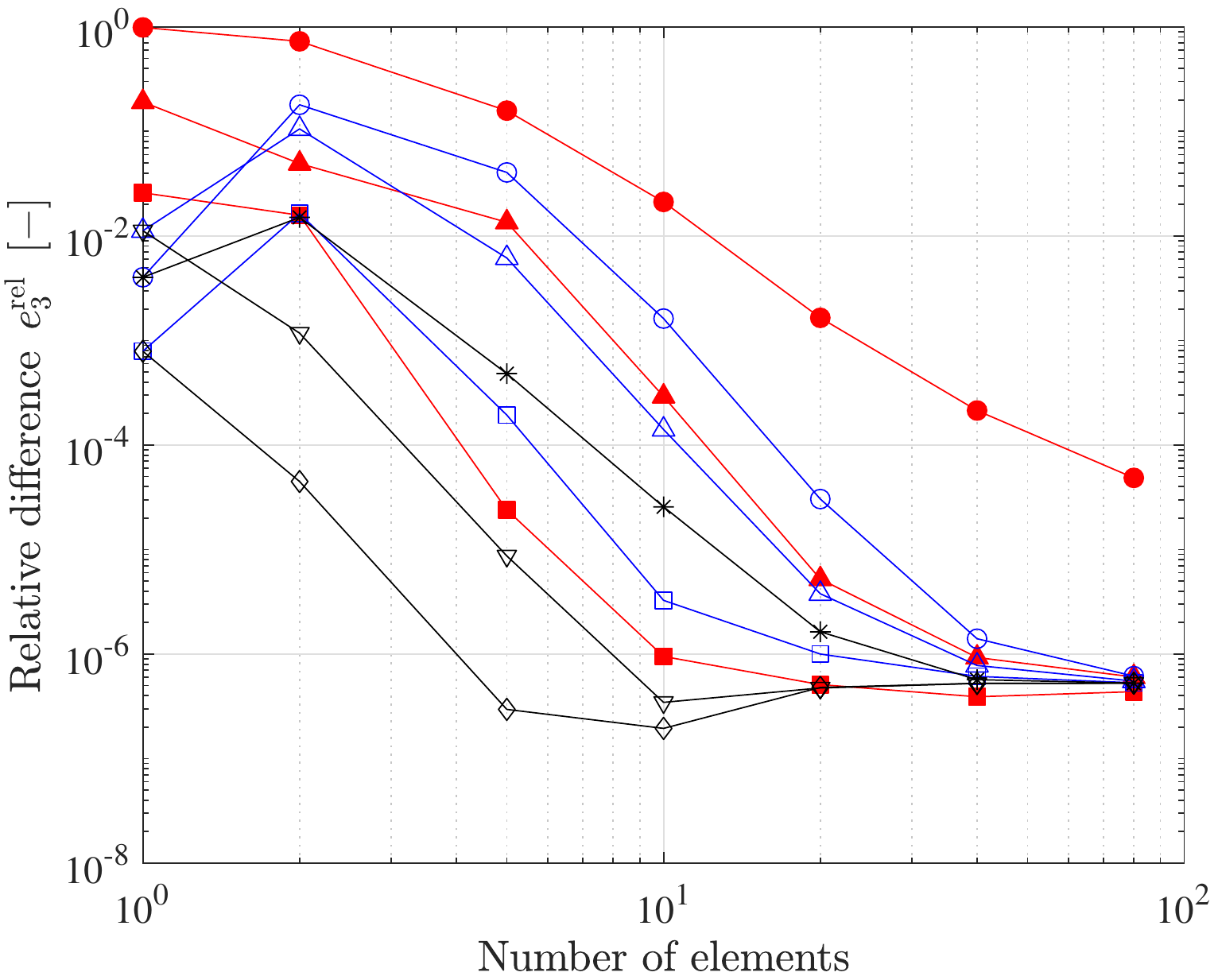}	
		\caption{$Z$-displacement (per element accuracy)}
	\end{subfigure}
	\begin{subfigure}[b] {0.4455\textwidth} \centering
		\includegraphics[width=\linewidth]{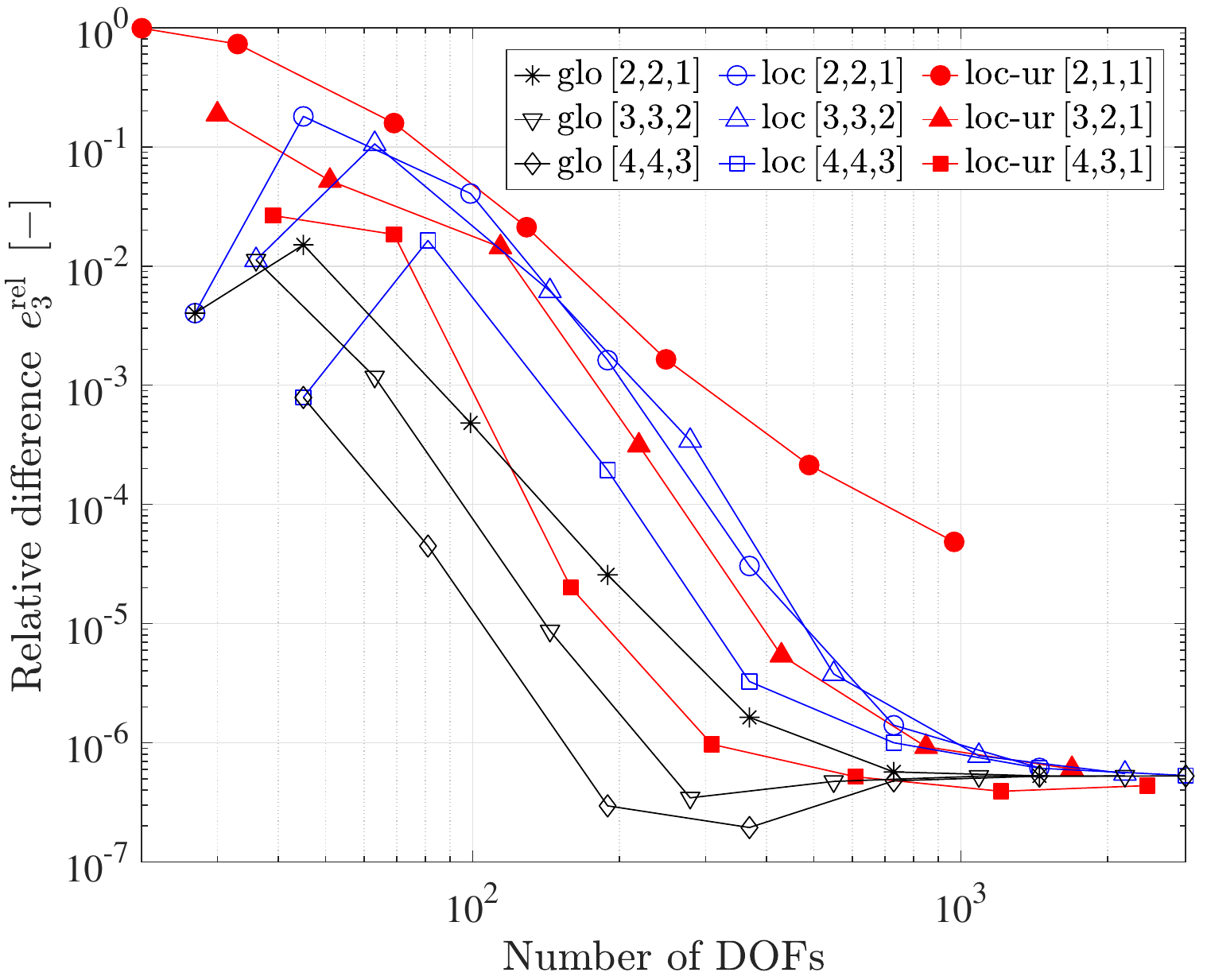}	
		\caption{$Z$-displacement (per DOF accuracy)}
	\end{subfigure}		
	\caption{$45^\circ$-arc cantilever beam (Case 1, \textcolor{black}{$R/d=10^2$}): Relative difference of the displacement at the tip between beam and brick element solutions. The latter is for a single quadratic B-spline element in the cross-section (i.e.,\,``IGA,\,$\mathrm{deg.}\!=\!(3,2,2)$, and $n_\mathrm{el}=320\times1\times1$''). The numbers in the bracket $[\bullet]$ denote the degrees, $[p,p_\mathrm{d},{p_\mathrm{p}}]$. All the beam solutions are obtained by the IGA-based mixed formulation. In the local approach ``loc-ur'', for $p=3,4$, we use $p_\mathrm{p}=p-n_\mathrm{el}$, if $n_\mathrm{el}<{p-1}$. Otherwise, we use $p_\mathrm{p}=1$.}
	\label{45deg_conv_disp}
\end{figure}

\begin{figure}[H]	
	\centering
	\begin{subfigure}[b] {0.329\textwidth} \centering
		\includegraphics[width=\linewidth]{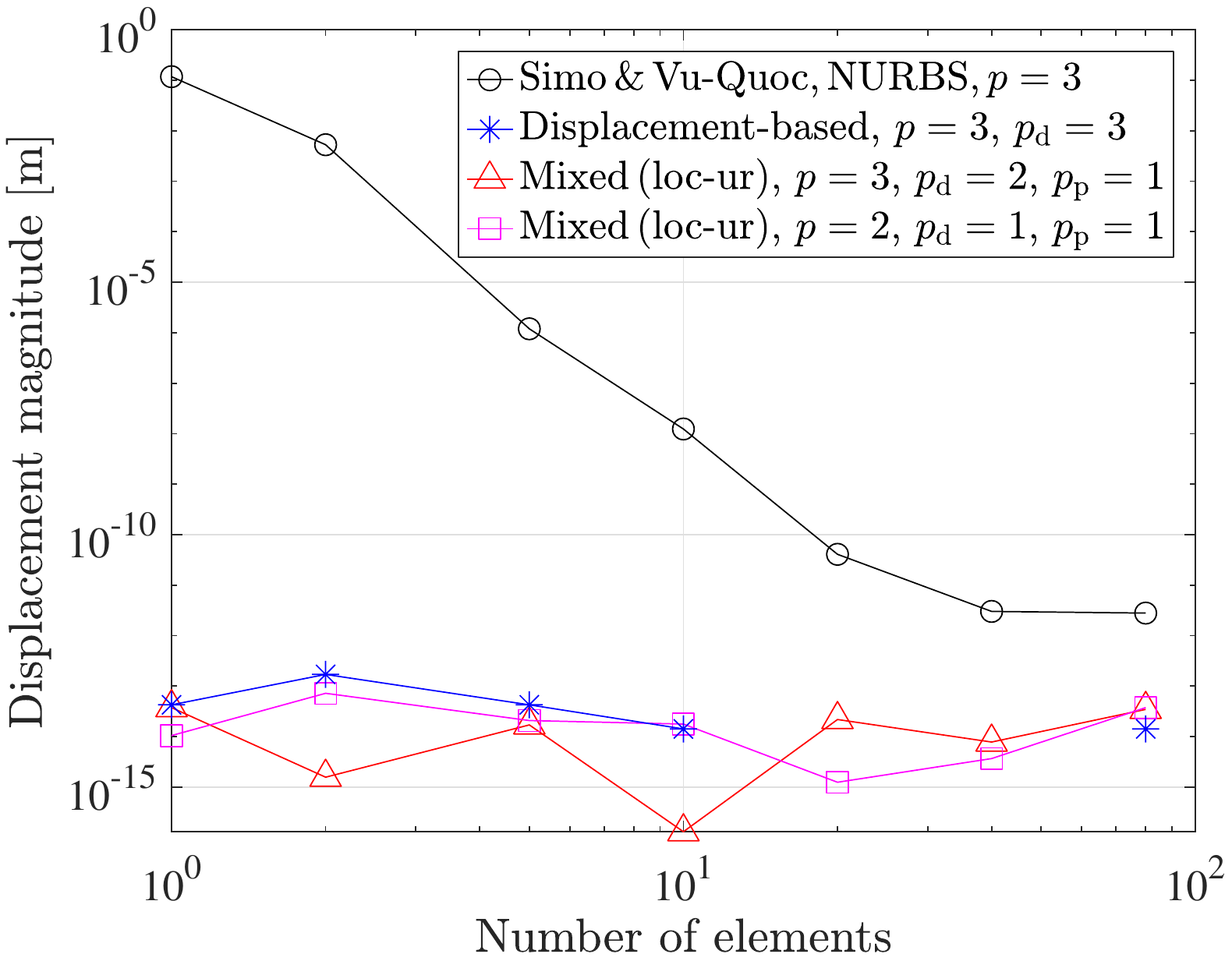}	
		\caption{$X$-displacement}		
	\end{subfigure}		
	\begin{subfigure}[b] {0.329\textwidth} \centering
		\includegraphics[width=\linewidth]{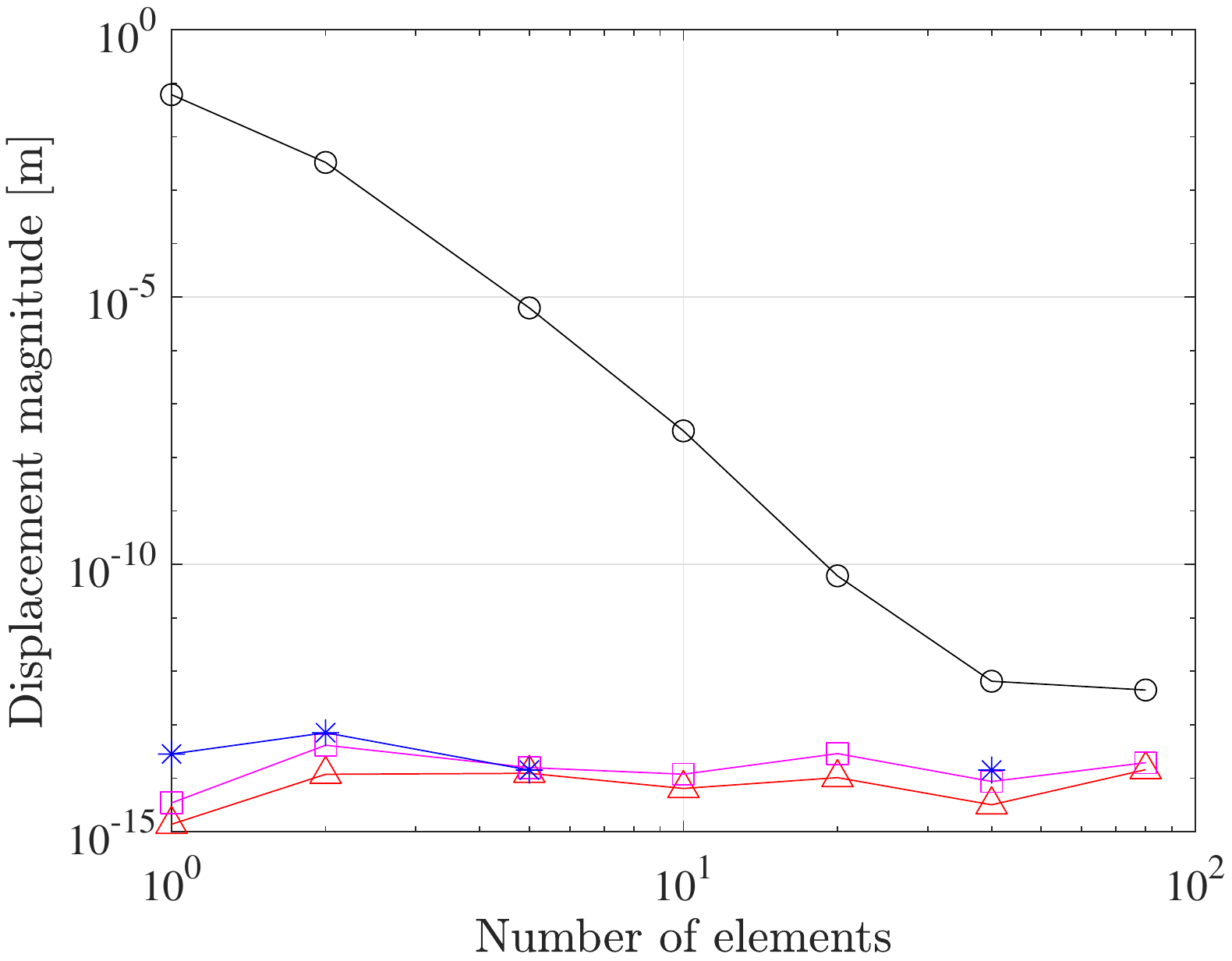}	
		\caption{$Y$-displacement}		
	\end{subfigure}		
	\begin{subfigure}[b] {0.329\textwidth} \centering
		\includegraphics[width=\linewidth]{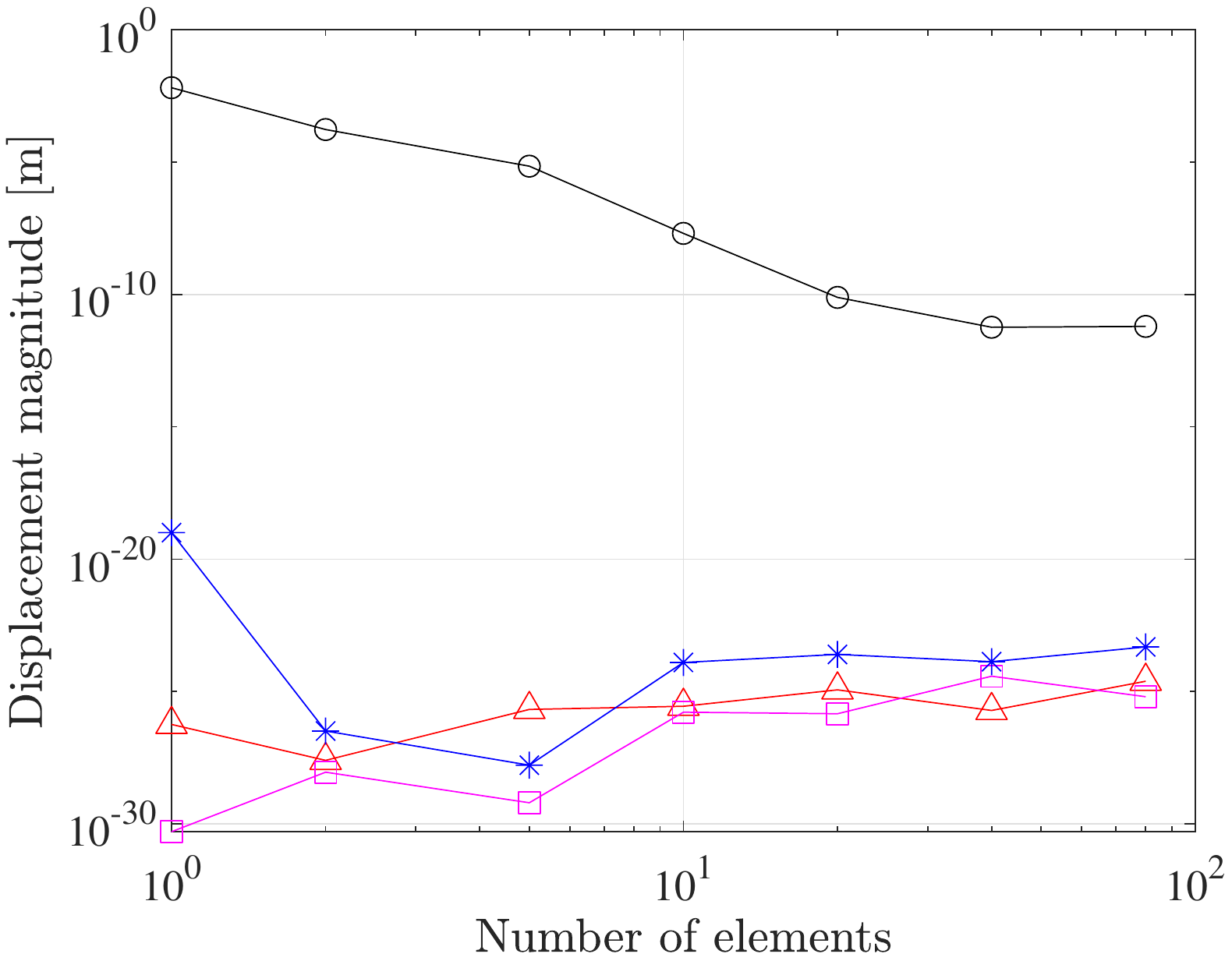}	
		\caption{$Z$-displacement}		
	\end{subfigure}			
	\caption{$45^\circ$-arc cantilever beam (Case 1\textcolor{black}{, $R/d=10^2$}): Displacement at the tip (point A) after unloading. All results are obtained by IGA. The missing data points in the blue curves in (a) and (b) are due to their values being exact zeros.}
	\label{45deg_disp_unload}
\end{figure}

\noindent Further, Fig.\,\ref{45deg_conv_disp} compares the convergence behavior between different finite element approximations of the independent solution fields, i.e., the global or local approach, and different degree of bases $p$, $p_\mathrm{d}$, and $p_\mathrm{p}$. We utilize the relative difference in the tip displacement, 
\begin{equation}
	\label{cant45_rel_error_disp_tip}
	{e^\mathrm{rel}_i} \coloneqq \left| \dfrac{{u^\mathrm{beam}_i - {u^\mathrm{brick}_i}}}{{u^\mathrm{brick}_i}} \right|,\,\,i\in\left\{1,2,3\right\},
\end{equation}
where $u^\mathrm{beam}_i$ and $u^\mathrm{brick}_i$ denote the displacement component $u_i\coloneqq\boldsymbol{u}\cdot{\boldsymbol{e}_i}$ of the beam and brick elements, respectively. The black curves show the results from the global approach of the IGA-based mixed formulation (``glo''), which exhibits the highest level of accuracy among all the shown results, but is computationally very expensive. The blue curves show the results from the local approach of the IGA-based mixed formulation using $p_\mathrm{p}=p-1$ (``loc''), with $p_\mathrm{d}=p$. It is the same formulation as the global approach for a single element ($n_\mathrm{el}=1$), so their results coincide then. However, as we increase the number of elements to $n_\mathrm{el}=2$, which has an interface between elements, ``loc'' suffers from a serious locking due to the $C^{p-1}$ continuity condition in the displacement field, see the drastic increase of the relative difference in all components $e^\mathrm{rel}_1$, $e^\mathrm{rel}_2$, and $e^\mathrm{rel}_3$, as we increase $n_\mathrm{el}$ from 1 to 2. The red curves show the beam solutions with the IGA-based mixed formulation using ${p_\mathrm{p}}=1$ (``loc-ur''), with ${p_\mathrm{d}}=p-1$. 
\begin{observation}	
	\label{cant45_observ_rigid}
	\textcolor{black}{	
	It should be noted that, in this example, we need at least quadratic NURBS to exactly represent the initial geometry of the center axis (circular arc), i.e., $p_\mathrm{g}=2$, so that we need $p_\mathrm{d}\ge{p_\mathrm{g}=2}$ to represent rigid body rotations exactly, see Remark \ref{rem_red_pd_curv}. The inability to represent rigid body rotations in case of $p_\mathrm{d}=1<p_\mathrm{g}$, which is seen in Table \ref{cant45deg_case1_eval_p234}, leads to artificial bending stiffness. This error can be alleviated by increasing the number of elements; however the convergence rate is significantly lowered, see the results of red curves with circular markers (``loc-ur'') in Fig.\,\ref{45deg_conv_disp}, compared with the results from using $p_\mathrm{d}=2$ (blue curves with circular markers).}
\end{observation}
\noindent\textcolor{black}{In cases of degrees $p=3,4$, with $p_\mathrm{d}=p-1$, the rigid body rotation can be represented exactly (see Table \ref{cant45deg_case1_eval_p234}), and it gives significantly improved per DOF accuracy (red curves with triangle and square markers) in Fig.\,\ref{45deg_conv_disp}, showing much lower error, compared with the local approach (``loc'').} 
\begin{table}[H]
	\centering	
	\caption{\textcolor{black}{$45^\circ$-arc cantilever beam (Case 1, $R/d=10^2$): Smallest six eigenvalues, ${\omega_i}^2$ [$\mathrm{rad}^2/\mathrm{s}^2$], $i=1,2,\cdots,6$. The bold-faced eigenvalues represent inability to represent rigid body rotations for $p_\mathrm{d}=p-1=1$, while all the other cases of $p$ and $p_\mathrm{d}$ correctly show six vanishing eigenvalues. We use $p_\mathrm{p}=1$, i.e., the approach ``IGA, loc-ur'', with $n_\mathrm{el}=10$.}}
	\label{cant45deg_case1_eval_p234}	
	\begin{tabular}{cccccccc}
		\toprule
		$p$   & $p_\mathrm{d}$ & ${\omega_1}^2$ & ${\omega_2}^2$ & ${\omega_3}^2$ & ${\omega_4}^2$ & ${\omega_5}^2$ & ${\omega_6}^2$ \\
		\midrule
		2     & 2     & -6.5E-12 & -3.4E-12 & 9.0E-14 & 8.7E-12 & 9.3E-12 & 6.4E-11 \\
		2     & 1     & -3.1E-12 & -1.3E-12 & 5.0E-12 & \textbf{6.2E-06} & \textbf{8.7E-05} & \textbf{3.7E-03} \\
		3     & 2     & -9.5E-12 & -9.1E-13 & 4.3E-12 & 6.0E-12 & 2.3E-11 & 3.4E-11 \\
		4     & 3     & -1.3E-10 & -2.7E-11 & -1.3E-11 & -2.0E-12 & 1.1E-11 & 3.0E-11 \\
		\bottomrule
	\end{tabular}%
\end{table}
\textcolor{black}{In order to verify the alleviation of locking, we investigate the tip displacement in $Z$-direction, as increasing the beam's slenderness ratio, $R/d$. The applied end force is inversely proportional to the bending stiffness, according to Eq.\,(\ref{cant45deg_app_end_f}), see Table \ref{app_num_ex_cant45_ext_f} for the force in each case of the slenderness ratio. Table \ref{cant45deg_allev_lock_compare_glo_loc} compares the $Z$-displacement at the tip, from two different approaches in IGA, the global (``glo''), and local approaches (``loc''), both with $p_\mathrm{d}=p$, and $p_\mathrm{p} = p-1$. In the results from the global approach, for every degree $p$, it is seen that the tip displacement converges, as we increase the slenderness ratio. In contrast, in the results of the local approach ``loc'' with lower degrees $p=2,3,4$ exhibit significantly decreased displacements, as increasing the slenderness ratio, due to locking. For higher degree $p$, the locking is much less pronounced. 
\begin{table}[H]
	\centering
	\caption{$45^\circ$-arc cantilever beam (Case 1): \textcolor{black}{$Z$-displacement (unit:\,$\mathrm{m}$) at the tip (point A). In all cases, we use a single load step, and $n_\mathrm{el}=10$. In the results of ``IGA, loc'', the significant decrease in the bold-faced values, as increasing the slenderness ratio, indicates the locking.}}
	\label{cant45deg_allev_lock_compare_glo_loc}		
	\begin{tabular}{ccccccrcccc}
	\toprule
	\multirow{2}[4]{*}{$p$} & \multirow{2}[4]{*}{$p_\mathrm{d}$} & \multicolumn{4}{c}{IGA, glo ($p_\mathrm{d}=p$, $p_\mathrm{p}=p-1$)} &       & \multicolumn{4}{c}{IGA, loc ($p_\mathrm{d}=p$, $p_\mathrm{p}=p-1$)} \\
	\cmidrule{3-6}\cmidrule{8-11}      &       & $R/d=10^2$ & $10^3$ & $10^4$ & $10^5$ &       & $10^2$ & $10^3$ & $10^4$ & $10^5$ \\
	\midrule
	2     & 2     & 53.4759 & 53.4426 & 53.4425 & 53.4425 &       & \textbf{53.3906} & \textbf{50.5364} & \textbf{42.9073} & \textbf{33.4118} \\
	3     & 3     & 53.4773 & 53.4687 & 53.4686 & 53.4686 &       & 53.4590 & \textbf{52.4263} & \textbf{45.1947} & \textbf{42.1914} \\
	4     & 4     & 53.4773 & 53.4687 & 53.4686 & 53.4686 &       & 53.4771 & 53.4645 & \textbf{53.2287} & \textbf{50.6505} \\
	5     & 5     & 53.4773 & 53.4687 & 53.4686 & 53.4686 &       & 53.4772 & 53.4686 & 53.4663 & \textbf{53.3805} \\
	6     & 6     & 53.4773 & 53.4687 & 53.4686 & 53.4686 &       & 53.4772 & 53.4687 & 53.4685 & 53.4668 \\
	7     & 7     & 53.4773 & 53.4687 & 53.4686 & 53.4686 &       & 53.4772 & 53.4687 & 53.4686 & 53.4685 \\
	8     & 8     & 53.4773 & 53.4687 & 53.4686 & 53.4686 &       & 53.4773 & 53.4687 & 53.4686 & 53.4686 \\
	9     & 9     & 53.4773 & 53.4687 & 53.4686 & 53.4686 &       & 53.4773 & 53.4687 & 53.4686 & 53.4686 \\
	10    & 10     & 53.4773 & 53.4687 & 53.4686 & 53.4686 &       & 53.4773 & 53.4687 & 53.4686 & 53.4686 \\
	\bottomrule
	\end{tabular}%
\end{table}
\noindent Table \ref{cant45deg_allev_lock_compare_ur_sr} shows the results from the reduction of $p_\mathrm{d}$, and $p_\mathrm{p}$. As mentioned in Observation \ref{cant45_observ_rigid}, the selection of degree $p_\mathrm{d}=p-1=1$ for $p=2$ suffers from the inability to represent rigid body rotations, which accentuates the locking. Thus, in the results of local approaches ``loc-ur'' and ``loc-sr'', using $p=2$ and $p_\mathrm{d}=p-1=1$ gives much smaller displacements than those from using $p_\mathrm{d}=p=2$. This is more pronounced, as the slenderness ratio increases. It is also seen that the local approach ``loc-ur'' with $p_\mathrm{p}=1$ (left-side of Table \ref{cant45deg_allev_lock_compare_ur_sr}) further alleviates locking, compared with the results from using $p_\mathrm{p}=p-1$ (``loc'', right-side of Table \ref{cant45deg_allev_lock_compare_glo_loc}). However, the decrease in the displacement, as increasing the slenderness ratio, is still clearly seen for lower degrees $p$ (bold-faced values). Further, as we discussed in Section \ref{ex_end_mnt_kref_vs_p}, numerical instability occurs for very high degrees $p$ in the approach ``loc-ur'' with $p_\mathrm{p}=1$, e.g., Newton-Raphson solution process even diverges for $p=10$. It is noticeable that the local approach ``loc-sr'' with selectively reduced degrees $p_\mathrm{p}$, given by Table \ref{tab_selec_deg_r}, effectively alleviates numerical instability and locking. 
\begin{table}[H]
	\centering
	\caption{$45^\circ$-arc cantilever beam (Case 1): \textcolor{black}{$Z$-displacement (unit:\,$\mathrm{m}$) at the tip (point A). In all cases, we use a single load step, and $n_\mathrm{el}=10$. In the results from using ``loc-ur'', the significant decrease in the bold-faced values, as increasing the slenderness ratio, indicates the locking. The Newton-Raphson iteration using the local approach (``loc-ur'') diverges for $p=10$, due to numerical instability.}}
	\label{cant45deg_allev_lock_compare_ur_sr}		
	\begin{tabular}{ccccccrcccc}
		\toprule
		\multirow{2}[4]{*}{$p$} & \multirow{2}[4]{*}{$p_\mathrm{d}$} & \multicolumn{4}{c}{IGA, loc-ur ($p_\mathrm{p}=1$)} &       & \multicolumn{4}{c}{IGA, loc-sr ($p_\mathrm{p}$ is given by Table \ref{tab_selec_deg_r})} \\
		\cmidrule{3-6}\cmidrule{8-11}      &       & $R/d=10^2$ & $10^3$ & $10^4$ & $10^5$ &       & $10^2$ & $10^3$ & $10^4$ & $10^5$ \\
		\midrule
		2     & 2     & \textbf{53.3906} & \textbf{50.5364} & \textbf{42.9073} & \textbf{33.4118} &       & 53.4758 & 53.4671 & 53.4671 & 53.4671 \\
		2     & 1     & \textbf{52.3462} & \textbf{42.4174} & \textbf{16.0015} & \textbf{0.7206} &       & 53.2293 & 52.3571 & 52.1100 & 52.1074 \\
		3     & 2     & 53.4606 & \textbf{52.4450} & \textbf{46.3224} & \textbf{44.8302} &       & 53.4848 & 53.4762 & 53.4761 & 53.4761 \\
		4     & 3     & 53.4772 & 53.4682 & \textbf{53.4612} & \textbf{53.3084} &       & 53.4818 & 53.4731 & 53.4730 & 53.4730 \\
		5     & 4     & 53.4773 & 53.4687 & 53.4683 & \textbf{53.4558} &       & 53.4824 & 53.4733 & 53.4730 & 53.4696 \\
		6     & 5     & 53.4774 & 53.4687 & 53.4686 & 53.4684 &       & 53.4829 & 53.4733 & 53.4730 & 53.4724 \\
		7     & 6     & 53.4776 & 53.4688 & 53.4687 & 53.4687 &       & 53.4838 & 53.4743 & 53.4733 & 53.4731 \\
		8     & 7     & 53.4780 & 53.4706 & 53.4688 & 53.4687 &       & 53.4774 & 53.4687 & 53.4686 & 53.4686 \\
		9     & 8     & 53.4774 & 53.4913 & 53.4755 & 53.4689 &       & 53.4773 & 53.4687 & 53.4686 & 53.4686 \\
		10    & 9     & -     & -     & -     & -     &       & 53.4773 & 53.4687 & 53.4686 & 53.4686 \\
		\bottomrule
	\end{tabular}%
\end{table}
}
\textcolor{black}{
The numerical instability in the formulation ``loc-ur'', observed in the left-side of Table\,\ref{cant45deg_allev_lock_compare_ur_sr}, is indicated very well from the results of the generalized eigenvalue analysis at the initial (undeformed) configuration, which is shown in the left-side of Table \ref{cant45deg_eval_loc_ur_sr_p2_10}. It is seen that, in the results from using ``loc-ur'' with $p_\mathrm{p}=1$, the eigenvalues for $p=8,9,10$ are much smaller than those from using other degrees, see bold-faced numbers in the left-side of Table \ref{cant45deg_eval_loc_ur_sr_p2_10}. This turns out to be alleviated by a selective adjustment of the degree $p_\mathrm{p}$ by Table \ref{tab_selec_deg_r}, see the eigenvalue analysis results in the right-side of Table \ref{cant45deg_eval_loc_ur_sr_p2_10}. However, the selection of $p_\mathrm{p}$ in Table \ref{tab_selec_deg_r} may not be optimal, and we may still observe spuriously decreased eigenvalues, e.g., $\omega_{10}$ in the result from using $p=9$. A further mathematical investigation to determine the optimal $p_\mathrm{p}$ remains future work. We also observe that the eigenvalues in the result from using the local approach ``loc-ur'' with $p_\mathrm{d}=p=2$ are spuriously higher than the others, due to locking. This turns out to be alleviated by the selective reduction of $p_\mathrm{p}$, i.e., the local approach ``loc-sr''. Further, eigenvalues from using $p_\mathrm{d}=p-1=1$ ($p=2$) are spuriously higher, in both local approaches ``loc-ur'', and ``loc-sr'', due to the inability to represent rigid body rotations (see Observation \ref{cant45_observ_rigid}).
}
\begin{table}[H]
	\centering
	\caption{\textcolor{black}{$45^\circ$-arc cantilever beam (Case 1, $R/d=10^2$): The smallest four nonzero natural (angular) frequencies, $\omega_i$ [$\mathrm{rad}/\mathrm{s}$], $i=7,8,9,10$. The bold-faced ones represent spurious eigenvalues due to numerical instability.
In all cases, $n_\mathrm{el}=10$.}}
	\label{cant45deg_eval_loc_ur_sr_p2_10}
	\begin{tabular}{ccccccrcccc}
		\toprule
		\multirow{2}[4]{*}{$p$} & \multirow{2}[4]{*}{$p_\mathrm{d}$} & \multicolumn{4}{c}{IGA, loc-ur, $p_\mathrm{p}=1$} &       & \multicolumn{4}{c}{IGA, loc-sr, $p_\mathrm{p}$ is given in Table \ref{tab_selec_deg_r}} \\
		\cmidrule{3-6}\cmidrule{8-11}      &       & $\omega_7$ & $\omega_8$ & $\omega_9$ & $\omega_{10}$ &       & $\omega_7$ & $\omega_8$ & $\omega_9$ & $\omega_{10}$ \\
		\midrule
		2     & 2     & 3.2514 & 9.1246 & 9.3017 & 15.4063 &       & 3.2377 & 8.8688 & 9.0037 & 14.5618 \\
		2     & 1     & 3.2631 & 9.0594 & 9.2475 & 14.9950 &       & 3.2569 & 9.0450 & 9.1958 & 14.9650 \\
		3     & 2     & 3.2387 & 8.8851 & 9.0238 & 14.6197 &       & 3.2381 & 8.8809 & 9.0159 & 14.6109 \\
		4     & 3     & 3.2379 & 8.8786 & 9.0115 & 14.5997 &       & 3.2380 & 8.8784 & 9.0136 & 14.6024 \\
		5     & 4     & 3.2379 & 8.8778 & 9.0106 & 14.5990 &       & 3.2380 & 8.8785 & 9.0137 & 14.6028 \\
		6     & 5     & 3.2377 & 8.8750 & 9.0079 & 14.5960 &       & 3.2380 & 8.8780 & 9.0137 & 14.6009 \\
		7     & 6     & 3.2368 & 8.8608 & 8.9944 & 14.5824 &       & 3.2380 & 8.8664 & 9.0137 & 14.5846 \\
		8     & 7     & \textbf{3.2301} & \textbf{8.7626} & \textbf{8.9097} & \textbf{14.4326} &       & 3.2380 & 8.8806 & 9.0136 & 14.6056 \\
		9     & 8     & \textbf{3.1713} & \textbf{6.6820} & \textbf{6.8214} & \textbf{7.2754} &       & 3.2379 & 8.8808 & 9.0007 & \textbf{10.0067} \\
		10    & 9     & \textbf{1.1305} & \textbf{1.1416} & \textbf{1.5804} & \textbf{1.6422} &       & 3.2380 & 8.8809 & 9.0137 & 14.6063 \\
		\bottomrule
	\end{tabular}%
\end{table}
Table \ref{cant45deg_comp_nlstep_niter} shows that the presented mixed formulation allows much larger load steps, and eventually requires much fewer number of iterations, compared to various reference formulations, including the displacement-based one (i.e.,\,``IGA, displacement-based''). 
\begin{table}[H]
	\centering
	\caption{$45^\circ$-arc cantilever beam (Case 1\textcolor{black}{, $R/d=10^2$}): Comparison of the number of load steps and Newton-Raphson iterations. The reference data from literature are obtained from Table 2 of \citet{frischkorn2013solid}.}
	\label{cant45deg_comp_nlstep_niter}
	\begin{tabular}{lcc}
		\toprule
		& \#load steps & \#iterations \\
		\midrule		
		IGA, mixed (loc-ur, EAS), $p=4$, $n_\mathrm{el}=80$ & 1     & 8 \\
		IGA, displacement-based, $p=4$, $n_\mathrm{el}=80$ & 1     & 18 \\
		IGA, brick, $\mathrm{deg.}\!=\!(3,1,1)$, $n_\mathrm{el}\!=\!320\!\times\!1\!\times\!1$ & 1     & 18 \\
		IGA, brick, $\mathrm{deg.}\!=\!(3,2,2)$, $n_\mathrm{el}\!=\!320\!\times\!1\!\times\!1$ & 1     & 18 \\
		\citet{frischkorn2013solid} & 1     & 19 \\
		\citet{wackerfuss2009mixed} & 1     & 9 \\
		\citet{rhim1998vectorial} & 1     & 13 \\
		\citet{crisfield1990consistent} & 3     & 17 \\
		\citet{cardona1988beam}   & 6     & 47 \\
		\citet{dvorkin1988non}  & 10    & 34 \\
		\citet{simo1986three}    & 3     & 27 \\
		\citet{bathe1979large}    & 60    & - \\
		\bottomrule
	\end{tabular}%
\end{table}
We further show the path-independence of the proposed beam formulation employing the finite element approximation of the total displacements in the extensible directors. We investigate if the displacement at the tip (point A) vanishes after unloading. \textcolor{black}{Here, we consider the slenderness ratio $R/d=10^2$.} We also compare the results with those from the beam formulation of \citet{simo1986three}. In the present paper, we approximate their field variables, displacement of the center axis, and the incremental rotations by NURBS basis functions, in the framework of IGA, see \citet{choi2019isogeometric} for more details. Fig.\,\ref{45deg_disp_unload} shows the $X$-, $Y$- and $Z$-directional displacements at the tip after unloading. The beam formulation based on \citet{simo1986three} exhibits significant displacement error, which decreases with element number (black curves). In contrast, the proposed beam formulation yields vanishing displacements to machine precision, in both the displacement-based and the mixed formulations. 
\subsubsection{Case 2: Nonlinear material with Poisson effect}
\begin{figure}[H]	
	\centering
	\begin{subfigure}[b] {0.455\textwidth} \centering
		\includegraphics[width=\linewidth]{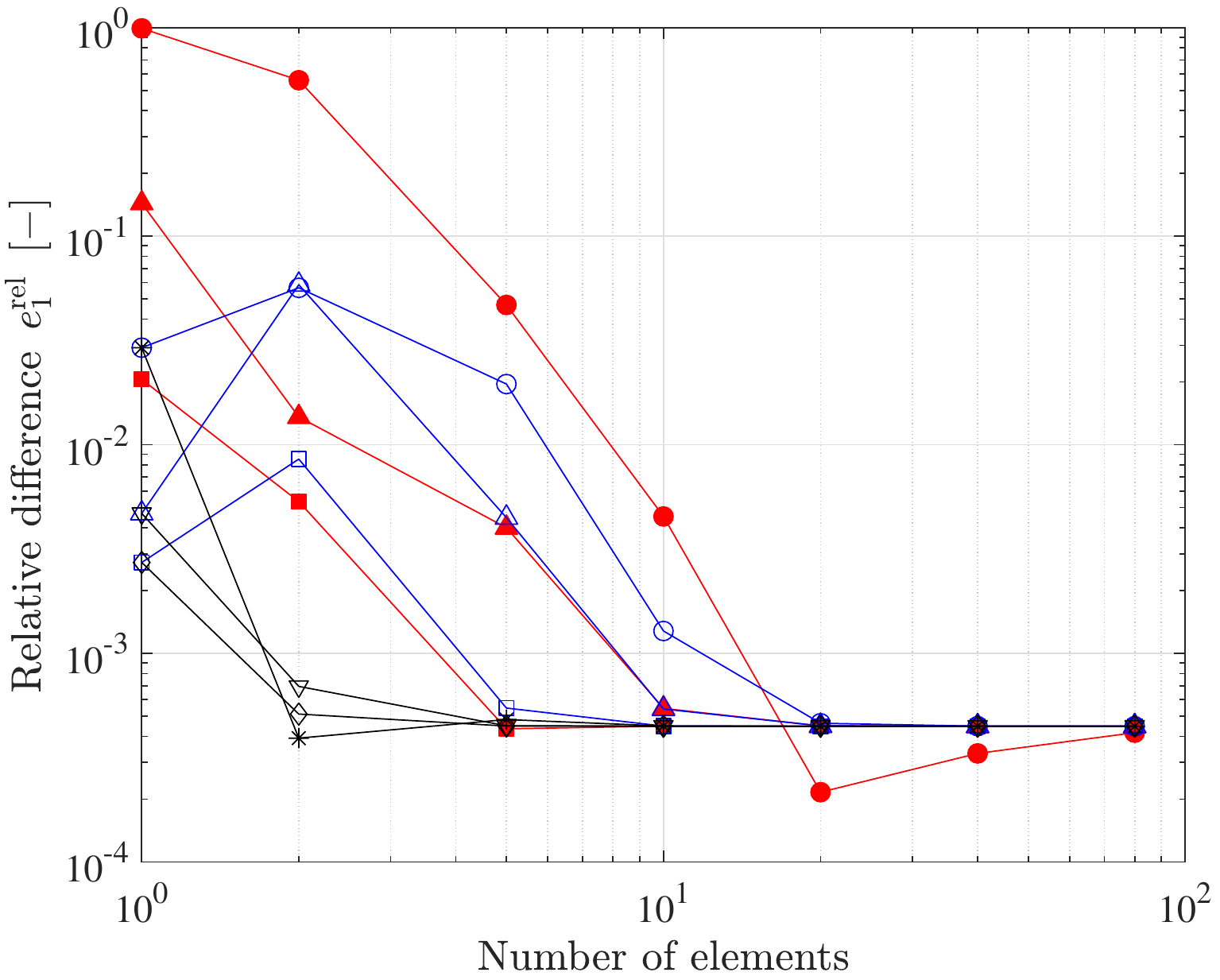}	
		\caption{$X$-displacement (per element accuracy)}
	\end{subfigure}
	\begin{subfigure}[b] {0.455\textwidth} \centering
		\includegraphics[width=\linewidth]{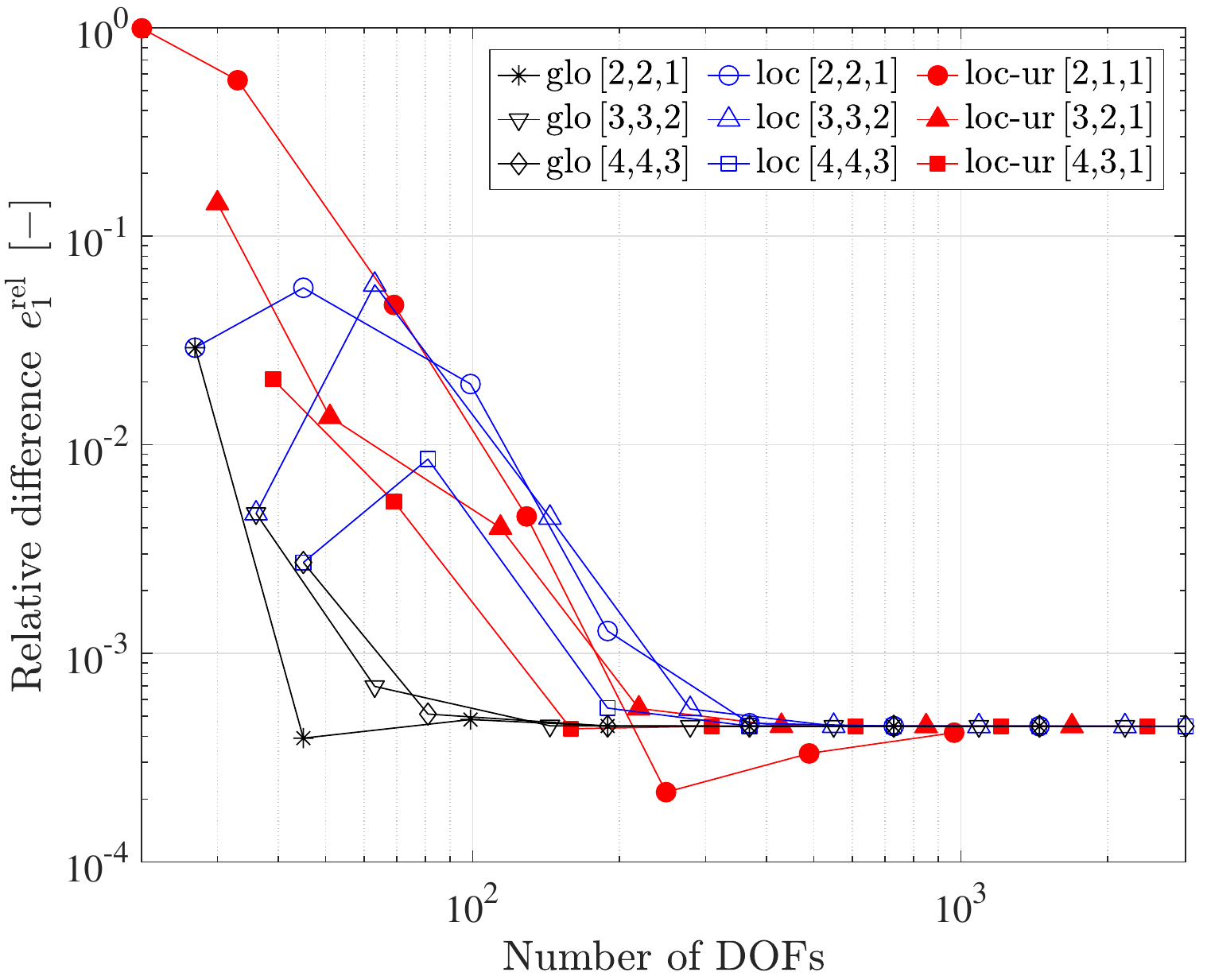}	
		\caption{$X$-displacement (per DOF accuracy)}
	\end{subfigure}
	\begin{subfigure}[b] {0.455\textwidth} \centering
		\includegraphics[width=\linewidth]{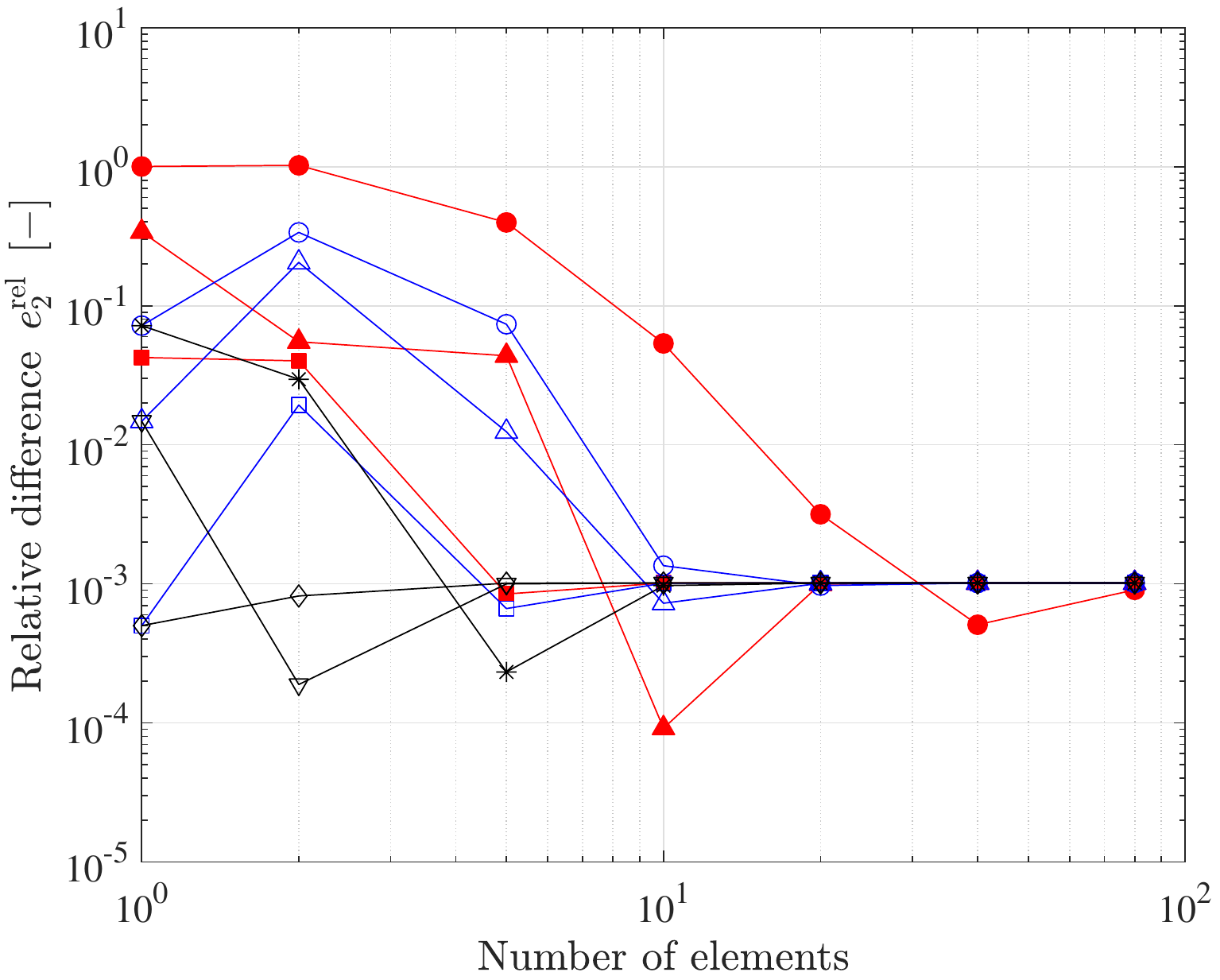}	
		\caption{$Y$-displacement (per element accuracy)}
	\end{subfigure}
	\begin{subfigure}[b] {0.455\textwidth} \centering
		\includegraphics[width=\linewidth]{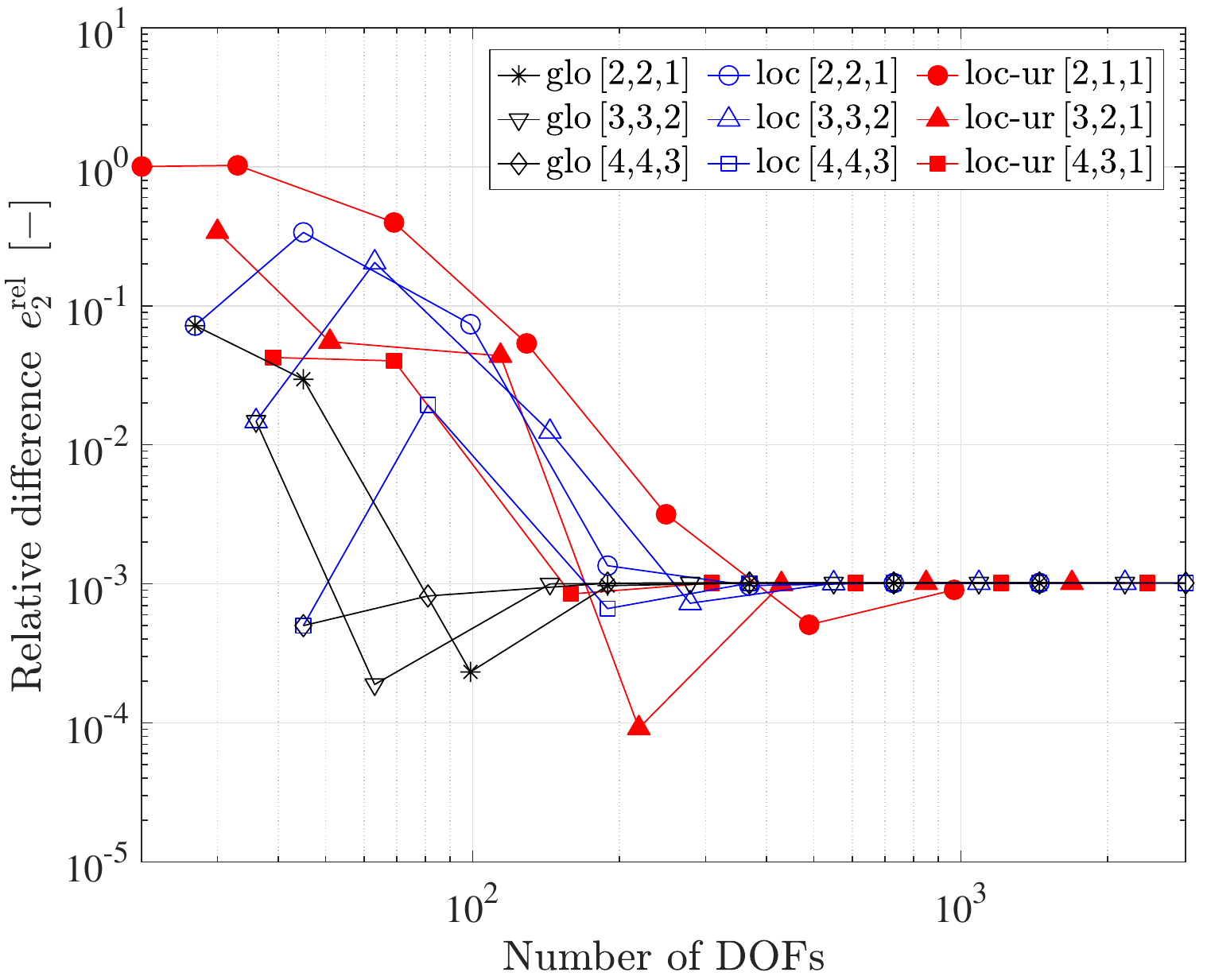}	
		\caption{$Y$-displacement (per DOF accuracy)}
	\end{subfigure}		
	\begin{subfigure}[b] {0.455\textwidth} \centering
		\includegraphics[width=\linewidth]{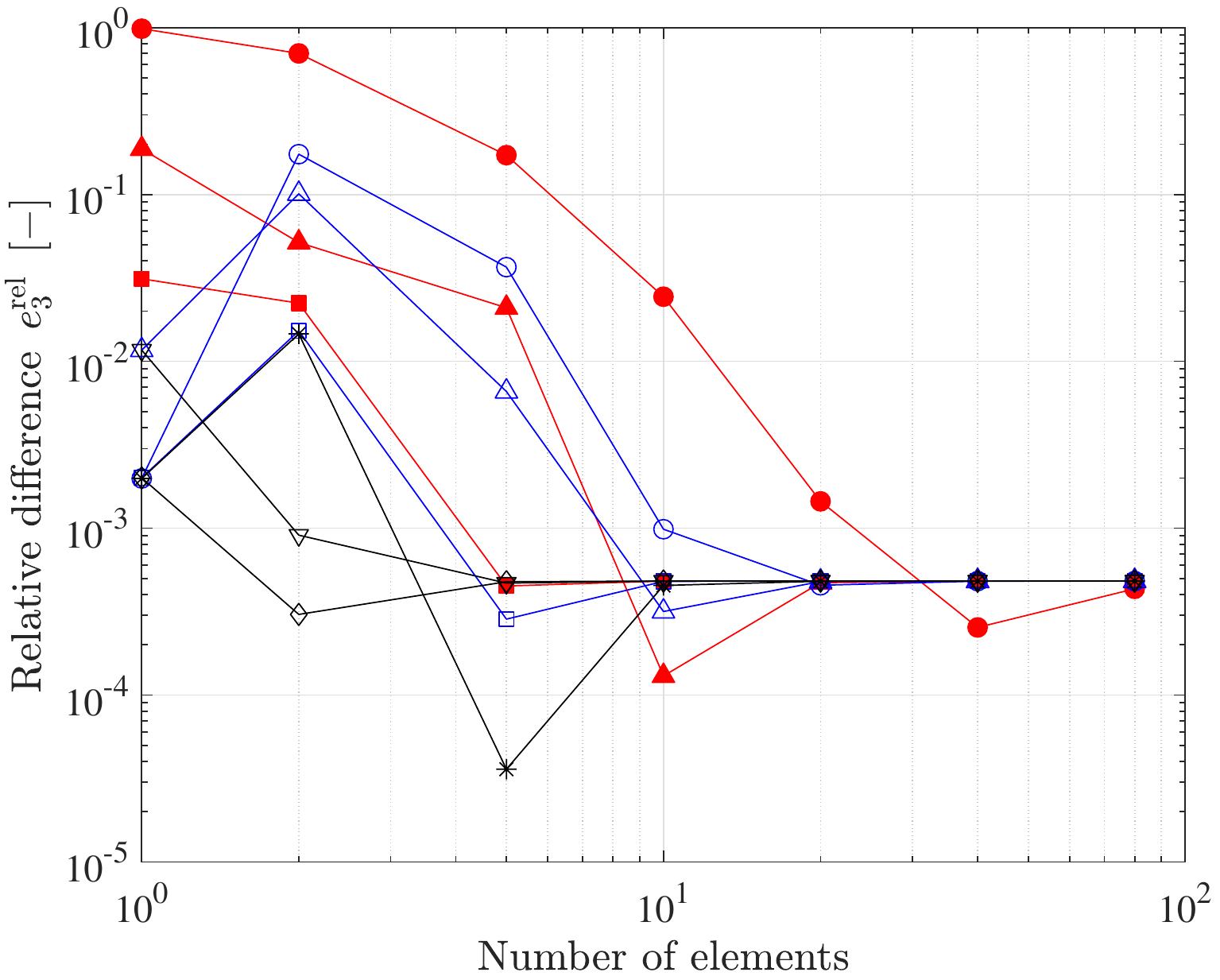}	
		\caption{$Z$-displacement (per element accuracy)}
	\end{subfigure}
	\begin{subfigure}[b] {0.455\textwidth} \centering
		\includegraphics[width=\linewidth]{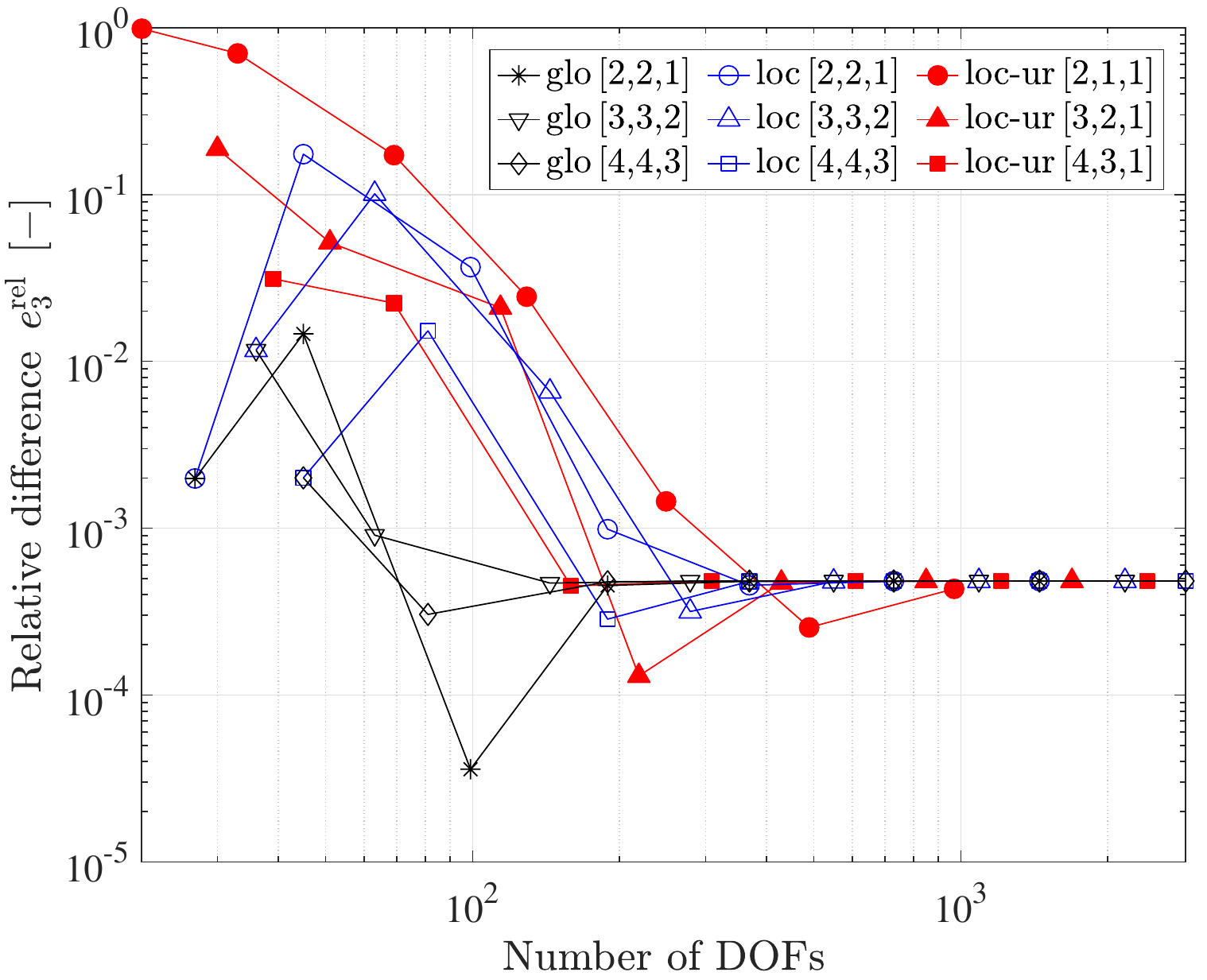}	
		\caption{$Z$-displacement (per DOF accuracy)}
	\end{subfigure}		
	\caption{$45^\circ$-arc cantilever beam (Case 2): Relative difference of the displacement at the tip between beam and brick element solution. The latter is for a single quadratic B-spline element in the cross-section (i.e.\,``IGA,\,$\mathrm{deg.}\!=\!(3,2,2)$, and $n_\mathrm{el}=320\times1\times1$''). The numbers in the bracket $[\bullet]$ denotes the degrees, $[p,p_\mathrm{d},{p_\mathrm{p}}]$. All the beam solutions are obtained by the IGA-based mixed formulation. In the local approach ``loc-ur'', for $p=3,4$, we use $p_\mathrm{p}=p-n_\mathrm{el}$, if $n_\mathrm{el}<{p-1}$. Otherwise, we use $p_\mathrm{p}=1$.}
	\label{45deg_case2_conv_disp}
\end{figure}

In Case 2, we consider a compressible Neo-Hookean type material with Poisson's ratio $\nu=0.3$, and Young's modulus $E\!=\!{10}\,\mathrm{MPa}$. The same boundary conditions as in Case 1 are considered. Table \ref{cant45deg_case2_comp_tip_disp_brick} compares the displacement at the tip (point A in Fig.\,\ref{45deg_init_config}), with various brick element solutions. Severe artificial stiffening is observed in the brick element solution with a single linear quadrilateral element in the cross-section (``IGA, brick, $\mathrm{deg.}\!=(3,1,1)$, $n_\mathrm{el}\!=\!320\times{1}\times{1}$''),
due to the lack of linear transverse normal strain, coupled with the bending deformation. The beam solution without any enrichment of the cross-sectional strains (``IGA, mixed, loc-ur'') suffers from even larger stiffening, since it lacks not only the quadratic but also the bilinear in-plane displacement field (i.e.,\,trapezoidal cross-section deformations). After enriching the in-plane cross-sectional strains by the EAS method, the beam solution (``IGA, mixed, loc-ur, EAS'') shows a good agreement with the brick element solution based on a single quadratic element in the cross-section. Its deviation from the brick element solution using a single cubic element in the cross-section (``IGA, brick, $\mathrm{deg.}\!=\!(3,3,3)$, $n_\mathrm{el}\!=\!320\times{1}\times{1}$'') shows the significance of the cross-sectional warping, which our current beam formulation does not account for. This remains future work. 
\begin{table}[H]
	\centering
	\caption{$45^\circ$-arc cantilever beam (Case 2): Comparison of the tip displacements.}
	\label{cant45deg_case2_comp_tip_disp_brick}
	\begin{tabular}{lccc}
		\toprule
		& $u_1\,[\mathrm{m}]$ & $u_2\,[\mathrm{m}]$ & $u_3\,[\mathrm{m}]$ \\
		\midrule
		IGA, brick, $\mathrm{deg.}\!=(3,1,1)$, $n_\mathrm{el} = 320\times{1}\times{1}$ & 12.5816 & -21.6881 & 51.5348 \\
		{IGA, brick, $\mathrm{deg.}\!=(3,2,2)$, $n_\mathrm{el} = 320\times{1}\times{1}$} & {13.8087} & {-23.9547} & {53.6642} \\
		IGA, brick, $\mathrm{deg.}\!=(3,3,3)$, $n_\mathrm{el} = 320\times{1}\times{1}$ & 13.9947 & -24.2610 & 53.8291 \\
		IGA, mixed (loc-ur), $p=4$, ${p_\mathrm{d}}=3$, $p_\mathrm{p} = 1$, $n_\mathrm{el}=80$ & 11.3090 & -19.3285 & 49.1188 \\
		{IGA, mixed (loc-ur, EAS), $p=4$, ${p_\mathrm{d}=3}$, $p_\mathrm{p} = 1$, $n_\mathrm{el}=80$} & {13.8148} & {-23.9790} & {53.6901} \\
		\bottomrule
	\end{tabular}%
\end{table}
\noindent In Fig.\,\ref{45deg_case2_conv_disp}, we compare the convergence of the relative difference in the tip displacements, $e^\mathrm{rel}_i$ in Eq.\,(\ref{cant45_rel_error_disp_tip}), between several different finite element approximations for independent solution fields. As observed in Case 1 (Observation \ref{cant45_observ_rigid}), the results from using $p=2$, and $p_\mathrm{d}=p-1=1$ exhibit severe locking due to the inability to represent the rigid body rotations (red curves with circular markers), in comparison to the results from the same degree $p=2$ with higher $p_\mathrm{d}=p=2$ (blue curves with circular markers). However, in the cases of higher $p=3,4$ and $p_\mathrm{d}=p-1$ (red curves with triangle and square markers), the accuracy improvement is noticeable.  

\subsection{L-shape frame}
\begin{figure}[H]
	\centering
	\begin{subfigure}[b] {0.4\textwidth} \centering
		\includegraphics[width=\linewidth]{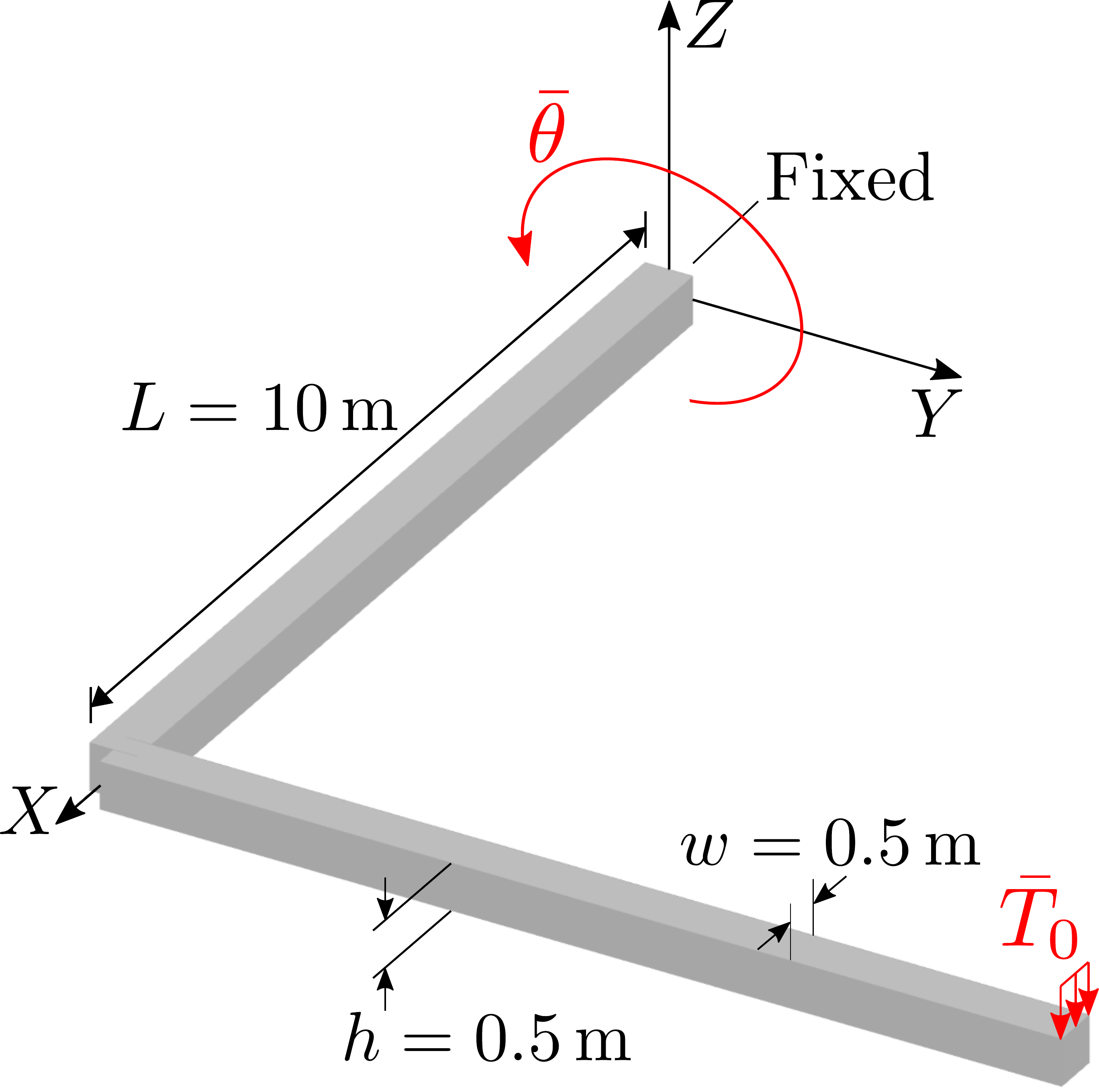}	
		\caption{}
	\end{subfigure}		
	\begin{subfigure}[b] {0.4\textwidth} \centering
		\includegraphics[width=\linewidth]{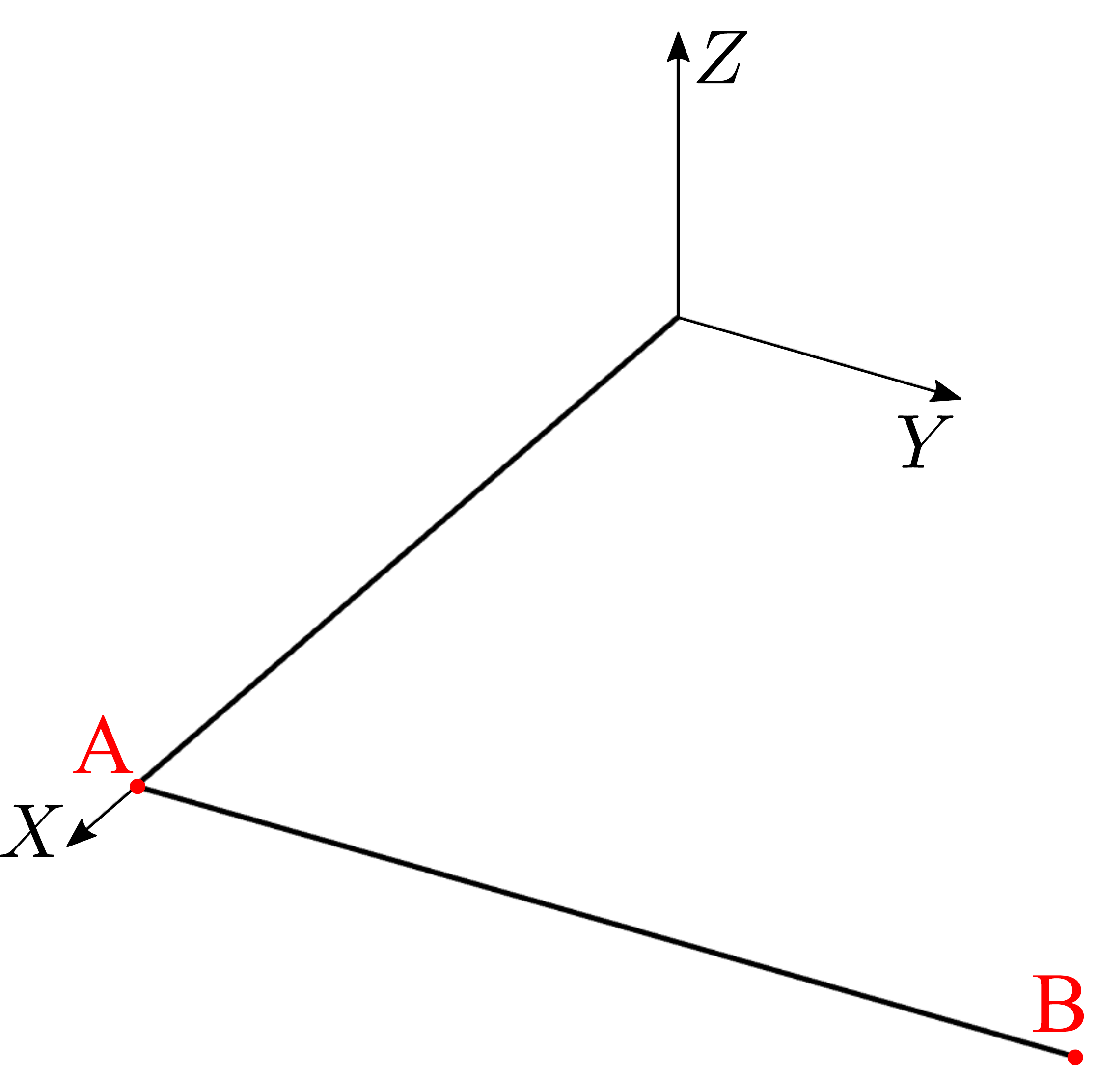}	
		\caption{}
	\end{subfigure}	
	\caption{L-shape frame. (a) Undeformed configuration and boundary conditions. (b) Center axes of two connected beams, where A and B indicate selected points where displacements are evaluated to test path-independence.}
	\label{Lshape_init_config}
\end{figure}

In this example, we show the path-independence of the presented mixed isogeometric beam formulation, considering two beams connected by a rotational continuity condition. We consider a L-shape frame with two straight beams having the initial length $L=10\,\mathrm{m}$ and a square cross-section of dimension $h=w=0.5\,\mathrm{m}$. We select compressible Neo-Hookean type material with Young's modulus $E=10\,\mathrm{MPa}$ and Poisson's ratio $\nu=0.3$. A distributed force $\bar{T}_0 = -200\mathrm{N}/\mathrm{m}$ is applied in $Z$-direction on the upper edge of the end face, and the other end is fixed, see Fig.\,\ref{Lshape_init_config}. The load is applied in five load steps with uniform increments. In order to test the path-independence of the beam formulation, we impose a prescribed rotation at the fixed end, keeping the applied load fixed, such that the force maintains the negative $Z$-direction, and does not rotate along $\bar \theta$. This procedure is inspired by a similar test in \citet{ibrahimbegovic2002role}. The rotation is prescribed about the $X$-axis, and we choose the total rotation angle $\bar {\theta}=200\pi$, i.e., 100 turns around the $X$-axis. In the additional prescribed rotation, we use $10^3$ load steps in total, such that each full ($360^\circ$) turn is uniformly divided by 10 increments, i.e., $36^\circ$ rotation is prescribed in each load step. Here, all the results are obtained by the proposed isogeometric mixed formulation (IGA, ``loc-ur'') with $p=3$, $p_\mathrm{d}=2$, $p_\mathrm{p}=1$. We also use $n_\mathrm{el}=5$ for each beam, so that we have 10 elements in total. 

\begin{figure}[H]
	\centering
	\begin{subfigure}[b] {0.4875\textwidth} \centering
		\includegraphics[width=\linewidth]{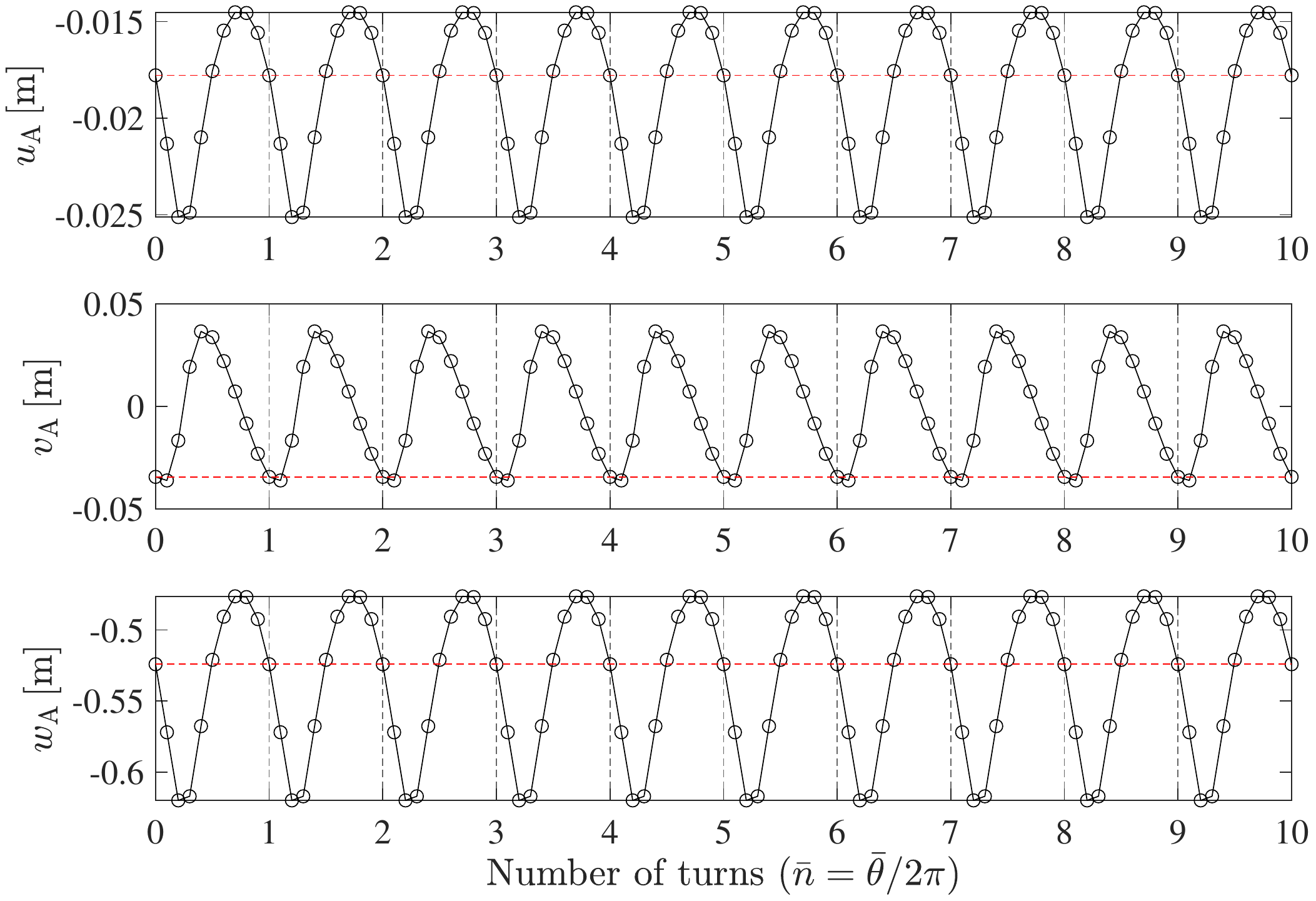}	
		\caption{Displacement at A}
	\end{subfigure}		
	\begin{subfigure}[b] {0.475\textwidth} \centering
		\includegraphics[width=\linewidth]{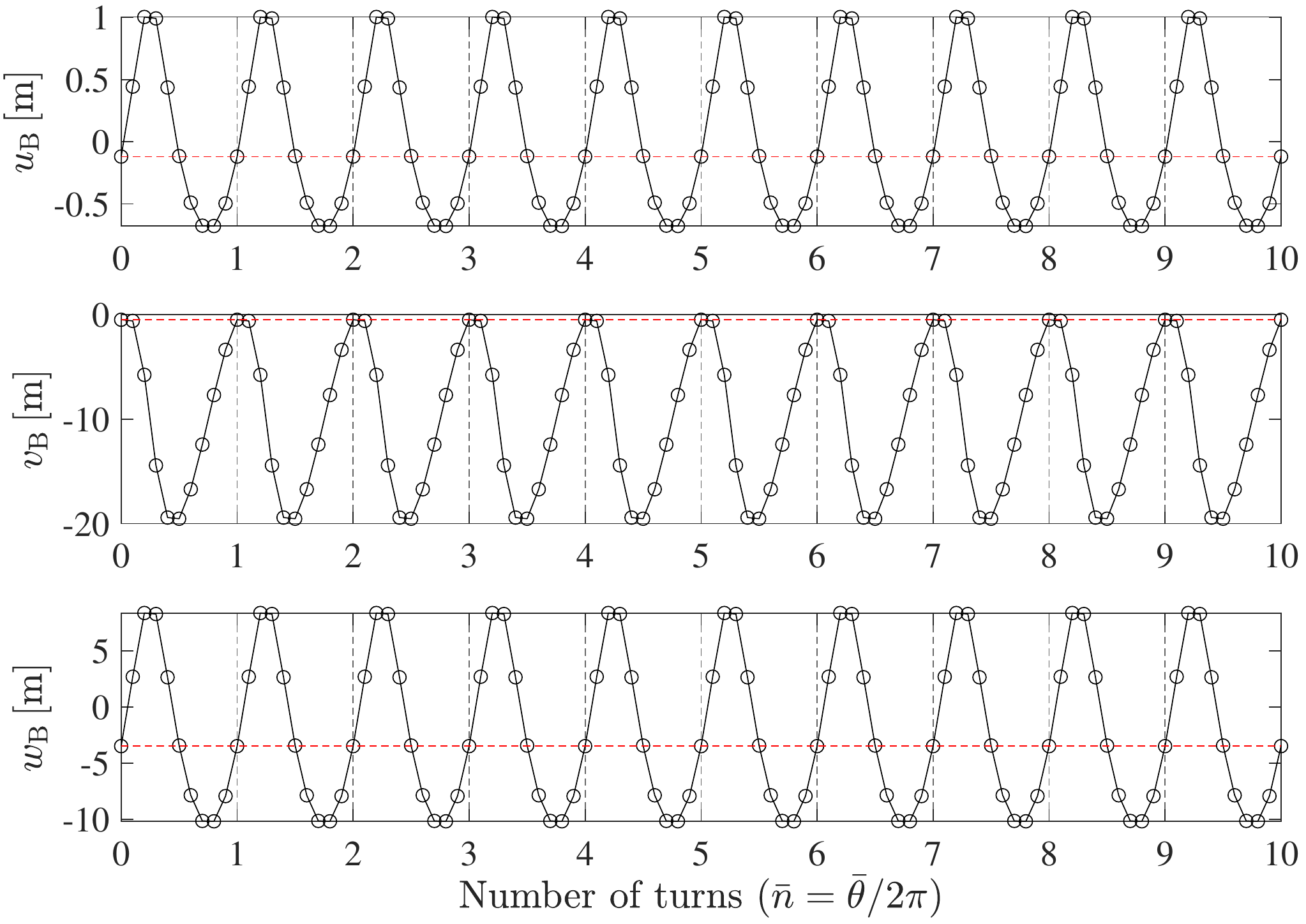}	
		\caption{Displacement at B}
	\end{subfigure}	
	\caption{L-shape frame: Displacements at points A and B due to the additional prescribed rotation. The horizontal dashed line marks the displacements at the deformed configuration by the applied force ${\bar T}_0$ and $\bar \theta=0$, and the vertical dashed lines mark every full turn, i.e. ${\bar n}\coloneqq{\bar \theta}/{2\pi}=1,2,...,10$.}
	\label{Lshape_disp_add_rot_pt_AB}
\end{figure}

\begin{figure}[H]
	\centering
	\begin{subfigure}[b] {0.4875\textwidth} \centering
		\includegraphics[width=\linewidth]{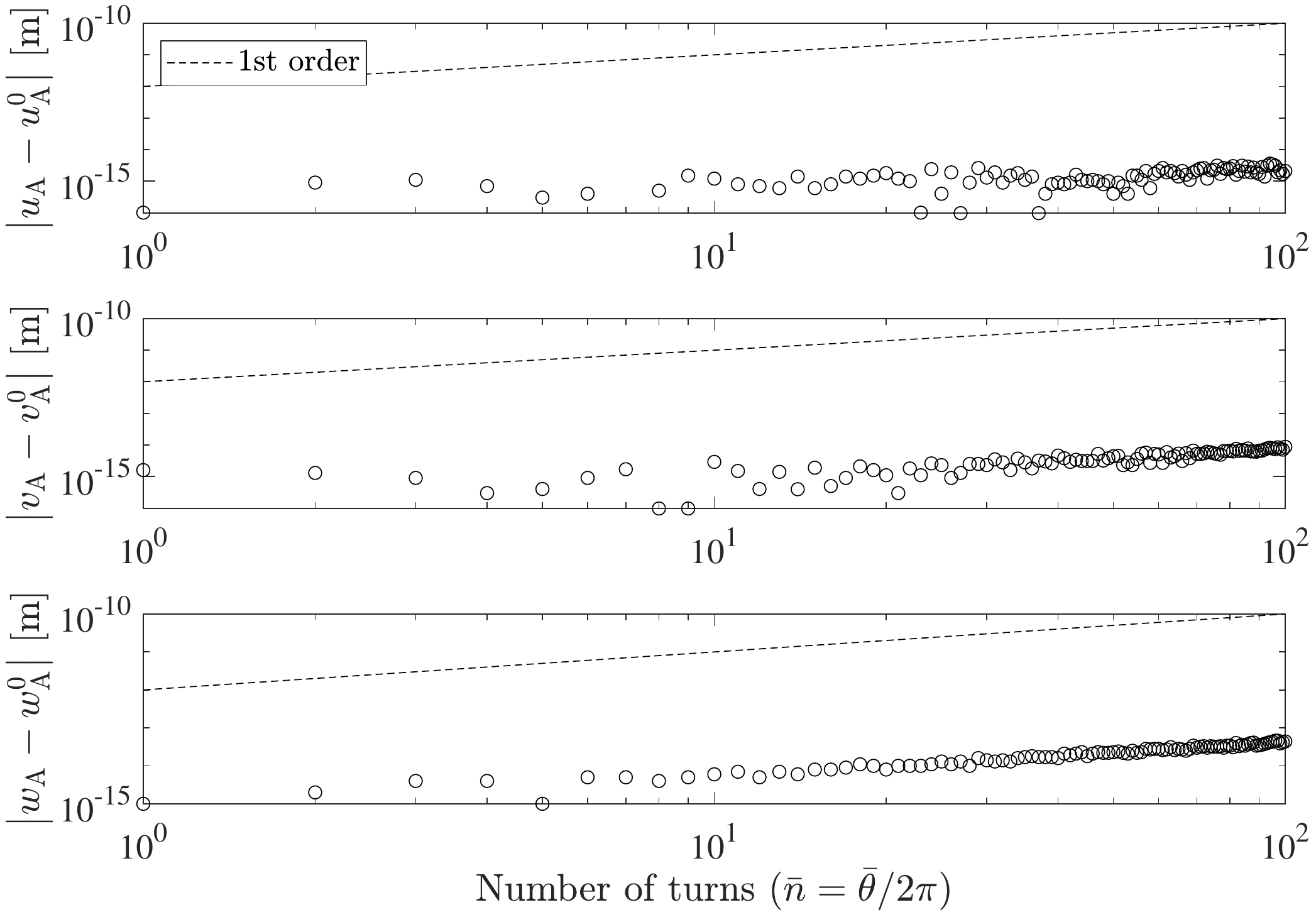}	
		\caption{Displacements at A}
	\end{subfigure}		
	\begin{subfigure}[b] {0.4875\textwidth} \centering
		\includegraphics[width=\linewidth]{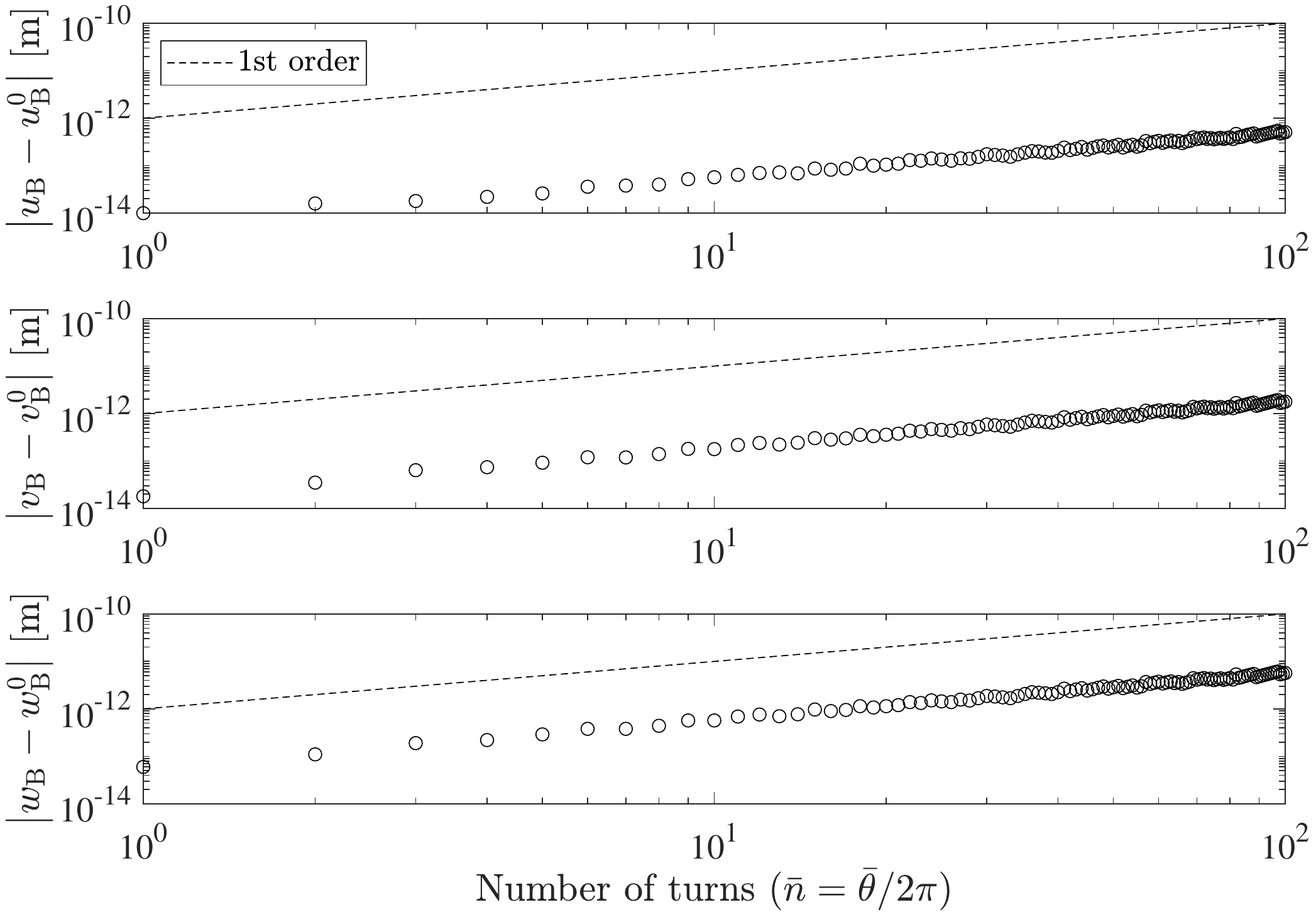}	
		\caption{Displacements at B}
	\end{subfigure}	
	\caption{L-shape frame. Difference of displacements at points A and B due to the additional prescribed rotations $\bar \theta = 2{\bar n}\pi$ with $\bar n=1,2,...,100$. The subscripts A and B denote the displacements at points A and B, respectively, and the superscript 0 denotes the displacements at the turn number ${\bar n}=0$.}
	\label{Lshape_disp_add_rot_pt_AB_1}
\end{figure}

Fig.\,\ref{Lshape_disp_add_rot_pt_AB} shows a periodic change of the displacements at the points A and B along the turns. Here, $u$, $v$, and $w$ indicate the $X$-, $Y$-, and $Z$-displacements, respectively. It is seen that the displacement at every full turn exhibits the same value, which demonstrates the path-independence of our beam formulation. For a more thorough investigation, in Fig.\,\ref{Lshape_disp_add_rot_pt_AB_1}, we plot the change of the displacements at points A and B along the 100 full turns. $u^0_\mathrm{A}$, $v^0_\mathrm{A}$, and $w^0_\mathrm{A}$ denote the displacements due to the applied force ${\bar T}_0$ at $\bar \theta = 0$. For all displacement components, the difference vanishes up to machine precision. A linearly increasing small error is due to the accumulated numerical error by the limited machine precision.

\section{Conclusion}
\label{conclu}
This paper presents an isogeometric mixed finite element formulation for nonlinear beams with extensible directors. We particularly investigate a selective reduction of the degree of bases for the independent solution fields of axis displacement, directors, and additional stress and strain fields. This approach effectively alleviates additional locking phenomena due to the higher order continuity of the displacement field in IGA. Based on this, we show the superior per DOF accuracy of IGA over conventional FEA. Further the developed method is computationally efficient, due to its local (element-wise) static condensation process, and the symmetry of the tangent stiffness matrix under conservative loads, and much less internal DOFs due to the reduced degree of bases. It is also shown that the mixed formulation yields improved convergence in the Newton-Raphson iteration in the thin beam limit, compared to the displacement-based one. Important extensions of this work are expected in the following directions:
\begin{itemize}
	\item Further investigation on the numerical stability in the proposed selection of degree $p_\mathrm{p}$, based on a generalized inf-sup test for multi-field variational principle, e.g., see \citet{krischok2019generalized}, 
	\item An optimal selection of degrees $p_\mathrm{p}$ for each component of the stress resultant and strain,
	\item Further enhancement of the out-of-plane cross-sectional strains including torsion-induced warping.
\end{itemize}

\nopagebreak[4] 
\appendix
\gdef\thesection{\Alph{section}} 
\makeatletter
\counterwithin{figure}{section}
\counterwithin{table}{section}
\renewcommand\@seccntformat[1]{\csname the#1\endcsname.\hspace{0.5em}}
\makeatother
\numberwithin{equation}{section}
\section{Appendix to the beam formulation}
\label{app_hypelas_conv_test}
\subsection{Operator expressions in the spatial discretization}
As we separately arrange the control coefficients for the center axis position, and the director displacement parts, some operator expressions need to be modified from those in \citet{choi2021isogeometric}, as follows.
\subsubsection{Discretization of the internal virtual work}
\label{app_oper_intvw}
In Eq.\,(\ref{submat_intf_disp}), we use the following operator
\begin{equation}
{\Bbb{B}}_{{\rm{total}}}^e\coloneqq {\left[ {\begin{array}{*{20}{c}}
				{{\Bbb{B}}_{\varepsilon \varphi }^1}& \cdots &{{\Bbb{B}}_{\varepsilon \varphi }^{{n_e}}}&\vline& {}&{{{\bf{0}}_{1 \times 6n_e^\mathrm{d}}}}&{}\\
				{{\Bbb{B}}_{\rho \varphi }^1}& \cdots &{{\Bbb{B}}_{\rho \varphi }^{{n_e}}}&\vline& {{\Bbb{B}}_{\rho\mathrm{d}}^1}& \cdots &{{\Bbb{B}}_{\rho\mathrm{d}}^{n_e^\mathrm{d}}}\\
				{}&{{{\bf{0}}_{3 \times 3{n_e}}}}&{}&\vline& {{\Bbb{B}}_{\kappa\mathrm{d}}^1}& \cdots &{{\Bbb{B}}_{\kappa\mathrm{d}}^{n_e^\mathrm{d}}}\\
				{{\Bbb{B}}_{\delta \varphi }^1}& \cdots &{{\Bbb{B}}_{\delta \varphi }^{{n_e}}}&\vline& {{\Bbb{B}}_{\delta\mathrm{d}}^1}& \cdots &{{\Bbb{B}}_{\delta\mathrm{d}}^{n_e^\mathrm{d}}}\\
				{}&{{{\bf{0}}_{4 \times 3{n_e}}}}&{}&\vline& {{\Bbb{B}}_{\gamma\mathrm{d}}^1}& \cdots &{{\Bbb{B}}_{\gamma\mathrm{d}}^{n_e^\mathrm{d}}}\\
				{}&{{{\bf{0}}_{3 \times 3{n_e}}}}&{}&\vline& {{\Bbb{B}}_{\chi\mathrm{d}}^1}& \cdots &{{\Bbb{B}}_{\chi\mathrm{d}}^{n_e^\mathrm{d}}}
		\end{array}} \right]_{15 \times {m_e}}},
\end{equation}
such that
\begin{subequations}	
	\allowdisplaybreaks
	\begin{alignat}{3}
		\delta {\varepsilon ^h} &= {\left[ {{\boldsymbol{\varphi }}_{,s}^{\rm{T}}N_{I,s}^p} \right]}\delta {{\boldsymbol{\varphi }}^e_I} &&\eqqcolon {\Bbb{B}}_{\varepsilon \varphi }^I\delta {{\boldsymbol{\varphi }}^e_I},\\
		\left\{ {\begin{array}{*{20}{c}}
				{\delta \rho _1^h}\\
				{\delta \rho _2^h}
		\end{array}} \right\} &= {\left[ {\begin{array}{*{20}{c}}
					{{\boldsymbol{d}}_{1,s}^{\rm{T}}N_{I,s}^p}\\
					{{\boldsymbol{d}}_{2,s}^{\rm{T}}N_{I,s}^p}
			\end{array}} \right]}\delta {{\boldsymbol{\varphi }}^e_I} + {\left[ {\begin{array}{*{20}{c}}
					{{\boldsymbol{\varphi }}_{,s}^{\rm{T}}N_{J,s}^{{p_\mathrm{d}}}}&{{{\bf{0}}_{1 \times 3}}}\\
					{{{\bf{0}}_{1 \times 3}}}&{{\boldsymbol{\varphi }}_{,s}^{\rm{T}}N_{J,s}^{{p_\mathrm{d}}}}
			\end{array}} \right]}\left\{ {\begin{array}{*{20}{c}}
				{\delta {{\bf{d}}^e_{1J}}}\\
				{\delta {{\bf{d}}^e_{2J}}}
		\end{array}} \right\} &&\eqqcolon {\Bbb{B}}_{\rho \varphi }^I\delta {{\boldsymbol{\varphi}}^e_I} + {\Bbb{B}}_{\rho \mathrm{d}}^J\delta {{\bf{d}}^e_J},\\
		\left\{ {\begin{array}{*{20}{c}}
				{\delta \kappa _{11}^h}\\
				{\delta \kappa _{22}^h}\\
				{2\delta \kappa _{12}^h}
		\end{array}} \right\} &= \left[ {\begin{array}{*{20}{c}}
				{{\boldsymbol{d}}_{1,s}^{\rm{T}}N_{J,s}^{{p_{\rm{d}}}}}&{{{\bf{0}}_{1 \times 3}}}\\
				{{{\bf{0}}_{1 \times 3}}}&{{\boldsymbol{d}}_{2,s}^{\rm{T}}N_{J,s}^{{p_{\rm{d}}}}}\\
				{{\boldsymbol{d}}_{2,s}^{\rm{T}}N_{J,s}^{{p_{\rm{d}}}}}&{{\boldsymbol{d}}_{1,s}^{\rm{T}}N_{J,s}^{{p_{\rm{d}}}}}
		\end{array}} \right]\left\{ {\begin{array}{*{20}{c}}
				{\delta {{\bf{d}}^e_{1J}}}\\
				{\delta {{\bf{d}}^e_{2J}}}
		\end{array}} \right\} &&\eqqcolon {\Bbb{B}}_{\kappa \mathrm{d}}^J\delta {{\bf{d}}^e_J},\\
		\left\{ {\begin{array}{*{20}{c}}
				{\delta \delta _1^h}\\
				{\delta \delta _2^h}
		\end{array}} \right\} &= \left[ {\begin{array}{*{20}{c}}
				{{\boldsymbol{d}}_1^{\rm{T}}N_{I,s}^p}\\
				{{\boldsymbol{d}}_2^{\rm{T}}N_{I,s}^p}
		\end{array}} \right]\delta {{\boldsymbol{\varphi }}^e_I} + \left[ {\begin{array}{*{20}{c}}
				{{\boldsymbol{\varphi }}_{,s}^{\rm{T}}N_J^{{p_{\rm{d}}}}}&{{{\bf{0}}_{1 \times 3}}}\\
				{{{\bf{0}}_{1 \times 3}}}&{{\boldsymbol{\varphi }}_{,s}^{\rm{T}}N_J^{{p_{\rm{d}}}}}
		\end{array}} \right]\left\{ {\begin{array}{*{20}{c}}
				{\delta {{\bf{d}}^e_{1J}}}\\
				{\delta {{\bf{d}}^e_{2J}}}
		\end{array}} \right\} &&\eqqcolon {\Bbb{B}}_{\delta \varphi }^I\delta {{\boldsymbol{\varphi }}^e_I} + {\Bbb{B}}_{\delta {\rm{d}}}^J\delta {{\bf{d}}^e_J},\\
		\left\{ {\begin{array}{*{20}{c}}
				{\delta \gamma _{11}^h}\\
				{\delta \gamma _{12}^h}\\
				{\delta \gamma _{21}^h}\\
				{\delta \gamma _{22}^h}
		\end{array}} \right\} &= \left[ {\begin{array}{*{20}{c}}
				{{\boldsymbol{d}}_1^{\rm{T}}N_{J,s}^{{p_{\rm{d}}}} + {\boldsymbol{d}}_{1,s}^{\rm{T}}N_J^{{p_{\rm{d}}}}}&{{{\bf{0}}_{1 \times 3}}}\\
				{{\boldsymbol{d}}_{2,s}^{\rm{T}}N_J^{{p_{\rm{d}}}}}&{{\boldsymbol{d}}_1^{\rm{T}}N_{J,s}^{{p_{\rm{d}}}}}\\
				{{\boldsymbol{d}}_2^{\rm{T}}N_{J,s}^{{p_{\rm{d}}}}}&{{\boldsymbol{d}}_{1,s}^{\rm{T}}N_J^{{p_{\rm{d}}}}}\\
				{{{\bf{0}}_{1 \times 3}}}&{{\boldsymbol{d}}_2^{\rm{T}}N_{J,s}^{{p_{\rm{d}}}} + {\boldsymbol{d}}_{2,s}^{\rm{T}}N_J^{{p_{\rm{d}}}}}
		\end{array}} \right]\left\{ {\begin{array}{*{20}{c}}
				{\delta {{\bf{d}}^e_{1I}}}\\
				{\delta {{\bf{d}}^e_{2I}}}
		\end{array}} \right\} &&\eqqcolon {\Bbb{B}}_{\gamma {\rm{d}}}^J\delta {{\bf{d}}^e_J},\\
		\left\{ {\begin{array}{*{20}{c}}
				{\delta \chi _{11}^h}\\
				{\delta \chi _{22}^h}\\
				{2\delta \chi _{12}^h}
		\end{array}} \right\} &= \left[ {\begin{array}{*{20}{c}}
				{{\boldsymbol{d}}_1^{\rm{T}}N_J^{{p_{\rm{d}}}}}&{{{\bf{0}}_{1 \times 3}}}\\
				{{{\bf{0}}_{1 \times 3}}}&{{\boldsymbol{d}}_2^{\rm{T}}N_J^{{p_{\rm{d}}}}}\\
				{{\boldsymbol{d}}_2^{\rm{T}}N_J^{{p_{\rm{d}}}}}&{{\boldsymbol{d}}_1^{\rm{T}}N_J^{{p_{\rm{d}}}}}
		\end{array}} \right]\left\{ {\begin{array}{*{20}{c}}
				{\delta {{\bf{d}}^e_{1J}}}\\
				{\delta {{\bf{d}}^e_{2J}}}
		\end{array}} \right\} &&\eqqcolon {\Bbb{B}}_{\chi {\rm{d}}}^J\delta {{\bf{d}}^e_J},
	\end{alignat}
\end{subequations}
where the repeated indices $I\in\left\{1,2,\cdots,n_e\right\}$ and $J\in\left\{1,2,\cdots,n^\mathrm{d}_e\right\}$ imply summations.
\subsubsection{Geometric part of the tangent stiffness matrix}
\label{app_oper_gtan}
In Eq.\,(\ref{gtan_uti_from_prev}), we utilize the following operators from \citet{choi2021isogeometric} with the original stress resultants replaced by the physical (independent) stress resultants.
\begin{equation}
	\label{geom_stiff_mat_k_g}
	{{\boldsymbol{k}}_{\mathrm{G}}} \coloneqq \left[ {\begin{array}{*{20}{c}}
			\begin{aligned}
				{{\boldsymbol{k}}_\varepsilon }\,\,\,\,\,\,\,{{\boldsymbol{k}}_\rho }\,\,\,\,\,{{\boldsymbol{k}}_\delta }\\
				{\,\,}\,\,\,\,\,\,\,{{\boldsymbol{k}}_\kappa }\,\,\,\,\,{{\boldsymbol{k}}_\gamma }\\
				{{\rm{sym}}{\rm{.}}}\,\,\,\,\,{\,\,\,\,\,}\,\,\,\,\,{{\boldsymbol{k}}_\chi }
			\end{aligned}
	\end{array}} \right]_{15\times15},
\end{equation}
with
\begingroup
\allowdisplaybreaks
\begin{subequations}
	\begin{alignat}{3}
		{\boldsymbol{k}}_\varepsilon  &\coloneqq {\tilde n}_\mathrm{p}{{\boldsymbol{1}}_{3\times3}},\\
		{\boldsymbol{k}}_\rho  &\coloneqq \left[ {\begin{array}{*{20}{c}}
				{{{\tilde m}^1_\mathrm{p}}{{\boldsymbol{1}}_{3\times3}}}&{{{\tilde m}^2_\mathrm{p}}{{\boldsymbol{1}}_{3\times3}}}
		\end{array}} \right],\\
		{\boldsymbol{k}}_\delta  &\coloneqq \left[ {\begin{array}{*{20}{c}}
				{{{\tilde q}^1_\mathrm{p}}{{\boldsymbol{1}}_{3\times3}}}&{{{\tilde q}^2_\mathrm{p}}{{\boldsymbol{1}}_{3\times3}}}
		\end{array}} \right],\\
		{\boldsymbol{k}}_\kappa  &\coloneqq \left[ {\renewcommand{\arraystretch}{1.5}\begin{array}{*{20}{c}}
				{{{\tilde h}^{11}_\mathrm{p}}{{\boldsymbol{1}}_{3\times3}}}&{{{\tilde h}^{12}_\mathrm{p}}{{\boldsymbol{1}}_{3\times3}}}\\
				{{\rm{sym}}{\rm{.}}}&{{{\tilde h}^{22}_\mathrm{p}}{{\boldsymbol{1}}_{3\times3}}}
		\end{array}} \right],\\
		{\boldsymbol{k}}_\gamma  &\coloneqq \left[ {\renewcommand{\arraystretch}{1.5}\begin{array}{*{20}{c}}
				{{{\tilde m}^{11}_\mathrm{p}}{{\boldsymbol{1}}_{3\times3}}}&{{{\tilde m}^{21}_\mathrm{p}}{{\boldsymbol{1}}_{3\times3}}}\\
				{{{\tilde m}^{12}_\mathrm{p}}{{\boldsymbol{1}}_{3\times3}}}&{{{\tilde m}^{22}_\mathrm{p}}{{\boldsymbol{1}}_{3\times3}}}
		\end{array}} \right],\\
		{\boldsymbol{k}}_\chi  &\coloneqq \left[ {\renewcommand{\arraystretch}{1.5}\begin{array}{*{20}{c}}
				{{{\tilde l}^{11}_\mathrm{p}}{{\boldsymbol{1}}_{3\times3}}}&{{{\tilde l}^{12}_\mathrm{p}}{{\boldsymbol{1}}_{3\times3}}}\\
				{{\rm{sym}}{\rm{.}}}&{{{\tilde l}^{22}_\mathrm{p}}{{\boldsymbol{1}}_{3\times3}}}
		\end{array}} \right],
	\end{alignat}
\end{subequations}
\endgroup
and
\begin{equation}	
	{{\Bbb{Y}}_e} \coloneqq {\left[ {\begin{array}{*{20}{c}}
				{N_{1,s}^p{{\bf{1}}_{3 \times 3}}}& \cdots &{N_{{n_e},s}^p{{\bf{1}}_{3 \times 3}}}&\vline& {}&{{{\bf{0}}_{3 \times 6n_e^{\rm{d}}}}}&{}\\
				{}&{{{\bf{0}}_{6 \times 3{n_e}}}}&{}&\vline& {N_{1,s}^{{p_{\rm{d}}}}{{\bf{1}}_{6 \times 6}}}& \cdots &{N_{n_e^{\rm{d}},s}^{{p_{\rm{d}}}}{{\bf{1}}_{6 \times 6}}}\\
				{}&{{{\bf{0}}_{6 \times 3{n_e}}}}&{}&\vline& {N_1^{{p_{\rm{d}}}}{{\bf{1}}_{6 \times 6}}}& \cdots &{N_{n_e^{\rm{d}}}^{{p_{\rm{d}}}}{{\bf{1}}_{6 \times 6}}}
		\end{array}} \right]},
\end{equation}
such that
\begin{equation}	
	\left\{ {\begin{array}{*{20}{c}}
			{\Delta {\boldsymbol{\varphi }}_{,s}^h}\\
			{\Delta {\boldsymbol{d}}_{1,s}^h}\\
			{\Delta {\boldsymbol{d}}_{2,s}^h}\\
			{\Delta {\boldsymbol{d}}_1^h}\\
			{\Delta {\boldsymbol{d}}_2^h}
	\end{array}} \right\} = {{\Bbb{Y}}_e}\delta {{\bf{y}}^e},\,\,\mathrm{with}\,\,\,\delta {{\bf{y}}^e} \coloneqq \left\{ {\begin{array}{*{20}{c}}
	{\delta {\boldsymbol{\varphi }}_1^e}\\
	\vdots \\
	{\delta {\boldsymbol{\varphi }}_{{n_e}}^e}\\
	{\delta {\bf{d}}_1^e}\\
	\vdots \\
	{\delta {\bf{d}}_{n_e^{\rm{d}}}^e}
\end{array}} \right\}.
\end{equation}
\subsection{Selectively reduced degree $p_\mathrm{p}$ in each case of $p=2,3,...,10$}
\textcolor{black}{Table \ref{tab_selec_deg_r} shows a list of selective reduced degrees $p_\mathrm{p}$ for physical stress resultants and strains, employed in the local approach ``loc-sr''. Note that $p_\mathrm{p}$ is always the same for the physical stress resultant and strain in every work-conjugate pair, so that ${\bf{k}}^e_{\mathrm{r}\varepsilon}$ of Eq.\,(\ref{la_submat_k_re}) is always a square (invertible) matrix.}
\begin{table}[H]
	\centering
	\caption{\textcolor{black}{Selection of degree $p_\mathrm{p}$ in the approach ``IGA, loc-sr''. $p^\mathrm{b}_\mathrm{p}$ denotes the degree at the first and last elements, which is shown here only if $p_\mathrm{p}\ne{p^\mathrm{b}_\mathrm{p}}$. Note that $\alpha,\beta\in\left\{1,2\right\}$.}}
	\label{tab_selec_deg_r}
	\begin{tabular}{lllllllllll}
		\toprule
		\multicolumn{1}{l}{\multirow{2}[4]{*}{\makecell{Physical\\strains}}} & \multicolumn{1}{l}{\multirow{2}[4]{*}{\makecell{Physical stress\\resultants}}} & \multicolumn{9}{c}{\multirow{1}[2]{*}{$p_\mathrm{p}$ ($p_\mathrm{p}^\mathrm{b}$)}} \\
		\cmidrule{3-11}      &       & $p=2$ & 3     & 4     & 5     & 6     & 7     & 8     & 9     & 10 \\
		\midrule
		$\varepsilon_\mathrm{p}$ & ${\tilde n}_\mathrm{p}$ & 0 (1) & 0 (1) & 0 (1) & 1     & 1     & 1     & 1     & 1     & 2 \\
		$\rho^\mathrm{p}_\alpha$ & ${\tilde m}^\alpha_\mathrm{p}$ & 1     & 1     & 2     & 2     & 2     & 2     & 2     & 2     & 3 \\
		$\kappa^\mathrm{p}_{\alpha\beta}$ & ${\tilde h}^{\alpha\beta}_\mathrm{p}$ & 1     & 1     & 2     & 2     & 2     & 2     & 2     & 2     & 3 \\
		$\delta^\mathrm{p}_\alpha$ & ${\tilde q}^\alpha_\mathrm{p}$ & 0 (1) & 0 (1) & 0 (2) & 1 (2) & 1 (2) & 1 (2) & 1 (2) & 1 (2) & 2 \\
		$\gamma^\mathrm{p}_{\alpha\beta}$ & ${\tilde m}^{\alpha\beta}_\mathrm{p}$ & 0 (1) & 0 (1) & 0 (1) & 0 (1) & 0 (1) & 0 (1) & 1 & 1 (2) & 2 \\
		$\chi^\mathrm{p}_{\alpha\beta}$ & ${\tilde \ell}^{\alpha\beta}_\mathrm{p}$ & 0     & 0     & 0 (1) & 0     & 1     & 1     & 1     & 1     & 2 \\
		\bottomrule
	\end{tabular}%
\end{table}
\subsection{Imposition of rotational continuity between beams}
\label{app_jct_cond_index_corr}
Table \ref{tab_corr_indic_jct} shows the global control coefficient index $K$ for the director displacement field, corresponding to the given local coefficient index $I$ in $e$th element, in two different cases of the location of junction. 
\begin{table}[H]
	\centering
	\caption{Two cases of the joint condition: Correspondence between the local coefficient index $I$ in $e$th element, and the global coefficient index $K$.}
	\begin{tabular}{cccc}
		\toprule
		& Element number $e$ & Local index $I$ & Global index $K$ \\
		\midrule
		Case 1 & 1     & 1     & 1 \\
		Case 2 & $n_\mathrm{el}$ & $n^\mathrm{d}_e$ & $n^\mathrm{d}_\mathrm{cp}$ \\
		\bottomrule
	\end{tabular}%
	\label{tab_corr_indic_jct}%
\end{table}%
\section{Appendix: Consistent mass matrix}
\label{app_cons_mass_mat}
We have the kinetic energy bilinear form for a three-dimensional body, as
\begin{equation}
	{d}({{\boldsymbol{x}}}_{,tt},\delta {{\boldsymbol{x}}}) = \int_\mathcal{B} {{\rho _0}\,\delta {\boldsymbol{x}}\cdot{{\boldsymbol{x}}_{,tt}}\,{j_0}\,{\mathrm{d}}\mathcal{B}},
\end{equation}
where $(\bullet)_{,tt}$ denotes the second order derivative with respect to time $t$. Applying the beam kinematics in Eq.\,(\ref{bkin_x_d_zta_ph}), we have
\begin{equation}
	d({\boldsymbol{y}}_{,tt},\delta {\boldsymbol{y}}) = \int_0^L {\delta {{\boldsymbol{y}}^{\rm{T}}}{{\bf{I}}_\rho }\,{\boldsymbol{y}}_{,tt}\,{\mathrm{d}}s},
\end{equation}
where ${\bf{I}}_\rho$ is an inertia matrix, defined by
\begin{equation}
	{{\bf{I}}_\rho } \coloneqq \left[ {\renewcommand\arraystretch{1.5}\begin{array}{*{20}{c}}
			{{\rho_A}{{\bf{1}}_{3\times3}}}&{I_\rho^{1}{{\bf{1}}_{3\times3}}}&{I_\rho^{2}{{\bf{1}}_{3\times3}}}\\
			{}&{I_\rho^{11}{{\bf{1}}_{3\times3}}}&{I_\rho^{12}{{\bf{1}}_{3\times3}}}\\
			{{\rm{sym}}{\rm{.}}}&{}&{I_\rho ^{22}{{\bf{1}}_{3\times3}}}
	\end{array}} \right],
\end{equation}
with the initial line density (mass per unit undeformed length)
\begin{equation}
	\rho_A\coloneqq \int_\mathcal{A} {{\rho _0}\,{j_0}\,{\rm{d}}\mathcal{A}},
\end{equation}
and the first moment of inertia components
\begin{equation}
	I_\rho ^{\gamma} \coloneqq \int_\mathcal{A} {{\rho _0}\,{\zeta ^\gamma}{j_0}\,{\mathrm{d}}\mathcal{A}},\,\,\gamma\in\left\{1,2\right\},
\end{equation}
and the second moment of inertia components
\begin{equation}
	I_\rho ^{\gamma \delta } \coloneqq \int_\mathcal{A} {{\rho _0}\,{\zeta ^\gamma }{\zeta ^\delta }{j_0}\,{\mathrm{d}}\mathcal{A}},\,\,\gamma,\delta\in\left\{1,2\right\}.
\end{equation}
It should be noted that the inertia matrix is constant during deformation, and it depends only on the initial geometry and the initial mass density distribution. Applying the spatial discretization using NURBS basis functions, we have
\begin{equation}
	d({{\boldsymbol{y}}_{,tt}},\delta {{\boldsymbol{y}}}) \approx \delta {{\bf{y}}^{\rm{T}}}{\bf{M}}\,{\bf{y}}_{,tt}\,,\,\,\mathrm{with}\,\,{\bf{M}} \coloneqq \mathop {\mathop {\mathlarger{\mathlarger{\bf{A}}}}}\limits_{e = 1}^{{n_{{\rm{el}}}}} {{\bf{m}}^e},
\end{equation}
where the element mass matrix is obtained as
\begin{equation}
	{{\bf{m}}^e} \coloneqq \int_{{\Xi ^e}} {{{{\Bbb{N}}_e}{{(\xi )}^{\rm{T}}}{{\bf{I}}_\rho\, }{{\Bbb{N}}_e}(\xi )\,\tilde j}\,{\rm{d}}\xi}.
\end{equation}
Note that the mass matrix is constant in deformation.
\section{Appendix to numerical examples}
\subsection{Cantilever beam under bending moment}
\textcolor{black}{Table \ref{app_tab_cant_hist_disp-b_fea_uri} shows the Newton-Raphson iteration history for the displacement-based formulation, with the increased load step number, 100, from that in the FEA result of Table \ref{ex_end_bend_mnt_hist_h1e-4}.} 
\begin{table}[H]
	\centering
	\caption{Cantilever beam under bending moment: Newton-Raphson iteration history for the displacement-based formulation using FEA (URI) with $p=p_\mathrm{d}=2$, and $n_\mathrm{el}=10$. Total number of load steps is 100, and here we present only the history in the first ($n=1$) and last ($n=100$) load steps.}
	\label{app_tab_cant_hist_disp-b_fea_uri}
	\begin{tabular}{llllll}
		\toprule
		\multirow{3}[3]{*}{Iteration\#} & \multicolumn{2}{c}{$n=1$} &       & \multicolumn{2}{c}{$n=100$} \\		
		\cmidrule{2-3}\cmidrule{5-6}
		 								& {\makecell{Euclidean norm\\of residual}} & {\makecell{Energy\\norm}} &    & \makecell{Euclidean norm\\of residual} & \makecell{Energy\\norm} \\
		\midrule
		1     & 6.3E-09 & 3.9E-10 &       & {6.3E-09} & {3.9E-10} \\
		2     & 3.4E+00 & 1.9E-02 &       & {3.7E+00} & {1.8E-02} \\
		3     & 8.7E-03 & 9.0E-08 &       & {9.1E-03} & {8.6E-08} \\
		4     & 6.8E-08 & 9.5E-16 &       & {6.6E-08} & {8.6E-15} \\
		5     & 1.9E-06 & 2.6E-15 &       & {1.3E-05} & {1.4E-13} \\
		6     & 2.0E-09 & 1.9E-16 &       & {3.3E-10} & {1.4E-15} \\
		7     & 1.6E-06 & 4.3E-15 &       & {1.3E-05} & {2.5E-13} \\
		8     & 2.4E-12 & 4.7E-19 &       & {1.4E-11} & {3.4E-20} \\
		9     & 4.2E-09 & 3.4E-20 &       & {3.3E-10} & {1.6E-22} \\
		10    & 9.2E-13 & 8.1E-28 &       &       	  &  \\
		\bottomrule
	\end{tabular}
\end{table}
\subsection{$45^\circ$-arc cantilever beam}
\textcolor{black}{Table \ref{app_num_ex_cant45_ext_f} shows the applied external force in each case of the slenderness ratio $R/d$.}
\begin{table}[H]
	\centering
	\caption{$45^\circ$-arc cantilever beam: \textcolor{black}{Applied force in each case of the slenderness ratio $R/d$.}}
	\label{app_num_ex_cant45_ext_f}	
	\begin{tabular}{lllcll}		
		\toprule
		$R\,[\mathrm{m}]$ & {$d\,[\mathrm{m}]$} & {$A\,[\mathrm{m}^2]$} & {$R/d\,[-]$} & $F\,[\mathrm{N}/\mathrm{m}^2]$ & $F\cdot{A}\,[\mathrm{N}]$ \\
		\midrule				
		\multirow{4}[1]{*}{$10^2$} 	& {$10^0$} 		& {$10^{0}$} 	& {$10^{2}$}  & $6\times{10^2}$		& $6\times{10^2}$  		\\
									& {$10^{-1}$} 	& {$10^{-2}$} 	& {$10^{3}$}  & $6\times{10^0}$		& $6\times{10^{-2}}$  	\\
									& {$10^{-2}$} 	& {$10^{-4}$} 	& {$10^{4}$}  & $6\times{10^{-2}}$ 	& $6\times{10^{-6}}$ 	\\
									& {$10^{-3}$} 	& {$10^{-6}$} 	& {$10^{5}$}  & $6\times{10^{-4}}$  & $6\times{10^{-10}}$	\\
		\bottomrule
	\end{tabular}
\end{table}
\noindent Table \ref{app_num_ex_cant45_c1_conv} presents the convergence test results of the brick element solution of the tip displacement for Case 1 considering $R/d=10^2$.
\begin{table}[H]
	\centering
	\caption{$45^\circ$-arc cantilever beam (Case 1, $R/d=10^2$): Convergence of the brick element solution for $d=1\,\mathrm{m}$. In all cases, we use IGA.}
	\begin{tabular}{cccrlll}
		\Xhline{3\arrayrulewidth}
		\multicolumn{3}{l}{Brick, $\mathrm{deg.}=(3,3,3)$} &       & \multicolumn{3}{l}{Tip displacements} \\
		\cline{5-7}
		\multicolumn{3}{l}{$n_\mathrm{el}$} &       & $u_1\,[\mathrm{m}]$ & $u_2\,[\mathrm{m}]$ & $u_3\,[\mathrm{m}]$ \\
		\Xhline{3\arrayrulewidth}
		\multicolumn{3}{l}{$40\times{1}\times{1}$} &       & 1.3729E+01 & -2.3822E+01 & 5.3607E+01 \\
		\multicolumn{3}{l}{$80\times{1}\times{1}$} &       & 1.3730E+01 & -2.3825E+01 & 5.3609E+01 \\
		\multicolumn{3}{l}{$160\times{1}\times{1}$} &       & 1.3731E+01 & -2.3826E+01 & 5.3610E+01 \\
		\multicolumn{3}{l}{$240\times{5}\times{5}$} &       & 1.3731E+01 & -2.3826E+01 & 5.3610E+01 \\
		\multicolumn{3}{l}{$240\times{8}\times{8}$} &       & 1.3731E+01 & -2.3826E+01 & 5.3610E+01 \\
		\multicolumn{3}{l}{$320\times{8}\times{8}$} &       & 1.3731E+01 & -2.3826E+01 & 5.3610E+01 \\
		\Xhline{3\arrayrulewidth}
	\end{tabular}%
	\label{app_num_ex_cant45_c1_conv}%
\end{table}%

\section*{Acknowledgement}
M.-J. Choi would like to gratefully acknowledge the financial support of a postdoctoral research fellowship from the Alexander von Humboldt Foundation in Germany.


\bibliography{mybibfile}

\end{document}